\documentclass[12pt]{article}
\usepackage[utf8]{inputenc}
\usepackage[russian]{babel}
\usepackage{amsmath,amssymb}
\usepackage{amsthm}
\usepackage{amscd}

\usepackage{array}
\usepackage[matrix,arrow,curve]{xy}
\usepackage{color}
\usepackage{bbm}
\usepackage{hyperref}

\date{}

\def \nlb {\nolinebreak}
\newcommand{\udot}{{\:\raisebox{3pt}{\text{\circle*{1.5}}}}}
\newcommand\Mod{\operatorname{Mod}}
\newcommand\Z{{\mathbb Z}}
\newcommand\Pb{{\mathbb P}}
\newcommand\Ker{\operatorname{Ker}}
\newcommand\Coker{\operatorname{Coker}}
\newcommand\Hom{\operatorname{Hom}}
\newcommand\Map{\operatorname{Map}}
\newcommand\Ext{\operatorname{Ext}}
\newcommand\Aut{\operatorname{Aut}}
\newcommand\Supp{\operatorname{Supp}}
\newcommand\SSupp{\mathcal{S}upp}
\newcommand\GL{\operatorname{GL}}
\newcommand\Aff{\operatorname{Aff}}
\newcommand\Sp{\operatorname{Sp}}
\newcommand\Ort{\operatorname{O}}
\newcommand\SL{\operatorname{SL}}
\newcommand\PGL{\operatorname{PGL}}

\newcommand\res{\operatorname{res}}
\newcommand\End{\operatorname{End}}
\newcommand\Mat{\operatorname{M}}
\newcommand\Br{\operatorname{Br}}
\newcommand\Pic{\operatorname{Pic}}
\newcommand\NS{\operatorname{NS}}
\newcommand\Nrd{\operatorname{Nrd}}
\newcommand\Gal{\operatorname{Gal}}
\newcommand\KM{\operatorname{K}}
\newcommand\Mcal{{\mathcal M}}
\newcommand\Fc{{\mathcal F}}
\newcommand\OO{{\mathcal O}}
\newcommand\LL{{\mathcal L}}
\newcommand\QQ{{\mathcal Q}}

\newcommand\Xc{{\mathcal X}}
\newcommand{\Cb}{{\mathbb C}}
\newcommand{\Q}{{\mathbb Q}}
\renewcommand{\P}{{\mathbb P}}
\newcommand{\A}{{\mathbb A}}
\newcommand{\Rb}{{\mathbb R}}
\newcommand{\Gb}{{\mathbb G}}
\newcommand{\F}{{\mathbb F}}
\newcommand\N{{\mathrm N}}
\def\G {\mathrm G}
\newcommand\trdeg{{\mathrm{tr.deg}}}
\newcommand\Nm{{\mathrm{Nm}}}
\newcommand\Spec{\mathrm{Spec}}
\newcommand\Div{\mathrm{Div}}
\newcommand\Sing{\mathrm{Sing}}
\newcommand\Tr{\operatorname{Tr}}
\newcommand\m{{\mathfrak m}}
\newcommand\p{{\mathfrak p}}
\newcommand\cl{{\mathrm{cl}}}
\newcommand\mult{{\mathrm{mult}}}
\newcommand\disc{{\mathrm d}}
\newcommand\ord{\mathrm{ord}}
\newcommand\codim{{\mathrm{codim}}}
\newcommand{\Ab}{{\mathbb A}}
\newcommand{\Fb}{{\mathbb F}}

\makeatletter\@addtoreset{equation}{section} \makeatother

\newcounter{probcounter}[section]

\makeatletter\@addtoreset{probcounter}{subsection} \makeatother

\theoremstyle{remark}
\newtheorem{examp}[probcounter]{Пример}

\theoremstyle{definition}
\newtheorem{defin}[probcounter]{Определение}
\newtheorem{prob}[probcounter]{Упражнение}

\title{Неразветвлённая группа Брауэра и её приложения}
\author{С.\,О.\,Горчинский, К.\,А.\,Шрамов}

\begin{document}

\maketitle
\thispagestyle{empty}

\newpage
\tableofcontents

\newpage
\section*{Предисловие}
\addcontentsline{toc}{section}{Предисловие}
\phantom{a}
\vskip 1cm
Эта книга представляет собой переработанные
материалы учебного семинара ``Арифметические
методы в алгебраической геометрии'', который авторы вели в НОЦ МИАН
весной~2011 года. Цель семинара состояла в том, чтобы дать введение в теорию
неразветвлённых групп Брауэра и стабильной рациональности, начиная с самых
первых определений теории когомологий групп. В соответствии с этой целью
мы проигнорировали многие популярные и хорошо освещённые в учебниках темы
(относящиеся к когомологиям Галуа, группам Брауэра и т.\,п.),
которые можно было бы затронуть по пути.
С другой стороны, достаточно много внимания уделено подробному разбору
примеров доказательства стабильной нерациональности с помощью неразветвлённой
группы Брауэра и примера нерационального стабильно рационального
многообразия, так как эти темы, насколько нам известно, в учебной литературе
не затрагивались, а для того, чтобы проследить за доказательствами
в статьях-первоисточниках, требуются некоторые усилия.
Также стиль семинара определил и форму подачи материала в задачах. При этом
мы постарались разбить доказательства почти всех нужных нам фактов
на как можно более простые шаги (которые, собственно,
и представлены в качестве
задач) и снабдить все нетривиальные места указаниями.
Поэтому мы надеемся, что
изучать этот текст будет
не сложнее (или по крайней мере не сильно сложнее),
чем читать обычные учебники (и тем более статьи).
Для понимания большей части текста требуется знание основ алгебры,
теории Галуа и начальных понятий алгебраической геометрии.

В главе~\ref{section:group-cohomology} мы собрали необходимые для дальнейшего
определения и утверждения, относящиеся к когомологиям абстрактных групп.
В главе~\ref{section:Galois-cohomology} то же сделано
для когомологий Галуа, а в главе~\ref{section:Brauer} --- для группы Брауэра.
Поскольку значение этих тем далеко не ограничивается их применениями
к стабильной рациональности, и знакомство с ними давно стало частью
общей математической культуры, мы рекомендуем заинтересовавшемуся
читателю продолжить их изучение по подробным классическим учебникам.
Для изучения когомологий групп
подойдёт глава~IV в книге~\cite{CasselsFrolich},
для изучения когомологий Галуа ---
книга Ж.-П.\,Серра~\cite{Serre-Galois}
и глава~V  в~\cite{CasselsFrolich},
для знакомства с группой Брауэра --- глава~X книги
\mbox{Ж.-П.\,Серра~\cite{Serre-corps-locaux}} и глава~VIII
книги~\cite{Bourbaki-algebra}.
Далее, в главе~\ref{section:unramified-Brauer} мы переходим к более
специальным вопросам, связанным с группой Брауэра полей дискретного
нормирования, и, в частности, определяем
неразветвлённую группу Брауэра поля. В качестве дополнительного чтения
по этим темам можно порекомендовать книгу~\cite{Serre-corps-locaux} и~\S1
главы~VI в книге~\cite{CasselsFrolich}.
Кроме того, почти весь материал из
глав~\ref{section:group-cohomology}--\ref{section:unramified-Brauer}
содержится в существенно более подробной форме в учебнике~\cite{GS}.
Также о когомологиях Галуа можно почитать в доступной форме
в~\cite{Ber}.

В главе~\ref{section:primer-Bogomolova} разбирается пример многообразия
вида~\mbox{$X=V/G$}, где~$G$~--- конечная группа, а~$V$~--- её представление,
определённое над алгебраически замкнутым полем $k$ характеристики нуль,
для которого при помощи введённых
в предыдущих главах понятий удаётся
установить нерациональность (и даже стабильную нерациональность).
Она следует из нетривиальности неразветвлённой группы Брауэра
поля $k(X)$, то есть поля инвариантов $k(V)^G$.
Примеры такого сорта впервые появились в статьях Д.\,Солтмена~\cite{Saltman}
и Ф.\,А.\,Богомолова~\cite{Bogomolov}, однако источником для нашего изложения
послужила статья И.\,Р.\,Шафаревича~\cite{Shafarevich}, в которой подход
Солтмена и Богомолова был значительно упрощён.
Многообразие $X$ имеет довольно большую размерность,
поэтому представляет интерес вопрос о том, можно ли применить
неразветвлённую группу Брауэра для доказательства нерациональности многообразий
меньшей размерности. Оказывается, что это можно сделать для некоторых
многообразий размерности~$3$ (которые
строятся совершенно другим образом, чем многообразия из
главы~\ref{section:primer-Bogomolova}).
Соответствующий пример разбирается
в главе~\ref{section:Artin-Mumford} --- это хорошо известный пример
нерационального двойного накрытия проективного
пространства~$\P^3$ с ветвлением в (особой) квартике. Первоисточником для этого
примера является статья М.\,Артина и Д.\,Мамфорда~\cite{ArtinMumford},
однако наше изложение основывается на более простом подходе
М.\,Гросса (см.~\cite[Appendix]{Aspinwall-i-td}). Перед тем, как
привести эту конструкцию, мы уделяем в главе~\ref{section:2-dim-quadrics}
некоторое время введению инварианта
Клиффорда и доказательству
необходимых вспомогательных утверждений о квадриках над
алгебраически незамкнутыми полями. Подробное изложение
теории квадрик над незамкнутыми полями можно найти в книге~\cite{Karpenko}.
В главе~\ref{section:Artin-Mumford} в числе прочего приведена конструкция
унирациональности двойного накрытия проективного
пространства~$\P^3$ с ветвлением в квартике
(наиболее существенная часть этой конструкции изложена в соответствии
с доказательством
теоремы~IV.7.7 в книге Ю.\,И.\,Манина~\cite{Manin-KubFormy}).

В главе~\ref{section:Weil} вводится понятие ограничения скаляров по
Вейлю и обсуждаются его основные свойства, а также
устанавливаются некоторые нужные в дальнейшем
в главе~\ref{section:Polietilen} свойства алгебраических торов.
Более подробно об алгебраических торах можно прочесть в книге
В.\,Е.\,Воскресенского~\cite{Voskresenskii-rus}.
Также в главе~\ref{section:Weil}
объясняется понятие универсального торсора,
и устанавливаются простейшие свойства поверхностей Шатле.
После того, как у нас
имеются примеры многообразий, для которых мы умеем строить препятствия
к стабильной рациональности, естественно задаться вопросом, не сводится
ли это понятие к обычному понятию рациональности. Оказывается, что не
сводится, однако привести разделяющий эти два класса многообразий пример
не так просто. Это делается в главе~\ref{section:Polietilen},
следуя статье А.\,Бовиля, Ж.-Л.\,Кольо-Телена,
Ж.-Ж.\,Сансюка и П.\,Свиннертон-Дайера~\cite{CT-i-drugie}.
В конце главы~\ref{section:Polietilen}
приводится рассуждение Н.\,Шепард-Баррона из~\cite{ShB},
несколько усиливающего конструкцию из~\cite{CT-i-drugie}.
На протяжении всей главы~\ref{section:Polietilen} мы старались использовать
геометрический язык и как можно дольше
не прибегать к введению координат и рассмотрению явных уравнений, что, как
мы надеемся, может сделать изложение чуть более прозрачным,
чем в~\cite{CT-i-drugie} и~\cite{ShB}.

Главы~\ref{section:Min-Hasse}
и~\ref{section:Br-Manin} посвящены ещё одному приложению
неразветвлённой группы Брауэра --- построению препятствия Брауэра--Манина.
Глава~\ref{section:Min-Hasse} является в какой-то мере мотивировочной:
в ней мы излагаем доказательство
классической теоремы Минковского--Хассе для квадрик
(точнее говоря, мы выводим эту теорему из некоторых фундаментальных
фактов теории полей классов).
В целом следуя
главе~IV в книге Ж.-П.\,Серра~\cite{Serr-kurs-arifmetiki}, мы стараемся
в большей степени использовать геометрический язык при сведении многомерного
случая к одномерному.
В главе~\ref{section:Br-Manin} определяется препятствие
Брауэра--Манина и при помощи него строится контрпример к утверждению,
аналогичному теореме Минковского--Хассе, для кривых рода~$1$.
Подробнее про препятствие Брауэра--Манина написано, например, в~\cite{Skorobogatov} (см. также краткое изложение в~\cite[5.2.3]{PanchishkinManin}).

В приложении~\ref{section:etale} содержится подборка
ссылок на основные утверждения
об этальных когомологиях, которые надо знать, чтобы в терминах
этальных когомологий проинтерпретировать неразветвлённую группу Брауэра
(более подробно об этом можно почитать, например,
в~\cite{Dan} или~\cite{Mil}).
Любителям читать первоисточники можно в числе прочего порекомендовать
текст А.\,Гротендика~\cite{Grothendieck66b},
где именно таким образом
была впервые введена неразветвлённая группа Брауэра.

\smallskip
Разумеется, многим связанным с неразветвлённой группой
Брауэра темам мы не смогли уделить внимание.
Подробный обзор вопросов, связанных со стабильной рациональностью
фактормногообразий (в частности, множество конструкций
стабильной рациональности для таких многообразий) и построением препятствий
к ней при помощи неразветвлённой группы Брауэра,
содержится в~\cite{CT-Sansuc}. Там же приведена большая подборка
ссылок на оригинальные работы по этой теме.
Эффективные методы вычисления неразветвлённой группы Брауэра
поля инвариантов $k(V)^G$ в терминах бициклических подгрупп
группы~$G$ описаны в~\cite{Bogomolov}
(см. также~\mbox{\cite[\S7]{CT-Sansuc}}).

Неразветвлённая группа Брауэра не является единственным известным
препятствием к стабильной рациональности: существуют обобщающие её
старшие неразветвлённые когомологические группы,
причём вторые когомологии как
раз соответствуют неразветвлённой группе Брауэра.
Про них можно почитать, например, в~\cite{CT-O}
и~\cite{Peyre-unramified}, где также содержатся примеры
препятствий к стабильной рациональности, строящихся с помощью
старших когомологий.
Особенно подробно изучена группа третьих неразветвлённых когомологий,
и даже существует
аналогичный разбиравшемуся
в главе~\ref{section:primer-Bogomolova}
явный пример многообразия вида~\mbox{$X=V/G$}, где~$V$~---
представление конечной группы $G$ над алгебраически
замкнутым полем характеристики нуль (см.~\cite{Peyre-gruppa}).
При этом нерациональность многообразия~$X$ нельзя доказать при помощи
неразветвлённой группы Брауэра,
но можно при помощи третьих неразветвлённых
когомологий.

Пример, разобранный в главе~\ref{section:primer-Bogomolova},
не является оптимальным в том смысле, что группа $G$ в этом примере имеет
порядок~$p^9$, где~\mbox{$p>2$}~--- простое число,
при том что на самом деле порядок такой группы можно понизить до $p^6$,
а также построить аналогичный пример с $p=2$ (см.~\cite{Bogomolov}).
Более того, при $p>2$ можно привести такой пример даже для группы
порядка $p^5$ (см.~\cite{KangHoshi}), а вот для групп~$G$ порядка~$p^4$ все
факторы $V/G$ представлений $G$, определённых над алгебраически
замкнутым полем характеристики нуль, уже оказываются рациональными
(см.~\cite{ChuKang}).
Отметим также, что над некоторыми алгебраически
незамкнутыми полями (в частности, над~$\Q$) можно построить
нерациональное многообразие вида~\mbox{$V\slash\Gamma$},
где~$V$~--- представление циклической группы $\Gamma\cong\Z/n\Z$
(см.,~например,~\cite{Swan}).

Мы не смогли уделить внимания
сравнительно новой и популярной теме, имеющей непосредственное
отношение к стабильной рациональности~---
универсальной тривиальности групп Чжоу. Заинтересованный читатель может обратиться
к статьям \cite{Voisin} и~\cite{CT-P}, а также к обзору~\cite{Pirutka}.

\medskip
В заключение мы бы хотели выразить благодарность
за конструктивную и стимулирующую атмосферу семинара всем его участникам.
Отдельно хотелось бы отметить, что
на разных этапах создания этих записок нам помогали полезными советами
А.\,Вишик, В.\,Вологодский, В.\,Жгун, Ж.-Л.\,Кольо-Телен,
А.\,Кузнецов, И.\,Нетай, Ю.\,Прохоров, С.\,Рыбаков,
А.\,Скоробогатов, М.\,Тёмкин, Д.\,Теста, А.\,Трепалин, А.\,Фонарёв,
Т.\,Шабалин и Е.\,Шиндер. Во время работы над книгой авторы получали поддержку от
фонда Д.\,Зимина ``Династия'', из средств государственной поддержки ведущих университетов Российской Федерации ``5-100'', а также от Лаборатории зеркальной симметрии НИУ ВШЭ, грант правительства РФ Договор № 14.641.31.0001.


\newpage

\section*{Список обозначений}
\addcontentsline{toc}{section}{Список обозначений}

\noindent $\Z$ --- кольцо целых чисел

\medskip

\noindent $\Q$ --- поле рациональных чисел

\medskip

\noindent $\Rb$ --- поле вещественных чисел

\medskip

\noindent $\Cb$ --- поле комплексных чисел

\medskip

\noindent $\F_q$ --- конечное поле из $q$ элементов

\medskip

\noindent $\Z_p$ --- кольцо целых $p$-адических чисел

\medskip

\noindent $\Q_p$ --- поле рациональных $p$-адических чисел

\medskip

\noindent $\widehat{\Z}$ --- проконечное пополнение группы $\Z$ (по сложению)

\medskip

\noindent $k(t_1,\ldots,t_n)$ ---
поле рациональных функций от формальных переменных $t_1,\ldots,t_n$
с~коэффициентами в поле~$k$

\medskip

\noindent $k((t))$ --- поле рядов Лорана от формальной переменной~$t$
с коэффициентами в поле~$k$

\medskip

\noindent $\Hom(X,Y)$ --- множество морфизмов между
объектами $X$ и $Y$ (в категории, которая обычно очевидна из контекста)

\medskip

\noindent $\Aut(X)$ --- группа автоморфизмов
объекта $X$

\bigskip

\bigskip

\noindent $\mathrm{S}_n$ --- группа перестановок множества из $n$ элементов

\medskip

\noindent $\mathrm{Stab}_G(x)$ ---
стабилизатор в группе $G$ элемента $x$ множества $X$ с действием $G$

\medskip

\noindent $G/H$ --- множество левых смежных классов группы $G$ по подгруппе
$H\subset G$

\medskip

\noindent $A_n$ --- $n$-кручение в абелевой группе $A$ для натурального
числа $n$

\medskip

\noindent $\Z[S]$ --- свободная абелева группа, порождённая множеством $S$

\medskip

\noindent $\langle S\rangle$ --- подгруппа,
порождённая подмножеством $S$, в группе или двусторонний идеал,
порождённый подмножеством $S$, в ассоциативной алгебре

\medskip

\noindent $M^G$ --- группа $G$-инвариантных элементов в $G$-модуле $M$

\medskip

\noindent $\Hom_G(M,M')$ --- группа морфизмов $G$-модулей из $M$ в $M'$

\medskip

\noindent $H^i(G,M)$ --- группа $i$-ых когомологий группы $G$ с коэффициентами в $G$-модуле $M$

\medskip

\noindent $Z/\Gamma$ --- множество орбит группы $\Gamma$, действующей
на множестве $Z$

\bigskip

\bigskip

\noindent ${\rm char}(R)$ --- характеристика кольца $R$

\medskip

\noindent $R^*$ --- мультипликативная группа обратимых элементов кольца $R$

\medskip

\noindent $\KM_2(K)$ --- вторая $\KM$-группа Милнора поля $K$

\medskip

\noindent $\bar K$ --- алгебраическое замыкание поля $K$

\medskip

\noindent $K^{sep}$ --- сепарабельное замыкание поля $K$

\medskip

\noindent $\mu_n$ --- группа корней $n$-ой степени из $1$ в $K^{sep}$, где $n$ взаимно просто с ${\rm char}(K)$

\medskip

\noindent $\dim_K(V)$ --- размерность векторного пространства $V$ над полем $K$

\medskip

\noindent $[L:K]=\dim_K(L)$ --- степень конечного расширения полей $K\subset L$

\medskip

\noindent $\Nm_{L/K}\colon L^*\to K^*$ --- норма Галуа для
конечного сепарабельного расширения полей~\mbox{$K\subset L$}

\medskip

\noindent $\Gal(L/K)$ --- группа Галуа для расширения Галуа $K\subset L$

\medskip

\noindent $G_K=\Gal(K^{sep}/K)$ --- абсолютная группа Галуа поля $K$

\bigskip

\bigskip

\noindent $\GL_n$ --- группа обратимых матриц порядка $n$

\medskip

\noindent $\Gb_m=\GL_1$ --- мультипликативная групповая схема

\medskip

\noindent $\SL_n$ --- группа обратимых матриц порядка $n$ с тривиальным определителем

\medskip

\noindent $\PGL_n=\GL_n/\Gb_m$ ---
группа обратимых матриц порядка $n$ с точностью до скалярных матриц, то есть
группа автоморфизмов $(n-1)$-мерного проективного пространства

\medskip

\noindent $\Pb^n$ --- проективное пространство размерности $n$

\medskip

\noindent $\OO_{\Pb^n}(r)=\OO_{\Pb^n}(1)^{\otimes r}$ --- $r$-ая тензорная
степень линейного расслоения $\OO_{\Pb^n}(1)$ на $\Pb^n$,
двойственного к тавтологическому расслоению~$\OO_{\Pb^n}(-1)$,
где $r$ --- целое число

\medskip

\noindent $V_L=L\otimes_K V$ --- расширение скаляров векторного пространства $V$ над полем~$K$ для расширения полей~\mbox{$K\subset L$}

\medskip

\noindent $\GL(V)$ --- группа $K$-линейных автоморфизмов векторного пространства $V$ над $K$

\medskip

\noindent $\Pb(V)\cong \Pb^{n-1}$ --- проективизация
$n$-мерного векторного пространства $V$ над полем $K$
(множество $K$-точек на $\Pb(V)$
равно множеству одномерных подпространств в $V$)

\bigskip

\bigskip

\noindent $\Spec(R)$ --- спектр коммутативного кольца $R$ с единицей

\medskip

\noindent $k[X]$ --- кольцо регулярных функций на многообразии~$X$
над полем~$k$

\medskip

\noindent $k(X)$ --- поле рациональных функций
на неприводимом многообразии $X$ над полем~$k$

\medskip

\noindent $\dim(X)$, $\dim_k(X)$ --- размерность
многообразия $X$ над полем $k$

\medskip

\noindent $\Div(X)$ --- группа дивизоров на гладком многообразии $X$

\medskip

\noindent $\Pic(X)$ --- группа Пикара многообразия $X$

\medskip

\noindent ${\mathbb T}_x(X)\subset \Pb^n$ --- вложенное
проективное касательное
пространство в точке $x$ к проективному многообразию $X\subset \Pb^n$

\medskip

\noindent $X_K$ --- расширение скаляров многообразия $X$, определённого над полем~$k$, для расширения полей~\mbox{$k\subset K$}

\medskip

\noindent $X(K)$ --- множество $K$-точек на многообразии $X$, определённом над полем~$k$, для расширения полей~\mbox{$k\subset K$}

\bigskip

\bigskip

\noindent $\Mat_n(A)$ --- алгебра $(n\times n)$-матриц с коэффициентами в ассоциативной алгебре $A$

\medskip

\noindent $\Br(K)$ --- группа Брауэра поля $K$

\medskip

\noindent $\Br(L/K)=\Ker\big(\Br(K)\to\Br(L)\big)$ --- относительная группа Брауэра для расширения полей $K\subset L$

\medskip

\noindent $b(X)\in \Br(K)$ --- класс многообразия Севери--Брауэра $X$, определённого над полем~$K$

\medskip

\noindent $\res\colon\Br(K)\to\Hom(G_{\kappa},\Q/\Z)$ --- отображение вычета для полного поля дискретного нормирования $K$ с совершенным полем вычетов $\kappa$

\medskip

\noindent $\OO_K$ --- кольцо нормирования в поле дискретного нормирования $K$ или кольцо целых в глобальном поле $K$ нулевой характеристики

\medskip

\noindent $K_v$ --- пополнение поля $K$ по дискретному (или, более
общо, мультипликативному) нормированию~$v$

\medskip

\noindent $\res_v\colon \Br(K)\to \Hom(G_{\kappa_v},\Q/\Z)$ --- вычет, соответствующий дискретному нормированию~$v$ поля~$K$, где~$\kappa_v$ обозначает поле вычетов для $K_v$, и $\kappa_v$ совершенно

\medskip

\noindent $\res_D\colon \Br\big(k(X)\big)\to \Hom(G_{k(D)},\Q/\Z)$ ---
вычет, соответствующий дискретному нормированию,
заданному неприводимым приведённым дивизором $D$
на неприводимом (нормальном)
многообразии $X$ над полем нулевой характеристики~$k$

\medskip

\noindent $\Br^{nr}\big(K)$ --- неразветвлённая группа Брауэра конечно порождённого поля $K$ над полем нулевой характеристики $k$

\medskip

\noindent $\Br^{nr}(X)$ --- неразветленная группа Браэура неприводимого (нормального) многообразия $X$ над полем нулевой характеристики $k$

\medskip

\noindent $\Br(X)$ --- группа Брауэра многообразия или схемы $X$

\bigskip

\bigskip

\noindent $\disc(Q)\in K^*/(K^*)^2$ --- дискриминант чётномерной квадрики $Q$ над полем $K$, где~\mbox{${\rm char}(K)\ne 2$}

\medskip

\noindent $\cl(Q)\in \Br(K)$ --- инвариант Клиффорда квадрики $Q$ над $K$ с тривиальным $\disc(Q)$, где ${\rm char}(K)\ne 2$

\medskip

\noindent $R_{K/k}(Y)$ --- ограничение скаляров по Вейлю для многообразия $Y$ над~$K$ и конечного сепарабельного расширения полей $k\subset K$

\medskip

\noindent $K^*=R_{K/k}(\Gb_m)$ --- алгебраический тор над полем $k$
для конечного сепарабельного расширения полей $k\subset K$,
множество $k$-точек которого равно $K\setminus\{0\}$

\medskip

\noindent $T^{\vee}$ ---
двойственный $G_k$-модуль к алгебраическому тору $T$ над полем $k$

\medskip

\noindent $M^{\vee}$ --- двойственный алгебраический тор
над $k$ к $G_k$-модулю $M$,
являющемуся свободной конечно порождённой абелевой группой

\newpage
\section{Когомологии групп}
\label{section:group-cohomology}

\subsection{Определение и общие свойства}
\label{subsection:cohom-generalities}

Для произвольных множеств $X$ и $Y$ через $\Map(X,Y)$ будем обозначать
множество всех отображений из $X$ в $Y$.
Для произвольных групп $G$ и $H$ через $\Hom(G,H)$ будем обозначать
множество всех гомоморфизмов из $G$ в $H$.
В частности, для абелевой группы $A$
множество $\Hom(G,A)$ является абелевой группой.

Для произвольной группы $G$ через~$e$
будем обозначать нейтральный элемент, а через $\Mod_G$~--- категорию
(левых) {\it $G$-модулей}, то есть абелевых групп $M$ с действием группы
$G$, коммутирующим со сложением в~$M$. Результат действия элемента
$g\in G$ на элементе $G$-модуля $m\in M$ будем обозначать~${}^gm$. В
частности, имеется равенство
$${}^g({}^h{m})={}^{gh}m\,.$$
Через $\Hom_G(M,M')$
обозначим абелеву
группу морфизмов в категории $\Mod_G$ между $G$-модулями $M$ и $M'$.
Через $M^G$ обозначим подгруппу в $M$, состоящую из \mbox{$G$-инвариантных}
элементов.
Будем говорить, что $G$-модуль $M$ {\it тривиальный}, если $G$
действует тождественным образом на $M$. По умолчанию группа $\Z$
рассматривается с тривиальной структурой $G$-модуля.

Для произвольного множества $S$ через $\Z[S]$ мы будем обозначать решётку
с базисом, занумерованным элементами из $S$. В частности, если на
$S$ действует группа~$G$, то $\Z[S]$ приобретает структуру $G$-модуля.
Кроме того, $\Z[G]$ имеет структуру кольца, и называется
\emph{групповым кольцом} группы~$G$, а естественное действие
группы~$G$ на себе левыми сдвигами определяет на $\Z[G]$ структуру $G$-модуля.

\bigskip

\begin{prob}{\bf Категория $G$-модулей}
\hspace{0cm}\nopagebreak
\begin{itemize}
\item[(i)] Покажите, что категория $\Mod_G$ является абелевой категорией и эквивалентна категории
левых модулей над групповым кольцом $\Z[G]$.
\item[(ii)] Покажите, что имеется канонический изоморфизм
$$\Hom_G(\Z,M)\cong M^G\,,$$
и взятие $G$-инвариантов является точным
слева функтором.
\end{itemize}
\end{prob}

\bigskip
На протяжении всей книги мы будем активно использовать
когомологии~\mbox{$H^i(G,M)$} группы $G$ с коэффициентами в модуле $M$.
Определение~\ref{defin:group-cohomology} использует понятие функтора
$\Ext$ в абелевой категории.
Читатель, незнакомый с производными функторами, может
считать упражнение~\ref{exer-stcoml}(iv), (v) определением групп
когомологий $H^i(G,M)$.

\begin{defin}\label{defin:group-cohomology}
{\it Группы когомологий $G$ с коэффициентами в $G$-модуле $M$}
задаются по формуле
$$
H^i(G,M)=\Ext_G^i(\Z,M),\quad i\geqslant 0\,,
$$
где $\Ext$ рассматривается в абелевой категории $G$-модулей.
\end{defin}

Таким образом, $H^i(G,M)$ являются правыми производными функторами к
взятию $G$-инвариантов.
В частности, морфизм $G$-модулей
$$f\colon M\to M'$$
определяет
гомоморфизмы групп когомологий
$$f\colon H^i(G,M)\to H^i(G,M')\,,$$
а точная тройка \mbox{$G$-модулей}
$$
0\to M'\to M\to M''\to 0
$$
определяет длинную точную последовательность абелевых групп
\begin{multline*}
0\to H^0(G,M')\to H^0(G,M)\to H^0(G,M'')\to{}\\
{}\to H^1(G,M')\to H^1(G,M)\to H^1(G,M'')\to\ldots
\end{multline*}

\bigskip

\begin{prob}{\bf Стандартный комплекс}\label{exer-stcoml}
\nopagebreak
\hspace{0cm}
\begin{itemize}
\item[(i)] Докажите, что $G$-модуль
$$
\mbox{$\Z[\underbrace{G\times\ldots \times G}_{i}]=\left\{\sum
n_{g_1\ldots g_i}\cdot(g_1,\ldots,g_i)\mid\,
 n_{g_1\ldots g_i}\in\Z\right\}$}
$$
с покоординатным умножением слева на $G$ при любом $i\geqslant 1$
является проективным (и
даже свободным) объектом в $\Mod_G$. (Указание: рассмотрите элементы
вида~\mbox{$(e,g_2,\ldots,g_i)$} в $\Z[G^{\times i}]$.)
\item[(ii)]
Рассмотрим линейные отображения
$$\partial\colon \Z[G^{\times (i+1)}]\to
\Z[G^{\times i}], \quad i\geqslant 0\,,$$
заданные по формуле
$$\partial(g_1,\ldots,g_{i+1})=
\sum\limits_{j=1}^{i+1}(-1)^{j+1}(g_1,\ldots,\hat{g}_j,\ldots,g_{i+1})\,,$$
где символ $\hat{g}$ обозначает пропускание элемента $g$; в
частности, при $i=0$ рассматривается отображение
$$\Z[G]\to \Z,\quad \sum n_g g\mapsto \sum n_g\,.$$
Докажите, что $\partial$ является
морфизмом $G$-модулей и что $\partial\circ \partial=0$. Это
определяет комплекс $G$-модулей $\Z[G^{\udot}]$
(для краткости мы пишем $G^{\udot}$ вместо $G^{\times\udot}$).
\item[(iii)]
Докажите, что комплекс $\Z[G^{\udot}]$ точен, то есть
является проективной резольвентой для тривиального $G$-модуля~$\Z$.
(Указание:
используйте отображения
$$G^i\to G^{i+1}, \quad (g_1,\ldots,g_i)\mapsto (e,g_1,\ldots,g_i)\,,$$
для построения стягивающей гомотопии.)
\item[(iv)]
Для $G$-модуля $M$ рассмотрим комплекс $K^{\udot}$ с членами
$$K^i=\Hom_G\big(\Z[G^{\times{(i+1)}}], M\big), \quad i\geqslant 0\,.$$
Покажите,
что
$$H^i(G,M)\cong H^i(K^{\udot})$$
при $i\geqslant 0$.
\item[(v)]
Вычислив явно члены комплекса $K^{\udot}$ из предыдущего пункта,
докажите, что группа $H^i(G,M)$ канонически изоморфна $i$-ой
группе когомологий комплекса
\begin{equation}\label{eq-stcompl}
0\to M\stackrel{d}\to
\Map(G,M)\stackrel{d}\to\ldots\stackrel{d}\to\Map(G^{\times
i},M)\stackrel{d}\to\ldots,
\end{equation}
где
дифференциал $d$ задаётся по формуле
\begin{multline*}
(d\varphi)(g_1,\ldots,g_{i+1})={}^{g_1}\big(\varphi(g_2,\ldots,g_{i+1})\big)+\\
+\sum_{j=1}^i(-1)^j\varphi(g_1,\ldots,g_jg_{j+1},
\ldots,g_{i+1})+(-1)^{i+1}\varphi(g_1,\ldots,g_i)\,.
\end{multline*}
В частности, в степени нуль дифференциал имеет вид
$$(dm)(g)={}^gm-m, \quad m\in\nolinebreak M\,.$$
\end{itemize}
\end{prob}

\bigskip

\begin{prob}{\bf Группы когомологий $H^0$, $H^1$ и $H^2$}\label{exer-012}
\hspace{0cm}
\begin{itemize}
\item[(o)]
Покажите, что $H^0(G,M)\cong M^G$.
\item[(i)]
Покажите, что $1$-коциклы в комплексе~\eqref{eq-stcompl} являются
{\it скрещёнными гомоморфизмами} из $G$ в $M$, то есть отображениями
$\varphi\colon G\to M$, удовлетворяющими условию
$$\varphi(gh)=\varphi(g)+{}^g\big(\varphi(h)\big)$$
для всех $g,h\in G$. При
этом $1$-коцикл $\varphi$ является кограницей тогда и только тогда,
когда существует такой элемент $m\in M$, что
$$\varphi(g)={}^gm-m$$
для всех
$g\in G$. В частности, для тривиального модуля $M$ первая группа
когомологий $H^1(G,M)$ изоморфна группе гомоморфизмов $\Hom(G,M)$.
\item[(ii)]
Докажите, что группа $H^2(G,M)$ биективна множеству классов
изоморфизма расширений вида
$$
0\to M\to \widetilde{G}\stackrel{\pi}\to G\to 0\,,
$$
где изоморфизм между двумя расширениями --- это изоморфизм между
средними членами соответствующих последовательностей, индуцирующий
тождественные отображения на крайних членах. (Указание:
соответствующий расширению $2$-коцикл задаётся по формуле
$$\varphi(g,h)=\tilde{g}\tilde{h}(\widetilde{gh})^{-1}\,,$$
где
$g\mapsto \tilde{g}$ --- произвольное теоретико-множественное
сечение гомоморфизма~$\pi$.
При этом структура $G$-модуля на $M$ задаётся сопряжением, то есть
$${}^gm=\tilde{g}m\tilde{g}^{-1},\quad m\in M, g\in G\,.$$
В обратную сторону, для коцикла $\varphi$ группа $\widetilde{G}$
определяется как множество~\mbox{$M\times G$} с групповой
операцией
$$(m, g)\cdot (n,h)=\big(m+{}^gn+\varphi(g,h),gh\big), \quad m,n\in M,
g,h\in G\,.$$
Условие коцикла равносильно ассоциативности этой
операции.)
\item[(iii)]
Пусть даны морфизм $G$-модулей $f\colon M\to M'$ и элемент
$\alpha\in H^2(G,M)$, соответствующий расширению
$$
0\to M\stackrel{\iota}\to\widetilde{G}\stackrel{\pi}\to G\to 0\,.
$$
Заметьте, что гомоморфизм $\pi$ задаёт структуру $\widetilde{G}$-модуля
на $M'$, и покажите, что образ вложения
$$(f,\iota)\colon M\to M'\rtimes\widetilde{G}$$
является нормальной подгруппой.
Проверьте, что элемент
$$f(\alpha)\in H^2(G,M')$$
соответствует расширению
$$
0\to M'\to \widetilde{G}'\to G\to 0\,,
$$
где
$$\widetilde{G}'=\big(M'\rtimes\widetilde{G}\big)/(f,\iota)(M)\,.
$$
\item[(iv)]
Докажите, что для простого числа $p$ имеется изоморфизм
$$
H^2(\Z/p\Z,\,\Z/p\Z)\cong\Z/p\Z\,.
$$
(Указание: используйте классификацию групп порядка $p^2$. На самом деле, данное утверждение верно и после замены $p$ на любое целое число, ср. с упражнением~\ref{exer-cyclgrp}(i) ниже.)
\end{itemize}
\end{prob}

\bigskip

\begin{prob}{\bf Когомологии циклических групп}\label{exer-cyclgrp}
\hspace{0cm}
\begin{itemize}
\item[(i)]
Пусть $G$ --- конечная циклическая группа порядка $n$ с образующей
$s$. Докажите, что комплекс
$$
\ldots\to\Z[G]\to\ldots\to\Z[G]\stackrel{s-1}\longrightarrow\Z[G]\stackrel{\N}\longrightarrow\Z[G]
\stackrel{s-1}\longrightarrow\Z[G]\to \Z\to 0
$$
является проективной резольвентой для $\Z$, где $\N$ ---
\emph{отображение нормы}, то есть
$$\N=\sum\limits_{g\in G}g\in \Z[G]\,.$$
Применив функтор $\Hom_G(-,M)$, выведите отсюда,
что группа $H^i(G,M)$ изоморфна~\mbox{$i$-ой} группе когомологий комплекса
$$0\to M\stackrel{s-1}\longrightarrow
M\stackrel{\N}\longrightarrow
M\stackrel{s-1}\to{}\ldots{}\to
M\to\ldots
$$
Таким образом, имеются изоморфизмы:
$$H^i(G,M)\cong M^G/{\N}(M)$$
для чётных $i>0$, и
$$H^i(G,M)\cong \Ker(\N)/(s-1)M$$
для нечётных $i$. В частности, такая формула для $H^i(G,M)$ позволяет
по-другому решить упражнение~\ref{exer-012}(iv).
Проверьте, что эти изоморфизмы функториальны
по $M$.
\item[(ii)]
Пусть $M=\Z$ является тривиальным $G$-модулем. Покажите,
что класс элемента $1\in\Z$ в факторгруппе
$$M^G/\N(M)\cong\Z/n\Z$$
при изоморфизме
$$H^2(G,M)\cong M^G/{\N}(M)$$
соответствует расширению
$$0\longrightarrow\Z\stackrel{n}\longrightarrow
\Z\stackrel{\lambda}\longrightarrow G\to 0\,,$$
где $\lambda(1)=s$. В частности, это расширение
зависит от выбора образующей~$s$.
Точно так же от выбора $s$ зависят изоморфизмы из пункта~(i).
(В какой момент этот выбор использовался при их построении?)
\item[(iii)]
Рассмотрим произвольный $G$-модуль $M$ и элемент $m\in M^G$.
Докажите, что класс элемента~$m$ в факторгруппе
$$M^G/\N(M)\cong H^2(G,M)$$
соответствует (см. пункт~(i)) расширению
$$
0\to M\to \widetilde{G}\stackrel{\pi}\to\Z/n\Z\to 0\,,
$$
где
$$\widetilde{G}=\big(M\rtimes\Z\big)/\langle (m,n)\rangle\,,$$
и $\Z$ действует на $M$ через гомоморфизм $\lambda$ из пункта~(ii).
(Указание: сначала рассмотрите случай тривиального модуля
$M=\Z$ и элемента~\mbox{$m=1\in\Z$}, и воспользуйтесь пунктом~(ii).
Для произвольных $M$ и $m$ рассмотрите морфизм $G$-модулей
$$\Z\to M,\quad l\mapsto l\cdot m\,,$$
и воспользуйтесь
упражнением~\ref{exer-012}(iii).)
\item[(iv)]
В обозначениях из пункта~(iii) покажите, что имеется
изоморфизм групп
$$\widetilde{G}\cong M\times\{1,\sigma,\ldots,\sigma^{n-1}\}\,,$$
где групповой закон на правой части
определяется соотношениями $\sigma^n=m$ и
$$\sigma\cdot m'={}^s(m')\cdot \sigma,\quad m'\in M\,.$$
При этом гомоморфизм
$$\pi\colon\widetilde{G}\to\Z/n\Z\cong G$$
из пункта~(iii) переводит $\sigma$ в $s^{-1}$.
\end{itemize}
\end{prob}

\bigskip

\begin{prob}{\bf Кограничные отображения}\label{exer-cobound-map}

Рассмотрим точную тройку $G$-модулей
$$
0\to M'\to M\stackrel{f}\to M''\to 0
$$
и соответствующую длинную точную последовательность когомологий.
Обозначим через
$$\delta\colon H^i(G,M'')\to H^{i+1}(G,M')$$
кограничные отображения в ней.
\begin{itemize}
\item[(i)]
Покажите, что для элемента
$$m''\in (M'')^G\cong H^0(G,M'')$$
имеет место равенство
$$(\delta m'')(g)=-m+{}^gm\,,$$
где $m\in M$ --- любой прообраз элемента
$m''$.
\item[(ii)]
Покажите, что для $1$-коцикла $\varphi\colon G\to M''$ имеет место равенство
$$
(\delta\varphi)(g,h)=\tilde\varphi(g)+{}^{g}\big(\tilde\varphi(h)\big)-
\tilde\varphi(gh)\,,
$$
где $\tilde\varphi(g)\in M$ --- любой прообраз элемента $\varphi(g)$
для каждого $g\in G$.
\item[(iii)]
Предположим, что $M'$, $M$ и $M''$ --- тривиальные $G$-модули.
Рассмотрим гомоморфизм $\varphi\colon G\to M''$. Покажите,
что элемент $\delta\varphi\in H^2(G, M')$ соответствует расширению
(см. упражнение~\ref{exer-012})
$$0\to M'\to \widetilde{G}\to G\to 1\,,$$
где
$$\widetilde{G}=M\times_{M''} G=\left\{(m,g)\in M\times G\mid
f(m)=\varphi(g)\right\}\,.$$
\item[(iv)]
Предположим, что $G$ --- конечная циклическая группа порядка $n$ с
образующей~$s$.
Покажите, что при отождествлениях из упражнения~\ref{exer-cyclgrp}(i)
кограничные отображения соответствуют отображениям
\begin{align*}
&\Ker(\N\colon M''\to M'')/(s-1)M''\longrightarrow (M')^G/\N(M'),\quad
m''\mapsto \N(m)\,,\\
&(M'')^G/\N(M'')\longrightarrow
\Ker(\N\colon M'\to M')/(s-1)M',\quad
m''\mapsto (s-1)m\,,
\end{align*}
где $m\in M$ --- произвольный элемент, для которого $f(m)=m''$.
(Указание: воспользуйтесь резольвентой из упражнения~\ref{exer-cyclgrp}(i).)
\end{itemize}
\end{prob}

\bigskip

\begin{prob}{\bf Когомологии с неабелевыми
коэффициентами}\label{exer-nonabcohom}

Пусть $G$ --- группа, действующая автоморфизмами на (возможно,
неабелевой) группе $\Gamma$. Пусть множество
$Z^1(G,\Gamma)$ состоит из отображений
множеств
$\varphi\colon G\to
\Gamma$, для которых равенство
$$\varphi(gh)=\varphi(g)\cdot\,{}^{g}\big(\varphi(h)\big)$$
выполняется при всех
$g,h\in G$.
\begin{itemize}
\item[(o)]
Покажите, что для любого $\varphi\in Z^1(G,\Gamma)$
имеется равенство $\varphi(e)=e$.
\item[(i)]
Докажите, что формула
$$
({}^{\gamma}\varphi)(g)=
\gamma\cdot\varphi(g)\cdot{}\,{}^g\big(\gamma^{-1}\big),\quad\gamma\in\Gamma,
g\in G\,,
$$
определяет действие группы $\Gamma$ на множестве $Z^1(G,\Gamma)$.
Положим
$$
H^0(G,\Gamma)=\Gamma^G,\quad H^1(G,\Gamma)=Z^1(G,\Gamma)/\Gamma\,.
$$
Заметим, что $H^1(G,\Gamma)$ является множеством с отмеченным
элементом, задаваемым отображением $\varphi(g)\equiv
e$. Покажите, что если $G$ действует тождественно
на $\Gamma$, то множество $H^1(G,\Gamma)$ канонически биективно
фактору множества~\mbox{$\Hom(G,\Gamma)$}
по действию группы $\Gamma$ сопряжениями.
\item[(ii)]
Пусть дана точная тройка (возможно, неабелевых) $G$-модулей
$$
1\to A\to \widetilde{\Gamma}\to \Gamma\to 1\,.
$$
Покажите, что формула
$$(\delta\gamma)(g)=\tilde\gamma^{-1}\cdot\,{}^g\tilde\gamma$$
определяет отображение множеств
$$\delta\colon H^0(G,\Gamma)\to H^1(G,A)\,,$$
где
$\gamma\in \Gamma^G\cong H^0(G,\Gamma)$, а
$\tilde{\gamma}\in\widetilde{\Gamma}$ --- любой прообраз элемента
$\gamma$.
Докажите, что имеется точная
последовательность множеств с отмеченным элементом
\begin{multline*}
1\to H^0(G,A)\to H^0(G,\widetilde{\Gamma})\to
H^0(G,\Gamma)\stackrel{\delta}\to{}\\
{}\stackrel{\delta}\to H^1(G,A)\to H^1(G,\widetilde{\Gamma})\to
H^1(G,\Gamma)\,,
\end{multline*}
то есть в каждом члене образ входящего отображения совпадает с прообразом
отмеченного элемента относительно исходящего отображения.
\item[(iii)]
В обозначениях пункта (ii) предположим,
что $A$ является центральной (в частности, абелевой) подгруппой
в $\widetilde{\Gamma}$.
Покажите, что
$$\delta\colon H^0(G,\Gamma)\to H^1(G,A)$$
является гомоморфизмом групп. Покажите, что формула
$$
(\delta\varphi)(g,h)=\tilde\varphi(g)\cdot\,{}^{g}\big(\tilde\varphi(h)\big)\cdot\tilde\varphi(gh)^{-1}
$$
определяет отображение множеств с отмеченным элементом
$$\delta\colon H^1(G,\Gamma)\to H^2(G,A)\,,$$
где $\varphi\in Z^1(G,\Gamma)$,
а $\tilde\varphi(g)\in\widetilde{\Gamma}$ --- любой прообраз
элемента~$\varphi(g)$, $g\in G$.
Докажите, что имеется точная
последовательность множеств с отмеченным элементом
\begin{multline*}
1\to H^0(G,A)\to H^0(G,\widetilde{\Gamma})\to
H^0(G,\Gamma)\stackrel{\delta}\to{}\\
{}\stackrel{\delta}\to H^1(G,A)\to H^1(G,\widetilde{\Gamma})\to
H^1(G,\Gamma)\stackrel{\delta}\to H^2(G,A)\,.
\end{multline*}
\item[(iv)]
В обозначениях пункта (iii) докажите, что абелева группа $H^1(G,A)$
действует на множестве $H^1(G,\widetilde{\Gamma})$ поточечным
умножением $1$-коциклов, и ядро этого действия равно образу
$\delta\big(H^0(G,\Gamma)\big)$ группы $H^0(G,\Gamma)$ при кограничном
отображении~$\delta$. Кроме того,
естественное отображение
$$H^1(G, \widetilde{\Gamma})\to H^1(G, \Gamma)$$
индуцирует вложение
$$H^1(G, \widetilde{\Gamma})/H^1(G,A)\hookrightarrow H^1(G, \Gamma)\,.$$
\item[(v)]
Пусть $G$ --- конечная циклическая группа порядка $n$ с образующей
$s$. Покажите, что множество $Z^1(G,\Gamma)$
биективно множеству $Z(G,\Gamma)$,
состоящим из таких элементов $\alpha\in \Gamma$, что
$$
\alpha\cdot\,{}^s\alpha\cdot\,{}^{s^2}\alpha\cdot\ldots\cdot{}^{s^{n-1}}\alpha
=e\,.
$$
При этом определённое в пункте (i) действие группы $\Gamma$ на $Z^1(G,\Gamma)$
соответствует действию $\Gamma$ на
$Z(G,\Gamma)$, заданному по формуле
$${}^{\gamma}\alpha=\gamma\cdot\alpha\cdot{}\,{}^s\big(\gamma^{-1}\big)\,.$$
(Указание: сопоставьте коциклу $\varphi$
элемент $\alpha=\varphi(s)\in\Gamma$.) Таким образом,
и для неабелева $G$-модуля $\Gamma$ первые когомологии
можно вычислять так, как описано в упражнении~\ref{exer-cyclgrp}.
\end{itemize}
\end{prob}

Имеется следующее продолжение упражнения~\ref{exer-nonabcohom}(iii). Если группа $\Gamma$ абелева, то $\delta\colon H^1(G,\Gamma)\to H^2(G,A)$ является отображением между абелевыми группами. Однако в общем случае $\delta$ не является гомоморфизмом абелевых групп. Более точно, для любых $\varphi,\varphi'\in H^1(G,\Gamma)$ выполняется равенство
$$
\delta(\varphi\varphi')=\delta(\varphi)+\delta(\varphi')+c(\varphi,\varphi')
$$
в $H^2(G,A)$, где отображение
$$
c\colon H^1(G,\Gamma)\times H^1(G,\Gamma)\to H^2(G,A)
$$
индуцировано коммутаторным спариванием~\mbox{${\Gamma\times\Gamma\to A}$}. Подробности можно найти в работе~\cite{Zarhin}.

\bigskip

\begin{prob}{\bf Ацикличная резольвента}\label{exer-acycl}
\hspace{0cm}
\begin{itemize}
\item[(i)]
Для абелевой группы $A$ положим $A_+=\Map(G,A)$. Определим на
$A_+$ структуру $G$-модуля по формуле
$$({}^g\varphi)(g')=\varphi(g'g),\quad \varphi\in A_+, g,g'\in G$$
(ср. с упражнением~\ref{exer-ind}(o)). Докажите, что
$$H^{i}(G,A_+)=0$$
для всех $i>0$. (Указание: посчитайте
$H^{i}(G,A_+)$ при помощи стандартного
комплекса из упражнения~\ref{exer-stcoml}(iv),
используя при этом существование изоморфизма
$$\Hom_G(M,A_+)\cong\Hom(M,A)$$
для любого $G$-модуля~$M$. Для этого вычисления надо также использовать упражнение~\ref{exer-stcoml}(iii) и то, что
все абелевы группы $\Z[G^{\times i}]$ свободны.)
\item[(ii)]
Для любого $G$-модуля $M$ рассмотрим вложение $G$-модулей
$M\hookrightarrow M_+$, заданное по формуле
$$m\mapsto\nolinebreak\{g\mapsto {}^g m\}\,.$$
Положим
$$M^0=M_+, \quad M^1=(M_+/M)_+,\quad M^2=\big(M^1/(M_+/M)\big)_+\,,$$
и продолжим индуктивно этот процесс дальше. Таким образом, определим
\mbox{$G$-модули}
$M^i=(C^{i-1})_+$\,, $i\geqslant 2$, где
$$C^{i-1}={\rm Coker}(M^{i-2}\to M^{i-1})\,,$$
а отображение
$$M^{i-1}\to M^i$$
определяется как композиция сюръекции $M^{i-1}\twoheadrightarrow C^{i-1}$ и
естественного вложения
$$C^{i-1}\hookrightarrow (C^{i-1})_+=M_i\,.$$
Покажите, что это
определяет комплекс $G$-модулей $M^{\udot}$, причём
$H^0(M^{\udot})\cong M$ и $H^{i}(M^{\udot})=0$ при $i>0$.
Такой комплекс называют
\emph{ацикличной резольвентой} модуля $M$, поскольку
$H^i(G,M^j)=0$ при $i>0$ и $j\geqslant 0$.
\item[(iii)]
Докажите, что
$$H^i(G,M)\cong H^i\big((M^{\udot})^G\big)$$
при $i\geqslant 0$.
(Указание: используйте ацикличность резольвенты $M^{\udot}$
модуля~$M$.)
\item[(iv)]
Покажите, что точная последовательность $G$-модулей
$$
0\to M'\to M\to M''\to 0
$$
индуцирует точную последовательность ацикличных резольвент
$$
0\to M'^{\udot}\to M^{\udot}\to M''^{\udot}\to 0
$$
\end{itemize}
\end{prob}

\subsection{Поведение при замене группы}

\begin{defin}
Пусть $f\colon H\to G$ --- гомоморфизм групп. Тогда $f$ определяет функтор
$$f^*\colon\Mod_G\to\Mod_H\,,$$
который сопоставляет $G$-модулю $M$ ту же
абелеву группу $M$ с действием $H$, пропускающимся через $f$. Обычно
$f^*M$ обозначается просто через $M$, если ясно, что рассматривается
действие группы $H$, а не $G$.
\end{defin}

\begin{prob}{\bf Обратный образ на
когомологиях}\label{exer-invimage}
\hspace{0cm}
\begin{itemize}
\item[(i)]
Проверьте, что $f^*$ является точным функтором.
\item[(ii)]
Покажите, что функтор  $f^*$ определяет отображения
$$f^*\colon H^i(G,M)\to
H^i(H,M),\quad i\geqslant 0$$
для любого $G$-модуля $M$. (Указание:
используйте точность функтора $f^*$.)
\item[(iii)]
Рассмотрим естественные отображения
$$f^*\colon\Map(G^{\times i},M)\to\Map(H^{\times i},M)\,.$$
Докажите, что они задают морфизм
комплексов (см. упражнение~\ref{exer-stcoml}(v))
$$f^*\colon\Map(G^{\times \udot},M)\to\Map(H^{\times
\udot},M)\,,$$
а соответствующие отображения на когомологиях
совпадают с отображениями, определёнными в предыдущем пункте.
\item[(iv)] Пусть $G$ является конечной циклической группой
с образующей $s$, а $f\colon H\hookrightarrow G$ является вложением подгруппы
индекса $r$.
Зафиксируем образующую
$$t=s^r$$
группы~$H$. Пусть $\N_G\colon M\to M$
и $\N_H\colon M\to M$
обозначают отображения нормы,
соответствующие группам $G$ и $H$
(см. упражнение~\ref{exer-cyclgrp}(i)).
Покажите, что гомоморфизм
$$f^*\colon H^i(G,M)\to H^i(H,M)$$
при отождествлениях
из упражнения~\ref{exer-cyclgrp}(i) относительно образующих~$s$ и~$t$
при чётном~$i$ соответствует отображению
$$M^G/\N_G(M)\longrightarrow M^H/\N_H(M)\,,$$
индуцированному тождественным отображением из $M$ в себя.
Покажите, что при нечётном $i$ гомоморфизм $f^*$ соответствует
отображению
$$
\Ker(\N_G)/(s-1)M\longrightarrow \Ker(\N_H)/(t-1)M\,,
$$
индуцированному отображением
$$
M\to M,\quad m\mapsto\sum\limits_{l=0}^{r-1}s^l(m)\,.
$$
(Указание: постройте морфизм между комплексами из
упражнения~\ref{exer-cyclgrp}(i) для групп $G$ и $H$, продолжающий
тождественное отображение из $M$ в себя в нулевом члене.)
\end{itemize}
\end{prob}

\bigskip

\begin{defin}\label{defin-coind}
Пусть $f\colon H\to G$ --- гомоморфизм групп. Тогда $f$ определяет функтор
$$
f_*\colon\Mod_H\to\Mod_G\,,\quad N\mapsto \Hom_H(\Z[G],N)\,,
$$
называемый {\it коиндуцированием}, где $N$ является $H$-модулем,
действие $H$ на $\Z[G]$ задаётся умножением слева, а действие $G$ на
$\Hom_H(\Z[G],N)$ задаётся по формуле
$$\big({}^g\varphi\big)(g')=\varphi(g'g),\quad \varphi\in\Hom_H(\Z[G],N),
\quad g,g'\in G\,.$$
\end{defin}

Заметим, что в обозначениях определения~\ref{defin-coind}
имеется канонический изоморфизм $G$-модулей
$$\Hom_H(\Z[G], N)\cong \Map_H(G,N)\,,$$
где $\Map_H$ обозначает множество $H$-эквивариантных отображений
между множествами с действием группы $H$.

\bigskip
\bigskip
\bigskip
\bigskip

\begin{prob}{\bf Коиндуцирование}\label{exer-ind}
\hspace{0cm}
\begin{itemize}
\item[(o)]
Покажите, что $\iota_*A\cong A_+$ для любой абелевой группы $A$, см.~упражнение~\ref{exer-acycl}, где
$\iota\colon\{e\}\hookrightarrow G$ является вложением тривиальной подгруппы.
\item[(i)]
Проверьте, что определение~\ref{defin-coind} корректно, то есть на
абелевой группе
$$f_*N=\Hom_H(\Z[G],N)$$
корректно определено действие группы $G$.
\item[(ii)] Покажите, что функтор $f_*$ точен слева.
(Указание: воспользуйтесь точностью слева функторов $\Hom(\Z[G], -)$ и взятия
$H$-инвариантов.)
\item[(iii)]
Докажите, что $(f^*,f_*)$ является парой сопряжённых функторов, то есть
для любого $G$-модуля $M$ и $H$-модуля $N$ имеется
канонический изоморфизм
$$
\Hom_H(M,N)\cong\Hom_G(M,f_*N)\,.
$$
(Указание: действуйте по аналогии с подобным утверждением о модулях
над коммутативными кольцами. В частности, отображению $\psi\in
\Hom_H(M,N)$ сопоставьте отображение
$$m\mapsto\{g\mapsto {}^g(\psi(m))\}\,,$$
а отображению $\xi\in \Hom_G(M,f_*N)$ сопоставьте
отображение
$$m\mapsto \xi(m)(e)\,,$$
где $e\in \Z[G]$ --- единица
группы и группового кольца.)
\item[(iv)]
Покажите, что для любой композиции гомоморфизмов групп $f\circ g$
имеется канонический изоморфизм функторов
$$f_*\circ g_*\cong (f\circ g)_*\,.$$
\end{itemize}
\end{prob}

\bigskip

\begin{prob}{\bf Коиндуцирование для инъективного гомоморфизма}
\label{exer-ind-inj}

Пусть $f\colon H\hookrightarrow G$ --- инъективный гомоморфизм групп, и $N$ ---
произвольный $H$-модуль.

\begin{itemize}
\item[(i)]
Заметьте, что существует (неканонический) изоморфизм абелевых групп
$$\mbox{$f_*N\cong \prod\limits_{H\setminus G}N\,,$}$$
и выведите отсюда, что функтор $f_*$ точен.
\item[(ii)] Докажите, что любой проективный $G$-модуль $M$
также проективен как $H$-модуль. (Указание: воспользуйтесь
пунктом (i), а также упражнением~\ref{exer-ind}(iii).)
\item[(iii)]
Докажите {\it лемму Шапиро}: для любого
$H$-модуля $N$ имеется канонический изоморфизм
$$
H^i(H,N)\cong H^i(G,f_*N),\quad i\geqslant 0\,.
$$
(Указание: воспользуйтесь
пунктом (ii), а также тем, что когомологии групп можно вычислять
при помощи любой проективной резольвенты для модуля~$\Z$.
После этого требуемый изоморфизм возникает из упражнения~\ref{exer-ind}(iii).
Другой способ: используйте упражнение~\ref{exer-acycl}(iii),
изоморфизм
ацикличных резольвент
$$f_*(N^{\udot})\cong (f_*N)^{\udot}\,,$$
вытекающий из упражнения~\ref{exer-ind}(o), (iv),
а также изоморфизм
$$(f_*N)^G\cong N^H$$
из упражнения~\ref{exer-ind}(iii).)
\item[(iv)]
Предположим, что подгруппа $H$ нормальна в $G$,
и выберем пред\-ста\-ви\-те\-лей~\mbox{$\{g_i\}$}
из правых смежных классов $H\setminus G$. Возникают
автоморфизмы
$$g_i\colon H\to H, \quad h\mapsto g_ihg_i^{-1}\,.$$
Сопоставив элементу $\varphi\in\Hom_H(\Z[G], N)$
набор
$$\mbox{$\big\{\varphi(g_i)\big\}\in \prod\limits_i g_i^*N\,,$}$$
проверьте, что имеется изоморфизм $H$-модулей
$$\mbox{$f^*f_*N\cong \prod\limits_i g_i^*N\,.$}$$
\item[(v)]
Покажите, что
$$H^i(H, N_+)=0$$
при $i>0$, где
$$N_+=\Map(G,N)\,.$$
(Указание: воспользуйтесь пунктом (iv), упражнением~\ref{exer-ind}(o), (iv),
а также упражнением~\ref{exer-acycl}(i).)
\end{itemize}
\end{prob}

\bigskip

Читатель, незнакомый со спектральными последовательностями,
может пропустить все пункты следующего упражнения, кроме пункта~(i).
В дальнейшем их результаты практически не будут использоваться.

\begin{prob}{\bf Коиндуцирование для сюръективного гомоморфизма}
\label{exer-ind-sur}

Пусть $f\colon H\twoheadrightarrow G$ --- сюръективный гомоморфизм групп
с ядром $I$, и $N$ --- произвольный $H$-модуль.

\begin{itemize}
\item[(i)]
Покажите, что имеется изоморфизм $G$-модулей
$$f_*N\cong N^I\,.$$
\item[(ii)]
Покажите, что
$$H^i\big(G,(N_+)^I\big)=0$$
при $i>0$, где
$$N_+=\Map(H,N)\,.$$
(Указание:
воспользуйтесь пунктом (i), и затем, как в упражнении~\ref{exer-ind-inj}(v),
примените упражнение~\ref{exer-ind}(o), (iv) и упражнение~\ref{exer-acycl}(i).)
\item[(iii)]
Покажите, что на группах $H^j(I,N)$ есть каноническая структура $G$-модулей,
и докажите, что имеется спектральная последовательность со вторым членом
\begin{equation}\label{eq:spectral-sequence-coind-sur}
E^{ij}_2=H^i\big(G,H^j(I,N)\big)\,,
\end{equation}
сходящаяся к $H^{i+j}(H,N)$.
(Указание: сначала рассмотрите ацикличную резольвенту
$N^{\udot}$ для $N$ относительно группы $H$ (упражнение~\ref{exer-acycl}
при~\mbox{$G=H$} и~\mbox{$M=N$}) и комплекс $G$-модулей
$(N^{\udot})^I$.
По упражнению~\ref{exer-ind-inj}(v)
его когомологии изоморфны $H^j(I,N)$, что задаёт действие $G$ на $H^j(I,N)$.
Затем рассмотрите бикомплекс
$$C^{\udot,\udot}=\big((N^{\udot})^I\big)^{\udot}\,,$$
полученный при
помощи применения ацикличной резольвенты для членов комплекса
$(N^{\udot})^I$ относительно группы $G$.
Далее используйте две спектральные
последовательности, связанные с
бикомплексом $\big(C^{\udot,\udot}\big)^{G}$.
По пункту~(ii)
одна из этих спектральных последовательностей выродится к одной строке,
равной~$(N^{\udot})^G$,
а другая спектральная последовательность
будет иметь вид~\eqref{eq:spectral-sequence-coind-sur}
по упражнениям~\ref{exer-ind-inj}(v) и~\ref{exer-acycl}(iv).
Отметим,
что это частный случай спектральной последовательности,
связанной с композицией производных функторов,
называемой иногда спектральной последовательностью Лере.)
\item[(iv)]
Покажите, что имеется точная последовательность
$$
0\to H^1(G,N^I)\to H^1(H,N)\to H^1(I,N)^G\to H^2(G,N^I)\to H^2(H,N)\,.
$$
(Указание: используйте спектральную последовательность из
предыдущего пункта.)
\end{itemize}
\end{prob}

\bigskip

\begin{prob}{\bf Прямой образ на когомологиях}\label{exer-trace}

Пусть $f\colon H\hookrightarrow G$ --- такое вложение групп, что $H$ имеет
конечный индекс в $G$.
\begin{itemize}
\item[(i)]
Докажите, что $(f_*,f^*)$
является парой сопряжённых функторов, то есть для любого $G$-модуля $M$
и $H$-модуля $N$ имеется канонический изоморфизм
$$
\Hom_G(f_*N,M)\cong\Hom_H(N,M)\,.
$$
(Указание: отображению $\psi\in \Hom_H(N,M)$ сопоставьте отображение
$$\varphi\mapsto\sum\limits_{[g]\in
G/H}{}^g\big(\psi(\varphi(g^{-1}))\big),$$
где сумма берётся по
представителям в классах смежности.
Отображению~\mbox{$\xi\in\Hom_G(f_*N,M)$} сопоставьте отображение
$$n\mapsto\xi(\chi_n)\,,$$
где $\chi_n\in\Hom_H(\Z[G],N)$ определяется по формуле
$$\chi_n(h)={}^h n$$
для $h\in H$, и $\chi_n(g)=0$ для $g\notin H$. Воспользуйтесь тем,
что для любого гомоморфизма $\varphi\in\Hom_H(\Z[G],N)$ имеется равенство
$$\varphi=\sum\limits_{[g]\in G/H}{}^g\big(\chi_{\varphi(g^{-1})}\big)\,,$$
где сумма берётся по представителям в классах смежности.)
Таким образом, в предположениях этого упражнения функтор $f^*$ является и левым, и правым сопряжённым
для функтора $f_*$, ср. с упражнением~\ref{exer-ind}(iii).
\item[(ii)]
Пользуясь пунктом~(i), постройте для любого  $G$-модуля $M$ канонический
морфизм $G$-модулей
$$f_*f^*M\to M$$
и определите прямой образ на
когомологиях
$$f_*\colon H^i(H,M)\to H^i(G,M)\,,$$
используя лемму Шапиро
(упражнение~\ref{exer-ind-inj}(iii)).
\item[(iii)]
Докажите, что композиция
$$f_*\circ f^*\colon H^i(G,M)\to H^i(G,M)$$
равна
умножению на число классов смежности $|G/H|$.
(Указание: сделайте
явное вычисление для $i=\nolinebreak 0$,
а затем выведите из этого общий случай
при помощи индукции по $i$ и используя длинную точную
последовательность когомологий, связанную с точной последовательностью
(см. упражнение~\ref{exer-acycl})
$$
0\to M\to M_+\to M_+/M\to 0\,.
$$
При этом проверьте согласованность $f_*$ и $f^*$ с кограничными
отображениями.)
\end{itemize}
\end{prob}

\bigskip
Отметим, что иногда используется следующая терминология.
Если дано вложение групп
$f\colon H\hookrightarrow G$, то отображение обратного образа
$$f^*\colon H^i(G,M)\to H^i(H,M)$$
обозначается символом $\rm Res$ и
называется {\it ограничением}. Если к тому же индекс~$H$ в~$G$
конечен, то отображение прямого образа
$$f_*\colon H^i(H,M)\to H^i(G,M)$$
обозначается символом $\rm Cor$ и называется {\it коограничением}. В этом
случае
$${\rm Cor}\circ{\rm Res}=[G:H]$$
по упражнению~\ref{exer-trace}(iii). Если гомоморфизм $f\colon H\to G$
сюръективен и имеет ядро~$I$, то композиция отображения обратного образа
$$f^*\colon H^i(G,M^I)\to H^i(H,M^I)$$
с естественным отображением
$H^i(H,M^I)\to H^i(H,M)$ называется {\it инфляцией} и обозначается
символом $\rm Inf$. В частности, имеется точная последовательность
$$
0\to H^1(G,N^I)\stackrel{\rm Inf}\longrightarrow
H^1(H,N)\stackrel{\rm Res}\longrightarrow H^1(I,N)
$$
(см. упражнение~\ref{exer-ind-sur}(iv)).

\subsection{Когомологии конечных групп}

\begin{prob}{\bf Конечность групп когомологий}\label{exer-cohomfingrp}

Пусть $G$ --- конечная группа.

\begin{itemize}
\item[(i)]
Докажите, что для любого $G$-модуля $M$
когомологии $H^i(G,M)$ аннулируются умножением на порядок группы
$|G|$ при $i>0$. (Указание: используйте
упражнение~\ref{exer-trace}(iii) для отображения
$\{e\}\hookrightarrow G$.)
\item[(ii)]
Докажите, что для любого конечно порождённого
$G$-модуля $M$ когомологии~\mbox{$H^i(G,M)$} конечны при $i>0$.
(Указание:
используя стандартный комплекс, покажите сначала, что группы
$H^i(G,M)$ конечно порождены, а затем используйте предыдущий пункт).
\end{itemize}
\end{prob}

\bigskip
\begin{prob}{\bf Ограничение на силовские подгруппы}\label{exer-Sylow}

Пусть $G$ --- конечная группа, $M$ --- конечно порождённый $G$-модуль,
и $i>0$.
Для каждого простого $p$ пусть
$G_p$ обозначает силовскую $p$-подгруппу в $G$ (в частности, если $p$ не
делит порядок группы $G$, то подгруппа $G_p$ тривиальна).

\begin{itemize}
\item[(i)]
Пусть $p$ --- простое число. Пусть $H_p$
обозначает наибольшую $p$-примарную подгруппу в группе $H^i(G,M)$ (это определение
корректно, так как группа~${H^i(G,M)}$ является конечной абелевой группой по
упражнению~\ref{exer-cohomfingrp}(ii)).
Покажите, что отображение ограничения
$$H_p\to H^i(G_p, M)$$
инъективно. (Указание:
рассмотрите композицию
$$H_p\to H^i(G_p,M)\to H^i(G,M)$$
и воспользуйтесь упражнением~\ref{exer-trace}(iii) для вложения
$G_p\hookrightarrow G$.)
\item[(ii)]
Предположим, что для любого простого числа $p$ и для любой силовской подгруппы $G_p\subset G$
имеется равенство $H^i(G_p, M)=0$.
Докажите, что $H^i(G,M)=0$.
(Указание: воспользуйтесь пунктом~(i).)
\item[(iii)] Верно ли утверждение, обратное к пункту~(ii)?
(Указание: рассмотрите группу~\mbox{$G=\mathrm{S}_3$} и
тривиальный $G$-модуль $M=\Z/3\Z$.)
\end{itemize}
\end{prob}

\bigskip

\begin{prob}{\bf Когомологии с коэффициентами в $\Q/\Z$}
\label{exercise:division-module}

Пусть $G$ --- конечная группа.
\begin{itemize}
\item[(i)]
Пусть $M$ --- произвольный $G$-модуль, являющийся также векторным пространством над $\Q$. Докажите,
что
$$H^i(G,M)=0$$
при всех $i>0$. (Указание: используйте
упражнение~\ref{exer-cohomfingrp}(i).)
\item[(ii)]
Докажите, что
$H^1(G,\Z)=0$, и имеется канонический изоморфизм
$$H^2(G,\Z)\cong\Hom(G,\Q/\Z)\,,$$
где $\Z$ и $\Q/\Z$ рассматриваются со
структурой тривиальных $G$-модулей. (Указание: воспользуйтесь точной
последовательностью $G$-модулей
$$0\to\Z\to\Q\to\Q/\Z\to 0$$
и пунктом~(i).)
\item[(iii)]
Пусть $G$ --- циклическая группа порядка $n$ с образующей $s$.
Пусть
$$\theta\colon\Z/n\Z\stackrel{{}\sim{}}\longrightarrow\Hom(G,\Q/\Z)$$
является
композицией изоморфизмов из пункта~(ii) и упражнения~\ref{exer-cyclgrp}(i).
Покажите, что $\theta$ переводит элемент $1\in\Z/n\Z$ в
гомоморфизм
$$G\to\Q/\Z,\quad s\mapsto [1/n]\,.$$
(Указание: воспользуйтесь упражнениями~\ref{exer-cyclgrp}(ii)
и~\ref{exer-cobound-map}(iii).)
\end{itemize}
\end{prob}

\subsection{Пермутационные и стабильно пермутационные модули}

\begin{defin}\label{defin:perm-module}
Пусть $G$ --- произвольная группа. Рассмотрим $G$-модуль $M$, свободный как
$\Z$-модуль. Говорят, что $M$ является {\it пермутационным}
$G$-модулем, если в $M$ существует базис над $\Z$, на котором
группа $G$ действует перестановками.
\end{defin}

\bigskip
\begin{prob}{\bf
Когомологии с коэффициентами в пермутационном модуле}
\label{exercise:permutation-module} \hspace{0cm}

Пусть $G$ --- конечная группа,
$M$ --- пермутационный $G$-модуль,
$B$ --- базис в $M$, на котором
группа $G$ действует перестановками,
и $\mathbb{O}$ --- множество $G$-орбит в $B$.
\begin{itemize}
\item[(i)] Покажите, что $M$ имеет вид
$$M\cong
\mbox{$\bigoplus\limits_{O\in\mathbb{O}}
\Z\big[G/\mathrm{Stab}_G(x_O)\big]$}\,,$$
где $x_O$ ---
(какой-нибудь) элемент в орбите $O$, подгруппа $\mathrm{Stab}_{G}(x_O)$ ---
стабилизатор элемента $x_O$ в группе $G$,
а $G/\mathrm{Stab}_G(x_O)$ обозначает множество классов смежности
по этой подгруппе.
\item[(ii)] Докажите, что $H^1(G,M)=0$ и
$$
H^2(G,M)\cong\mbox{$\bigoplus\limits_{O\in
\mathbb{O}}\Hom\big(\mathrm{Stab}_{G}(x_O),\Q/\Z\big)$}\,.
$$
(Указание:
воспользуйтесь пунктом~(i) и тем, что для вложения произвольной подгруппы
$i\colon H\hookrightarrow G$ имеется изоморфизм $G$-модулей
$$i_*\Z\cong\Z[G/H]\,.$$
После этого примените лемму Шапиро (см. упражнение~\ref{exer-ind-inj}(iii))
и упражнение~\ref{exercise:division-module}(ii).)
\end{itemize}
\end{prob}

\bigskip
\begin{defin}
\label{defin:stab-perm-module}
Пусть $G$ --- произвольная группа, а $M$ --- некоторый \mbox{$G$-модуль}.
Говорят, что $M$ является {\it стабильно пермутационным}
$G$-модулем, если для некоторых пермутационных
$G$-модулей $N_1$ и $N_2$
имеется изоморфизм $G$-модулей
$$M\oplus N_1\cong N_2\,.$$
\end{defin}

\bigskip
\begin{prob}{\bf Ограничение стабильно пермутационного модуля}
\label{exercise:stab-permutation-module-restriction}

Пусть $G$ --- произвольная группа, $M$ --- (стабильно) пермутационный
$G$-модуль.
Пусть $f\colon H\to G$ --- гомоморфизм групп.
Покажите, что $f^*M$ является (стабильно) пермутационным $H$-модулем.
В частности, ограничение
(стабильно) пермутационного модуля на подгруппу снова является (стабильно)
пермутационным модулем.
\end{prob}

\bigskip
\begin{prob}{\bf Когомологии с коэффициентами в стабильно пермутационном модуле}
\label{exercise:stab-permutation-module}

Пусть $G$ --- конечная группа, $M$ --- стабильно пермутационный
$G$-модуль.

\begin{itemize}
\item[(i)]
Докажите, что
$$H^1(G,M)=0\,.$$
(Указание: воспользуйтесь упражнением~\ref{exercise:permutation-module}(ii).)
\item[(ii)] Используя пункт (i) и
упражнение~\ref{exercise:stab-permutation-module-restriction}, покажите,
что модуль $M$ удовлетворяет следующему условию:
\begin{equation}\label{eq:Endo-Miyata}
\mbox{\ для любой подгруппы $H\subset G$
выполнено равенство $H^1(H,M)=0$.}
\end{equation}
\item[(iii)] Пусть теперь $N$ --- произвольный $G$-модуль,
который при ограничении на каждую силовскую подгруппу в $G$
становится стабильно пермутационным.
Докажите, что $N$ удовлетворяет условию~\eqref{eq:Endo-Miyata}.
(Указание: воспользуйтесь пунктом~(ii) и
упражнением~\ref{exer-Sylow}(ii).)
\end{itemize}
\end{prob}

\newpage

\section{Когомологии Галуа}
\label{section:Galois-cohomology}

Пусть конечное расширение полей $K\subset L$ является расширением
Галуа. Когомологии группы $\Gal(L/K)$ называются {\it
когомологиями Галуа}. Для краткости положим
$$G={\Gal}(L/K)\,.$$
Эти обозначения будут использоваться
в разделах~\ref{subsection:spusk-dlya-kategorij}
и~\ref{subsection:forms-and-H1}.

\subsection{Спуск для расслоенных категорий}
\label{subsection:spusk-dlya-kategorij}

Когомологии Галуа тесно связаны со следующими вопросами. Пусть дан
некоторый алгебро-геометрический объект $Y$, определённый над полем $L$ (см. пример~\ref{examp-fibrcat}).
Какие условия гарантируют то, что на самом деле $Y$ определён над меньшим полем~$K$, то есть существует
аналогичный объект $X$ над $K$, и при этом $Y$
получается из $X$ расширением скаляров с $K$ на $L$?
Как описать все $X$, решающие данную задачу? Следующие несколько определений, примеров и упражнений
посвящены формализации этих вопросов и объяснению того,
как они описываются в терминах когомологий Галуа.

\bigskip

Пусть для каждого поля $E$ определена категория $\Mcal(E)$, а для
каждого расширения полей~\mbox{$E\subset F$} дан функтор
$$\Mcal(E)\to \Mcal(F),\quad X\mapsto X_F\,.$$
Этот функтор обычно называют \emph{расширением скаляров}. Кроме того, пусть
для каждой башни полей $E\subset F\subset M$ и каждого объекта $X$
из $\Mcal(E)$ зафиксирован функториальный по $X$ изоморфизм
$$(X_F)_M\cong X_{M}\,,$$
правильно ведущий себя в башнях вида
$$E\subset F\subset M\subset N\,.$$
Такой набор $\Mcal$ категорий $\Mcal(E)$ называется {\it
категорией, расслоенной над полями}.
В этой главе для краткости мы будем
также называть его просто \emph{расслоенной категорией}.
В дальнейшем~$\Mcal$
обозначает расслоенную категорию.

\bigskip

\begin{examp}{\it Категории, расслоенные над полями}\label{examp-fibrcat}

Следующие примеры категорий $\Mcal(E)$ определяют категории,
расслоенные над полями, причём функторы $\Mcal(E)\to\Mcal(F)$
определяются при помощи тензорного произведения на $F$ над $E$:
\begin{itemize}
\item[(i)]
категория векторных пространств (в том числе бесконечномерных)
над полем~$E$;
\item[(ii)]
категория векторных пространств с квадратичной формой над $E$;
\item[(iii)]
категория ассоциативных алгебр над $E$;
\item[(iv)]
категория аффинных многообразий над $E$;
\item[(v)]
категория квазипроективных многообразий над $E$;
\item[(vi)]
категория квазипроективных многообразий над $E$ с
отмеченной $E$-точкой;
\item[(vii)]
категория алгебраических групп над $E$.
\end{itemize}
\end{examp}

\bigskip

Напомним, что $K\subset L$
обозначает конечное расширение Галуа с группой Галуа~$G$.
Каждый элемент $g\in G$ определяет функтор
$$g_*\colon\Mcal(L)\to \Mcal(L)\,,$$
возникающий из ``расширения'' полей $g\colon L\stackrel{\sim}\longrightarrow
L$. Таким образом, для каждого объекта $Y$ из $\Mcal(L)$ определён
объект $g_*Y$ из $\Mcal(L)$, а для каждого морфизма $\theta\colon Y\to Y'$
определён морфизм
$$g_*\theta\colon g_*Y\to g_*Y'\,.$$
Можно сказать, что
группа $G$ действует на категории~$\Mcal(L)$.

Если $Y\cong X_L$ для некоторого объекта $X$ из $\Mcal(K)$, то из свойства
расширений скаляров в башне полей
$$K\hookrightarrow L\stackrel{g}\hookrightarrow L$$
возникает канонический изоморфизм
$$g_*X_L\stackrel{\sim}\longrightarrow X_L\,.$$
При этом изоморфизм
$$(gh)_*X_L\stackrel{\sim}\longrightarrow X_L$$
совпадает с композицией
изоморфизмов $g_*(h_*X_L\stackrel{\sim}\longrightarrow X_L)$ и
$g_*X_L\stackrel{\sim}\longrightarrow X_L$ для любых $g,h\in G$. Для
любых двух объектов~$X$ и~$X'$ из~$\Mcal(K)$ группа $G$ действует на
множестве морфизмов из~$X_L$ в~$X'_L$ в категории~$\Mcal(L)$.
А именно, для любого морфизма
$$\theta\colon X_L\to X'_L$$
результат действия~$g$ на~$\theta$, обозначаемый, как обычно,
символом ${}^g\theta$, определяется коммутативностью следующей диаграммы:
\begin{equation}\label{eq-actaut}
\begin{CD}
g_*X_L@>\sim>>X_L\\
@VVg_*\theta V@VV{}^g\theta V\\
g_*X'_L@>\sim>>X'_L
\end{CD}
\end{equation}
В частности, группа
автоморфизмов $\Aut(X_L)$ объекта $X_L$ в категории $\Mcal(L)$
является (возможно, неабелевым) $G$-мо\-ду\-лем.

Будем говорить, что на объекте $Y$ из категории $\Mcal(L)$ определены
{\it данные спуска}, если для любого элемента $g\in G$ задан
изоморфизм $\rho(g)\colon g_*Y\stackrel{\sim}\longrightarrow Y$,
удовлетворяющий приведённому выше условию на композицию элементов из
$G$, то есть
$$\rho(gh)=\rho(g)\circ g_*(\rho(h))$$
для любых $g,h\in G$.
Можно сказать, что $Y$ является $G$-эквивариантным объектом
в~$\Mcal(L)$.

В частности, как показано выше, для каждого объекта~$X$ из~$\Mcal(K)$
имеются канонические данные спуска на объекте~$X_L$, которые мы
будем обозначать
$$\rho_X(g)\colon g_*X\stackrel{\sim}\longrightarrow X\,.$$
Морфизмы между объектами с данными спуска $(Y,\rho)$ и $(Y',\rho')$
определяются естественным образом: это
морфизмы $\theta\colon Y\to Y'$, коммутирующие с ``действием'' $G$, то есть
такие, что для любого $g\in G$ имеется равенство
$$\theta\circ\rho(g)=\rho'(g)\circ g_*\theta$$
между морфизмами из $g_*Y$ в~$Y'$.
Обозначим через $\Mcal(L)^G$ категорию объектов из~$\Mcal(L)$ с данными
спуска.

Изоморфизм объектов с данными спуска~---
это морфизм объектов с данными спуска, у которого есть обратный.
В частности, данные спуска $\rho$ и $\rho'$ на одном и том же объекте $Y$
считаются изоморфными, если существует автоморфизм $\sigma$ объекта~$Y$,
задающий изоморфизм объектов с данными спуска
$$\sigma\colon (Y,\rho)\to (Y,\rho')\,.$$

Будем говорить, что $\Mcal$ {\it удовлетворяет условию спуска}, если
для любого конечного расширения Галуа $K\subset L$
функтор
$$X\mapsto (X_L, \rho_X)$$
из $\Mcal(K)$ в $\Mcal(L)^G$ является эквивалентностью категорий.
Другими словами, выполнены два условия.
Во-первых, для любых двух объектов~$X$ и~$X'$ из~$\Mcal(K)$ естественное
отображение между множествами морфизмов объектов в соответствующих
категориях
$$\Hom(X,X')\to\Hom\big((X_L,\rho_X), (X'_L, \rho_{X'})\big)$$
является
биекцией. Во-вторых, для любого объекта~\mbox{$(Y,\rho)$} из
$\Mcal(L)^G$ существует такой объект $X$ из
$\Mcal(K)$, что имеется изоморфизм объектов с данными спуска
$$(Y,\rho)\cong (X_L,\rho_X)\,.$$
В частности, объект $X$ определён однозначно
с точностью до изоморфизма в $\Mcal(K)$. Такой объект~$X$
мы будем называть \emph{спуском} объекта~$(Y,\rho)$ с~$L$ на~$K$.
Для краткости мы иногда будем говорить, что $X$ является спуском $Y$.
Также довольно часто мы будем ограничиваться проверкой лишь второго
из двух выписанных выше условий, так как выполнение первого условия почти всегда очевидно.

Для объекта $X$ из $\Mcal(K)$ его {\it $L$-формой} называется такой
объект $X'$ из $\Mcal(K)$, что существует изоморфизм $X_L\cong X'_L$ в
категории $\Mcal(L)$. Две $L$-формы объекта $X$ называются
изоморфными, если они изоморфны в категории $\Mcal(K)$. {\it Формой}
объекта $X$ называется $L$-форма для некоторого конечного
расширения Галуа $K\subset L$. Множество классов
изоморфизма $L$-форм объекта $X$ обозначим через $\Phi(X,L)$,
а множество классов
изоморфизма форм объекта $X$ над всевозможными полями~--- через $\Phi(X)$.

\bigskip
\begin{examp}{\it Данные спуска на расслоенных категориях}\label{examp-spusk}
\nopagebreak
\hspace{0cm}
\begin{itemize}
\item[(i)]
Пусть $\Mcal$ --- это расслоенная категория векторных пространств,
$U$ --- векторное пространство над полем $L$. Тогда $g_*U$
является той же абелевой группой~$U$, но с новым действием поля $L$,
заданным по формуле
$$\lambda\cdot u=\big(g^{-1}(\lambda)\big) u$$
для
$\lambda\in L$, $u\in U$. Если $V$ --- векторное пространство над
$K$, то канонический изоморфизм
$g_*V_L\stackrel{\sim}\longrightarrow V_L$ (напомним, что
$V_L=L\otimes_K V$) задаётся по формуле
$$\lambda\otimes v\mapsto g(\lambda)\otimes v\,.$$
После
выбора базиса в $V$ над $K$ изоморфизм
$g_*V_L\stackrel{\sim}\longrightarrow V_L$ задаётся просто применением $g$
к координатам векторов. При этом действие~$G$ на группе автоморфизмов
$\Aut(V_L)\cong\GL(V_L)$
задаётся применением $g$ к матричным элементам \mbox{$L$-линейных}
операторов. Определить данные спуска на векторном пространстве~$U$
над $L$ --- то же самое, что задать \emph{полулинейное действие} $G$ на~$U$,
то есть такое действие, что
$$g(\lambda u)=g(\lambda)g(u)$$
для любых $g\in G$, $\lambda\in L$ и $u\in U$.
Например, такая ситуация возникает, если задано $K$-линейное представление группы $G$ на векторном
пространстве $V$ над~$K$, и $U=V_L$ (ср. с упражнением~\ref{exer-spusk}(i)
ниже).
\item[(ii)]
Пусть $\Mcal$ --- это расслоенная категория векторных пространств c
квадратичной формой, $U$ --- векторное пространство над
полем $L$, а $q\colon U\to L$~--- квадратичная форма. Тогда пара
$g_*(U,q)$ состоит из векторного пространства~\mbox{$g_*U$} и квадратичной
формы, заданной по формуле
$$(g_*q)(u)=g\big(q(u)\big)$$
для $u\in U$. Для
пространства с квадратичной формой $(V,p)$ над $K$ изоморфизмы
$$g_*(V_L,p_L)\stackrel{\sim}\longrightarrow (V_L,p_L)$$
и действие
группы $G$ на $\Aut(V_L,p_L)$ определяются так же, как в пункте~(i).
Заметим, что действие группы $G$ на матричные элементы $L$-линейных
операторов из $V_L$ в себя сохраняет ортогональную подгруппу
$\Ort(V_L,p_L)$ в $\GL(V_L)$. Определить данные спуска на объекте
$(U,q)$ из $\Mcal(L)$~---
то же самое, что задать полулинейное действие $G$ на $U$,
для которого
$$g\big(q(u)\big)=q\big(g(u)\big)$$
при всех $g\in G$ и $u\in U$.
\item[(iii)]
Пусть $\Mcal$ --- это расслоенная категория ассоциативных алгебр,
$B$~--- ассоциативная алгебра над полем $L$. Тогда $g_*B$
имеет структуру векторного пространства, описанную в пункте (i). В
частности, $g_*B$ совпадает с $B$ как абелева группа. Отображение
умножения
$$g_*B\times g_*B\to g_*B$$
то же самое, что и для $B$. Для
алгебры $A$ над полем $K$ изоморфизмы
$$g_*A_L\stackrel{\sim}\longrightarrow A_L$$
и действие $G$ на $A_L$
определяются так же, как в пункте~(i). Определить данные спуска на
алгебре $B$ \mbox{над~$L$~---}
то же самое, что задать полулинейное действие~$G$ на~$B$,
коммутирующее с умножением в $B$.
\item[(iv)]
Пусть $\Mcal$ --- это расслоенная категория аффинных
многообразий, $Y$ --- аффинное многообразие над полем $L$. Пусть $B$ является $L$-алгеброй
регулярных функций на~$Y$, то есть $Y\cong\Spec(B)$. Тогда
алгебра регулярных функций на многообразии $g_*Y$ определяется как
$g_*B$, и всё делается так же, как в пункте (iii).
Заметим, что~$g_*Y$ канонически изоморфно над $L$ многообразию,
уравнения которого получаются применением $g$ к уравнениям исходного
многообразия $Y$.
\item[(v)]
Пусть $\Mcal$ --- это расслоенная категория квазипроективных
многообразий, $Y$~--- квазипроективное
многообразие над полем $L$. Представим $Y$ в виде $Y_1\setminus Y_2$, где~$Y_i$ являются проективными многообразиями в некотором проективном пространстве $\Pb^n$. Тогда $g_*Y=g_*Y_1\setminus g_*Y_2$, и каждое $g_*Y_i$ является проективным многообразием в $\Pb^n$, уравнения которого получены применением $g$ к уравнениям многообразия $Y_i$.
При этом для квазипроективного многообразия $X$ над~$K$
изоморфизм
$$g_*X_L\stackrel{\sim}\longrightarrow X_L$$
соответствует
тождественному отображению из $X$ в себя. Действие группы~$G$ на автоморфизм из~$\Aut(X_L)$ задаётся
применением $g$ к коэффициентам однородных многочленов на $\Pb^n$, задающих (локально по $X_L$) данный автоморфизм.
\item[(vi)]
Пусть $\Mcal$ --- это расслоенная категория квазипроективных
многообразий с отмеченной точкой, $(Y,y)$ --- квазипроективное многообразие над~$L$ с точкой
$y\in Y$ с координатами из $L$. Тогда $g_*(Y,y)$ --- это многообразие $g_*Y$ с точкой $g_*y$,
полученной применением $g$ к координатам исходной точки $y$. Остальное определяется подобно пункту (v).
\item[(vii)]
Пусть $\Mcal$ --- это расслоенная категория алгебраических групп,
$H$~--- алгебраическая группа над~$L$. Тогда $g_*H$ имеет
структуру многообразия над~$L$, описанную в пункте (v), единица в
$g_*H$ определяется как в пункте~(vi), а морфизмы умножения и взятия
обратного получаются применением~$g$ к коэффициентам соответствующих
морфизмов для исходной группы $H$. Остальное определяется подобно
пункту (v).
\end{itemize}
\end{examp}

\bigskip
Многие часто встречающиеся алгебро-геометрические объекты образуют
расслоенные категории, удовлетворяющие условию спуска. В частности,
в упражнениях~\ref{exer-spuskvp} и~\ref{exer-generalspusk} мы
покажем, что таковы все расслоенные категории из
примера~\ref{examp-fibrcat}. Во всех этих случаях единственность
объекта $X$, решающего задачу спуска для данных спуска $(Y,\rho)$,
очевидна, и главное --- проверить существование объекта $X$.

\begin{prob}{\bf Спуск для векторных пространств}\label{exer-spuskvp}
\hspace{0cm}
\begin{itemize}
\item[(i)]
Пусть $(U,\rho)$ --- векторное пространство над $L$ с данными спуска, то есть
с полулинейным действием $\rho$
группы $G$ на $U$ (см. пример~\ref{examp-spusk}(i)).
Положим~\mbox{$V=U^G$}.
Постройте каноническое отображение
векторных пространств над $L$ с данными спуска
$$(V_L,\rho_V)\to (U,\rho)\,,$$
то есть
$L$-линейное отображение $\alpha\colon V_L\to U$, коммутирующее с
(полулинейным) действием группы $G$.
\item[(ii)] Пусть дано произвольное расширение полей $K\subset E$. Рассмотрим $L$-алгебру
$$A=E\otimes_K L$$
с действием группы $G$ через правый сомножитель $L$.
Определим также $A$-модуль
$$W=A\otimes_L U\cong E\otimes_K U$$
с $A$-полулинейным действием группы $G$ через правый сомножитель $U$.
Проверьте, что $A^G=E\otimes 1\cong E$. Покажите, что $L$-линейное отображение
$$
\alpha\colon L\otimes_K U^G\to U,
$$
определённое в пункте (i),
является изоморфизмом тогда и только тогда, когда естественный гомоморфизм $A$-модулей
$$
\beta\colon A\otimes_E W^G\to W
$$
является изоморфизмом. (Указание: воспользуйтесь
изоморфизмом~\mbox{$W^G\cong E\otimes_K U^G$}.)
\item[(iii)] Пусть $E=L$ или $E=K^{sep}$. Докажите, что тогда имеет место
\mbox{$G$-эк}\-ви\-ва\-риантный изоморфизм колец
$$
A\cong \prod_G E\,,
$$
где действие группы $G$ на $\prod_G E$ заключается в перестановке сомножителей. (Указание: пусть $L=K(\theta)$, а $f\in K[T]$ является минимальным многочленом элемента $\theta\in L$. Тогда многочлен $f$ разлагается в произведение линейных сомножителей над полем $E$.)
\item[(iv)] Сохраняя предположение пункта (iii), докажите, что любой
\mbox{$G$-эк}\-ви\-ва\-риантный $A$-модуль имеет вид $\bigoplus_G P$,
где $P$ является $E$-векторным пространством, действие группы $G$
заключается в перестановке прямых слагаемых,
а структура $A$-модуля естественно определяется при помощи
изоморфизма~\mbox{ $A\cong \prod_G E$} из пункта (iii).
\item[(v)] Сохраняя предположение пункта (iii), покажите, что для любого
\mbox{$G$-эк}\-ви\-ва\-риантного $A$-модуля $M$
естественный гомоморфизм $A$-модулей
$$A\otimes_E M^G\to M$$
является изоморфизмом. (Указание: воспользуйтесь пунктом (iv)). Применяя пункт (ii), выведите из этого, что отображение $\alpha$ является изоморфизмом.
Таким образом, расслоенная категория векторных пространств удовлетворяет
условию спуска.
\end{itemize}
\end{prob}

Более явно, упражнение~\ref{exer-spuskvp}
означает, что в любом $L$-векторном пространстве с
полулинейным действием группы $G$ найдётся $G$-инвариантный
базис.

\bigskip

\begin{prob}{\bf Другие примеры спуска}\label{exer-generalspusk}
\hspace{0cm}
\begin{itemize}
\item[(i)]
Докажите, что расслоенная категория векторных пространств c
квадратичной формой удовлетворяет условию спуска. (Указание: для
векторного пространства с квадратичной формой $(U,q)$ над $L$, на котором определены
данные спуска (см. пример~\ref{examp-spusk}(ii)),
сначала постройте векторное пространство~$V$ над $K$
как в упражнении~\ref{exer-spuskvp}. Затем определите квадратичную форму~$p$
на~$V$ как ограничение $q$ на
$V\subset V_L\cong U$.)
\item[(ii)]
Докажите, что расслоенная категория ассоциативных алгебр
удовлетворяет условию спуска. (Указание: действуйте как в пункте (i).)
\item[(iii)]
Докажите, что расслоенная категория аффинных
многообразий удовлетворяет условию спуска.
(Указание: это является частным случаем пункта (ii).)
\item[(iv)]
Докажите, что расслоенная категория квазипроективных
многообразий удовлетворяет условию спуска. (Указание: многообразие $Y$, определённое над $L$,
можно рассматривать как многообразие над полем $K$,
поскольку расширение~\mbox{$K\subset L$}
конечно, и конечно порождённые алгебры над
полем~$L$ также конечно порождены и над $K$.
Для данных спуска $\rho$ на $Y$ (см. пример~\ref{examp-spusk}(v))
рассмотрим действие группы $G$ на многообразии
$Y$ над $K$, заданное композицией~\mbox{$\rho(g)\circ g$},
где
$$g\colon Y\to g_*Y$$
действует применением $g$
к координатам точек на~$Y$. Поскольку $Y$ квазипроективно,
существует покрытие многообразия~$Y$ над~$K$
открытыми аффинными $G$-инвариантными подмножествами.
Это позволяет всё свести к пункту (iii).
Таким образом, спуск многообразия $Y$ является фактором
$Y/G$ относительно определённого выше действия группы $G$,
где $Y$ рассматривается как многообразие над полем~$K$.)
\item[(v)]
Докажите, что расслоенная категория квазипроективных
многообразий с отмеченной точкой удовлетворяет условию спуска. (Указание: для
квазипроективного многообразия с отмеченной точкой $(Y,y)$ над $L$, на котором определены
данные спуска (см. пример~\ref{examp-spusk}(vi)), сначала постройте квазипроективное многообразие
$X$ над $K$ как в пункте (iv).
Затем заметьте, что точка $y$ соответствует точке на многообразии $X_L\cong Y$,
координаты которой инвариантны относительно действия $G$, что определяет $K$-точку $x$ на $X$.)
\item[(vi)]
Докажите, что расслоенная категория алгебраических групп
удовлетворяет условию спуска.
(Указание: воспользуйтесь тем, что алгебраические группы
квазипроективны, и действуйте как в пункте (v).)
\end{itemize}
\end{prob}

\bigskip
\begin{prob}{\textbf{Спуск для проективных многообразий}}
\label{prob:spusk-projective}

Пусть $Y$ --- многообразие над полем $L$ с данными спуска $\rho$,
а $X$ --- спуск многообразия~$Y$ с~$L$ на~$K$.
Пусть $U$ является конечномерным $L$-векторном пространством
с данными спуска, то есть на $U$ задано полулинейное
действие группы $G$ (см. пример~\ref{examp-spusk}(i)).

\begin{itemize}
\item[(o)]
Покажите, что спуск проективного пространства
$\mathbb{P}(U)$ над $L$ является проективным
пространством $\mathbb{P}(U^G)$ над $K$.
(Указание: пользуясь упражнением~\ref{exer-spuskvp}(ii),(v),
постройте изоморфизм
$$
\mathbb{P}(U^G)_L\cong \mathbb{P}(U)
$$
многообразий с данными спуска над $L$.)
\item[(i)]
Пусть $\dim(U)=m+1$. Предположим, что замкнутое вложение
$$\iota\colon Y\hookrightarrow\P^m=\P(U)$$
над $L$ является морфизмом многообразий
с данными спуска.
Покажите, что имеется замкнутое вложение
$$X\hookrightarrow\P^m=\P(U^G)$$
над $K$.
(Указание: воспользуйтесь пунктом~(o).)
\item[(ii)]
В предположениях пункта~(i) пусть $K\subset F$ является
произвольным расширением полей. Выберем ненулевой вектор $v\in (U^G)_F$
и рассмотрим соответствующую ему $F$-точку
$$[v]\in \P(U^G)(F)\,.$$
Пусть
$$\phi\colon (U^G)_F\to U_F$$
является расширением скаляров
с $K$ на $F$ естественного вложения $K$-векторных пространств
$U^G\to U$.
Покажите, что $[v]\in X(F)$ тогда и только тогда,
когда
$$[\phi(v)]\in Y(F)\,.$$
(Указание: это непосредственно следует
из изоморфизма $X_K\cong Y$.)
\item[(iii)]
Пусть $Y\subset\P^n=\P(W)$ --- замкнутое вложение над $L$,
где $W$ является \mbox{$L$-векторным} пространством размерности $n+1$.
Докажите, что имеется замкнутое вложение
$$X\subset\Pb^{(n+1)^d-1}$$
над $K$, где $d=[L:K]$. Таким образом,
спуск проективного многообразия~--- проективное многообразие.
(Указание:
рассмотрим векторное пространство (см. пример~\ref{examp-spusk}(i))
$$U=\mbox{$\bigotimes\limits_{g\in G} g_*W\,.$}$$
Определим данные спуска
на $U$, задав (полулинейное) действие группы $G$ на разложимых
тензорах по формуле
$$h\colon {}\otimes_{g}w_g\mapsto {}\otimes_{g}w_{h^{-1}g}\,.$$
Определим вложение $\iota\colon Y\hookrightarrow\Pb(U)$
как композицию вложения
$$\mbox{$\prod\limits_g\rho(g)^{-1}\colon Y\hookrightarrow
\prod\limits_{g\in G} g_*Y$}\,,$$
естественного отображения
$$\mbox{$\prod\limits_{g\in G} g_*Y$}\hookrightarrow
\mbox{$\prod\limits_{g\in G} \P(g_*W)$}$$
и вложения Сегре
$$\mbox{$\prod\limits_{g\in G} \P(g_*W)$}\hookrightarrow\P(U)\,.$$
Покажите, что $\iota$ является
морфизмом объектов с данными спуска. Наконец, воспользуйтесь пунктом~(i).)
\end{itemize}
\end{prob}

\subsection{Формы и первые когомологии Галуа}
\label{subsection:forms-and-H1}

В этом разделе $\Mcal$ обозначает категорию, расслоенную над полями.

\begin{prob}{\bf Данные спуска, $1$-коциклы и формы}\label{exer-spusk}
\hspace{0cm}
\begin{itemize}
\item[(i)]
Покажите, что для любого объекта $X$ из $\Mcal(K)$ существует
каноническая биекция между $H^1\big(G,\Aut(X_L)\big)$ (см.
упражнение~\ref{exer-nonabcohom}) и множеством классов изоморфизма
данных спуска на $X_L$, которая сопоставляет отмеченному элементу из
$H^1\big(G,\Aut(X_L)\big)$ канонические данные спуска $\rho_X$ на $X_L$.
(Указание: используйте соответствие
$$\rho(g)=\varphi(g)\circ\rho_X(g)$$
между данными спуска $\rho$ и
\mbox{$1$-коциклами} $\varphi$. Также вспомните, что
$${}^g\theta\circ\rho_X(g)=\rho_X(g)\circ g_*\theta$$
для любого морфизма
$\theta\colon X_L\to X_L$, см. диаграмму~\eqref{eq-actaut}.)
В частности, любой гомоморфизм $G\to\Aut(X)$ определяет
данные спуска на $X_L$ (см. упражнение~\ref{exer-012}(i)).
\item[(ii)]
Для любого объекта $X$ из $\Mcal(K)$ его $L$-форма $X'$ определяет
данные спуска~$\rho_{X'}$ на объекте $X_L\cong X'_L$ (возможно,
отличные от данных спуска $\rho_X$). Докажите, что это определяет
отображение из множества $\Phi(X,L)$ классов изоморфизма $L$-форм
объекта $X$ в множество классов изоморфизма данных спуска на $X_L$.
\item[(iii)]
Покажите, что если $\Mcal$ удовлетворяет условию спуска, то имеется
каноническая биекция между $\Phi(X,L)$ и множеством классов
изоморфизма данных спуска на $X_L$. Таким образом, в этом случае
из пункта~(i) возникает каноническая биекция
$$\Phi(X,L)\cong H^1\big(G,\Aut(X_L)\big)\,.$$
В частности, $1$-коцикл
$$\varphi\in\nolinebreak Z^1\big(G,\Aut(X_L)\big)$$
определяет форму многообразия~$X$, которая называется
{\it скруткой $X$ при помощи $\varphi$}.
\item[(iv)]
Опишите все формы многообразия $X$ над $K$, состоящего из двух
точек.
\item[(v)]
Пусть $R\subset L$ --- множество корней неприводимого
многочлена, задающего расширение $K\subset L$.
Покажите, что $\Spec(L)$ является $L$-формой многообразия
$$\mbox{$X=\coprod\limits_R\Spec(K)$}\,,$$
и соответствующий $1$-коцикл задаётся
(ср. с пунктом~(i))
естественным гомоморфизмом
$$G\to\Aut(X)=\Aut(R)\,.$$
\item[(vi)] Докажите \emph{теорему Гильберта 90}:
$$H^1\big(G,\GL_n(L)\big)=\{1\}\,.$$
(Указание: по упражнению~\ref{exer-spuskvp}
категория (конечномерных)
векторных пространств удовлетворяет условию спуска.
Поскольку все векторные пространства одинаковой размерности изоморфны,
остаётся применить пункт~(iii).)
В частности,
при~\mbox{$n=1$} получаем
$$H^1(G,L^*)=\{1\}\,.$$
\end{itemize}
\end{prob}

\bigskip

Существует (по сути равносильное)
явное алгебраическое доказательство теоремы Гильберта~90.

\begin{prob}{\bf Теорема Гильберта 90 (алгебраическое доказательство)}
\hspace{0cm}
\begin{itemize}
\item[(i)]
Докажите \emph{лемму Артина о независимости характеров}:
пусть задан набор
$$\chi_1,\ldots,\chi_n\colon H\to L^*$$
попарно различных гомоморфизмов
из произвольной группы $H$ в $L^*$. Тогда гомоморфизмы $\chi_1,\ldots,\chi_n$
линейно независимы над $L$, то есть
если~\mbox{$\lambda_1,\ldots,\lambda_n\in L$} таковы, что
$$\sum\limits_{i=1}^n\lambda_i\chi_i(h)=0$$
для любого $h\in H$, то
все $\lambda_i$ равны нулю.
(Указание: рассмотрите минимальную линейную зависимость
$$
\sum\limits_{i=1}^n\lambda_i\chi_i=0
$$
и равенство
$\sum_{i=1}^n\lambda_i\chi_i(gh)=0$
для произвольных $g,h\in H$.)
\item[(ii)]
Рассмотрим $1$-коцикл
$$\varphi\in Z^1\big(G,\GL_n(L)\big)\,.$$
Для каждого вектора-столбца
$u$ длины $n$ с компонентами из $L$ определим вектор-столбец
$$
\Psi(u)=\sum\limits_{g\in G}\varphi(g)\cdot\,{}^gu\,.
$$
Пусть $v\in L^n$ --- такой вектор-строка, что
$v\cdot \Psi(\lambda u)=0$ для любого элемента~\mbox{$\lambda\in L$}.
Докажите, что тогда $v\cdot u=0$.
(Указание: воспользуйтесь равенством
$$
v\cdot \Psi(\lambda u)=
\sum_{g\in G}(v\cdot\varphi(g)\cdot\,{}^gu)\,\,{}^g\lambda
$$
и леммой Артина, применённой к характерам
$\lambda\mapsto {}^g\lambda$ из $L^*$ в себя.)
\item[(iii)]
Для каждой
матрицы $A$ размера $n\times n$ с коэффициентами
из $L$ положим
$$
\Psi(A)=\sum\limits_{g\in G}\varphi(g)\cdot\,{}^gA\,.
$$
Покажите, что существует такая матрица $A$,
что матрица $\Psi(A)$ обратима. (Указание:
воспользуйтесь тем, что
$$\Psi(u_1,\ldots,u_n)=\big(\Psi(u_1),\ldots,\Psi(u_n)\big)$$
для векторов-столбцов $u_1,\ldots,u_n$, а также тем, что по пункту~(ii)
всевозможные векторы $\Psi(u)$ порождают пространство $L^n$.)
\item[(iv)]
Докажите, что
$${}^g \Psi(A)=\varphi(g)^{-1}\Psi(A)$$
для любого $g\in G$.
(Указание: вычислите явно
${}^g \Psi(A)$ и воспользуйтесь условием $1$-коцикла для $\varphi$.)
Воспользовавшись пунктом (iii), выведите отсюда
тривиальность \mbox{$1$-коцикла}~$\varphi$.
\end{itemize}
\end{prob}

\bigskip
\begin{prob}{\bf Поведение при замене поля}\label{exer-changefield}

Пусть $X$ --- объект из $\Mcal(K)$, и пусть дано конечное расширение
Галуа~\mbox{$K\subset F$}, такое что имеется
композиция конечных расширений Галуа
$$K\subset L\subset F\,.$$
Положим $N=\Gal(F/L)$ и $H=\Gal(F/K)$, так что имеется точная
последовательность
$$
1\to N\to H\to G\to 1\,.
$$
\begin{itemize}
\item[(i)]
Покажите, что следующая диаграмма коммутативна:
$$
\begin{CD}
\Phi(X,L)@>>>\Phi(X,F)\\
@VVV@VVV\\
H^1\big(G,\Aut(X_L)\big)@>>>H^1\big(H,\Aut(X_F)\big)
\end{CD}
$$
Здесь верхняя горизонтальная стрелка определяется естественным
образом, вертикальные стрелки являются отображениями из
упражнения~\ref{exer-spusk}(ii),(iii), а нижняя горизонтальная стрелка
является композицией отображения обратного образа
$$H^1\big(G,\Aut(X_L)\big)\to H^1\big(H,\Aut(X_L)\big)\,,$$
определённого гомоморфизмом
групп $H\twoheadrightarrow G$, и естественного отображения
$$H^1\big(H,\Aut(X_L)\big)\to H^1\big(H,\Aut(X_F)\big)\,.$$
\item[(ii)]
Покажите, что следующая диаграмма коммутативна:
$$
\begin{CD}
\Phi(X,F)@>>>\Phi(X_L,F)\\
@VVV@VVV\\
H^1\big(H,\Aut(X_F)\big)@>>>H^1\big(N,\Aut(X_F)\big)
\end{CD}
$$
Здесь верхняя горизонтальная стрелка задаётся расширением скаляров для
форм с $K$ на $L$, вертикальные стрелки являются отображениями из
упражнения~\ref{exer-spusk}(ii),(iii), а нижняя горизонтальная стрелка
является отображением обратного образа, определённого гомоморфизмом
групп $N\hookrightarrow H$. Заметьте, что данная коммутативная
диаграмма также имеет место, если $L$ над $K$ не является
расширением Галуа.
\end{itemize}
\end{prob}

\bigskip

\begin{prob}{\bf Когомологии с аддитивными коэффициентами}\label{exer-cohomadd}
\hspace{0cm}
\begin{itemize}
\item[(o)]
Чему равна группа $H^0(G,L)$?
\item[(i)]
Покажите, что $G$-модуль $K[G]$ изоморфен $\iota_*(K)$, см.
определение~\ref{defin-coind}, где
$\iota\colon\{e\}\hookrightarrow G$ является вложением тривиальной подгруппы.
\item[(ii)]
Докажите, что
$$H^i\big(G,K[G]\big)=0$$
при $i>0$. (Указание: используйте лемму
Шапиро (см. упражнение~\ref{exer-ind-inj}(iii))
для вложения $\{e\}\hookrightarrow G$.)
\item[(iii)]
Докажите, что
$$H^i(G,L)=0$$
при $i>0$. (Указание: воспользуйтесь
\emph{теоремой о нормальном базисе}, согласно которой в $L$ над $K$
существует такой базис
$$\{e_g\mid g\in G\}\,,$$
что $g(e_h)=e_{gh}$ для
любых $g,h\in G$; затем используйте пункт~(ii).)
\end{itemize}
\end{prob}

\bigskip

\begin{prob}{\bf Первые когомологии Галуа некоторых групп}
\label{exer-1st-Galois}
\hspace{0cm}
\begin{itemize}
\item[(i)]
Покажите, что
$$H^1\big(G,\SL_n(L)\big)=\{1\}\,.$$
(Указание: воспользуйтесь точной последовательностью
$$1\to \SL_n(L)\to \GL_n(L)\to L^*\to 1\,,$$
теоремой Гильберта~90 и упражнением~\ref{exer-nonabcohom}(ii).)
\item[(ii)]
Покажите, что
$$H^1\big(G,\Aff_n(L)\big)=\{1\}\,,$$
где $\Aff_n$ обозначает группу аффинных автоморфизмов
$n$-мерного аффинного пространства.
(Указание: воспользуйтесь точной последовательностью
$$1\to L^n\to \Aff_n(L)\to \GL_n(L)\to 1\,,$$
теоремой Гильберта~90 и упражнением~\ref{exer-nonabcohom}(ii).)
\item[(iii)]
Покажите, что
$$H^1\big(G,\Sp_{2n}(L)\big)=\{1\}\,,$$
если ${\rm char}(K)\ne 2$,
где $\Sp_{2n}$ обозначает симплектическую группу матриц
размера~\mbox{$2n\times 2n$}. (Указание:
вспомните, как устроены невырожденные кососимметрические формы над $K$.)
\item[(iv)]
Пусть множество $S$ состоит из таких комплексных матриц $M$
размера~\mbox{$n\times n$}, что
$$\bar{M}=M^T=M^{-1}\,,$$
где $\bar M$ обозначает комплексное сопряжение матрицы~$M$,
а $M^T$~--- её транспонирование.
Рассмотрим действие ортогональной группы $\Ort_n(\Cb)$ на $S$,
заданное по формуле
$$A\colon M\mapsto AM (\bar{A})^{-1}$$
для $A\in \Ort_n(\Cb)$, $M\in S$.
Сколько элементов в фактормножестве $S/\Ort_n(\Cb)$?
(Указание: используя упражнение~\ref{exer-nonabcohom}(v),
покажите, что имеется биекция
$$
H^1\big(\Gal(\Cb/\Rb),\Ort_n(\Cb)\big)\cong S/\Ort_n(\Cb).
$$
Далее вспомните, как устроены невырожденные квадратичные формы
над~$\Rb$.)
\end{itemize}
\end{prob}

\bigskip

\begin{defin}\label{defin:Severi-Brauer}
{\it Многообразием Севери--Брауэра} называется форма проективного пространства.
\end{defin}

\begin{prob}{\bf Формы $\Pb^1$}

Докажите, что множество $\Phi(\Pb^1)$ находится во взаимно
однозначном соответствии с множеством классов проективной
эквивалентности гладких коник в~$\Pb^2$.
(Указание: рассмотрите вложение формы $X$ в $\Pb^2$
при помощи очень обильного антиканонического пучка~$\omega_X^{-1}$.)
\end{prob}

\bigskip

\begin{prob}{\bf Многообразия Севери--Брауэра с точкой}\label{exer-SBforms}
\begin{itemize}
\item[(o)]
Покажите, что множество классов изоморфизма многообразий Севери--Брауэра,
являющихся $L$-формами $\Pb^{n-1}$,
находится в канонической биекции с множеством $H^1\big(G,\PGL_n(L)\big)$.
\item[(i)]
Докажите, что
многообразие Севери--Брауэра имеет точку над полем $K$ тогда и только
тогда, когда оно изоморфно проективному пространству
над~$K$. (Указание: какая группа автоморфизмов у проективного
пространства с отмеченной точкой?)
\item[(ii)]
Пусть $X$ является многообразием Севери--Брауэра над полем $K$.
Докажите, что $X_{K(X)}$ изоморфно проективному
пространству над полем $K(X)$, где~$K(X)$ обозначает поле
рациональных функций на $X$. (Указание: постройте \mbox{$K(X)$-точку}
на $X_{K(X)}$ и воспользуйтесь пунктом (i).)
\end{itemize}
\end{prob}

\subsection{Когомологии проконечных групп}
\label{subsection:cohom-profinite}

Часто бывает удобнее использовать абсолютную группу Галуа
$$G_K=\Gal(K^{sep}/K)\,,$$
не фиксируя конечного расширения $L$ над $K$. Поскольку $G_K$ является проконечной группой,
это мотивирует рассмотрение когомологий для таких групп.

\bigskip

Пусть $(I,\leqslant)$ --- направленное
частично упорядоченное множество, то есть такое частично упорядоченное
(возможно, бесконечное)
множество, что для любых элементов $i_1,i_2\in I$ существует элемент
$j\in I$, для которого $i_1\leqslant j$ и $i_2\leqslant j$. Пусть для
каждого~$i\in I$ задана конечная группа $G_i$, а для каждой пары
сравнимых элементов~\mbox{$i\leqslant j$} в~$I$ задан такой гомоморфизм групп
$$\varphi_{ji}\colon G_j\to G_i\,,$$
что $\varphi_{ii}={\rm id}_{G_i}$ и
$$\varphi_{ji}\circ\varphi_{kj}= \varphi_{ki}$$
для любых $i\leqslant
j\leqslant k$. Набор $\{G_i\}$, $i\in I$, называется {\it проективной системой конечных групп}.

Напомним, что на (возможно,
бесконечном) произведении конечных дискретных групп $\prod_{i\in
I}G_i$ определена \emph{топология Тихонова}:
база окрестностей единицы в
ней~--- это подмножества вида
$$\mbox{$U_S=\prod\limits_{i\in S}\{e\}\times \prod\limits_{i\in
G\backslash S}G_i$}\,,$$
где $S\subset I$~--- конечное подмножество.
\emph{Обратный предел} $\varprojlim G_i$ состоит из таких наборов
$$\mbox{$(g_i)\in\prod\limits_{i\in I} G_i\,,$}$$
что $\varphi_{ji}(g_j)=g_i$ для любых
$i\leqslant j$. Будем рассматривать на множестве $\varprojlim G_i$ топологию,
являющуюся ограничением с $\prod_{i\in I}G_i$ топологии Тихонова.
Легко проверить, что обратный предел проективной системы групп является группой.

\begin{defin}\label{defin-profinitegrp}
{\it Проконечной группой}
называется топологическая группа, изоморфная обратному пределу
$\varprojlim G_i$ для некоторой проективной
системы конечных групп $\{G_i\}$, $i\in I$.
\end{defin}

\begin{examp}
Следующие примеры топологических групп являются проконечными группами:
\begin{itemize}
\item[(i)]
группа $p$-адических чисел по сложению $\Z_p=\varprojlim\Z/p^n\Z$;
\item[(ii)]
группа $\widehat\Z=\varprojlim\Z/n\Z$;
\item[(iii)]
{\it проконечное пополнение} произвольной группы $G$, которое определяется
как
$$\widehat G=\varprojlim G/N\,,$$
где предел берётся по всем
нормальным подгруппам $N$ конечного индекса в группе~$G$;
\item[(iv)]
{\it про-$p$-пополнение} произвольной группы $G$, которое определяется
как
$$\widehat G_p=\varprojlim G/N\,,$$
где предел берётся по всем
нормальным подгруппам~$N$, индекс которых в группе~$G$ является
степенью числа~$p$.
\end{itemize}
\end{examp}

\bigskip

\begin{prob}{\bf Проконечные группы}\label{exer-profgrp}
\nopagebreak \hspace{0cm}
\begin{itemize}
\item[(o)]
Покажите, что проекции из произведения $\prod_i G_i$ в $G_j$
определяют гомоморфизмы
$$\pi_j\colon \varprojlim G_i\to G_j\,.$$
Докажите,
что $\varprojlim G_i$ обладает следующим универсальным свойством:
если для группы $H$ заданы такие гомоморфизмы
$$\xi_i\colon H\to G_i\,,$$
что $\varphi_{ji}\circ\xi_j=\xi_i$ для любых $i\leqslant j$, то
существует единственный гомоморфизм
$$\xi\colon H\to \varprojlim G_i\,,$$
для которого
$\pi_j\circ \xi=\xi_j$ при всех $j\in I$.
\item[(i)]
Покажите, что любая проконечная группа компактна. (Указание: сначала
докажите, что произведение $\prod_{i\in I}G_i$ компактно в топологии
Тихонова, а затем покажите, что $\varprojlim G_i$ является замкнутым
подмножеством в $\prod_{i\in I}G_i$.)
\item[(ii)]
Покажите, что топология на $\varprojlim G_i$ --- это слабейшая
топология, в которой все естественные отображения проекции
$$\pi_j\colon\varprojlim G_i \to G_j$$
непрерывны. Покажите, что подгруппа
$$\Ker(\pi_j)\subset \varprojlim G_i$$
открыта и замкнута.
\item[(iii)]
Докажите, что любая проконечная группа изоморфна своему
проконечному пополнению.
\item[(iv)]
Докажите, что любая открытая подгруппа в проконечной группе имеет
конечный индекс. (Указание: покажите, что фактортопология на
множестве смежных классов дискретна и компактна.)
\end{itemize}
\end{prob}

\bigskip

\begin{prob}{\bf Проконечная теория Галуа}

Пусть $L\subset K$~--- расширение Галуа,
то есть сепарабельное нормальное алгебраическое расширение полей
(возможно, бесконечное).
\begin{itemize}
\item[(i)]
Постройте изоморфизм
$$\Gal(L/K)\cong \varprojlim {\Gal}(L_i/K)\,,$$
где $L_i$ пробегает все подполя в $L$, являющиеся конечными
расширениями Галуа поля $K$. Таким образом, группа Галуа
$\Gal(L/K)$ является
проконечной группой.
\item[(ii)]
Покажите, что база открытых окрестностей единицы в $\Gal(L/K)$
образована стабилизаторами элементов из $L$.
\item[(iii)]
Докажите \emph{основную теорему теории Галуа}: существует естественная биекция
между множествами замкнутых подгрупп $H\subset \Gal(L/K)$ и всех
промежуточных расширений $K\subset E\subset L$, задаваемая по
обычной формуле
$$E=L^H, \quad H=\Gal(L/E)\,.$$
При этом подгруппа $H$
нормальна в $\Gal(L/K)$
тогда и только тогда, когда поле~$E$ нормально над~$K$.
(Указание: чтобы установить замкнутость подгруппы~$\Gal(L/E)$
для заданного поля $E$ используйте равенство
$$\Gal(L/E)=\bigcap_i \Gal(L/E_i)\,,$$
где $E_i$ пробегает все
подполя $E_i\subset E$, конечные над $K$. Также воспользуйтесь тем, что
$\Gal(L/E_i)$ является замкнутой подгруппой в $\Gal(L/K)$: это
доказывается применением упражнения~\ref{exer-profgrp}(ii) к
гомоморфизмам
$$\Gal(L/K)\to\Gal(F_i/K)\,,$$
где $F_i$
обозначает нормальное замыкание $E_i$ над $K$ в $L$. Для заданной
замкнутой подгруппы
$$H\subset \Gal(L/K)$$
покажите, что $H$ также является
замкнутой подгруппой в группе~$\Gal(L/E)$, где $E=L^H$. При помощи обычной
теории Галуа покажите, что для любой нормальной подгруппы конечного
индекса
$$N\subset \Gal(L/E)$$
индуцированное отображение
$$H\to\Gal(L/E)/N$$
сюръективно. Используя замкнутость подгруппы $H$ в~$\Gal(L/E)$,
выведите отсюда, что $H=\Gal(L/E)$.)
\end{itemize}
\end{prob}

\bigskip

Например, для конечного поля $\F_p$ имеется изоморфизм
$$\Gal(\bar\F_p/\F_p)\cong \widehat{\Z}\,.$$
Эта проконечная группа топологически
порождена отображением Фробениуса
$$\mathrm{Fr}\colon\bar{\F}_p\to\bar{\F}_p,\quad x\mapsto x^p\,.$$

\bigskip

\begin{prob}{\bf Дискретные модули}
\nopagebreak
\hspace{0cm}
\begin{itemize}
\item[(i)]
Пусть $G$ --- проконечная группа, $M$ --- некоторый
$G$-модуль с дискретной топологией.
Покажите,
что отображение действия
$$G\times M\to M$$
непрерывно тогда и только тогда,
когда для любого элемента
из $M$ его стабилизатор открыт в~$G$
(напомним, что по упражнению~\ref{exer-profgrp}(iv)
открытая подгруппа проконечной группы имеет конечный индекс).
Такие $G$-модули называются {\it дискретными}.
\item[(ii)]
Покажите, что для любого дискретного $G$-модуля $M$
имеется равенство
$$M=\bigcup_U M^U\,,$$
где $U$ пробегает все открытые нормальные подгруппы в $G$.
\end{itemize}
\end{prob}

\bigskip
Пусть, как и раньше, $(I,\leqslant)$ --- направленное частично упорядоченное
множество. Пусть $\{X_i\}$, $i\in I$ ---
\emph{индуктивная система множеств}, то есть
для каждой пары сравнимых элементов
$i\leqslant j$ задано отображение
$$\psi_{ij}\colon X_i\to X_j\,,$$
причём $\psi_{ii}={\rm id}_{X_i}$ и
$\psi_{jk}\circ\psi_{ij}= \psi_{ik}$ для любых $i\leqslant
j\leqslant k$.
\emph{Прямой предел}~\mbox{$\varinjlim X_i$}
определяется как фактормножество
$$\mbox{$\big(\coprod\limits_{i\in I}X_i\big)/
\big(\varphi_{ij}(x_i)\sim x_j\big)\,.$}$$
Также прямой предел можно определять при помощи
универсального свойства, ср.~с~упражнением~\ref{exer-profgrp}(o).
При этом прямой предел индуктивной системы
абелевых групп является абелевой группой.

\begin{defin}\label{defin:pryamoj-predel}
Для дискретного $G$-модуля $M$ и открытых нормальных подгрупп
$U_2\subset U_1$ в $G$ рассмотрим отображение инфляции
$$H^i\big(G/U_1,M^{U_1}\big)\to H^i\big(G/U_2,M^{U_2}\big),
\quad i\geqslant 0\,,$$
то есть композицию отображения обратного образа
$$H^i\big(G/U_1,M^{U_1}\big)\to H^i\big(G/U_2,M^{U_1}\big)$$
и естественного отображения
$$H^i\big(G/U_2,M^{U_1}\big)\to H^i\big(G/U_2,M^{U_2}\big)\,.$$
Это задаёт индуктивную систему абелевых групп.
{\it Когомологии
проконечной группы $G$ с коэффициентами в дискретном модуле~$M$}
определяются по формуле
$$
H^i(G,M)=\varinjlim H^i\big(G/U,M^{U}\big)\,,
$$
где прямой предел берётся по всем открытым нормальным подгруппам $U\subset G$.
Аналогично определяются когомологии $H^0(G,\Gamma)$ и $H^1(G,\Gamma)$
для неабелева дискретного $G$-модуля $\Gamma$
(ср. с упражнением~\ref{exer-nonabcohom}).
\end{defin}

Отметим, что когомологии проконечной группы~$G$ с коэффициентами в дискретном
модуле~$M$ могут отличаться от когомологий $G$ как абстрактной группы
с коэффициентами в $M$ (см. определение~\ref{defin:group-cohomology}).
В то же время для обеих этих групп когомологий мы используем
одинаковые обозначения. В дальнейшем если группа~$G$
проконечна, и $G$-модуль $M$ дискретен, то
$H^i(G,M)$ обозначает когомологии в смысле
определения~\ref{defin:pryamoj-predel}.

\bigskip
Следующее упражнение показывает, что для проконечных групп когомологии обладают
теми же основными свойствами, что и для конечных групп.

\begin{prob}{\bf Когомологии проконечных
групп}\label{exer-profcohom}

Пусть $G$ --- проконечная группа, $M$ --- дискретный $G$-модуль.
\begin{itemize}
\item[(o)]
Покажите, что $H^0(G,M)\cong M^G$.
\item[(i)]
Рассмотрим комплекс
$\Map_c(G^{\udot},M)$, члены которого состоят из непрерывных отображений из
$G^i$ в $M$, $i\geqslant 0$, а дифференциал определяется так же,
как в комплексе~\eqref{eq-stcompl} из
упражнения~\ref{exer-stcoml}(v). Докажите, что
$$H^i(G,M)\cong H^i\big(\Map_c(G^{\udot},M)\big)\,.$$
(Указание: покажите, что
$$\Map_c\big(G^{\udot},M\big)\cong
\varinjlim \Map\big((G/U)^{\udot},M^U\big)\,,$$
где прямой предел берётся по всем
открытым нормальным подгруппам \mbox{$U\subset G$}.
Далее воспользуйтесь тем, что
когомологии коммутируют с направленным прямым пределом.)
\item[(ii)] Покажите, что если $M$ --- тривиальный $G$-модуль, то
$$H^1(G,M)\cong\Hom(G,M)\,,$$
где $\Hom(G,M)$ обозначает группу
непрерывных гомоморфизмов, то есть гомоморфизмов с конечным образом.
\item[(iii)] Покажите, что
$H^1(G,\Z)=0$, и имеется канонический изоморфизм
$$H^2(G,\Z)\cong\Hom(G,\Q/\Z)\,,$$
где $\Z$ и $\Q/\Z$ рассматриваются со
структурой тривиальных $G$-модулей.
(Указание: воспользуйтесь упражнением~\ref{exercise:division-module}(ii).)
\item[(iv)]
Покажите, что $H^i(G,M)$ являются группами кручения при $i>0$.
(Указание: воспользуйтесь
упражнением~\ref{exer-cohomfingrp}(i).)
\item[(v)]
Покажите, что точная тройка дискретных $G$-модулей
$$
0\to M'\to M\to M''\to 0
$$
определяет длинную точную последовательность когомологий.
\item[(vi)]
Докажите, что
$$H^i(G,M)\cong\Ext^i_G(\Z,M)\,,$$
где $\Ext$ рассматривается в абелевой категории
дискретных $G$-модулей.
\item[(vii)]
Пусть
$G$ --- проконечная группа, $H$ --- конечная группа, и
$f\colon G\twoheadrightarrow H$ --- сюръективный гомоморфизм
с ядром $I$.
Пусть $M$ --- дискретный $G$-модуль.
(Проверьте, что $I$ --- проконечная группа,
и $M$ --- дискретный $I$-модуль!)
Покажите, что имеется точная последовательность:
$$
0\to H^1(H,M^I)\to H^1(G,M)\to H^1(I,M)^H\to H^2(H,M^I)\to H^2(G,M)\,.
$$
(Указание:
воспользуйтесь упражнением~\ref{exer-ind-sur}(iv) и перейдите к пределу.)
На самом деле, то же утверждение верно и в случае,
если не требовать конечности группы $H$, а требовать только замкнутость
проконечной подгруппы $I\subset G$.
\end{itemize}
\end{prob}

\bigskip
\begin{prob}{\bf Первые когомологии группы $\widehat{\Z}$}
\label{prob:H1-ot-Z-s-kryshkoj}

Пусть $M$ --- дискретный модуль над
группой $\widehat{\Z}$,
каждый элемент которого имеет конечный порядок
(как элемент абелевой группы $M$).
Покажите, что
$$H^1(\widehat{\Z},M)\cong M/(s-1)M\,,$$
где $s\in\widehat{\Z}$
является произвольной топологической образующей группы
$\widehat{\Z}$.
(Указание: докажите,
что для любого элемента $m\in M$ существует натуральное число
$n$, для которого выполнено равенство
$$\sum\limits_{i=0}^{n-1}s^i(m)=0\,.$$
После этого
воспользуйтесь явным видом когомологий для конечных
циклических групп, полученным
в упражнении~\ref{exer-cyclgrp}(i), и перейдите к прямому пределу
из определения~\ref{defin:pryamoj-predel}.)
\end{prob}

\bigskip
\begin{prob}{\bf Конечность первых когомологий}
\label{prob:H1-finite}

Пусть $G$ --- проконечная группа, $M$ --- дискретный $G$-модуль,
конечно порождённый как абелева группа.
\begin{itemize}
\item[(i)]
Покажите, что ядро
$I$ действия группы $G$ на модуле $M$ имеет конечный индекс в $G$.
(Указание: выберите какой-нибудь набор порождающих
$$m_1,\ldots, m_r\in M$$
модуля $M$ как абелевой группы,
рассмотрите пересечение $I\subset G$ стабилизаторов
элементов $m_i$ в группе $G$, и воспользуйтесь
упражнением~\ref{exer-profgrp}(iv).)
\item[(ii)]
Предположим, что модуль $M$ свободен
как абелева группа.
Проверьте, что
$$H^1(I, M)=0\,.$$
(Указание: воспользуйтесь упражнением~\ref{exer-profcohom}(ii).)
\item[(iii)]
Предположим, что модуль $M$ свободен
как абелева группа. Пусть $H=G/I$.
Тогда $M$ можно рассматривать как дискретный $H$-модуль.
Докажите, что
$$H^1(G, M)\cong H^1(H,M)\,.$$
(Указание: воспользуйтесь пунктом~(ii)
и упражнением~\ref{exer-profcohom}(vii).)
В частности, группа $H^1(G, M)$ конечна
по упражнению~\ref{exer-cohomfingrp}(ii).
\item[(iv)] Что изменится, если в пункте~(iii) не предполагать,
что модуль $M$ свободен как абелева группа?
(Указание: в этом случае утверждение пункта~(iii) не выполняется
даже для конечной группы $G$.)
\end{itemize}
\end{prob}

\bigskip
\begin{prob}{\bf Когомологии проконечных
групп с коэффициентами в пермутационном модуле}
\label{prob:stab-perm-profinite}

Пусть $G$ --- проконечная группа, $M$ --- дискретный $G$-модуль.
\begin{itemize}
\item[(i)]
Пусть $M$ --- пермутационный $G$-модуль,
$B$ --- базис в $M$, на котором
группа~$G$ действует перестановками,
и $\mathbb{O}$ --- множество $G$-орбит в $B$.
Докажите, что $H^1(G,M)=0$ и
$$
H^2(G,M)\cong\mbox{$\bigoplus\limits_{O\in
\mathbb{O}}\Hom\big(\mathrm{Stab}_{G}(x_O),\Q/\Z\big)$}\,,
$$
где $x_O$ ---
(какой-нибудь) элемент в орбите $O$, подгруппа $\mathrm{Stab}_{G}(x_O)$ ---
стабилизатор элемента $x_O$ в группе $G$,
а $G/\mathrm{Stab}_G(x_O)$ обозначает множество классов смежности
по этой подгруппе.
(Указание: воспользуйтесь
упражнением~\ref{exercise:permutation-module}(ii).)
\item[(ii)]
Пусть $M$ --- стабильно пермутационный
$G$-модуль. Докажите, что
$$H^1(G,M)=0\,.$$
(Указание: воспользуйтесь
упражнением~\ref{exercise:stab-permutation-module}(i).)
\end{itemize}
\end{prob}

\subsection{Когомологии абсолютной группы Галуа}

В дальнейшем мы продолжим обозначать
абсолютную группу Галуа поля $K$ через~$G_K$,
то есть $G_K=\Gal(K^{sep}/K)$, где $K^{sep}$~--- сепарабельное
замыкание поля $K$.

\begin{prob}{\bf Первые когомологии и формы}
\begin{itemize}
\item[(i)]
Пусть $\Mcal$~--- категория, расслоенная над полями, удовлетворяющая условию спуска, $X$~---
объект из $\Mcal(K)$. Покажите, что $\Aut(X_{K^{sep}})$
является дискретным $G_K$-модулем.
\item[(ii)]
Покажите, что имеется каноническая биекция
$$\Phi(X)\cong H^1\big(G_K,\Aut(X_{K^{sep}})\big)\,.$$
(Указание: воспользуйтесь
упражнениями~\ref{exer-spusk}(iii)
и~\ref{exer-changefield}.)
\end{itemize}
\end{prob}

\bigskip

\begin{prob}{\bf Теория Куммера}\label{exer-Kummer}

Пусть $n$ --- натуральное число, взаимно простое с характеристикой поля $K$.
\begin{itemize}
\item[(o)]
Покажите, что
$$H^1\big(G_K,\GL_d(K^{sep})\big)=\{1\}
$$
для любого $d\geqslant 1$.
(Указание: воспользуйтесь теоремой Гильберта~90 и перейдите к пределу.)
В~частности,
$H^1\big(G_K,(K^{sep})^*\big)=\{1\}$.
\item[(i)]
Покажите, что существует короткая точная последовательность
дискретных $G_K$-модулей
$$
1\to \mu_n\longrightarrow (K^{sep})^*\stackrel{n}\longrightarrow (K^{sep})^*\to 1\,,
$$
где $\mu_n$ обозначает группу корней $n$-ой степени
из единицы в поле $K^{sep}$, а правый гомоморфизм является
возведением в $n$-ую степень.
\item[(ii)]
Предположим, что поле $K$ содержит все элементы из $\mu_n$.
Постройте изоморфизм абелевых групп
$$\Hom(G_K,\mu_n)\cong K^*/(K^*)^n\,.$$
(Указание: воспользуйтесь длинной точной последовательностью
когомологий, связанной с точной последовательностью $G_K$-модулей
из пункта~(i).)
\item[(iii)]
В предположении из пункта (ii) докажите, что любое конечное
расширение Галуа~\mbox{$K\subset L$}
с циклической группой Галуа, порядок которой делит $n$,
имеет вид
$$L=K\big(\sqrt[n]{a}\big)$$
для некоторого элемента $a\in K^*$.
(Указание: по каждому элементу $a\in K^*$ можно построить отображение
$$G_K\to\mu_n, \quad g\mapsto g\big(\sqrt[n]{a}\big)/\sqrt[n]{a}\,.$$
Сравните это отображение
с образом элемента
$$a\in K^*\cong H^0\big(G_K,(K^{sep})^*\big)$$ при
кограничном отображении,
связанном с точной последовательностью из пункта (i),
см. упражнение~\ref{exer-cobound-map}(i).)
\end{itemize}
\end{prob}

\bigskip

\begin{prob}{\bf Теория Артина--Шрейера}

Пусть $K$ --- поле характеристики $p>0$.
\begin{itemize}
\item[(o)]
Покажите, что
$$H^i(G_K,K^{sep})=0$$
при $i>0$.
(Указание: воспользуйтесь упражнением~\ref{exer-cohomadd}(iii)
и перейдите к пределу.)
Отметим, что данное утверждение выполнено также и для поля нулевой
характеристики.
\item[(i)]
Покажите, что существует короткая точная последовательность дискретных $G_K$-модулей
$$
0\to \Z/p\Z\longrightarrow K^{sep}\stackrel{\mathrm{Fr}-1}\longrightarrow K^{sep}\to 0\,,
$$
где $\mathrm{Fr}$ обозначает морфизм Фробениуса $\mathrm{Fr}(x)=x^p$.
\item[(ii)]
Постройте изоморфизм абелевых групп
$$\Hom(G_K,\Z/p\Z)\cong K/(\mathrm{Fr}-1)K\,.$$
(Указание: воспользуйтесь длинной точной последовательностью
когомологий, связанной с точной последовательностью $G_K$-модулей
из пункта~(i).)
\item[(iii)]
Докажите, что любое конечное расширение Галуа
$K\subset L$ с группой Галуа $\Z/p\Z$
имеет вид~\mbox{$L=K(x)$}, где элемент
$x\in K^{sep}$ удовлетворяет уравнению
$$x^p-x=a$$
для некоторого элемента $a\in K$.
(Указание: по каждому элементу $a\in K$ можно
построить отображение
$$G_K\to\Z/p\Z, \quad g\mapsto g(x)-x\,,$$
где $x\in K^{sep}$ удовлетворяет уравнению $x^p-x=a$.
Сравните это отображение
с образом элемента
$$a\in K\cong H^0(G_K, K^{sep})$$
при кограничном отображении,
связанном с точной последовательностью из пункта (i),
см. упражнение~\ref{exer-cobound-map}(i).)
\end{itemize}
\end{prob}

\subsection{Группа Пикара как стабильно пермутационный модуль}

Теперь мы приведём возникающие из алгебраической геометрии примеры дискретных
стабильно пермутационных модулей Галуа. Пусть $X$ --- гладкое алгебраическое
многообразие над полем $K$. Группа дивизоров $\Div(X_{K^{sep}})$
и группа Пикара~$\Pic(X_{K^{sep}})$ являются дискретными
$G_K$-модулями.

\bigskip
\begin{prob}{\bf Когомологии группы дивизоров}\label{divisorscohom}
\hspace{0cm}
\begin{itemize}
\item[(i)]
Покажите, что группа $\Div(X_{K^{sep}})$
является пермутационным $G_K$-модулем.
\item[(ii)]
Пусть дано (возможно, бесконечное) расширение полей $K\subset F$,
и $K_F$ является сепарабельным замыканием поля $K$ в поле $F$.
Предположим, что степень~\mbox{$n=[K_F:K]$} конечна.
Докажите, что тогда имеется изоморфизм \mbox{$K^{sep}$-алгебр}
$$
K^{sep}\otimes_K F\cong\prod_{i=1}^n F_i\,,
$$
где $F_i$ являются полями, содержащими $K^{sep}$. При этом группа Галуа $G_K$ действует транзитивно на множестве
полей $\{F_i\}$. Пусть $F_{i_0}$ является одним из этих полей.
Покажите, что стабилизатор поля $F_{i_0}$ соответствует подполю в~$K^{sep}$, совпадающему с образом поля $K_F$ относительно композиции
$$
K_F\subset F\subset F_{i_0}\,,
$$
где второе вложение задаётся проекцией $K^{sep}\otimes_K F\to F_{i_0}$.
Не забудьте проверить, что данный образ содержится в $K^{sep}\subset F_{i_0}$.
(Указание: воспользуйтесь изоморфизмом
$$
K^{sep}\otimes_K F\cong (K^{sep}\otimes_K K_F)\otimes_{K_F}F\,,
$$
изоморфизмом
$$
K^{sep}\otimes_K K_F\cong \prod_{i=1}^n K^{sep}
$$
а также тем, что $K^{sep}\otimes_{K_F}F$ является областью целостности,
так как $K_F$ сепарабельно замкнуто в $F$, см.~\cite[Theorem~IV.21(2)]{Jac}.)
\item[(iii)]
Докажите, что имеются канонические изоморфизмы
$$
H^0\big(G_K,\Div(X_{K^{sep}})\big)\cong\Div(X),\quad H^1\big(G_K,\Div(X_{K^{sep}})\big)=0\,,
$$
а также
$$
H^2\big(G_K,\Div(X_{K^{sep}})\big)\cong\mbox{$\bigoplus\limits_{D\subset X}
\Hom\big(G_{K_D},\Q/\Z\big)$}\,,
$$
где $D$ пробегает все неприводимые приведённые дивизоры на $X$,
а $K_D$ является сепарабельным замыканием поля $K$
в поле $K(D)$ рациональных функций на $D$.
(Указание: выбор одной из неприводимых компонент $E$
дивизора~\mbox{$D_{K^{sep}}\subset X_{K^{sep}}$}
над $K^{sep}$
фиксирует вложение
$$K(D)\subset K^{sep}(E)\,,$$
которое определяет также вложение $K_D\subset K^{sep}$.
Из пункта~(ii) следует, что возникающая подгруппа Галуа
$G_{K_D}\subset G_{K}$ является стабилизатором компоненты $E$.
Далее воспользуйтесь пунктом~(i) и упражнением~\ref{exer-profcohom}.)
\end{itemize}
\end{prob}

\bigskip

Предположим дальше, что поле $K$ имеет нулевую характеристику.
В частности, имеется равенство $K^{sep}=\bar{K}$.

\begin{prob}{\bf Стабильно пермутационные модули и группа Пикара}
\label{exercise:stab-permutation-module-from-Pic} \hspace{0cm}

\begin{itemize}
\item[(i)] Пусть $\pi\colon\widetilde{X}\to X$ --- собственный бирациональный
морфизм гладких многообразий над~$K$. Покажите, что отображения
прямого и обратного образа,
$$\pi_*\colon\Pic(\widetilde{X}_{\bar{K}})\to
\Pic(X_{\bar{K}}), \quad
\pi^*\colon\Pic(X_{\bar{K}})
\to\Pic(\widetilde{X}_{\bar{K}})\,,$$
определяют изоморфизм $G_K$-модулей
$$\Pic(\widetilde{X}_{\bar{K}})\cong \Pic(X_{\bar{K}})\oplus N\,,$$
где $N=\Ker(\pi_*)$.
(Указание:
воспользуйтесь тем, что
$$
\pi_*\circ\pi^*\colon \Pic(X_{\bar{K}})\to\Pic(X_{\bar{K}})
$$
является тождественным отображением.
Это следует из того, что множество точек неопределённости рационального отображения
$$\pi^{-1}\colon X\dasharrow\widetilde{X}$$
имеет коразмерность не меньше двух.)
\item[(ii)]
Пусть $\Sigma$ обозначает множество исключительных
дивизоров морфизма $\pi$, то есть таких неприводимых дивизоров
$D\subset\widetilde{X}$, что
$$\dim\big(\pi(D)\big)<\dim(D)\,.$$
Покажите, что имеется канонический изоморфизм $G_K$-модулей
$$\Z[\Sigma]\cong N\,.$$
(Указание: сначала постройте сюръективный морфизм $G_K$-модулей
$$\Z[\Sigma]\to N\,,$$
а потом докажите его инъективность, снова пользуясь тем,
что рациональное отображение~$\pi^{-1}$ не определено в коразмерности
не меньше двух.)
\item[(iii)] Используя пункт~(ii), покажите, что $N$ --- пермутационный
$G_K$-модуль.
\item[(iv)] Предположим, что многообразие $X$
рационально над $K$ и проективно.
Докажите, что $G_K$-модуль $\Pic(X_{\bar{K}})$ является
стабильно пермутационным.
(Указание: по теореме Хиронаки существуют гладкое проективное
многообразие $\widetilde{X}$ и собственные бирациональные морфизмы
$\widetilde{X}\to\Pb^n$ и $\widetilde{X}\to X$, где $n=\dim(X)$.
После этого можно
воспользоваться пунктами~(i) и~(iii).)
\item[(v)] Предположим, что многообразие
$X$ стабильно рационально над $K$ (см. определение~\ref{defin:stably-rational}
ниже) и проективно.
Докажите, что $G_K$-модуль $\Pic(X_{\bar{K}})$ является
стабильно пермутационным.
(Указание: воспользуйтесь изоморфизмом $G_K$-модулей
$$\Pic\big((X\times\P^n)_{\bar{K}}\big)\cong\Pic\big(X_{\bar{K}}\big)\oplus\Z$$
и пунктом~(iv).)
\end{itemize}
\end{prob}

\subsection{Торсоры}
\label{subsection:torsory-nad-tochkoj}

Пусть $U$ --- многообразие, определённое над полем $K$, а $\Gamma$ ---
(приведённая) алгебраическая группа над $K$. \emph{Действием} группы
$\Gamma$ на многообразии $U$ называется морфизм
$$m\colon \Gamma\times U\to U\,,$$
удовлетворяющий обычному условию ассоциативности действия.
Эквивалентно, для любой $K$-схемы $S$
морфизм $m$ определяет действие группы $\Gamma(S)$ на множестве $U(S)$.

\begin{defin}\label{defin-torsorpoint}
Действие $\Gamma$ на $U$ называется \emph{свободным
и транзитивным}, если $U\neq\varnothing$ (то есть у $U$ имеется
точка над $\bar{K}$),
и морфизм
$$(m,\mathrm{pr}_U)\colon \Gamma\times U\to U\times U$$
является изоморфизмом. Эквивалентно, для
любой $K$-схемы $S$ и точки $x\in\nolinebreak U(S)$ отображение
$$g\mapsto m(g,x)$$
определяет биекцию между $\Gamma(S)$ и
$U(S)$. Отметим, что последнее свойство достаточно
проверять для $S=\mathrm{Spec}(K^{sep})$.
Если действие~$\Gamma$ на~$U$ свободно и транзитивно,
то $U$ называется \emph{торс\'oром} над $\Gamma$
(а также \emph{$\Gamma$-торсором}, или \emph{главным однородным
пространством} над $\Gamma$). Если $U(K)\neq\varnothing$,
то $U$ и~$\Gamma$ изоморфны как многообразия над~$K$, и $U$ называется
\emph{тривиальным} $\Gamma$-торсором.
Понятие изоморфизма для торсоров
определяется естественным образом.
\end{defin}

\bigskip
\begin{prob}{\bf Торсоры и первые когомологии Галуа}
\label{prob:torsor-vs-H1}
\begin{itemize}
\item[(i)]
Докажите, что множество $\Gamma$-торсоров над полем $K$ с точностью до
изоморфизма находится во взаимно однозначном соответствии с
множеством~\mbox{$H^1\big(G_K,\Gamma(K^{sep})\big)$}. (Указание:
воспользуйтесь тем, что произвольный торсор является формой
тривиального торсора. Также покажите, что группа автоморфизмов
тривиального торсора над $K^{sep}$ изоморфна группе $\Gamma(K^{sep})$, и
воспользуйтесь упражнением~\ref{exer-spusk}(iii).)
\item[(ii)]
Пусть $\Gamma$ является одной из групп
$\GL_n$, $\SL_n$, $\Aff_n$, или $\Sp_{2n}$
(в последнем случае мы также дополнительно предположим, что
$\mathrm{char}(K)\neq 2$).
Докажите, что все $\Gamma$-торсоры тривиальны. (Указание:
воспользуйтесь пунктом~(i), а также теоремой Гильберта 90
и упражнением~\ref{exer-1st-Galois}.)
\end{itemize}
\end{prob}

Простым примером торсора над группой
$\mu_2\cong {\mathbb Z}/2{\mathbb Z}$
является подмногообразие в ${\mathbb A}^1$, заданное
уравнением
$$\{x^2=a\}\,,$$
где $a\in K^*$. Из определения~\ref{defin-torsorpoint}
легко увидеть, что данный торсор
тривиален тогда и только тогда, когда $a\in (K^*)^2$.
Согласно упражнению~\ref{prob:torsor-vs-H1}(i), это соответствует
изоморфизмам
$$H^1(G_K,\mu_2)\cong {\rm Hom}(G_K,\mu_2)\cong K^*/(K^*)^2\,,$$
см. упражнение~\ref{exer-Kummer}(ii).

\bigskip
\begin{defin}
\label{definition:torsorbase}
Пусть дан
морфизм $\phi\colon V\to X$ многообразий, определённых над полем $k$.
Предположим, что $\phi$ является плоским (см.~\cite[III.9]{Hartshorne})
и сюръективным (на множестве $\bar{k}$-точек).
Предположим, что (приведённая) алгебраическая
группа $\Gamma$ действует
на многообразии $V$, причём действие~$\Gamma$ сохраняет слои
морфизма $\phi$, и является свободным и транзитивным на этих
слоях. В этом случае $V$ называется \emph{$\Gamma$-торсором над
$X$}. (Иногда также употребляют и словосочетание
\emph{торсор над $\Gamma$}, если нет нужды упоминать базу и
если это не ведёт к путанице.) Понятие изоморфизма для торсоров над $X$
определяется естественным образом.
\end{defin}

\bigskip
\begin{prob}{\bf Бирациональные свойства торсоров}
\label{prob:torsor-birational}

Пусть дан $\Gamma$-торсор $\phi\colon V\to X$, где $X$ --- неприводимое
многообразие над полем~$k$. Рассмотрим соответствующий
$\Gamma$-торсор $U$
над полем функций $K=k(X)$.

\begin{itemize}
\item[(o)] Проверьте, что торсор $U$
тривиален тогда и только тогда, когда у морфизма~$\phi$ есть рациональное
сечение.
\item[(i)]
Пользуясь пунктом~(o), проверьте, что если торсор $U$
тривиален, то многообразия $V$ и
$X\times\Gamma$ бирационально эквивалентны (для ясности будем считать,
что~$\Gamma$ геометрически неприводима; из этого следует неприводимость
многообразия~\mbox{$X\times\Gamma$}).
\item[(ii)]
Предположим, что $H^1\big(G_K, \Gamma(K^{sep})\big)=\{1\}$
(например, $\Gamma$ является одной из групп,
перечисленных в
упражнении~\ref{prob:torsor-vs-H1}(ii)).
Докажите, что многообразия $V$
и $X\times\Gamma$ бирационально эквивалентны.
(Указание: воспользуйтесь пунктом~(i) и
упражнением~\ref{prob:torsor-vs-H1}(i).)
\end{itemize}
\end{prob}

\subsection{Когомологии обратного предела}

Пусть задана проективная система комплексов абелевых групп $\{C^\udot_i\}$,
где $i$ пробегает множество натуральных чисел с естественным порядком
(см. раздел~\ref{subsection:cohom-profinite}).
Таким образом, для всех пар натуральных чисел~\mbox{$i\leqslant j$}
заданы морфизмы комплексов
$$\varphi_{ji}\colon C^{\udot}_j\to C^{\udot}_i\,,$$
удовлетворяющие условиям согласованности.
Тогда для всех целых $p$ определены абелевы группы
$$C^p=\varprojlim C^p_i\,,$$
которые образуют комплекс
$$
C^{\udot}=\varprojlim C^{\udot}_i\,.
$$
Для каждого целого $p$ имеется каноническое отображение
$$
H^p(C^{\udot})\to \varprojlim H^p(C^{\udot}_i)\,.
$$

\begin{prob}{\bf Когомологии обратного предела комплексов}\label{invlimcomplex}
\hspace{0cm}
\begin{itemize}
\item[(i)]
Пусть дан сюръективный морфизм комплексов
$$\varphi\colon E^{\udot}\to D^{\udot}\,.$$
Предположим, что для коцикла
$\alpha\in D^p$ существует такой класс
$c\in H^p(E^{\udot})$, что выполняется равенство
$\varphi(c)=[\alpha]$ в $H^p(D^{\udot})$.
Покажите, что тогда существует такой коцикл $\beta\in E^p$,
что его класс в $H^p(E^{\udot})$ равен $c$,
и выполняется равенство $\varphi(\beta)=\alpha$ в $D^p$.
(Указание: рассмотрите сначала произвольный коцикл~\mbox{$\beta'\in E^{p}$},
для которого $[\beta']=c$. Для такого $\beta'$ существует
элемент~\mbox{$\zeta\in D^{p-1}$}, удовлетворяющий условию
$$\alpha=\varphi(\beta')+d(\zeta)\,.$$
Положите
$$\beta=\beta'+d(\xi)\,,$$
где $\xi\in E^{p-1}$ --- такой элемент, что $\varphi(\xi)=\zeta$.)
\item[(ii)]
Предположим, что все гомоморфизмы $\varphi_{ji}$ сюръективны.
Докажите, что тогда для каждого целого $p$ каноническое отображение
$$H^p(C^{\udot})\to \varprojlim H^p(C^{\udot}_i)$$
является изоморфизмом, то есть в этом случае взятие когомологий
коммутирует с обратным пределом. (Указание: примените пункт~(i).)
\end{itemize}
\end{prob}

\bigskip

\begin{prob}{\bf Когомологии обратного предела модулей}\label{invlimmods}
\hspace{0cm}

Пусть даны произвольная группа $G$ и $G$-модуль $M$.
Пусть на $M$ также задана убывающая
фильтрация $G$-подмодулями
$$
M=M^0\supset M^1\supset\ldots\supset M^i\supset M^{i+1}\supset\ldots
$$
Предположим, что естественное отображение
$$M\to\varprojlim M/M^i$$
является изоморфизмом.
\begin{itemize}
\item[(i)]
Покажите, что естественное отображение
$$H^p(G,M)\to\varprojlim H^p(G,M/M^i)$$
является изоморфизмом для любого $p\geqslant 0$.
(Указание: примените упражнение~\ref{invlimcomplex}(ii)
к обратному пределу комплексов
$$
C_i^{\udot}=\Map(G^{\times\udot},M/M^i)\,,
$$
см. упражнение~\ref{exer-stcoml}.)
\item[(ii)]
Предположим, что для некоторого $p\geqslant 0$
при всех $i\geqslant 1$ выполняются условия
$$
H^p(G,M^i/M^{i+1})=
H^{p+1}(G,M^i/M^{i+1})=0\,.
$$
Положим $\overline{M}=M/M^1$.
Докажите, что тогда естественное отображение
$$H^p(G,M)\to H^p(G,\overline{M})$$
является изоморфизмом. (Указание: докажите по индукции,
что для любого~\mbox{$i\geqslant 1$}
отображение
$$H^p(G,M/M^i)\to H^p(G,\overline{M})$$
является изоморфизмом, а потом воспользуйтесь пунктом~(i).)
\end{itemize}
\end{prob}

\bigskip

\begin{prob}{\textbf{Когомологии обратного предела неабелевых модулей}}
\label{kogominvlimit}

Пусть $G$ действует автоморфизмами
на (возможно, неабелевой) группе $\Gamma$.
Предположим, что на $\Gamma$ задана убывающая фильтрация
нормальными $G$-инвариантными подгруппами
$$
\Gamma\supset\Gamma^1\supset\ldots\supset\Gamma^i\supset\Gamma^{i+1}\supset
\ldots
$$
Предположим также, что естественное отображение
$$\Gamma\to\varprojlim\Gamma/\Gamma^i$$
является изоморфизмом.
\begin{itemize}
\item[(i)]
Покажите, что естественное отображение
$$H^1(G,\Gamma)\to\varprojlim H^1(G,\Gamma/\Gamma^i)$$
является биекцией. (Указание: модифицируйте рассуждение
из упражнеий~\ref{invlimmods}(i) и~\ref{invlimcomplex}.)
\item[(ii)]
В условиях пункта~(i) предположим также, что для любого $i\geqslant 1$
факторгруппа~$\Gamma^i/\Gamma^{i+1}$ абелева
и выполняются условия
$$
H^1(G,\Gamma^i/\Gamma^{i+1})=H^2(G,\Gamma^i/\Gamma^{i+1})=0\,.
$$
Положим $\overline{\Gamma}=\Gamma/\Gamma^1$.
Докажите, что тогда естественное отображение
$$H^1(G,\Gamma)\to H^1(G,\overline{\Gamma})$$
является биекцией. (Указание: действуйте так же,
как в упражнении~\ref{invlimmods}(ii).)
\end{itemize}
\end{prob}

\newpage

\section{Группа Брауэра I}
\label{section:Brauer}

\subsection{Определение и общие свойства}

Напомним некоторые понятия и факты, связанные с ассоциативными алгебрами
(подробности можно найти в
главе~X книги~\cite{Serre-corps-locaux}, в главе~VIII
книги~\cite{Bourbaki-algebra} или в главе 2 книги~\cite{GS}).
Пусть $K$~--- поле, $A$~---
ассоциативная алгебра с единицей над~$K$. Алгебра~$A$ называется
{\it центральной над $K$}, если центр $A$ совпадает с $K$.
Алгебра~$A$ называется {\it простой}, если в $A$ нет двусторонних идеалов,
кроме нулевого идеала и самой~$A$.
Класс конечномерных центральных простых алгебр замкнут
относительно тензорного произведения
(см.~\cite[VIII.1.2, VIII.7.4]{Bourbaki-algebra}).
Пусть $\Mat_n(K)$ обозначает алгебру
матриц размера $n\times n$ с коэффициентами из поля $K$.
Заметим, что~$\Mat_n(K)$ является
центральной простой алгеброй над $K$.

Пусть $A$ --- конечномерная центральная простая алгебра над $K$.
\emph{Теорема Веддерберна} (см.~\cite[VIII.5.4]{Bourbaki-algebra})
утверждает, что тогда существует такое
конечномерное центральное тело~$D$ над~$K$, что $A\cong\Mat_m(D)$
для некоторого натурального числа~$m$. Отметим, что
$$\Mat_m(D)\cong D\otimes_K\Mat_m(K)\,.$$
Кроме того, если для
конечномерных тел~$D_1$ и~$D_2$ над~$K$ и натурального числа $m$
имеет место изоморфизм
$$\Mat_m(D_1)\cong \Mat_m(D_2)\,,$$ то $D_1\cong
D_2$.

Над сепарабельно замкнутым полем не существует
центральных конечномерных тел, кроме самого поля
(см.~\cite[VIII.10.5]{Bourbaki-algebra}).
Следовательно, над сепарабельно замкнутым полем любая
конечномерная центральная простая алгебра изоморфна алгебре матриц.
С другой стороны, свойства конечномерных алгебр быть центральными и
простыми сохраняются при расширении скаляров
(см.~\cite[VIII.7.4]{Bourbaki-algebra}).
Поэтому
размерность конечномерной центральной простой алгебры над произвольным
полем является квадратом.

Будем говорить, что конечномерные центральные простые алгебры $A$ и $B$
эквивалентны (и будем писать в этом случае $A\sim B$),
если существует
конечномерное центральное тело $D$ над $K$ и натуральные числа $p$ и
$q$, для которых $A\cong \Mat_p(D)$ и $B\cong \Mat_q(D)$.
Другими словами, в качестве отношения эквивалентности
мы берём симметричное транзитивное замыкание отношения
$$A\sim A\otimes_K\Mat_m(K)\,.$$
В каждом классе эквивалентности можно выбрать единственного
представителя, являющегося центральным телом. Также, если $A\sim B$
и $\dim_K(A)=\dim_K(B)$, то~\mbox{$A\cong B$}. Через~\mbox{$\Br(K)$}
обозначим
множество классов эквивалентности конечномерных центральных простых
алгебр над $K$.
Для алгебры $A$ обозначим через $[A]$ её класс в~$\Br(K)$.
Рассмотрим структуру абелевой полугруппы на
$\Br(K)$, заданную тензорным произведением алгебр.
Отметим, что в этой полугруппе $[A]=0$ тогда и только тогда,
когда~$A$ является алгеброй матриц над полем $K$.

Для произвольной алгебры $A$ пусть $A^{op}$ обозначает то же
векторное пространство~$A$ с умножением, заданным в обратном порядке
по сравнению с умножением в~$A$. Заметим, что имеется изоморфизм
алгебр
$$A\otimes_K A^{op}\cong \End_K(A)\cong\Mat_n(K)\,,$$
где $n=\dim_K(A)$, задаваемый по формуле
$$a\otimes b\mapsto\{x\mapsto axb\}, \quad a,x\in A, b\in A^{op}\,.$$
Поэтому для любого элемента в
$\Br(K)$ существует обратный, и $\Br(K)$ является абелевой группой.
Эта группа~--- важнейший арифметический инвариант поля.

\begin{defin}
Группа $\Br(K)$ называется {\it группой Брауэра} поля $K$.
\end{defin}

Расширение полей $K\subset L$ определяет гомоморфизм групп Брауэра
$$\Br(K)\to\Br(L)\,,$$
задаваемый расширением скаляров для
алгебр $A\mapsto L\otimes_K A$. Положим
$$\Br(L/K)={\rm Ker}\big(\Br(K)\to \Br(L)\big)\,.$$
В частности, имеется равенство
$$
\Br(K)=\bigcup\limits_L \Br(L/K)\,,
$$
где $L$ пробегает все конечные расширения Галуа поля $K$, так как из сказанного выше следует, что $\Br(K^{sep})=0$.

Для центрального тела $D$ размерности $n^2$ над $K$ любое
максимальное под\-поле~\mbox{$L\subset D$} имеет размерность $n$ над $K$
(см.~\cite[VIII.10.3]{Bourbaki-algebra}).
Более того, соответствующий телу $D$ элемент
$[D]\in\Br(K)$ лежит в подгруппе
$$\Br(L/K)\subset\Br(K)$$
(см.~\cite[VIII.10.5]{Bourbaki-algebra}).

\bigskip

\begin{prob}{\bf Когомологическое определение группы
Брауэра}\label{exer-cohomBr}

Пусть $K\subset L$~--- конечное расширение Галуа с группой Галуа $G$.
\begin{itemize}
\item[(i)]
Постройте биекцию между множеством $\Phi\big(\Mat_n(K),L\big)$ классов
изоморфизма \mbox{$L$-форм} алгебры матриц $\Mat_n(K)$ и множеством
$H^1\big(G,\PGL_n(L)\big)$. (Указание: воспользуйтесь тем, что все
автоморфизмы алгебры матриц являются внутренними, то есть сопряжениями
при помощи обратимых матриц.)
\item[(ii)]
Покажите, что естественное отображение
$$H^1\big(G,\PGL_n(L)\big)\to \Br(L/K)\,,$$
возникающее из пункта (i), инъективно,
и при этом
$$\Br(L/K)=\mbox{$\bigcup\limits_n H^1\big(G,\PGL_n(L)\big)$}\,.$$
\item[(iii)]
Постройте каноническое вложение групп
$$\lambda\colon\Br(L/K)\hookrightarrow H^2(G,L^*)\,.$$
(Указание: покажите,
что кограничное отображение, связанное с точной последовательностью
$G$-модулей
$$1\to L^*\to \GL_n(L)\to\PGL_n(L)\to 1$$
определяет
отображение множеств
$$\lambda\colon\Br(L/K)\to H^2(G,L^*)\,,$$
причём \mbox{$\lambda^{-1}(0)=\{0\}$}. Также заметьте, что взятие
тензорного произведения алгебр соответствует отображению
$$
H^1\big(G,\PGL_m(L)\big)\times H^1\big(G,\PGL_n(L)\big)\to H^1\big(G,\PGL_{mn}(L)\big)\,,
$$
индуцированному тензорным произведением матриц
$$\PGL_m(L)\times\PGL_n(L)\to \PGL_{mn}(L)\,.$$
Далее, для
доказательства того, что $\lambda$~--- гомоморфизм групп,
используйте явную формулу для кограничного отображения из
упражнения~\ref{exer-nonabcohom}(iii). Наконец, установите
инъективность гомоморфизма~$\lambda$.)
\item[(iv)]
Положим $d=[L:K]$. Докажите, что кограничное отображение
$$H^1\big(G,\PGL_d(L)\big)\to H^2(G,L^*)$$
сюръективно. (Указание: рассмотрите
векторное пространство $U$ над $L$ с базисом
$\{\sigma_g\mid g\in G\}$, и для
произвольного $2$-коцикла
$$\psi\colon G\times G\to L^*$$
рассмотрите отображение
$$\tilde{\varphi}\colon G\to \Aut(U)\cong
\GL_d(L), \quad \tilde{\varphi}(g)\colon \sigma_h\mapsto
\psi(g,h)\sigma_{gh}\,,$$
а также
соответствующий $1$-коцикл
$$\varphi\colon G\to \PGL_d(L)$$
(ср. с упражнением~\ref{exer-nonabcohom}(iii)). Более инвариантно, по
упражнению~\ref{exer-012}(ii) элемент \mbox{$\alpha\in H^2(G,L^*)$}
соответствует расширению групп
$$
1\to L^*\stackrel{\imath}\to \widetilde{G}\to G\to 1\,.
$$
Рассмотрим скрученную групповую алгебру $L*\widetilde{G}$:
как $L$-векторное пространство она изоморфна обычной
групповой алгебре $L[\widetilde{G}]$,
а умножение в $L*\widetilde{G}$ однозначно определяется
правилом
$$g\cdot x=g(x) x\cdot g,\quad g\in \widetilde{G},\quad x\in L\,,$$
где действие $\widetilde{G}$ на $L$ пропускается через
естественное действие $G$ на $L$. В частности,
$L*\widetilde{G}$ явлется алгеброй над $K$, и $L$ не лежит в её центре.
Далее рассмотрим факторалгебру
$$
A= L*\widetilde{G} /\langle x\cdot e-1\cdot \iota(x) \rangle_{x\in L^*}\,,
$$
где $e\in \widetilde{G}$ обозначает единицу в группе. Тогда $A$~---
центральная простая алгебра над $K$,
и имеется равенство $\lambda([A])=\alpha$.)
\item[(v)]
Покажите, что
$$\lambda\colon\Br(L/K)\to H^2(G,L^*)$$
является изоморфизмом.
\end{itemize}
\end{prob}

\bigskip

\begin{prob}{\bf Поведение при замене поля}\label{exer-changefieldBr}

Пусть дана композиция конечных расширений Галуа
$$K\subset L\subset F\,,$$
причём $F$ явдяется расширением Галуа поля $K$.
Положим
$$G=\Gal(L/K),\quad N=\Gal(F/L)\,,$$
и $H=\Gal(F/K)$.
\begin{itemize}
\item[(i)]
Покажите, что следующая диаграмма коммутативна:
$$
\begin{CD}
\Br(L/K)@>>>\Br(F/K)\\
@VVV@VVV\\
H^2(G,L^*)@>>>H^2(H,F^*)
\end{CD}
$$
Здесь верхняя горизонтальная стрелка определяется естественным
образом, вертикальные стрелки являются изоморфизмами из
упражнения~\ref{exer-cohomBr}(v), а нижняя горизонтальная стрелка
является отображением инфляции, то есть композицией отображения
обратного образа
$$H^2(G,L^*)\to H^2(H,L^*)\,,$$
определённого
гомоморфизмом групп $H\twoheadrightarrow G$, и естественного
отображения
$$H^2(H,L^*)\to H^2(H,F^*)\,.$$
(Указание: воспользуйтесь
упражнениями~\ref{exer-changefield}(i) и~\ref{exer-cohomBr}(ii).)
\item[(ii)]
Постройте канонический изоморфизм групп
$$\Br(K)\cong H^2\big(G_K,(K^{sep})^*\big)\,,$$
где, как и раньше, $G_K$ обозначает абсолютную группу Галуа поля $K$.
(Указание: воспользуйтесь пунктом (i).)
\item[(iii)]
Покажите, что следующая диаграмма коммутативна:
$$
\begin{CD}
\Br(F/K)@>>>\Br(F/L)\\
@VVV@VVV\\
H^2(H,F^*)@>>>H^2(N,F^*)
\end{CD}
$$
Здесь верхняя горизонтальная стрелка задаётся расширением скаляров для
центральных простых алгебр с $K$ на $L$, вертикальные стрелки являются
изоморфизмами из
упражнения~\ref{exer-cohomBr}(v), а нижняя горизонтальная стрелка
является отображением обратного образа, определённым гомоморфизмом
групп~\mbox{$N\hookrightarrow H$}. (Указание: воспользуйтесь
упражнениями~\ref{exer-changefield}(ii) и~\ref{exer-cohomBr}(ii).)
Обратите внимание, что данная коммутативная диаграмма имеет
место, даже если $K\subset L$ не является расширением Галуа.
\item[(iv)]
Пользуясь упражнением~\ref{exer-cohomBr}(v),
докажите, что имеет место точная последовательность
$$
0\to H^2(G,L^*)\to H^2(H, F^*)\to H^2(N,F^*)\,.
$$
Также выведите её существование из упражнения~\ref{exer-ind-sur}(iii),
используя теорему Гильберта~90.
\end{itemize}
\end{prob}

\bigskip

\begin{prob}{\bf Случай $K=\Rb$}\label{prop-realBr}
\nopagebreak
\hspace{0cm}
\begin{itemize}
\item[(o)]
Приведите пример нетривиального элемента в группе Брауэра
$\Br(\Rb)$.
\item[(i)]
Покажите, что
$$\Br(\Rb)\cong\Z/2\Z\,.$$
(Указание: воспользуйтесь
упражнением~\ref{exer-cyclgrp}.) Сравните это с теоремой Фробениуса
о строении конечномерных ассоциативных тел над $\Rb$.
\end{itemize}
\end{prob}

\bigskip

\begin{prob}{\bf Кватернионные алгебры}
\label{exer-quaternion-alg}

Пусть характеристика поля $K$ не равна $2$.
Рассмотрим ненулевые элементы~\mbox{$a,b\in K$}. Обозначим через $A(a,b)$
факторалгебру
$$K\{i,j\}/\langle i^2-a,j^2-b,ij+ji\rangle$$
свободной
ассоциативной алгебры $K\{i,j\}$. Эта алгебра называется {\it
кватернионной алгеброй}, соответствующей элементам $a$ и $b$.
\begin{itemize}
\item[(i)]
Докажите, что $A(a,b)$ имеет базис $(1,i,j,ij)$ над $K$ и является
центральной простой алгеброй над $K$.
\item[(ii)]
Покажите, что $A(a,b)\cong \Mat_2(K)$
тогда и только тогда, когда
уравнение
$$u^2-av^2-bw^2=0$$
имеет ненулевое решение в $K$.
(Указание: сначала покажите, что изоморфизм~\mbox{$A(a,b)\cong\Mat_2(K)$}
имеет место тогда и
только тогда, когда в $A(a,b)$ есть необратимые элементы. Потом
вспомните, как вычисляются обратные элементы для кватернионов.
Наконец, покажите, что уравнение
$$x^2-ay^2-bz^2+ab t^2=0$$
имеет
ненулевое решение тогда и только тогда, когда уравнение
$$u^2-av^2-bw^2=0$$
имеет ненулевое решение, поскольку оба уравнения
равносильны тому, что $b$ лежит в образе нормы Галуа для расширения полей
$K\subset K(\sqrt{a})$. Сравните это рассуждение
с тем, что происходит ниже в упражнении~\ref{prob-section}(ii).)
\item[(iii)]
Покажите, что если $a$ или $b$ является квадратом в поле $K$, то
$A(a,b)\cong\nolinebreak\Mat_2(K)$. Также покажите, что если
$a+b=1$, то $A(a,b)\cong\Mat_2(K)$. (Указание: используйте пункт
(ii).)
\item[(iv)]
Докажите, что имеется изоморфизм
$$A(a,b)\otimes_K A(a',b)\cong A(aa',b)\otimes_K\Mat_2(K)\,.$$
(Указание: напишите естественный базис
в алгебре $A(a,b)\otimes_K A(a',b)$, а затем рассмотрите некоторые
его части, порождающие подалгебры $A(a,b^2)$ и $A(aa',b)$.) Выведите из
этого и из пункта (iii), что имеется изоморфизм
$$A(a,b)\otimes_K A(a,b)\cong\Mat_4(K)\,.$$
\end{itemize}
\end{prob}

\bigskip

Приведённая ниже конструкция обобщает упражнение~\ref{exer-quaternion-alg}.

\begin{prob}{\bf Циклические алгебры}
\label{exer-cyclic-alg}

Пусть $K\subset L$ является конечным расширением Галуа с циклической группой
Галуа~$G$. Зафиксируем образующую~$s$ группы~$G$.
\begin{itemize}
\item[(o)]
Покажите, что выбор образующей $s$ определяет изоморфизм
$$\Br(L/K)\cong K^*/\Nm_{L/K}(L^*)\,,$$
где $\Nm_{L/K}$ обозначает норму Галуа.
(Указание: воспользуйтесь упражнениями~\ref{exer-cohomBr}(v)
и~\ref{exer-cyclgrp}(i).)
\item[(i)]
Покажите, что при изоморфизме из пункта~(o) класс элемента
$a\in K^*$
в группе~\mbox{$K^*/\Nm_{L/K}(L^*)$}
соответствует факторалгебре
$$A=L_G[\sigma]/\langle \sigma^n-a\rangle\,,$$
где
$L_G[\sigma]$ как $L$-векторное пространство совпадает
с пространством многочленов~$L[\sigma]$, а
закон умножения в $L_G[\sigma]$ определяется соотношением
$$\sigma\cdot x={}^s x\cdot\sigma,\quad x\in L\,.$$
(Указание: воспользуйтесь упражнением~\ref{exer-cyclgrp}(iv) и
указанием к упражнению~\ref{exer-cohomBr}(iv).)
Такие центральные простые алгебры называются \emph{циклическими}.
Отметим, что алгебра $A$ зависит не только от
циклического расширения~\mbox{$K\subset L$}
и элемента $a\in K^*$,
но также и от выбора образующей $s\in G$.
\item[(ii)]
Предположим, что характеристика поля $K$ взаимно проста с $n$,
и в $K$ содержится первообразный корень $n$-ой степени из единицы
$\zeta\in K^*$. Тогда по теории Куммера
(см. упражнение~\ref{exer-Kummer}) имеется равенство
$L=K\big(\sqrt[n]{b}\big)$ для некоторого элемента $b\in K^*$, и
изоморфизм
$$G\cong \langle\zeta\rangle\subset K^*\,.$$
Выберем $s=\zeta$
в качестве образующей группы $G$. Покажите,
что имеется изоморфизм алгебр
$$A\cong A_{\zeta}(a,b)\,,$$
где~$A$~--- алгебра из пункта~(i),
и
$$A_{\zeta}(a,b)=K\{u,v\}/\langle u^n-a, v^n-b, uv-\zeta vu\rangle\,.$$
\item[(iii)]
В предположениях из пункта~(ii) покажите, что
имеют место равенства элементов группы Брауэра $\Br(K)$:
\begin{align*}
&[A_{\zeta}(a,b)]=-[A_{\zeta}(b,a)]=[A_{\zeta^{-1}}(b,a)]\,,\\
&[A_{\zeta}(a,b)]+[A_{\zeta}(a',b)]=[A_{\zeta}(aa',b)]\,.
\end{align*}
Кроме того, $[A_{\zeta}(a,b)]=0$ в $\Br(K)$ тогда и только тогда,
когда $a$ лежит в образе нормы Галуа для расширения полей
$K\subset K\big(\sqrt[n]{b}\big)$.
Докажите, что
$$[A_{\zeta}(a,1-a)]=0\,,$$
а также
$$[A_{\zeta}(a,-a)]=0\,.$$
(Указание: первое равенство непосредственно следует из
пункта~(ii). Бимультипликативность и утверждение про норму Галуа
следуют из пункта~(o).
Для доказательства двух последних равенств
покажите, что $1-a$ и $-a$ являются нормами в расширении
$K\subset K\big(\sqrt[n]{a}\big)$.)
\end{itemize}
\end{prob}

\bigskip
\begin{defin}\label{defin:K-group}
\emph{$\KM$-группа Милнора}
$\KM_2(K)$ поля $K$
определяется как фактор группы $K^*\otimes_\Z K^*$ по подгруппе,
порождённой всеми элементами вида $a\otimes (1-a)$,
где~\mbox{$a\in K$, $a\ne 0,1$}.
Образ элемента $a\otimes b$ в группе $\KM_2(K)$ обычно обозначается
символом~$\{a,b\}$. Тождество
$$\{a,1-a\}=0$$
в группе $\KM_2(K)$
называется \emph{соотношением Стейнберга}.
Гомоморфизмы из группы~$\KM_2(K)$ в произвольную (абелеву)
группу обычно называются \emph{символами}.
\end{defin}

\begin{prob}{\bf Соотношения в группе $\KM_2$}\label{exer-a-a}

Докажите, что в группе $\KM_2(K)$ выполнены соотношения
$$\{a,-a\}=0, \quad \{a,b\}=-\{b,a\}\,.$$
(Указание:
для доказательства первого соотношения вычислите значение
выражения
$$a^{-1}\otimes (1-a^{-1})+a\otimes (1-a)$$
в группе $K^*\otimes_{\Z} K^*$.
Для доказательства второго соотношения рассмотрите
элемент
$$\{ab,-ab\}\in\KM_2(K)$$
и воспользуйтесь первым соотношением.)
\end{prob}

\begin{prob}{\bf Символ норменного вычета}\label{exer-norm-symbol}

Предположим, что характеристика поля $K$ взаимно проста с $n$,
и в $K$ содержится первообразный корень $n$-ой степени из единицы
$\zeta\in K^*$. Покажите, что корректно определён
гомоморфизм групп
$$\nu_{\zeta}\colon\KM_2(K)/n\to\Br(K)_n, \quad
\{a,b\}\mapsto[A_{\zeta}(a,b)]\,.$$
(Указание: воспользуйтесь упражнением~\ref{exer-cyclic-alg}(iii).)
Это отображение называется \emph{символом норменного
вычета}.
\end{prob}

\emph{Теорема Меркурьева--Суслина} (см.~\cite{MeSu} или~\cite{Suslin})
утверждает, что $\nu_{\zeta}$ ---
изоморфизм. В частности, сюръективность отображения $\nu_{\zeta}$
означает, что любой элемент порядка $n$ в группе Брауэра можно
представить в виде тензорного произведения циклических алгебр
размерности $n^2$.

\bigskip

\begin{prob}{\bf Кручение в группе Брауэра}\label{exer-brtor}
\hspace{0cm}\nopagebreak
\begin{itemize}
\item[(o)]
Покажите, что $\Br(K)$ является группой кручения. (Указание:
используйте упражнение~\ref{exer-profcohom}(iv).)
\item[(i)]
Покажите, что для конечного сепарабельного расширения полей
$K\subset L$ (не обязательно нормального) все элементы в группе
Брауэра $\Br(L/K)$ аннулируются умножением на
степень расширения $[L:K]$. (Указание:
используйте упражнение~\ref{exer-trace}(iii) для вложения проконечных
групп
$$G_L\hookrightarrow G_K$$
и $G_K$-модуля $(K^{sep})^*$, а также
упражнение~\ref{exer-changefieldBr}(iii).)
\item[(ii)]
Пусть элемент $\alpha\in\Br(K)$ представим центральной простой
алгеброй над $K$ размерности $n^2$, где $n$ взаимно просто с
характеристикой поля $K$. Покажите, что в группе $\Br(K)$ выполнено
равенство $n\alpha=0$. (Указание:
воспользуйтесь упражнением~\ref{exer-cohomBr}(i), точной
последовательностью $G_K$-модулей
$$1\to \mu_n\to
\SL_n(K^{sep})\to\PGL_n(K^{sep})\to 1\,,$$
где $\mu_n$ обозначает
группу корней $n$-ой степени из $1$ в поле $K^{sep}$, а также
упражнением~\ref{exer-1st-Galois}(i).)
\item[(iii)] Пусть $n$ по-прежнему взаимно просто с характеристикой
поля $K$. Покажите, что имеется изоморфизм
$$H^2(G_K,\mu_n)\stackrel{\sim}\longrightarrow \Br(K)_n\,.$$
(Указание:
воспользуйтесь последовательностью Куммера
из упражнения~\ref{exer-Kummer}(i).)
\end{itemize}
\end{prob}

На самом деле, можно показать (см.~\cite[Exercise~X.5.2]{Serre-corps-locaux}),
что утверждение из
упражнения~\ref{exer-brtor}(i) верно для произвольного конечного
расширения $L$ поля $K$ (возможно, несепарабельного). Утверждение из
упражнения~\ref{exer-brtor}(ii) также верно для произвольного $n$
(возможно, делящегося на характеристику поля $K$). Это следует из обобщения
упражнения~\ref{exer-brtor}(i) на случай произвольного конечного
расширения, из того, что~\mbox{$\alpha\in\Br(L/K)$}, где $L$ является
максимальным подполем в теле $D$, представляющем элемент $\alpha$, а
также из того, что $[L:K]=n$, где $n^2$ является размерностью $D$
над~$K$.

Отметим, что обратное к утверждению из
упражнения~\ref{exer-brtor}(ii) неверно: произвольный элемент
порядка $n$ в группе Брауэра не обязательно представим алгеброй
размерности $n^2$. Например, можно показать (см.~\cite[Example~1.5.7]{GS}),
что для поля
$K=\Cb(t_1,t_2,t_3,t_4)$ алгебра
$$A=A(t_1,t_2)\otimes_{K}A(t_3,t_4)\,,$$
определяющая элемент порядка~$2$
в $\Br(K)$, является телом размерности~$16$, а значит, не
представима в виде алгебры размерности~$4$.
Ввиду теоремы Меркурьева--Суслина и упражнения~\ref{prob-quatconic}(ii)
на это можно смотреть как на утверждение о том, что
образ элемента
$$\{t_1,t_2\}+\{t_3,t_4\}$$
в группе $\KM_2(K)/2$ не является образом никакого разложимого элемента
из~$\KM_2(K)$, то есть элемента вида $\{a,b\}$.

\begin{defin}\label{defin-index}
{\it Индексом} элемента $\alpha\in\Br(K)$
называется квадратный корень из размерности центрального тела над
$K$, представляющего $\alpha$ (заметим, что индекс является
натуральным числом).
\end{defin}

Из упражнения~\ref{exer-brtor} и из сказанного перед
определением~\ref{defin-index} следует, что порядок произвольного
элемента в группе Брауэра делит его индекс. Отметим, что индекс
элемента $\alpha\in\Br(K)$ равен минимальной степени такого сепарабельного
расширения $L$ над $K$, что $\alpha\in\Br(L/K)$,
см.~\cite[Corollary~4.5.9]{GS}.
Кроме того, можно доказать
(см.~\mbox{\cite[Exercise~X.5.3b]{Serre-corps-locaux}}),
что у порядка и
индекса одинаковое множество простых делителей. Гипотеза
Кольо-Телена
утверждает, что если $K$ является полем рациональных
функций $k(X)$ на неприводимом многообразии $X$ размерности~$d$ над
сепарабельно замкнутым полем~$k$, то для любого элемента
$\alpha\in \Br(K)$ его индекс делит число~$m^{d-1}$, где~$m$~---
порядок элемента~$\alpha$.
Это равносильно тому, что любой элемент порядка~$m$ в группе~$\Br(K)$
можно представить центральной простой алгеброй размерности~$m^{2d-2}$.
Для случая $d=1$ это означает, что $\Br(K)=0$.
Этот частный случай хорошо известен и ниже (в упражнении~\ref{prob:Tzen}) мы
приведём его доказательство.
Для случая $d=2$ гипотеза Кольо-Телена означает, что
любой элемент порядка~$m$ в группе~$\Br(K)$ можно представить центральным
телом размерности $m^2$, или, что равносильно, для любого
центрального тела размерности $n^2$ над $K$ порядок его класса в
группе Брауэра~$\Br(K)$ равен~$n$. Последнее утверждение доказано де
Йонгом (см.~\cite{deJong-index}).

\bigskip

\begin{prob}{\bf Приведённая норма}\label{exer-rednorm}

Пусть $A$ является центральной простой алгеброй размерности $n^2$ над полем~$K$.
\begin{itemize}
\item[(i)]
Покажите, что определитель
$$\det\colon\Mat_n(K^{sep})\to K^{sep}$$
инвариантен относительно автоморфизмов алгебры $\Mat_n(K^{sep})$.
(Указание: вспомните, что все автоморфизмы алгебры $\Mat_n(K^{sep})$
внутренние.)
\item[(ii)]
Рассмотрим изоморфизм
$$\chi\colon K^{sep}\otimes_K
A\stackrel{\sim}\longrightarrow \Mat_n(K^{sep})\,.$$
Докажите, что
композиция
$$\det\circ\,\chi\colon K^{sep}\otimes_K A\to K^{sep}$$
коммутирует с действием группы Галуа~$G_K$ и не зависит от выбора
$\chi$. (Указание: действия группы Галуа на
$K^{sep}\otimes_K A$ и на $\Mat_n(K^{sep})$, а также различные выборы~$\chi$\,, отличаются на автоморфизм
алгебры~\mbox{$\Mat_n(K^{sep})$}.)
\item[(iii)]
Докажите, что существует единственное отображение степени $n$ из $A$
в $K$ (то есть отображение, задающееся полиномом степени $n$ от
координат элементов из $A$ в произвольном базисе в $A$ над $K$),
которое после расширения скаляров с $K$ на $K^{sep}$ соответствует
определителю
$$\det\colon\Mat_n(K^{sep})\to K^{sep}\,.$$
Это отображение
называется {\it приведённой нормой} и обозначается через
$$\Nrd\colon A\to K\,.$$
\item[(iv)]
Покажите, что для ненулевых элементов $a,b\in K$ и элемента
$$\alpha=x+yi+zj+tij$$
кватернионной алгебры $A(a,b)$ имеется равенство
$$\Nrd(\alpha)=x^2-ay^2-bz^2+ab t^2\,.$$
\item[(v)]
Покажите, что $A$ является телом тогда и только тогда, когда
гиперповерхность степени $n$ в $\Pb(A)\cong\Pb^{n^2-1}$, заданная
уравнением $\Nrd(\alpha)=0$, не имеет точек над~$K$.
\end{itemize}
\end{prob}

\subsection{Группа Брауэра и арифметические свойства полей}

\begin{defin}
Поле $K$ {\it имеет тип~$C_1$}, если любая гиперповерхность
$X_d\subset \Pb^{n}$ степени $d$, определённая над $K$, имеет
$K$-точку при $d\leqslant n$.
\end{defin}

Не вдаваясь в детали, связанные с особенностями, можно
сказать, что условие~$C_1$ означает, что любая гиперповерхность Фано
имеет $K$-точку.

\bigskip

\begin{prob}{\bf Группа Брауэра и свойство $C_1$}\label{exer-C1Br}

Докажите, что $\Br(K)=0$ для любого поля $K$ типа $C_1$. (Указание:
используйте упражнение~\ref{exer-rednorm}(v).)
\end{prob}

\bigskip

\begin{prob}{\bf Теорема Тзена}\label{prob:Tzen}
\hspace{0cm}
\begin{itemize}
\item[(i)]
Пусть $k$ --- алгебраически замкнутое поле, $t$ --- формальная
переменная. Докажите, что поле $K=k(t)$ имеет тип $C_1$. (Указание:
для однородного многочлена $f(x_0,\ldots,x_n)$ степени $d$ с
коэффициентами в $K$ можно считать, что его коэффициенты являются
многочленами от $t$ степени не выше $c$ для некоторого натурального
числа $c$. Для каждого натурального числа $r$ будем искать решение
уравнения $f=0$ в многочленах от $t$ степени меньше $r$ в виде
$$x_i=\sum\limits_{j=0}^{r-1} y_{ij}t^j\,.$$
Покажите, что тогда уравнение
$f=0$ над $K$ равносильно системе из не более, чем $c+dr$ однородных
уравнений от $(n+1)r$ переменных $y_{ij}$ над $k$, и рассмотрите
достаточно большое $r$.)
\item[(ii)]
Докажите, что любое конечное расширение $L\supset K$ поля $K$ типа
$C_1$ также имеет тип $C_1$. (Указание: выберите базис
$(e_1,\ldots,e_r)$ в $L$ над $K$ и для однородного многочлена
$f(x_0,\ldots,x_n)$ степени $d$ с коэффициентами в $L$ рассмотрите
многочлен
$$g(x_{01},\ldots,x_{0r},\ldots,
x_{n1},\ldots,x_{nr})={\rm Nm}_{L/K}\big(f(x_0,\ldots,x_n)\big)$$
степени
$dr$ от $(n+1)r$ переменных
с коэффициентами в $K$, где
$$x_i=\sum_{j=1}^r x_{ij}e_j\,,$$
а ${\rm Nm}_{L/K}\colon L\to K$ обозначает норму для расширения полей
$K\subset L$. Наконец, воспользуйтесь тем, что уравнение $f=0$ имеет
ненулевое решение над $L$ тогда и только тогда, когда уравнение
$g=0$ имеет ненулевое решение над $K$.)
\item[(iii)]
Докажите \emph{теорему Тзена}: для алгебраически замкнутого поля $k$
поле $K$ степени трансцендентности $1$ над $k$ имеет тип $C_1$. В
частности, $\Br(K)=0$ (см.~упражнение~\ref{exer-C1Br}).
\end{itemize}
\end{prob}

\bigskip

\begin{prob}{\bf Теорема Шевалле--Варнинга}\label{exer-Chevalley}

Пусть $K=\F_q$ является конечным полем из $q=p^s$ элементов.
\begin{itemize}
\item[(i)]
Покажите, что в поле $K$ выполнено равенство
$$\sum\limits_{x\in K}x^a=0$$
для любого целого $0\leqslant a< q-1$.
(Указание: при $a=0$ воспользуйтесь тем, что $q=0$ в поле $K$,
не забыв при этом, что $0^0=1$ по определению. При~\mbox{$a>0$}
выберите $t\in K^*$, для которого $t^a\neq 1$, и рассмотрите
сумму~\mbox{$\sum\limits_{x\in K}(tx)^a$}.\,)
\item[(ii)]
Докажите, что
$$\sum\limits_{x\in K^{n+1}}P(x)=0$$
для любого многочлена $P$ от $n+1$ переменной
над $K$, степень которого строго меньше,
чем $(n+1)(q-1)$. (Указание: если $P$ --- моном, то используйте пункт~(i).)
\item[(iii)]
Пусть $X\subset \Pb^{n}$ --- гиперповерхность степени $d$,
определённая над $K$, и~$d\leqslant n$. Докажите
\emph{теорему Шевалле--Варнинга}:
имеет место сравнение
$$|X(K)|\equiv 1\pmod{p}\,.$$
(Указание: воспользуйтесь тем, что
$$
|X(K)|(q-1)+1\equiv \sum_{x\in K^{n+1}}P(x)\pmod{p}\,,
$$
где
$$P(x)=1-f(x)^{q-1}\,,$$
а $f(x)=0$ является уравнением
гиперповерхности $X$, и примените пункт~(ii).)
В частности, конечное поле имеет тип $C_1$.
\item[(iv)]
Выведите из пункта (iii), что группа Брауэра конечного поля
тривиальна. Дайте другое доказательство этого факта, пользуясь тем,
что группа Галуа любого конечного расширения конечного поля
циклическая, теоремой Гильберта 90, а также
упражнением~\ref{exer-cyclgrp}(i).
\end{itemize}
\end{prob}

На самом деле, теорему Шевалле--Варнинга можно обобщить.
А именно, \emph{теорема Эн\'o} (см.~\cite{Esnault})
утверждает,
что количество точек на любом гладком многообразии Фано
(или, более общо, на гладком рационально связном многообразии)
над конечным полем $\mathbb{F}_q$ сравнимо с $1$ по модулю $q$.
В частности, такое многообразие имеет точку над полем $\mathbb{F}_q$.

\subsection{Группа Брауэра и многообразия Севери--Брауэра}
\label{subsection:Br-i-SB}

Заметим, что центральные простые алгебры размерности $n^2$ и
многообразия Севери--Брауэра
(см. определение~\ref{defin:Severi-Brauer})
размерности $n-1$ над полем $K$
канонически параметризуются одним и тем же множеством
$H^1\big(G_K,\PGL_n(K^{sep})\big)$, см. упражнения~\ref{exer-cohomBr}(i)
и~\ref{exer-SBforms}(o). Эту связь между алгебрами и многообразиями
можно также описать следующим образом.

\begin{prob}{\bf Многообразие Севери--Брауэра, построенное
по центральной простой алгебре}
\label{SBbyalgebra}
\hspace{0cm}
\begin{itemize}
\item[(i)]
Рассмотрим $n$-мерное векторное пространство $V$ над полем $K$.
Покажите, что множество $n$-мерных правых идеалов в $\End_K(V)$
канонически отождествляется с~$\P(V)$.
(Указание: каждой прямой в $V$ сопоставьте правый идеал
в~\mbox{$\End_K(V)$},
состоящий из всех эндоморфизмов, образ которых содержится в данной прямой.)
\item[(ii)]
Для центральной простой алгебры $A$ над $K$
размерности $n^2$ рассмотрим грассманиан ${\rm G}(A,n)$,
соответствующий линейным подпространствам размерности~$n$ в~$A$.
Покажите, что имеется замкнутое
подмногообразие $X\subset {\rm G}(A,n)$ над~$K$,
соответствующее тем подпространствам $I\subset A$ размерности $n$,
которые являются правыми идеалами.
При этом для любого расширения полей $K\subset L$ множество
$L$-точек~$X(L)$ канонически отождествляется с множеством правых
идеалов размерности~$n$ в~$A_L$.
(Указание: пусть $e_1,\ldots,e_{n^2}$ является базисом в $A$ над $K$,
и пусть
$$M_i\colon A\to A$$ является оператором умножения справа на $e_i$.
Тогда $n$-мерное подпространство~\mbox{$I\subset A$} является правым идеалом
тогда и только тогда, когда $M_i(I)\subset I$
для всех $1\leqslant i\leqslant n^2$.
То, что последнее условие задаёт замкнутое по Зарисскому подмножество,
легко проверить в локальных координатах на грассманиане.)
\item[(iii)]
Докажите, что $X$ является многообразием Севери--Брауэра,
соответствующим~$A$ так, как описано выше в терминах первых
когомологий. (Указание: воспользуйтесь пунктом~(i).)
\end{itemize}
\end{prob}

\bigskip

\begin{prob}{\bf Центральные простые алгебры
размерности $4$}\label{prob-quatconic}

Пусть характеристика поля $K$ не равна $2$.
\begin{itemize}
\item[(i)]
Покажите, что для кватернионной алгебры $A(a,b)$, $a,b\in K^*$,
соответствующее многообразие Севери--Брауэра изоморфно над $K$ конике,
заданной уравнением
$$u^2-av^2-bw^2=0\,.$$
(Указание: рассмотрим линейную функцию на $A(a,b)$, заданную формулой
$$x+yi+zj+tij\mapsto t\,.$$
Проверьте, что данная функция не равна тождественно
нулю ни на одном двумерном правом идеале $I$ в $A(a,b)$,
и в каждом таком $I$ существует единственный с точностью до
пропорциональности ненулевой элемент вида~\mbox{$x+yi+zj$}.
Данный элемент однозначно определяет идеал $I$, так
как в $A(a,b)$ нет одномерных правых идеалов.
Наоборот, каждой точке $(u:v:w)$ на конике сопоставьте двумерный правый
идеал в $A(a,b)$,
порождённый как $K$-векторное пространство элементами
$u+vi+wj$ и $av+ui-wij$, ср. с упражнением~\ref{prob-section}(ii).)
\item[(ii)]
Покажите, что любая центральная простая алгебра размерности $4$
является кватернионной. (Указание: рассмотрите соответствующее
многообразие Севери--Брауэра.)
\end{itemize}
\end{prob}

Последнее утверждение можно доказать и непосредственно
алгебраическим способом
(см.~\cite[VIII.11.2]{Bourbaki-algebra}).

\bigskip
\begin{defin}\label{defin-b}
\emph{Классом} многообразия Севери--Брауэра $X$, определённого над полем
$K$,
в группе Брауэра~$\Br(K)$ называется класс $b(X)$
соответствующей центральной простой
алгебры.
\end{defin}

\bigskip
\begin{prob}{\bf Тривиальность класса многообразия Севери--Брауэра}
\label{prob:SB-s-tochkoi}

Докажите, что многообразие Севери--Брауэра $X$,
определённое над полем $K$, имеет $K$-точку
тогда и только тогда, когда $b(X)=0$ в группе $\Br(K)$.
(Указание: воспользуйтесь
упражнениями~\ref{exer-SBforms}(o),(i) и~\ref{exer-cohomBr}(ii).)
\end{prob}

\bigskip

Дадим теперь более геометрическое описание класса $b(X)$.
Для этого сначала опишем классы в группе Брауэра,
связанные с инвариантами в группе Пикара.
Напомним, что многообразие $X$, определённое над полем $K$,
называется \emph{геометрически неприводимым},
если неприводимо многообразие $X_{\bar{K}}$
(и \emph{геометрически приводимым} в противном случае).

Пусть $X$ является
гладким проективным геометрически неприводимым
многообразием над полем $K$.
В частности,  для любого расширения полей $K\subset L$ корректно
определено поле рациональных функций на~$X_L$, которое мы будем
для краткости обозначать через $L(X)$. Отметим, что если $K\subset L$
является расширением Галуа, то возникает действие соответствующей
группы Галуа $G=\Gal(L/K)$ на $L(X)$,
причём для любой функции $f\in L(X)$
выполняется соотношение
$$g(f)=(g^{-1})^*(f)\,,$$
где в правой части равенства рассматривается действие группы~$G$
на схеме $X_L$. В частности, отображение взятия дивизора функции
$$L(X)^*\to \Div(X_L)$$
коммутирует с действием группы~$G$.

\begin{prob}{\bf Классы инвариантов в группе Пикара}\label{propb:Picarinvar}
\hspace{0cm}
\begin{itemize}
\item[(i)]
Покажите, что имеются каноническая точная последовательность
$$
H^1\big(G_K,K^{sep}(X)^*/(K^{sep})^*\big)\to \Br(K)\to \Br\big(K(X)\big)
$$
и естественный изоморфизм
$$
K(X)^*/K^*\cong \big(K^{sep}(X)^*/(K^{sep})^*\big)^{G_K}\,.
$$
(Указание: рассмотрите точную последовательность $G_K$-модулей
$$
1\to (K^{sep})^*\to K^{sep}(X)^*\to K^{sep}(X)^*/(K^{sep})^*\to 1\,,
$$
ср. с упражнением~\ref{prob:Br0X}(i).
Воспользуйтесь тем, что каноническое отображение
$$
H^2\big(G_K,K^{sep}(X)^*\big)\to \Br\big(K(X)\big)
$$
инъективно, а также упражнением~\ref{exer-spuskvp}(v).)
\item[(ii)]
Покажите, что имеется каноническая точная последовательность
$$
0\to\Pic(X)\to\Pic(X_{K^{sep}})^{G_K}\to
H^1\big(G_K,K^{sep}(X)^*/(K^{sep})^*\big)\to 0\,.
$$
(Указание: рассмотрите точную последовательность $G_K$-модулей
$$
1\to K^{sep}(X)^*/(K^{sep})^*\to\Div(X_{K^{sep}})\to
\Pic(X_{K^{sep}})\to 0\,,
$$
ср. с упражнением~\ref{prob:Br0X}(ii).
Воспользуйтесь пунктом~(i) и упражнением~\ref{divisorscohom}(iii).)
\item[(iii)]
Используя пункты~(i) и~(ii), постройте каноническую точную последовательность
$$
0\to\Pic(X)\longrightarrow\Pic(X_{K^{sep}})^{G_K}\stackrel{\xi}
\longrightarrow
\Br(K)\longrightarrow \Br\big(K(X)\big)\,.
$$
\item[(iv)]
Дадим теперь явное описание отображения $\xi$ из пункта~(iii).
Пусть элемент из группы $\Pic(X_{K^{sep}})^{G_K}$ определён
над конечным расширением Галуа $K\subset L$ с группой Галуа $G$,
то есть является классом $[H]$ в группе Пикара некоторого дивизора $H$ на $X_L$,
определённого над $L$. Покажите, что для любого элемента~\mbox{$g\in G$}
существует рациональная функция
$$f_g\in L(X)^*$$ на $X_L$,
дивизор которой равен $g(H)-H$.
Выведите отсюда, что для любых~\mbox{$g,h\in G$} дивизор рациональной функции
$$\omega(g,h)=f_g\cdot g(f_h)\cdot f_{gh}^{-1}$$
тривиален, то есть
$$\omega(g,h)\in L^*\subset L(X)^*\,.$$
Докажите, что $\omega$ является $2$-коциклом,
причём его класс в $H^2(G,L^*)$ не зависит от выбора дивизора $H$
и функций $f_g$ и равен $\xi([H])$.
(Указание: воспользуйтесь упражнением~\ref{exer-cobound-map}(i),(ii).)
\end{itemize}
\end{prob}

\bigskip
\begin{prob}{\textbf{Класс многообразия Севери--Брауэра}}\label{classSB}

Пусть $X$ является многообразием Севери--Брауэра размерности $n-1$.
\begin{itemize}
\item[(i)]
Пусть даны $n$-мерное векторное пространство $V$ над полем $K$
и автоморфизм~\mbox{$\gamma\in\PGL(V)$}
проективного пространства $\P^{n-1}=\P(V)$.
Пусть $\Pi\subset\P^{n-1}$ является гиперплоскостью над $K$,
заданной линейной функцией $l\in V^{\vee}$.
Проверьте, что дивизор $\gamma(\Pi)-\Pi$ на $\P(V)$
является дивизором рациональной функции
$$
\tilde{\gamma}(l)\cdot l^{-1}\,,
$$
где в качестве $\tilde{\gamma}\in\GL(V)$ можно взять
произвольное поднятие
элемента~\mbox{$\gamma\in\PGL(V)$},
и группа $\GL(V)$ действует на $V^{\vee}$ по формуле
\begin{equation}\label{eq:actionfunction}
\sigma(l)=(\sigma^{\vee})^{-1}(l),\quad \sigma\in\GL(V)\,.
\end{equation}
\item[(ii)]
Покажите, что группа Галуа $G_K$ действует тождественно на группе
$$\Pic(X_{K^{sep}})\cong\Z\,.$$
Таким образом,
гиперплоскость
$$\Pi\subset\P^{n-1}_{K^{sep}}$$
соответствует
дивизору $H\subset X_{K^{sep}}$, который определяет класс
$$[H]\in\Pic(X_{K^{sep}})^{G_K}\,.$$
Докажите, что выполняется равенство
$$\xi([H])=-b(X)\,,$$
где $\xi$ определяется
в упражнении~\ref{propb:Picarinvar}(iii).
(Указание: пусть $V$ является $n$-мерным векторным пространством над $K$.
Зафиксируем изоморфизм
$$\P(V)_L\cong X_L\,,$$
где $K\subset L$
является конечным расширением Галуа с группой Галуа $G$.
Пусть гиперплоскость $\Pi\subset \P(V)_L$ определена над $K$,
то есть задаётся линейной функцией~\mbox{$l\in V^{\vee}$}.
Пусть $1$-коцикл
$$\varphi\colon G\to\PGL(V_L)$$
задаёт форму $X$ проективного пространства $\P(V)$, см. упражнение~\ref{exer-spusk}, определённую рассмотренным выше изоморфизмом~\mbox{$\P(V)_L\cong X_L$}. Тогда дивизор~$g(H)$
соответствует гиперплоскости $\varphi(g)(\Pi)$.
Применим теперь явное вычисление класса~$\xi([H])$
из упражнения~\ref{propb:Picarinvar}(iv).
Из пункта~(i) следует,
что функция~$f_g$ соответствует рациональной функции на $\P^{n-1}_L$,
заданной отношением линейных функций
$$\tilde{\varphi}(g)(l)\cdot l^{-1}\,,$$
где $\tilde{\varphi}(g)\in\GL(V)$
является произвольным поднятием элемента $\varphi(g)\in\PGL(V)$.
Таким образом, имеется равенство
$$
\omega(g,h)=\tilde{\varphi}(g)(l)\cdot l^{-1}\cdot\, {}^g
\big(\tilde\varphi(h)\big)(l)\cdot\tilde\varphi(gh)(l)^{-1}\in L(X)^*\,,
$$
где рациональная функция $\omega(g,h)$ определена в упражнении~\ref{propb:Picarinvar}(iv).
Остаётся воспользоваться формулой~\eqref{eq:actionfunction},
тем, что оператор
$$
\tilde\varphi(g)\cdot\, ^{g}\big(\tilde\varphi(h)\big)\cdot
\tilde\varphi(gh)^{-1}\colon V\to V
$$
является скалярным, а также упражнением~\ref{exer-nonabcohom}(iii).)
\item[(iii)]
Покажите, что имеется точная последовательность
$$
0\to\Pic(X)\to\Z\to\Br(K)\to\Br\big(K(X)\big)\,,
$$
где первое отображение сопоставляет классу дивизора его степень
на
$$X_{K^{sep}}\cong\P^{n-1}_{K^{sep}}\,,$$
а второе отображение переводит $1$ в $-b(X)$.
(Указание: воспользуйтесь пунктом~(ii)
и упражнением~\ref{propb:Picarinvar}(iii).)
\item[(iv)]
Докажите, что порядок элемента $b(X)\in\Br(K)$ равен минимальной
положительной степени дивизоров на $X$, определённых над $K$.
(Указание: воспользуйтесь пунктом~(iii).)
\item[(v)]
Покажите, что класс $b(X)$ порождает ядро отображения
$$\Br(X)\to\Br\big(K(X)\big)\,.$$
В частности, если многообразия Севери--Брауэра $X$ и $Y$
бирационально эквивалентны над $K$, то $b(X)$ и $b(Y)$ порождают
одну и ту же подгруппу в~$\Br(K)$. (Указание: воспользуйтесь пунктом~(iii).)
\end{itemize}
\end{prob}

Гипотеза Амицура (ср. с~\cite{Amitsur}) утверждает, что верно утверждение,
обратное к упражнению~\ref{classSB}(v), а именно, если
классы $b(X)$ и $b(Y)$ многообразий Севери--Брауэра
над $K$ одинаковой размерности порождают
одну и ту же подгруппу в группе $\Br(K)$, то многообразия~$X$ и~$Y$
бирационально эквивалентны над $K$.

\bigskip

Предложим теперь непосредственное решение упражнения~\ref{classSB}(v).
Напомним, что любое рациональное отображение между проективными
пространствами единственным образом задаётся набором однородных форм
одинаковой степени без общих делителей. Степень этих форм называется
степенью соответствующего рационального отображения. Любое линейное
рациональное отображение, то есть рациональное отображение степени
один, является композицией линейного вложения и рациональной
линейной проекции.

\bigskip
\begin{prob}{\bf Морфизмы между многообразиями Севери--Брауэра}
\label{exer-morf-SB}
\hspace{0cm}
\begin{itemize}
\item[(i)]
Покажите, что при $m\leqslant n$
формы линейных вложений
$$\Pb^{\,m-1}\hookrightarrow\Pb^{\,n-1}\,,$$
а также формы рациональных линейных проекций
$$\Pb^{\,n-1}\dasharrow\Pb^{\,n-m-1}$$
параметризуются множеством
$H^1\big(G_K,\PGL_{n,m}(K^{sep})\big)$, где
$$\PGL_{n,m}(K^{sep})=\GL_{n,m}(K^{sep})/(K^{sep})^*\,,$$
а $\GL_{n,m}(K^{sep})$ состоит из матриц размера $n\times n$ с элементами
из поля $K^{sep}$, имеющих нулевой левый нижний
прямоугольник размера $(n-m)\times m$.
\item[(ii)]
Пусть $\varphi\colon X\dasharrow Y$ --- определённое над $K$ рациональное
отображение между многообразиями Севери--Брауэра, соответствующее
рациональному линейному отображению между проективными
пространствами после расширения скаляров с $K$ на $K^{sep}$.
Покажите, что
$$b(X)=b(Y)\,.$$
(Указание: для случая
линейного вложения или рациональной линейной проекции воспользуйтесь
пунктом~(i). При этом рассмотрите естественные отображения между
точными последовательностями $G_K$-модулей
$$
\xymatrix{
1\ar@{->}[r] & (K^{sep})^*
\ar@{->}[r]\ar@{->}[d]^{\mathrm{id}} &
\GL_{n,m}(K^{sep})\ar@{->}[r]\ar@{->}[d] &
\PGL_{n,m}(K^{sep})\ar@{->}[r]\ar@{->}[d] & 1\\
1\ar@{->}[r] & (K^{sep})^*\ar@{->}[r] & \GL_{r}(K^{sep})\ar@{->}[r]
& \PGL_{r}(K^{sep})\ar@{->}[r] & 1 }
$$
для $r=n,m$ и $n-m$, являющиеся тождественным отображением на
$(K^{sep})^*$. Затем воспользуйтесь
упражнением~\ref{exer-cohomBr}(iii)).
\item[(iii)]
Пусть $\varphi\colon X\dasharrow Y$ --- определённое над $K$ рациональное
отображение между многообразиями Севери--Брауэра, соответствующее
рациональному отображению степени $d$ между проективными
пространствами после расширения скаляров с $K$ на $K^{sep}$.
Покажите, что $$d\cdot b(X)=b(Y)\,.$$
(Указание: сначала разложите
произвольное рациональное отображение между проективными
пространствами в композицию вложения Веронезе и рационального
линейного отображения. Потом заметьте, что $d$-ая симметрическая
степень тавтологического представления группы $\GL_n(K^{sep})$
определяет отображение между точными последовательностями
$G_K$-модулей
$$
\xymatrix{ 1\ar@{->}[r] & (K^{sep})^* \ar@{->}[r]\ar@{->}[d]^{d} &
\GL_{n}(K^{sep})\ar@{->}[r]\ar@{->}[d] &
\PGL_{n}(K^{sep})\ar@{->}[r]\ar@{->}[d] & 1\\
1\ar@{->}[r] & (K^{sep})^*\ar@{->}[r] & \GL_{N}(K^{sep})\ar@{->}[r]
& \PGL_{N}(K^{sep})\ar@{->}[r] & 1 }
$$
Здесь $N$ равняется биномиальному коэффициенту
$$N=C_{n+d-1}^n={n+d-1\choose n}\,,$$
и на $(K^{sep})^*$ соответствующее отображение
является возведением в степень~$d$. Наконец, воспользуйтесь пунктом (ii).)
\item[(iv)]
Решите упражнение~\ref{classSB}(v), используя предыдущие пункты данного упражнения.
\end{itemize}
\end{prob}

\newpage

\section{Группа Брауэра II}
\label{section:unramified-Brauer}

\subsection{Полные поля дискретного нормирования}
\label{subsection:unramified-Brauer-generalities}

Сначала напомним несколько основных понятий и фактов о полных полях
дискретного нормирования
(подробности см. в~\cite{CasselsFrolich}, \cite{Serre-corps-locaux}
или~\cite{Lang-chisla}).
{\it Дискретным нормированием} поля~$K$
называется сюръективный гомоморфизм групп
$$v\colon K^*\to \Z\,,$$
для которого
выполнено неравенство
$$v(x+y)\geqslant {\rm min}\{v(x),v(y)\}.$$
при всех $x,y\in K$ (при
этом, по определению, $v(0)=+\infty$). {\it Кольцом нормирования}
называется подкольцо $\OO_K\subset K$, заданное условием
$v(x)\geqslant 0$. {\it Идеалом нормирования} называется идеал
$\m_K\subset K$, заданный условием $v(x)>0$. Кольцо $\OO_K$ является
локальным кольцом с главным максимальным идеалом $\m_K$, порождённым
любым {\it униформизующим} элементом, то есть таким элементом $x\in K$,
что $v(x)=1$.
В частности, из этого следует одно важное свойство дискретного нормирования:
если~\mbox{$v(x)\neq v(y)$}, то
$$v(x+y)=\min\{v(x),v(y)\}\,.$$
{\it Полем вычетов} называется фактор по максимальному
идеалу
$$\kappa=\OO_K/\m_K\,.$$

Если зафиксировать вещественное число $\alpha>1$, то дискретное
нормирование определяет метрику на поле $K$ по формуле
$$\rho(x,y)=\alpha^{-v(x-y)},\quad x,y\in K\,.$$
{\it Полное поле дискретного нормирования} --- это поле дискретного
нормирования, полное относительно соответствующей метрики. Полнота
$K$ не зависит от выбора $\alpha$ и равносильна тому, что
естественное отображение
$$\OO_K\to\varprojlim\OO_K/\m_K^i$$
является
изоморфизмом. Примерами полных полей дискретного нормирования
являются поле $p$-адических чисел $\Q_p$, а также поле рядов Лорана
$\kappa((u))$ с естественными нормированиями.

В дальнейшем $K$ обозначает полное поле дискретного нормирования. Одним из
основных арифметических свойств поля $K$ является выполнение {\it
леммы Гензеля} (см.~\cite[\S\,II.2]{Lang-chisla}):
пусть $f\in \OO_K[t]$ --- многочлен, редукция $\bar
f$ которого по модулю $\m_K$ имеет такой корень $\bar x$ в поле
вычетов $\kappa$, что $\bar f'(\bar x)\ne 0$. Тогда существует
единственный корень $x\in \OO_K$ многочлена $f$, для которого
$$x\equiv\bar x\pmod{\m_K}\,.$$
Если характеристики полей $K$ и $\kappa$
одинаковые, то у естественного гомоморфизма колец $\OO_K\to \kappa$
существует (неканоническое) сечение.
Если поле вычетов $\kappa$ имеет характеристику
$p>0$ и совершенно (то есть все неприводимые многочлены над $\kappa$
сепарабельны), то у естественного гомоморфизма групп
$$\OO_K^*\to\kappa^*$$
существует каноническое сечение $x\mapsto [x]$
(см.~\cite[Proposition~II.8]{Serre-corps-locaux}).
Элемент $[x]$ называется {\it представителем Тейхмюллера} элемента~$x$.

Пусть дано конечное расширение полей $K\subset L$. Тогда существует
единственное дискретное нормирование
$$w\colon L^*\to\Z\,,$$
для которого $w|_K$
пропорционально нормированию $v$, то есть
$$w|_K=e\cdot v$$
для некоторого натурального числа $e$
(см.~\cite[Теорема~II.10.1]{CasselsFrolich});
число $e$ называется {\it индексом ветвления} расширения $K\subset L$.
Кроме того, поле $L$
полно относительно дискретного нормирования $w$,
а кольцо $\OO_L$ является целым замыканием кольца $\OO_K$ в $L$
(см.~\cite[Proposition~II.3]{Serre-corps-locaux}).
Выполняется соотношение
$n=ef$,
где $n=[L:K]$, $f=[\lambda:\kappa]$, а
$\lambda=\OO_L/\m_L$ обозначает поле вычетов для $L$
(см.~\mbox{\cite[Предложение~I.5.3]{CasselsFrolich}}).

Расширение $K\subset L$ {\it неразветвлено}, если $e=1$ (или, что равносильно,
$n=f$) и расширение
$\kappa\subset\lambda$ сепарабельно.
В этом случае будем обозначать дискретное
нормирование поля~$L$ так же, как и соответствующее нормирование поля~$K$.
Поскольку $e=1$, это не
приводит к противоречию. Неразветвлённые
расширения сепарабельны (см.~\cite[Предложение~I.7.1]{CasselsFrolich}).
Категория неразветвлённых
расширений поля $K$ эквивалентна категории сепарабельных расширений
поля вычетов $\kappa$
(см.~\cite[Теорема~I.7.1]{CasselsFrolich}).
Примерами неразветвлённых расширений
являются расширение
$\Q_p\subset \Q_p(\sqrt{a})$, если $(a,p)=1$ и $p\neq 2$, а
также расширение $\kappa((u))\subset \lambda((u))$, если
$\kappa\subset\lambda$ --- сепарабельное конечное расширение полей.

Расширение $K\subset L$ {\it вполне разветвлено}, если $f=1$, то есть
$\kappa=\lambda$, или, что равносильно, $n=e$. Расширение $K\subset
L$  вполне разветвлено тогда и только тогда, когда $L=K(y)$, где $y$
является корнем {\it многочлена Эйзенштейна} $f(t)\in\OO_K[t]$, то есть
такого многочлена
$$f(t)=t^n+a_{n-1}t^{n-1}+\ldots+a_1 t+a_0\,,$$
что
$a_i\in\m_K$ для всех $i=0,\ldots,n-1$, и $v(a_0)=1$
(см.~\cite[Теорема~I.6.1]{CasselsFrolich}).
Примерами
вполне разветвлённых расширений являются расширение $\Q_p\subset
\Q_p(y)$, где $y$ --- корень уравнения
$$y^2+p^2y+p=0\,,$$
а также расширение $\kappa((u))\subset
\kappa((v))$, где $v^n=u$.

Для произвольного конечного расширения $K\subset L$ существует
наибольшее подполе $E\subset L$, содержащее $K$ и неразветвлённое
над ним
(см.~\cite[Теорема~I.7.2]{CasselsFrolich}).
При этом
поле вычетов для $E$
совпадает с сепарабельным замыканием поля $\kappa$ в поле
вычетов $\lambda$ для $L$. Если расширение $\kappa\subset\lambda$
сепарабельно, то
расширение $E\subset L$ вполне разветвлено, и
выполняются равенства
$$[E:K]=f\,,\quad [L:E]=e\,.$$

Пусть теперь $K\subset L$ --- конечное расширение Галуа с группой Галуа
$$G=\Gal(L/K)\,,$$
причём расширение полей вычетов $\kappa\subset
\lambda$ сепарабельно.
Тогда расширение $\kappa\subset
\lambda$ также является расширением Галуа
(см.~\cite[стр.~50]{CasselsFrolich}). Положим
$$G^{nr}=\Gal(\lambda/\kappa)\,.$$
Возникает канонический гомоморфизм
$$\psi\colon G\to\nolinebreak G^{nr}\,,$$
определяемый следующим
образом. Каждый автоморфизм $g\in G$ поля $L$ над~$K$ сохраняет
дискретное нормирование $w$ поля~$L$, совместимое с дискретным
нормированием~$v$ поля~$K$, в силу единственности нормирования~$w$
(см. выше). Поэтому~$g$ сохраняет подкольцо~$\OO_L$ в~$L$, а также
максимальный идеал $\m_L\subset \OO_L$. Можно рассмотреть~$g$ по
модулю идеала $\m_L$ и таким образом получить автоморфизм $\psi(g)$
поля~$\lambda$ над~$\kappa$. Гомоморфизм~$\psi$ всегда сюръективен.
Расширение $K\subset L$ неразветвлено тогда и только тогда, когда~$\psi$
--- изоморфизм. Наибольшее неразветвлённое над~$K$
подполе~\mbox{$E\subset L$}
является расширением Галуа, соответствующим по теории
Галуа подгруппе ${\rm Ker}(\psi)$. Имеется канонический изоморфизм
$$\Gal(E/K)\cong G^{nr}\,.$$
Пусть $K^{nr}\subset K^{sep}$
обозначает наибольшее неразветвлённое над $K$ подрасширение
в~$K^{sep}$. Имеется канонический изоморфизм
$$\Gal(K^{nr}/K)\cong G_{\kappa}\,,$$
где, как и раньше, $G_{\kappa}$ обозначает абсолютную группу
Галуа поля $\kappa$.

\bigskip

\begin{prob}{\bf Цикличность для
вполне разветвлённых расширений}\label{exer-cyclext}

Во введённых выше обозначениях
предположим, что поле $\kappa$ (а значит, и поле $K$) имеет характеристику
$0$, а также что $\kappa$ содержит все
корни из единицы. Обозначим через $\mu$ группу всех корней из единицы.

\begin{itemize}
\item[(i)]
Покажите, что поле $K$ также содержит все корни из единицы.
(Указание: воспользуйтесь леммой Гензеля.)
\item[(ii)]
Предположим, что расширение $K\subset L$ вполне разветвлено, и $L=K(y)$,
где $y$ является корнем многочлена Эйзенштейна
$$f(t)=t^n+a_{n-1}t^{n-1}+\ldots+a_1 t+a_0\in\OO_K[t].$$
Докажите, что
$$L=K\big(\sqrt[n]{-a_0}\big)\,.$$
(Указание: покажите, что элемент
$$z=\frac{a_{n-1}y^{n-1}+\ldots+a_1 y+a_0}{a_0}\in L$$
лежит в $\OO_L$ и
удовлетворяет сравнению
$$z\equiv 1\pmod{\m_L}\,.$$
Затем при
помощи леммы Гензеля покажите,
что $z$ является $n$-ой степенью в поле~$L$.)
\item[(iii)]
Пользуясь теорией Куммера (см. упражнение~\ref{exer-Kummer})
покажите, что любое вполне
разветвлённое конечное
расширение $K\subset L$ циклическое, а его группа Галуа
канонически вложена в $\mu$.
\item[(iv)]
Докажите, что для произвольного конечного расширения Галуа $K\subset L$
подгруппа
$$\Gal(L/E)\subset G$$
циклическая и центральная,
где $E$ является наибольшим подполем в $L$,
неразветвлённым над $K$. (Указание:
воспользуйтесь тем, что~$\mu$ является тривиальным $G$-модулем.)
\end{itemize}
\end{prob}

\bigskip

Нам также понадобится в дальнейшем следующий технический результат.
Пусть~\mbox{$K\subset L$} является конечным сепарабельным
расширением полей, $v$ является дискретным нормированием
поля $K$, а $K_v$ является пополнением поля $K$ относительно $v$
(в частности, в отличие от использованных выше обозначений,
$K$ не является \emph{полным} полем дискретного нормирования).
Пусть $W$ обозначает множество дискретных нормирований $w$ поля $L$,
продолжающих нормирование $v$, то есть таких, что
функция $w|_K$ пропорциональна нормированию $v$.
Для каждого $w\in W$ через $L_w$
обозначим пополнение поля $L$ относительно
нормирования $w$. Положим~\mbox{$A=L\otimes_K K_v$}.

\begin{prob}{\bf Продолжения дискретного нормирования}
\label{prob-val-extension}
\hspace{0cm}
\begin{itemize}
\item[(i)]
Проверьте, что имеется изоморфизм $K_v$-алгебр
$$\mbox{$A\cong \prod\limits_{i\in I}L_{i}$}\,,$$
где $I$ является некоторым конечным множеством,
а $L_{i}$, $i\in I$,
являются конечными расширениями поля $K_v$.
(Указание: по теореме о примитивном элементе~$L$ порождается над $K$ одним элементом.
Разложите его минимальный многочлен на неприводимые
сомножители над полем $K_v$.)
\item[(ii)]
Покажите, что для любого $w\in W$ однозначно определён
сюръективный гомоморфизм $K_v$-алгебр $A\to L_w$.
(Указание: требуемый гомоморфизм определяется
естественным вложением $L\subset L_w$,
а также индуцированным им вложением~\mbox{$K_v\subset L_w$}.)
Это задаёт отображение $W\to I$.
\item[(iii)]
Постройте обратное отображение $I\to W$.
(Указание: для каждого конечного расширения полей
$K_v\subset L_{i}$ существует единственное продолжение
нормирования $v$ с $K_v$ на $L_{i}$.
Положите $w$ равным его ограничению на $L\subset L_{i}$.)
Таким образом, имеются биекция $I\cong W$ и изоморфизм $K_v$-алгебр
$$
\mbox{$A\cong \prod\limits_{w\in W}L_{w}$}\,.
$$
\item[(iv)]
Предположим, что $K\subset L$ является расширением Галуа
с группой Галуа $G$.
Группа $G$ естественным образом
действует на множестве $W$.
Докажите, что это действие транзитивно.
(Указание: иначе $A\cong B\times C$, где $B$ и $C$ являются
$K_v$-алгебрами с действием группы $G$,
что противоречит равенству $A^G=K_v$.)
\item[(v)]
Докажите, что для любого $w\in W$ расширение $K_v\subset L_w$
является расширением Галуа, и естественный гомоморфизм
$\Gal(L_w/K_v)\to G$ индуцирует изоморфизм
$\Gal(L_w/K_v)\cong H$, где $H=\mathrm{Stab}_G(w)$.
(Указание: имеется естественный гомоморфизм $H\to \Gal(L_w/K_v)$.
Проверьте инъективность этого гомоморфизма,
а потом при помощи пункта~(iv) сравните $|H|$ и $[L_w:K_v]$.)
\end{itemize}
\end{prob}

\subsection{Группа Брауэра полного поля дискретного нормирования}
\label{subsection:res}

Пусть $K$ является полным полем дискретного нормирования
с кольцом нормирования $\OO_K$.
Предположим, что поле вычетов $\kappa$ совершенно.
В этом случае \emph{теорема Ленга}
(см.~\cite{Lang-quasi} или~\cite[X.6]{Serre-corps-locaux})
утверждает, что поле $K^{nr}$ имеет
тип $C_1$ (отметим, что из совершенности $\kappa$ следует
алгебраическая замкнутость поля вычетов поля~$K^{nr}$).
В частности,
$$\Br(K^{nr})=0$$
по упражнению~\ref{exer-C1Br}. Последний факт можно
доказать и другим способом,
используя некоторую технику работы с когомологиями конечных
групп и то, что для любой композиции сепарабельных расширений
$$K^{nr}\subset E\subset F$$
отображение нормы
$${\rm Nm}_{F/E}\colon F^*\to E^*$$
сюръективно (это нетривиальный арифметический факт,
см.~\cite[Proposition~X.11]{Serre-corps-locaux}
и~\cite[Proposition~V.7]{Serre-corps-locaux}).

\bigskip
\bigskip
\bigskip
\bigskip

\begin{prob}{\bf Отображение вычета на группе Брауэра}
\label{exercise:discretnoe-normirovanie-na-Br}
\hspace{0cm}
\begin{itemize}
\item[(i)]
Покажите, что естественное отображение
$$\Br(K^{nr}/K)\to\Br(K)$$
является изоморфизмом. (Указание: воспользуйтесь тем, что
$\Br(K^{nr})=0$.)
\item[(ii)]
Постройте каноническое отображение
$$\res\colon\Br(K)\to\Hom(G_{\kappa},\Q/\Z)\,.$$
Оно называется {\it отображением вычета}.
(Указание: воспользуйтесь изоморфизмом
$$\Br(K^{nr}/K)\cong\Br(K)$$
из пункта (i), отображением
$$H^2\big(G_{\kappa},(K^{nr})^*\big)\to H^2(G_{\kappa},\Z)\,,$$
которое возникает из отображения
$G_{\kappa}$-модулей $v\colon (K^{nr})^*\to\Z$, и
упражнением~\ref{exercise:division-module}(ii).)
\item[(iii)]
Пусть $K\subset L$ --- конечное неразветвлённое расширение
Галуа с группой Галуа~$G$, и $\alpha\in H^2(G,\OO_L^*)$.
Имеются естественные отображения (см. упражнение~\ref{exer-cohomBr}(v))
$$H^2(G,\OO_L^*)\to H^2(G, L^*)\stackrel{\sim}{\to}\Br(L/K)\hookrightarrow
\Br(K)\,.$$
Таким образом, определён вычет $\res(\alpha)$.
Покажите, что $\res(\alpha)=0$.
\end{itemize}
\end{prob}

Явным образом, утверждение из
упражнения~\ref{exercise:discretnoe-normirovanie-na-Br}(i)
заключается в том, что для любого центрального тела $D$ размерности
$n^2$ над $K$ существует такое неразветвлённое расширение $K\subset E$,
что
$$E\otimes_K D\cong \Mat_n(E)\,.$$
Это можно также доказать
непосредственно алгебраическим способом, показав, что среди
максимальных подполей в $D$ существует поле $E$, неразветвлённое над
$K$ (см.~раздел ``Добавление'' в~\mbox{\cite[VI.1]{CasselsFrolich}}).

Отметим, что элемент из $\Hom(G_{\kappa},\Q/\Z)$
соответствует циклическому расширению поля $\kappa$
с выбранной образующей группы Галуа.

\bigskip
\begin{defin}\label{defin-symbHilb}
Для ненулевых элементов $a,b\in K$ {\it символ Гильберта}
определяется по формуле
$$
(a,b)_v=(-1)^{v(a)v(b)}\,a^{-v(b)}b^{v(a)}\pmod{\m_K}\in\kappa^*.
$$
\end{defin}

\bigskip
\begin{prob}{\bf Вычет для циклических алгебр}
\label{exer-vychet-dlya-cycl}
\hspace{0cm}
\nopagebreak
\begin{itemize}
\item[(i)]
Покажите, что символ Гильберта определяет гомоморфизм
групп (см. определение~\ref{defin:K-group})
$$(-,-)_v\colon\KM_2(K)\to\kappa^*,\quad \{a,b\}\mapsto (a,b)_v\,.$$
\item[(ii)]
Пусть $K\subset L$ --- конечное неразветвлённое циклическое
расширение Галуа степени $n$ с группой Галуа $G$,
и $a\in K^*$.
Выберем образующую $s$ группы $G$. Пусть~$A$~---
циклическая алгебра над $K$, построенная
по этим данным (см. упражнение~\ref{exer-cyclic-alg}(i)).
Покажите, что если $v(a)=0$, то
$$\res([A])=0\,.$$
(Указание: воспользуйтесь
упражнением~\ref{exercise:discretnoe-normirovanie-na-Br}(iii).)
\item[(iii)] Предположим, что в обозначениях из пункта~(ii)
выполнено равенство~\mbox{$v(a)=1$}.
Докажите, что
$\res([A])$ соответствует циклическому расширению
$\kappa\subset\lambda$ с выбранной выше образующей
$$s\in G\cong\Gal(\lambda/\kappa)\,.$$
(Указание: воспользуйтесь функториальностью
изоморфизмов из упражнения~\ref{exer-cyclgrp}(i)
относительно морфизма $G$-модулей $v\colon L^*\to\Z$. Таким образом,
возникает коммутативная диаграмма
$$
\xymatrix{
K^*/\Nm_{L/K}(L^*)\ar@{->}[r]\ar@{->}[d]^{v} &
H^2(G,L^*)\ar@{->}[d]\\
\Z/n\Z\ar@{->}[r] & H^2(G,\Z)
}
$$
После этого воспользуйтесь
упражнением~\ref{exercise:division-module}(iii).)
\item[(iv)]
Пусть зафиксирован корень \mbox{$n$-ой} степени
из единицы $\zeta\in\kappa$, где $n$ взаимно просто с
характеристикой поля $\kappa$ (а значит, и поля $K$).
По теории Куммера (см.
упражнение~\ref{exer-Kummer}) возникает изоморфизм
$$\Hom(G_{\kappa},\Z/n\Z)\cong \kappa^*/(\kappa^*)^n\,.$$
Для простоты будем также обозначать через $\res$ композицию
$$\Br(K)_n\stackrel{\res}\longrightarrow
\Hom(G_{\kappa},\Z/n\Z)\cong \kappa^*/(\kappa^*)^n\,.$$
Докажите, что следующая диаграмма коммутативна
(см. упражнение~\ref{exer-norm-symbol})
$$
\xymatrix{
\KM_2(K)\ar@{->}[d]^{(-,-)_v}\ar@{->}[r]^{\nu_{\zeta}}&
\Br(K)_n\ar@{->}[d]^{\res}\\
\kappa^*\ar@{->}[r]& \kappa^*/(\kappa^*)^n
}
$$
Таким образом, вычет класса циклической алгебры
$A_{\zeta}(a,b)$ соответствует расширению Куммера
$$\kappa\subset\kappa\big(\sqrt[n]{(a,b)_v}\big)\,.$$
(Указание: напомним, что имеется (неканонический)
изоморфизм
$$K^*\cong \Z\times\OO_K^*\,.$$
Используя упражнение~\ref{exer-a-a},
выведите из этого, что группа $\KM_2(K)$ порождается
символами вида $\{a,b\}$, где $v(b)=0$, а $v(a)$ равно $0$ или $1$.
Дальше воспользуйтесь пунктами~(ii) и~(iii).)
\end{itemize}
\end{prob}

\bigskip
Из определения отображения $\res$ следует,
что для неразветвлённого расширения~\mbox{$K\subset L$}
и элемента $\alpha\in\Br(L/K)$ вычет $\res(\alpha)$
равен образу $\alpha$ относительно композиции
$$
H^2(G,L^*)\to H^2(G,\Z)\stackrel{\sim}\to \Hom(G,\Q/\Z)\cong
\Hom(G^{nr},\Q/\Z)\hookrightarrow \Hom(G_{\kappa},\Q/\Z)\,.
$$
Мы также будем обозначать возникающее отображение
$$\Br(L/K)\to\Hom(G^{nr},\Q/\Z)$$
символом $\res$.

\begin{prob}{\bf Группа Брауэра полного поля дискретного
нормирования}\label{exer-Brcompldv}

В пунктах (i) и (ii) предполагается, что расширение $K\subset L$
является конечным неразветвлённым расширением Галуа, и $G=\Gal(L/K)$.
\begin{itemize}
\item[(i)]
Покажите, что имеется изоморфизм
$$
\Ker\big(\Br(L/K)\stackrel{\res}\longrightarrow\Hom(G^{nr},\Q/\Z)\big)\cong
H^2(G,\OO_L^*)\,,
$$
и что выбор униформизующего элемента
$\pi\in\OO_K$ определяет расщепление
$$
\Br(L/K)\cong H^2(G,\OO_L^*)\oplus \Hom(G^{nr},\Q/\Z)\,.
$$
(Указание: воспользуйтесь тем, что $\pi$
определяет расщепление точной последовательности $G$-модулей
$$1\to\OO_L^*\to L^*\stackrel{v}\to \Z\to 0\,,$$
и тем, что $G\cong G^{nr}$.)
\item[(ii)]
Морфизм $G$-модулей
$$\OO_L^*\to\lambda^*,\quad x\mapsto x\pmod{\m_L}\,,$$
определяет отображение
$$H^2(G,\OO_L^*)\to H^2(G,\lambda^*)\,.$$
Покажите, что оно является изоморфизмом.
(Указание: рассмотрите на $\OO_L^*$ убывающую фильтрацию
подгруппами $1+\m_K^i$, $i\geqslant 1$.
Поскольку присоединённые факторы этой фильтрации изоморфны $\lambda$,
то по упражнению~\ref{exer-cohomadd}(iii) она удовлетворяет
всем условиям из упражнения~\ref{invlimmods}(ii).)
\item[(iii)]
Покажите, что имеет место точная последовательность
$$
0\to\Br(\kappa)\to\Br(K)\to\Hom(G_{\kappa},\Q/\Z)\to 0\,,
$$
причём выбор униформизующего элемента $\pi\in \OO_K$ задаёт её
расщепление. (Указание: воспользуйтесь пунктами (i) и (ii), перейдите
к пределу по всем неразветвлённым расширениям поля $K$ и примените
упражнение~\ref{exercise:discretnoe-normirovanie-na-Br}(i),(ii).)
\item[(iv)]
Покажите, что для произвольного конечного сепарабельного
расширения~\mbox{$K\subset\nlb L$}
(не обязательно нормального или неразветвлённого)
следующая диаграмма
коммутативна:
$$
\begin{CD}
\Br(K)@>\res>>\Hom(G_{\kappa},\Q/\Z)\\
@VVV@VVeV\\
\Br(L)@>\res>>\Hom(G_{\lambda},\Q/\Z)
\end{CD}
$$
Здесь правая вертикальная стрелка является композицией отображения
обратного образа
$$\Hom(G_{\kappa},\Q/\Z)\to \Hom(G_{\lambda},\Q/\Z)$$
и умножения на индекс ветвления~$e$.
\end{itemize}
\end{prob}

\bigskip

\begin{defin}\label{defin-local-field}
\emph{Неархимедовым локальным полем} называется
полное поле
дискретного нормирования с конечным полем вычетов
(ср. с определением локального поля из главы~\ref{section:Min-Hasse}).
\end{defin}

Можно показать, что любое неархимедово локальное поле
является либо конечным расширением поля $\Q_p$, где $p$ --- простое число,
либо полем $\F_q((T))$, где $q$ --- степень простого числа.

\bigskip
\begin{prob}{\bf Группа Брауэра неархимедова
локального поля}\label{exer-Brlocfield}

Пусть $K$ является неархимедовым локальным полем.

\begin{itemize}
\item[(i)]
Докажите, что имеется канонический изоморфизм
$$\res\colon\Br(K)\stackrel{\sim}\longrightarrow\Q/\Z\,.$$
(Указание:
воспользуйтесь упражнением~\ref{exer-Brcompldv}(iii), теоремой
Шевалле--Варнинга из упражнения~\ref{exer-Chevalley}, а также тем, что для
конечного поля группа Галуа канонически изоморфна~$\widehat{\Z}$.)
\item[(ii)]
Пусть $K\subset L$ --- конечное расширение неархимедовых локальных полей.
Покажите, что естественное отображение
$$\Br(K)\to \Br(L)$$
соответствует умножению на степень расширения полей
$$\Q/\Z\stackrel{n}\to\Q/\Z,\quad n=[L:K]\,.$$
(Указание: воспользуйтесь упражнением~\ref{exer-Brcompldv}(iv).)
\end{itemize}
\end{prob}

\bigskip

Пусть $\mu(K)\subset K^*$ --- группа корней из единицы в поле $K$
(мы не предполагаем, что $K$ содержит \emph{все} корни из единицы),
а $K\subset L$~--- конечное расширение Галуа с группой Галуа $G$, как и раньше.
Композицию естественного отображения
$$H^2\big(G,\mu(K)\big)\to H^2(G,L^*)\,,$$
возникающего из вложения $G$-модулей
$\mu(K)\hookrightarrow L^*$, и отображения вычета
$$\res\colon H^2(G,L^*)=\Br(L/K)\to\Hom(G_{\kappa},\Q/\Z)$$
для простоты будем также
обозначать через $\res$.

\begin{prob}{\bf Ядро отображения вычета}
\label{exercise:yadro-vycheta}

Пусть элемент
$\alpha\in H^2\big(G,\mu(K)\big)$ лежит в образе естественного отображения
$$H^2\big(G^{nr},\mu(K)\big)\to\nlb H^2\big(G,\mu(K)\big)\,.$$
Докажите, что $\res(\alpha)=0$.
(Указание: воспользуйтесь тем,
что отображение~$\res$ пропускается через изоморфизм
$$H^2(G^{nr}, E^*)\stackrel{\sim}{\longrightarrow}\nlb H^2(G, L^*)\,,$$
и тем, что $v\big(\mu(K)\big)=0$.)
\end{prob}

\subsection{Неразветвлённая группа Брауэра функционального поля}
\label{subsection:Br-function-field}

Пусть теперь $K$ --- конечно порождённое поле над полем нулевой
характеристики~$k$. Для каждого дискретного нормирования $v$ поля
$K$ через $K_v$ обозначим пополнение поля $K$ относительно $v$, а
через $\kappa_v$ обозначим поле вычетов для $K_v$.
Предположим, что нормирование $v$ тривиально на поле $k$,
то есть $v(k^*)=0$. Тогда поле $\kappa_v$ имеет нулевую характеристику;
в частности, оно совершенно.
Пусть $\res_v$
обозначает композицию естественного отображения
$$\Br(K)\to\Br(K_v)$$
и отображения вычета
(см.~упражнение~\ref{exercise:discretnoe-normirovanie-na-Br}(ii))
$$\Br(K_v)\to \Hom(G_{\kappa_v}, \Q/\Z)\,.$$
Поскольку поле
$\kappa_v$ совершенно, то вычет $\res_v$ корректно определён.
Если~\mbox{$\mathrm{char}(k)=p>0$}, то вычет становится корректно
определённым (и всё сказанное ниже остаётся верным)
после замены группы Брауэра на подгруппу элементов, порядки
которых взаимно просты с $p$.

Если $K\cong k(X)$, где $X$ --- нормальное неприводимое
многообразие над
полем~$k$, и~$D\subset\nlb X$~--- неприводимый приведённый дивизор, то
возникает дискретное нормирование
$$v_D(f)={\rm ord}_D(f),\quad f\in K^*\,,$$
для которого $\kappa_{v_D}\cong k(D)$.
Нормирования поля $K$, которые получаются таким образом для
некоторого выбора многообразия $X$ и дивизора $D$,
называются \emph{дивизориальными}.
В этом случае мы будем обозначать
отображение~$\res_{v_D}$ также через~$\res_D$.

Отметим, что не все нормирования являются дивизориальными.
А именно, нормирование $v$ дивизориально тогда и только тогда,
когда кольцо нормирования в~$K$ является
локализацией конечно порождённой $k$-алгебры
(см., например,~\mbox{\cite[Proposition~1]{Gabber}}). Последнее условие также
равносильно тому, что
$$\trdeg(K/k)-1=\trdeg(\kappa_v/k)\,,$$
где $\trdeg$ обозначает степень трансцендентности
расширения полей.

\bigskip
\begin{defin}\label{defin-unrBr}
{\it Неразветвлённой группой Брауэра} поля $K$ называется подгруппа
$$\Br^{nr}(K)=\bigcap_v \Ker(\res_v)$$
в группе $\Br(K)$, где $v$ пробегает
множество всех дискретных нормирований поля~$K$, тривиальных на $k$.
\end{defin}

На самом деле, при определении неразветвлённой группы
Брауэра достаточно в качестве $v$ брать только дивизориальные
нормирования (см.~\mbox{\cite[Proposition~2.1.8e]{CT}}
и~\mbox{\cite[\S2.2.2]{CT}}).

\bigskip

Следующее упражнение показывает,
что неразветвлённая группа Брауэра не меняется при
чисто трансцендентных расширениях.

\begin{prob}{\bf Теорема Фаддеева о
неразветвлённой группе Брауэра чисто трансцендентного расширения}
\label{exercise:Br-nr-is-stable}

Пусть $t$ --- формальная переменная, $\bar{K}$ обозначает
алгебраическое замыкание поля $K$, а $\Div(\A^1_{\bar{K}})$ ---
группу дивизоров на аффинной прямой $\A^1_{\bar{K}}$ над $\bar K$.
\begin{itemize}
\item[(i)]
Покажите, что естественное отображение
$$\Br\big(\bar{K}(t)/K(t)\big)\to\Br\big(K(t)\big)$$
является изоморфизмом. (Указание:
воспользуйтесь теоремой Тзена.)
\item[(ii)]
Докажите, что
$$
H^1\big(G_K, \Div(\A^1_{\bar{K}})\big)=0,
\quad H^2\big(G_K, \Div(\A^1_{\bar{K}})\big)\cong
\mbox{$\bigoplus\limits_{x\in \A^1_{K}}\Hom\big(G_{K(x)},\Q/\Z\big)$}\,,
$$
где прямая сумма берётся по всем схемным замкнутым точкам
$x\in\A^1_{K}$, а~\mbox{$K(x)$} обозначает поле вычетов в точке $x$.
(Указание: воспользуйтесь упражнением~\ref{divisorscohom}(iii).)
\item[(iii)]
Рассмотрим точную последовательность $G_K$-модулей
$$
1\to \bar{K}^*\to \bar{K}(t)^*\to \Div(\A^1_{\bar{K}})\to 0\,.
$$
Получите из неё точную последовательность
$$
0\to\Br(K)\stackrel{\theta}{\longrightarrow}
\Br\big(K(t)\big)\stackrel{\oplus_x \res_x}{\longrightarrow}
\mbox{$\bigoplus\limits_{x\in \A^1_{K}}\Hom\big(G_{K(x)},\Q/\Z\big)$}\,,
$$
где $\theta$ является естественным отображением, возникающим из
расширения полей~$K\subset K(t)$. (Указание: воспользуйтесь
изоморфизмом
$$G_K\cong\Gal\big(\bar{K}(t)/K(t)\big)\,,$$
а также пунктами~(i) и~(ii).)
\item[(iv)]
Покажите, что образ подгруппы
$$\Br^{nr}(K)\subset \Br(K)$$
при вложении $\theta$ содержится в подгруппе
$$\Br^{nr}\big(K(t)\big)\subset
\Br\big(K(t)\big)\,.$$
(Указание: дискретное
нормирование
$$w\colon K(t)^*\to\Z$$
можно ограничить на подполе $K\subset K(t)$.)
\item[(v)]
Докажите, что любое дискретное нормирование
$$v\colon K^*\to \Z$$
является
ограничением некоторого дискретного нормирования
$$w\colon K(t)^*\to\Z\,.$$
(Указание: рассмотрите кольцо нормирования $\OO\subset K$ и
униформизующий элемент $\varpi\in\OO$, связанные с $v$. Потом докажите,
что локализация кольца~$\OO[t]$ в простом идеале, порождённом $\varpi$,
является кольцом нормирования для некоторого дискретного
нормирования
$$w\colon K(t)^*\to\Z\,,$$
удовлетворяющего требуемому условию
$w|_K=v$.)
\item[(vi)]
Докажите, что имеет место вложение
$$\theta\big(\Br(K)\big)\cap
\Br^{nr}\big(K(t)\big)\subset \theta\big(\Br^{nr}(K)\big)\,.$$
(Указание: используйте пункт (v).)
\item[(vii)]
Используя пункты (iii), (iv) и (vi), докажите \emph{теорему
Фаддеева}: отображение~$\theta$ индуцирует
изоморфизм
$$\Br^{nr}(K)\stackrel{\sim}\longrightarrow\Br^{nr}\big(K(t)\big)\,.$$
\end{itemize}
\end{prob}

То, что поле $k$ (а значит, и поле $K$) имеет характеристику нуль,
используется для применения теоремы Тзена в
упражнении~\ref{exercise:Br-nr-is-stable}(i). Иначе надо требовать
совершенность поля $K$.

\bigskip
\begin{defin}\label{defin:stably-rational}
Неприводимые многообразия $X$ и $Y$ над полем $k$ называются
\emph{стабильно бирационально эквивалентными},
если для некоторых $m$ и $n$ многообразия
$X\times\P^m$ и $Y\times\P^n$ бирационально эквивалентны над $k$.
Многообразие, стабильно бирационально эквивалентное точке,
называется \emph{стабильно рациональным}.
\end{defin}

Упражнение~\ref{exercise:Br-nr-is-stable} показывает, что неразветвлённая
группа Брауэра $\Br^{nr}\big(k(X)\big)$ является инвариантом многообразия
$X$ относительно стабильно бирациональной эквивалентности. Это обстоятельство
будет нашим основным инструментом для доказательства стабильной
нерациональности многообразий
(и, в частности, их нерациональности в обычном смысле).

\subsection{Группа Брауэра многообразия}
\label{subsection:Br-of-a-variety}

Существует несколько способов обобщить определение группы Брауэра
с полей на многообразия (а также на схемы). Наиболее
непосредственный способ заключается в обобщении понятия
центральной простой алгебры.
{\it Алгеброй Адзумаи} над схемой~$X$ называется векторное
расслоение $A$ над $X$ с такой структурой алгебры,
что для любой (схемной) точки $x\in X$ слой $A\vert_x$ является центральной простой алгеброй над
полем вычетов~$k(x)$ (см. \cite[Предложение~IV.2.1]{Mil}).
Замыкая отношение
$$A\sim A\otimes_{\OO_X}\End(E)\,,$$
где $E$ --- произвольное векторное расслоение на $X$,
определим отношение эквивалентности на алгебрах Адзумаи.
Через $\Br(X)$ обозначим множество классов
эквивалентности алгебр Адзумаи над~$X$ (см. \cite[IV.2]{Mil}).
Тензорное произведение алгебр задаёт структуру полугруппы
на~$\Br(X)$. Как и для случая полей, имеется изоморфизм
$$A\otimes_{\OO_X}A^{op}\cong\End(A)\,.$$
Поэтому $\Br(X)$ является группой.

\begin{defin}\label{defin:Br-scheme}
Группа $\Br(X)$
называется {\it группой Брауэра} схемы $X$.
\end{defin}

С другой стороны, можно определить неразветвлённую группу
Брауэра неприводимого многообразия
над полем нулевой характеристики следующим образом.

\begin{defin}\label{defin:unramified-Br-X}
{\it Неразветвлённой группой Брауэра} неприводимого нормального
многообразия $X$ над полем $k$ характеристики нуль называется подгруппа
$$\Br^{nr}(X)=\bigcap_D \Ker(\res_D)$$
в группе $\Br\big(k(X)\big)$, где $D$ пробегает
множество всех неприводимых приведённых дивизоров $D\subset X$.
\end{defin}

Отметим, что определение~\ref{defin:unramified-Br-X}
можно обобщить и на случай ненормальных неприводимых многообразий,
однако нам это не понадобится.

Наконец, можно рассмотреть неразветвлённую группу
Брауэра $\Br^{nr}\big(k(X)\big)$ поля функций $k(X)$ на неприводимом
многообразии $X$ над полем характеристики нуль.

Опишем соотношения между тремя введёнными выше группами.
Имеются канонические гомоморфизмы групп
$$\Br(X)\stackrel{\phi}\longrightarrow
\Br^{nr}(X)\stackrel{\psi}\longleftarrow
\Br^{nr}\big(k(X)\big)\,.$$
Важный нетривиальный факт
заключается в том, что если многообразие $X$ гладкое,
то гомоморфизм $\phi$ является изоморфизмом.
Это утверждение получается из основного результата в~\cite{deJong}
и интерпретации
неразветвлённой группы Брауэра~$\Br^{nr}(X)$
в терминах этальных когомологий,
обобщающей изоморфизм из упражнения~\ref{exer-changefieldBr}(ii)
(подробности см.
в разделе~\ref{subsection:etale-Brauer} или в~\cite[\S3.4]{CT}).
Кроме того, если многообразие $X$ гладкое и полное, то
$\psi$ также является изоморфизмом
(см.~\cite[Proposition~4.2.3]{CT}).
В частности,
(неразветвлённые) группы Брауэра
у двух стабильно бирационально эквивалентных
гладких полных многообразий над полем нулевой
характеристики одинаковы
(см. замечание после определения~\ref{defin:stably-rational}).
Стабильно бирациональную инвариантность группы~$\Br^{nr}(X)$
можно доказать и более коротко, используя
только выражение этой группы через
этальные когомологии (см. разделы~\ref{subsection:etale-Kummer}
и~\ref{subsection:etale-Brauer}).

В разных случаях удобно использовать разные определения
группы Брауэра многообразия. Например,
для произвольного морфизма многообразий
$f\colon Y\to X$
легко определить гомоморфизм
обратного образа между группами Брауэра
$$f^*\colon\Br(X)\to\Br(Y)\,,$$
рассматривая обратный образ расслоений.
Если~$f$ является вложением точки в гладкую
кривую над полем характеристики нуль,
то построение~$f^*$ можно довольно просто провести и в терминах
неразветвлённых групп Брауэра.
Эта конструкция использует изоморфизм из упражнения~\ref{exer-Brcompldv}(i)
для дискретного нормирования~$v_Y$,
а также, в обозначениях
из этого упражнения, гомоморфизм
$$H^2(G,\OO_L^*)\to H^2(G,\lambda^*)\,.$$
Однако для произвольного морфизма $f$ непосредственное определение
обратного образа
в терминах неразветвлённых групп Брауэра
представляет техническую сложность (см.~\cite[\S12]{Rost}).

\bigskip
Как показывает следующее упражнение,
отображение вычета на группе Брауэра поля функций многообразия
обладает важным свойством: для каждого элемента
группы Брауэра его вычет равен нулю для почти всех дивизоров.

\begin{prob}{\textbf{Тривиальность вычета для почти
всех дивизоров}}
\label{prob:pochti-vse}

Пусть $X$ --- неприводимое нормальное многообразие над полем $k$
характеристики нуль. Положим $K=k(X)$.
Пусть $\alpha\in\Br(K)$.
\begin{itemize}
\item[(o)]
Покажите, что существует конечный морфизм
неприводимых нормальных многообразий $f\colon Y\to X$ над $k$, для которого
$\alpha\in\Br(L/K)$, и $K\subset L$ --- конечное расширение
Галуа, где $L=k(Y)$. Далее в этом упражнении
$G$ будет обозначать группу Галуа расширения $K\subset L$.
\item[(i)]
В обозначениях пункта~(o) рассмотрим неприводимый приведённый
дивизор~\mbox{$D\subset X$}
и произвольную неприводимую компоненту $E$
дивизора~\mbox{$f^{-1}(D)\subset Y$}.
Покажите, что дискретное нормирование $v_E$ является продолжением
дискретного нормирования $v_D$. Докажите, что все продолжения
дискретного нормирования $v_D$ имеют вид $v_E$ для некоторого $E$,
построенного как выше. (Указание: для доказательства первого утверждения
сделайте явное вычисление в локальных кольцах дивизоров
$D$ и $E$.
Для доказательства второго воспользуйтесь тем, что группа Галуа $G$
действует на множестве дивизоров $E$,
а также упражнением~\ref{prob-val-extension}(iv).)
\item[(ii)]
Пусть подгруппа $H\subset G$ состоит из таких
элементов $g\in G$, что $g(E)=E$. Пусть $K_D$ обозначает
пополнение поля $K$ относительно дискретного нормирования $v_D$,
а $L_E$ --- пополнение поля $L$ относительно дискретного
нормирования~$v_E$. Докажите, что $K_D\subset L_E$ является
расширением Галуа с группой Галуа~$H$. (Указание:
воспользуйтесь пунктом~(i),
а также упражнением~\ref{prob-val-extension}(v).)
\item[(iii)]
Рассмотрим $2$-коцикл $\varphi\colon G\times G\to L^*$,
соответствующий элементу
$$\alpha\in H^2(G,L^*)\cong\Br(L/K)\,.$$
Пусть $2$-коцикл $\omega$ задан композицией отображений
$$H\times H\to G\times G\stackrel{\varphi}\longrightarrow
L^*\hookrightarrow L_E^*\,.$$
Покажите, что образ $\alpha$ относительно естественного
отображения
$$\Br(K)\to\Br(K_D)$$
задаётся классом $2$-коцикла $\omega$.
\item[(iv)]
Докажите, что для почти всех
(то есть для всех кроме конечного числа)
неприводимых приведённых дивизоров $E\subset Y$ имеет место
равенство
$$v_E\big(\varphi(g_1,g_2)\big)=0$$
при любых $g_1, g_2\in G$.
(Указание: для любой
рациональной функции $\Phi\in k(Y)^*$ для почти всех
неприводимых приведённых дивизоров $E\subset Y$
выполнено равенство $v_E(\Phi)=0$.)
\item[(v)]
Докажите, что для почти всех
неприводимых приведённых дивизоров $D\subset X$ имеет место
равенство
$$\res_D(\alpha)=0\,,$$
где $\res_D$ обозначает вычет, связанный с
нормированием, определяемым дивизором~$D$.
(Указание:
почти все дивизоры~$D$ не содержатся в дивизоре
ветвления~\mbox{$\Delta\subset X$}
морфизма $f\colon Y\to X$, то есть расширение
$K_D\subset L_E$ неразветвлено.
Далее воспользуйтесь пунктами~(iii) и~(iv) и
упражнением~\ref{exercise:discretnoe-normirovanie-na-Br}(iii).)
\end{itemize}
\end{prob}

Из упражнения~\ref{prob:pochti-vse} следует, что для неприводимого
(нормального) многообразия $X$ имеется точная последовательность
(ср. с упражнением~\ref{exercise:Br-nr-is-stable})
$$
0\to\Br^{nr}(X)\to
\Br\big(k(X)\big)\stackrel{\oplus_D \res_D}\longrightarrow
\mbox{$\bigoplus\limits_{D\subset X} \Hom\big(G_{k(D)}, \Q/\Z\big)$}\,,
$$
где сумма берётся по всем неприводимым приведённым дивизорам
\mbox{$D\subset X$}.

\subsection{Геометрический смысл отображения вычета}
\label{subsection:models}

Пусть $B$ является целой схемой (например, $B$ является спектром
области целостности),
$K$ обозначает поле рациональных функций на $B$,
а $X$ является многообразием над $K$. {\it Моделью} многообразия $X$ над
$B$ называется схема $\mathcal X$ вместе с плоским
морфизмом $\pi\colon \mathcal{X}\to B$
(см.~\mbox{\cite[III.9]{Hartshorne}}) и изоморфизмом над $K$
между (схемным) общим слоем морфизма $\pi$ и многообразием $X$.
Если многообразие $X$ проективно,
то есть имеется замкнутое вложение
$X\subset\P^m$ над $K$,
то {\it проективной моделью} называется модель
$\mathcal X$ вместе с замкнутым
вложением~\mbox{$\mathcal X\subset\P^n_B$} над $B$.
Заметим, что после перехода к общему слою данное вложение
не обязано совпадать с исходным
проективным вложением $X\subset\P^m$ над $K$.
Если это выполнено (в частности, $m=n$), то мы будем говорить, что проективная
модель $\mathcal{X}\subset\P^n_B$ \emph{согласована
с проективным вложением}
многообразия $X$.
Если многообразие $X$ гладкое, то {\it гладкой моделью} называется
модель $\mathcal X$, для которой морфизм ${\mathcal X}\to B$ гладкий
(см.~\mbox{\cite[III.10]{Hartshorne}}).

Пусть $R$
является дедекиндовым кольцом,
то есть нётеровой целозамкнутой областью целостности, в которой
каждый ненулевой простой идеал максимален.
Например, $R$ является кольцом дискретного нормирования,
кольцом функций на регулярной аффинной кривой,
или кольцом целых в глобальном поле нулевой характеристики
(см. раздел~\ref{subsection:Min-Hasse-preliminary}).
Пусть $B=\Spec(R)$. Предположим, что $X\subset\P^n$ является проективным
многообразием над $K$. Тогда проективную модель
$\mathcal{X}$ можно явно задать как подсхему в проективном
пространстве $\P^n_R$ с однородными координатами~\mbox{$T_0,\ldots,T_n$}
следующим образом.
Рассмотрим идеал $I$ в кольце $R[T_0,\ldots,T_n]$,
состоящий из всех однородных многочленов,
обращающихся в нуль на $X$.
Этот идеал и задаёт подсхему $\mathcal{X}\subset\P^n_R$.
Другими словами, $\mathcal{X}$ задаётся уравнениями, соответствующими
образующим идеала $I$. В случае, если многообразие
$X$ геометрически неприводимо, гладкость модели $\mathcal{X}$
равносильна тому, что для любого ненулевого простого
идеала $\mathfrak{p}\subset R$ редукция уравнений схемы
$\mathcal{X}$ по модулю $\mathfrak{p}$ задаёт гладкое
многообразие над полем $R/\mathfrak{p}$.

\bigskip
Пусть теперь $R$ является кольцом нормирования
$\OO_K$ в полном поле $K$ дискретного нормирования $v$
с совершенным полем вычетов $\kappa$.
Предположим, что $X\subset \P^n$ является гиперповерхностью
степени $d$ над $K$.
Пусть $\mathcal X\subset \P^n_{\OO_K}$ является её проективной моделью,
согласованной с проективным вложением гиперповерхности~$X$
над~$K$. Тогда $\mathcal X$ задаётся в $\P^n_{\OO_K}$ однородным
уравнением $F=0$ степени $d$, все коэффициенты которого
из $\OO_K$, хотя бы один из них имеет нулевое нормирование,
причём над полем $K$ уравнение $F=0$ задаёт гиперповерхность,
проективно эквивалентную
гиперповерхности $X$. Это выводится из того,
что схема $\mathcal{X}$ плоская над $\OO_K$, а также из того,
любой дивизор на схеме $\P^n_{\OO_K}$
задаётся одним однородным уравнением.

\bigskip
Геометрический смысл отображения вычета заключается в следующем:
для многообразия Севери--Брауэра $X$ над $K$ существует
гладкая проективная модель над~$\OO_K$
тогда и только тогда, когда $\res\big(b(X)\big)=0$
(см. определение~\ref{defin-b}).

Доказательству этого факта посвящены следующие несколько упражнений.
Кроме того, ниже в упражнении~\ref{prob-conic-models}
мы разбираем отдельно случай коник элементарными методами.

\bigskip
\begin{prob}{\bf Тривиальность вычета}\label{trivialresidue}

Пусть $X$ является многообразием Севери--Брауэра над $K$,
для которого существует гладкая проективная модель $\mathcal{X}$ над $\OO_K$.
\begin{itemize}
\item[(i)]
Покажите, что схема $\mathcal X$ регулярна,
её замкнутый слой является неприводимым
главным дивизором, и отображение ограничения
$$\Pic(\mathcal{X})\to\Pic(X)$$
является изоморфизмом.
\item[(ii)]
Пусть $K\subset L$ является конечным неразветвлённым расширением Галуа
с группой Галуа $G$,
для которого многообразие $Y=X_L$ изоморфно
проективному пространству над $L$.
Пусть дивизор $H\subset Y$ соответствует при этом изоморфизме
гиперплоскости
в проективном пространстве.
Положим
$$
\mathcal{Y}=\Spec(\OO_L)\times_{\Spec(\OO_K)}\mathcal{X}\,.
$$
Пусть $\mathcal{H}\subset\mathcal{Y}$
является замыканием по Зарисскому дивизора $H$ в $Y\subset \mathcal{Y}$.
Группа $G$ естественным образом
действует на схеме $\mathcal Y$. Докажите, что для
любого элемента $g\in G$ существует
рациональная функция $f_g$ на схеме $\mathcal Y$, дивизор которой равен
$g(\mathcal{H})-\mathcal{H}$.
(Указание: схема $\mathcal Y$ обладает всеми свойствами
схемы~$\mathcal{X}$, перечисленными в пункте~(i).)
\item[(iii)] Докажите, что $\res\big(b(X)\big)=0$.
(Указание: примените конструкцию
из упражнения~\ref{propb:Picarinvar}(iv) к дивизору $\mathcal{H}$
и рациональным функциям $f_g$ на схеме $\mathcal Y$ из пункта~(ii).
Поскольку группа обратимых функций на схеме $\mathcal Y$
равна $\OO_L^*$, данная конструкция даёт класс из~\mbox{$H^2(G,\OO^*_L)$}.
Наконец, воспользуйтесь упражнением~\ref{classSB}(ii).)
\end{itemize}
\end{prob}

\bigskip
\begin{prob}{\textbf{Когомологии с целыми коэффициентами}}
\label{ingeralcohomology}

Пусть дано конечное неразветвлённое расширение Галуа
$K\subset L$ с группой Галуа~$G$,
а $\lambda$ обозначает поле вычетов кольца нормирования $\OO_L$ в поле $L$.
\begin{itemize}
\item[(i)]
Положим $\Gamma=\GL_n(\OO_L)$ и рассмотрим убывающую фильтрацию
$$\Gamma^i=I+\Mat_n\big(\m^i_L\big),\quad i\geqslant 1\,,$$
где $I$ обозначает единичную матрицу размера $n\times n$,
а $\m_L$ является максимальным идеалом в $\OO_L$.
Проверьте, что данная фильтрация удовлетворяет всем условиям
из упражнения~\ref{kogominvlimit}(ii).
(Указание: постройте изоморфизмы
$$\Gamma^i/\Gamma^{i+1}\cong \Mat_n(\lambda)\,,$$
где матрицы рассматриваются с аддитивным групповым законом,
и воспользуйтесь
упражнением~\ref{exer-cohomadd}(iii).)
\item[(ii)]
Покажите, что образ фильтрации из пункта~(i)
относительно гомоморфизма
$$\GL_n(\OO_L)\to\PGL_n(\OO_L)$$
также удовлетворяет всем условиям из упражнения~\ref{kogominvlimit}(ii).
(Указание: соответствующие присоединённые факторы
изоморфны фактору группы~$\Mat_n(\lambda)$ по подгруппе скалярных матриц.)
\item[(iii)]
Докажите, что
естественные отображения
$$H^1\big(G,\PGL_n(\OO_L)\big)\to H^1\big(G,\PGL_n(\lambda)\big)$$
при всех $n\geqslant 1$
являются биекциями, ср. с упражнением~\ref{exer-Brcompldv}(ii).
(Указание: воспользуйтесь пунктом~(ii) и упражнением~\ref{kogominvlimit}(ii).)
\item[(iv)]
Покажите, что имеется каноническая биекция между
множеством~\mbox{$H^1\big(G,\PGL_n(\OO_L)\big)$}
и множеством классов изоморфизма алгебр Адзумаи
ранга $n^2$
(см. раздел~\ref{subsection:Br-of-a-variety})
над $\OO_K$, которые становятся изоморфны матричной алгебре
над $\OO_L$. (Указание: рассмотрите расслоенную категорию
$\mathcal{M}$ над всеми конечными неразветвлёнными расширениями
$K\subset E$,
для которой~$\mathcal{M}(E)$ является
категорией алгебр Адзумаи над $\OO_E$.
Примените к ней теорию
спуска из разделов~\ref{subsection:spusk-dlya-kategorij}
и~\ref{subsection:forms-and-H1}.)
На самом деле, можно показать, что любая алгебра Адзумаи
над $\OO_K$ становится изоморфной матричной алгебре над $\OO_E$,
где~\mbox{$K\subset E$}~--- подходящее конечное неразветвлённое расширение.
\item[(v)]
Пусть дан элемент $\alpha\in H^2(G,\OO_L^*)$ и пусть
$\alpha_{\kappa}\in H^2(G,\lambda^*)$ обозначает его
образ относительно естественного отображения
$$H^2(G,\OO_L^*)\to H^2(G,\lambda^*)\,.$$
Пусть также дана конечномерная центральная простая алгебра
$\mathcal{A}_{\kappa}$ над $\kappa$,
для которой
$$[\mathcal{A}_{\kappa}]=\alpha_{\kappa}\in\Br(\kappa)\,.$$
Покажите, что существует однозначно определённая с точностью
до изоморфизма алгебра Адзумаи $\mathcal{A}$ над $\OO_K$,
редукция которой по модулю $\m_K$ изоморфна $\mathcal{A}_{\kappa}$.
При этом в $\Br(K)$ выполняется равенство
$$[\mathcal{A}_K]=\alpha_K\,,$$
где
$$\mathcal{A}_K=K\otimes_{\OO_K}\mathcal A\,,$$
а $\alpha_K$ является образом $\alpha$ относительно
естественного отображения
$$H^2(G,\OO_L^*)\to H^2(G,L^*)\,.$$
(Указание: воспользуйтесь пунктами~(iii) и~(iv).)
\end{itemize}
\end{prob}

\bigskip

\begin{prob}{\textbf{Построение гладкой модели}}\label{smoothmodel}
\hspace{0cm}
\begin{itemize}
\item[(i)]
Пусть дано многообразие Севери--Браэура $X$ над $K$.
Предположим, что существует алгебра Адзумаи
$\mathcal{A}$ над $\OO_K$,
для которой алгебра
$$\mathcal{A}_K=K\otimes_{\OO_K}\mathcal{A}$$
соответствует многообразию Севери--Брауэра $X$ так,
как это описано в начале
раздела~\ref{subsection:Br-i-SB}. Докажите, что тогда у $X$
существует гладкая проективная модель над $\OO_K$. (Указание:
примените аналог конструкции из начала раздела~\ref{subsection:Br-i-SB}
к алгебре $\mathcal{A}$ над $\OO_K$. А именно,
рассмотрите грассманиан $\G(\mathcal{A},n)$
над~$\OO_K$ и подсхему $\mathcal{X}$ над $\OO_K$ в нём,
параметризующую правые идеалы в $\mathcal{A}$. Тогда общий слой
$\mathcal{X}_K$ будет изоморфен $X$, а замкнутый слой будет
многообразием Севери--Брауэра, соответствующим конечномерной
центральной простой алгебре
$$\mathcal{A}_{\kappa}=\kappa\otimes_{\OO_K}\mathcal{A}$$
над $\kappa$. Таким образом, $\mathcal{X}$
будет гладкой проективной моделью многообразия~$X$ над~$\OO_K$.)
\item[(ii)]
Пусть существует алгебра Адзумаи $\mathcal{D}$ над $\OO_K$,
для которой $D=\mathcal{D}_K$ является телом.
Покажите, что тогда для любой конечномерной центральной простой
алгебры $A$ над $K$, удовлетворяющей условию
$b(A)=b(D)$, существует такая алгебра Адзумаи $\mathcal{A}$ над $\OO_K$,
что $A\cong \mathcal{A}_K$.
(Указание: из общей теории конечномерных центральных простых алгебр,
см. раздел~\ref{subsection:cohom-generalities},
следует, что $A\cong \Mat_p(D)$ для некоторого натурального числа $p$.)
\item[(iii)]
Пусть $\mathcal{D}$ является алгеброй Адзумаи над $\OO_K$,
для которой алгебра $\mathcal{D}_{\kappa}$ над~$\kappa$ является телом.
Докажите, что тогда алгебра $\mathcal{D}_K$ над $K$ также является телом.
(Указание: рассмотрите гиперповерхность $\mathcal H$ в $\P^n_{\OO_K}$,
заданную приведённой нормой алгебры $\mathcal{A}$,
см. упражнение~\ref{exer-rednorm}. Так как $\mathcal{D}_{\kappa}$ ---
тело, то у гиперповерхности $\mathcal{H}_{\kappa}$ нет
$\kappa$-точек по упражнению~\ref{exer-rednorm}(v).
Следовательно, у~$\mathcal{H}$ нет $\OO_K$-точек.
Из-за проективности $\mathcal{H}$ над $\OO_K$ это означает,
что и у гиперповерхности $\mathcal{H}_K$ над $K$ нет $K$-точек.
Используя вновь упражнение~\ref{exer-rednorm}(v),
мы получаем, что $\mathcal{D}_K$ является телом.)
\item[(iv)]
Пусть дано многообразие Севери--Браэура $X$ над $K$,
для которого выполнено равенство $\res\big(b(X)\big)=0$.
Докажите, что у $X$ существует гладкая модель над~$\OO_K$.
(Указание: пусть $K\subset L$ --- неразветвлённое расширение
Галуа с группой Галуа $G$,
для которого
$$b(X)\in\Br(L/K)\,.$$
По упражнению~\ref{exer-Brcompldv}(i)
из условия на вычет возникает элемент $\alpha\in H^2(G,\OO_L^*)$,
образ которого в $H^2(G,L^*)$ равен $b(X)$.
Рассмотрим его редукцию
$$\alpha_{\kappa}\in H^2(G,\lambda^*)\,.$$
Пусть $\mathcal{D}_{\kappa}$ является телом над $\kappa$,
класс которого в $\Br(\kappa)$ равен $\alpha_{\kappa}$.
Применяя упражнение~\ref{ingeralcohomology}(v),
мы получаем алгебру Адзумаи~$\mathcal{D}$ над~$\OO_K$.
Остаётся воспользоваться пунктами~(i),~(ii) и~(iii).)
\end{itemize}
\end{prob}

\bigskip

Рассмотрим теперь более подробно случай коник.
Для них равносильность тривиальности вычета с существованием
гладкой модели можно установить элементарными методами.
Для простоты предположим,
что характеристика поля $\kappa$ не равна~$2$
(а значит, и характеристика поля~$K$ также не равна $2$).
Пусть дана гладкая коника~\mbox{$X\subset \P^2$} над~$K$.
{\it Хорошей} моделью над $\OO_K$ коники $X$ будем называть проективную
модель $\mathcal{X}\subset\P^2_{\OO_K}$,
согласованную с вложением $X\subset\P^2_K$,
для которой схема $\mathcal X$ регулярна,
а замкнутый слой
$$\mathcal{X}_{\kappa}=\Spec(\kappa)\times_{\Spec(\OO_K)}\mathcal{X}$$
либо гладкий, либо является формой пары различных прямых в~$\P^2$
(ср. с упражнением~\ref{exer-Clifunr}).

\bigskip
\bigskip
\bigskip
\bigskip
\begin{prob}{\bf Модели коник}
\label{prob-conic-models}
\hspace{0cm}
\begin{itemize}
\item[(i)]
Покажите, что любая хорошая модель гладкой коники $X$
задаётся уравнением
\begin{equation}\label{eq:horoshaya-konika}
x^2-ay^2-bz^2=0\,,
\end{equation}
где $a,b\in\OO_K$,
причём хотя бы один из коэффициентов $a$ и $b$ обратим в $\OO_K$,
то есть имеет нулевое нормирование, а другой коэффициент
имеет нормирование~$0$ или~$1$. (Указание: поскольку $2$
обратимо в $\OO_K$, уравнение модели $\mathcal{X}$ можно
диагонализовать. Из условия на геометрический замкнутый слой следует,
что не менее двух коэффициентов у диагонального уравнения должны
быть обратимы. Из условия на регулярность
схемы $\mathcal X$ следует, что нормирование
оставшегося коэффициента не может быть больше $1$.)
\item[(ii)]
Докажите, что для любой гладкой коники $X$ над $K$
существует хорошая модель над $\OO_K$. (Указание:
уравнение коники можно привести к виду~\eqref{eq:horoshaya-konika}
умножением этого уравнения на подходящую степень униформизующего
элемента в~$\OO_K$ и диагональной заменой переменных.)
\item[(iii)]
Выберем хорошую модель $\mathcal X$ коники $X$.
Проверьте, что
имеется $G_{\kappa}$-эквивариантная биекция между множеством
неприводимых компонент многообразия
$\mathcal{X}_{\bar{\kappa}}$
и множеством $\big\{\pm\sqrt{(a,b)_v}\big\}$
(см. определение~\ref{defin-symbHilb}),
где $a$ и $b$ определяются как в пункте~(i).
Если модель $\mathcal X$
гладкая, то мы по определению считаем первое
действие тождественным действием на множестве из двух элементов.
(Указание: пользуясь пунктом~(i), разберите
по отдельности случай гладкого и особого замкнутого слоя.)
В частности, из упражнений~\ref{prob-quatconic}(i)
и~\ref{exer-vychet-dlya-cycl}(iv) следует, что
действие группы $G_{\kappa}$ на множестве неприводимых
компонент многообразия $\mathcal{X}_{\bar{\kappa}}$
не зависит от выбора хорошей модели.
\item[(iv)]
Предположим, что действие группы Галуа $G_{\kappa}$
на множестве
неприводимых компонент многообразия $\mathcal{X}_{\bar{\kappa}}$
тождественно.
Покажите, что тогда у $X$ существует гладкая проективная
модель над $\OO_K$. (Указание: в этом случае у коники~$X$
есть точка над $K$ и она проективно эквивалентна любой другой
конике с точкой, например, конике, заданной уравнением $x^2-y^2-z^2=0$.
Более геометрический подход заключается в том,
чтобы в случае особого замкнутого слоя
стянуть на схеме $\mathcal X$ одну из двух неприводимых компонент
замкнутого слоя~$X_{\bar{\kappa}}$,
каждая из которых по условию определена
над $\kappa$, и получить гладкую модель над~$\OO_K$.)
\item[(v)]
Покажите, что для любой гладкой проективной модели $\mathcal X$ коники $X$
существует замкнутое вложение $\mathcal X\subset\P^2_{\OO_K}$,
согласованное с вложением $X\subset\P^2_K$.
(Указание: убедитесь, что замкнутый
слой гладкой проективной модели~$\mathcal{X}$ является формой~$\P^1$,
так как относительный канонический пучок имеет степень~$-2$ на каждом слое.
Далее вложите схему~$\mathcal{X}$ при помощи относительного
антиканонического пучка в $\P^2_{\OO_K}$.)
\item[(vi)]
Докажите, что у коники $X$ существует гладкая проективная
модель тогда и только тогда, когда $\res\big(b(X)\big)=0$.
(Указание: воспользуйтесь предыдущими пунктами,
а также упражнением~\ref{exer-vychet-dlya-cycl}(iv).)
\end{itemize}
\end{prob}

Более подробно ознакомиться с вырождениями моделей
многообразий Севери--Брауэра можно по статье~\cite{Artin1982}.
Там доказано, что для любого многообразия Севери--Брауэра $X$ над $K$
существует проективная модель $\mathcal{X}$, которая является
регулярной схемой, и для которой замкнутый слой $\mathcal{X}_{\kappa}$
неприводим. При этом алгебраческое замыкание $\lambda$ поля $\kappa$
в поле рациональных функций $\kappa(\mathcal{X}_{\kappa})$
(ср. с упражнением~\ref{divisorscohom}(ii))
является циклическим расширением поля $\kappa$,
задаваемым вычетом
$$\res\big(b(X)\big)\in\Hom(G_{\kappa},\Q/\Z)\,.$$
Более того, многообразие $\mathcal{X}_{\bar\kappa}$ является объединением $r$
трансверсально пересекающихся копий одного гладкого проективного
рационального многообразия, где~$r$ равно порядку вычета
$\res\big(b(X)\big)$ в группе $\Hom(G_{\kappa},\Q/\Z)$,
и группа Галуа
$$\Gal(\lambda/\kappa)\cong\Z/r\Z$$
циклически переставляет неприводимые компоненты
многообразия $\mathcal{X}_{\bar\kappa}$. Это обобщает случай коник,
рассмотренный в упражнении~\ref{prob-conic-models}.

\newpage
\section{Пример унирационального нерационального многообразия}
\label{section:primer-Bogomolova}

В этой главе мы строим пример многообразия
вида~\mbox{$X=V/G$}, где~$G$~--- конечная группа, а~$V$~--- её представление,
определённое над алгебраически замкнутым полем~$k$ характеристики нуль,
для которого при помощи введённых
в предыдущих главах понятий удаётся установить нерациональность.
Примеры такого сорта впервые появились в статьях Д.\,Солтмена~\cite{Saltman}
и Ф.\,А.\,Богомолова~\cite{Bogomolov}, однако источником для нашего изложения
послужила статья И.\,Р.\,Шафаревича~\cite{Shafarevich}, в которой подход
Солтмена и Богомолова был значительно упрощён.

\subsection{Геометрические данные}

Пусть $G$ --- конечная подгруппа в группе автоморфизмов
гладкого неприводимого многообразия $V$, определённого над полем~$k$.
Положим $L=k(V)$ и
$K=k(V)^G$. Таким образом, расширение полей
$K\subset L$ является расширением Галуа с группой Галуа~$G$.

\bigskip

\begin{prob}{\bf Относительная группа Брауэра}
\label{exercise:G-na-Pic-V-0}

Пусть $\Pic(V)=0$ и $k[V]^*=k^*$.
Докажите, что имеется естественное вложение
$$H^2(G,k^*)\hookrightarrow\Br(L/K)\,.$$
(Указание: воспользуйтесь точной последовательностью $G$-модулей
$$1\to k^*\to k(V)^*\to\Div(V)\to 0$$
и упражнением~\ref{prob:stab-perm-profinite}(i).)
\end{prob}

\bigskip

Предположим теперь, что поле $k$ алгебраически замкнуто и имеет
характеристику нуль. Пусть $\mu$ обозначает группу корней из единицы
в $k$.

\begin{prob}{\bf Замена $k^*$ на $\mu$}
\label{exer:k-mu}

Покажите, что
$$H^i(G,\mu)\cong H^i(G,k^*)$$
при всех $i>0$.
(Указание: воспользуйтесь однозначной делимостью
в группе~$k^*/\mu$ и
упражнением~\ref{exercise:division-module}(i).)
\end{prob}

В дальнейшем $Z(\Gamma)$ обозначает центр группы $\Gamma$.

Предположим, что существует элемент $0\neq \alpha\in H^2(G,\mu)$
обладающий следующим свойством:

\begin{equation}\label{eq:Shafarevich}
\begin{minipage}[c]{0.8\textwidth}
для любой подгруппы $H\subset G$ и для любой центральной
циклической подгруппы $N\subset\nlb Z(H)$ в $H$ ограничение
$\alpha|_H$ элемента $\alpha$ на подгруппу $H$
содержится в образе естественного отображения $H^2(H/N,\mu)\to H^2(H,\mu)$.
\end{minipage}
\end{equation}

Пусть $\bar{\alpha}\in\nlb\Br(K)$~--- образ $\alpha$ при вложении
$$H^2(G, \mu)\hookrightarrow\Br(K)\,,$$
являющемся композицией
изоморфизма
$$H^2(G,\mu)\cong H^2(G,k^*)$$
из упражнения~\ref{exer:k-mu},
вложения
$$H^2(G,k^*)\hookrightarrow\Br(L/K)$$
из упражнения~\ref{exercise:G-na-Pic-V-0}, и естественного вложения
$\Br(L/K)\hookrightarrow\Br(K)$.

\begin{prob}{\bf Фактор по действию конечной группы}
\label{exercise:factor-po-hitrovyebnutoj-gruppe}

\begin{itemize}
\item[(i)]
Во введённых выше обозначениях докажите,
что для любого дискретного нормирования $v$ поля $K$,
тривиального на $k$, имеется равенство $\res_v(\bar{\alpha})=0$.
(Указание: пусть $w$ является произвольным продолжением на $L$
нормирования $v$. По упражнению~\ref{prob-val-extension}(v),
расширение соответствующих пополненных полей $K_v\subset L_w$
является расширением Галуа, причём его группа Галуа $H$
канонически вложена в $G$. При этом образ $\alpha$
относительно отображения $\Br(K)\to\Br(K_v)$ задаётся
ограничением $\alpha|_H\in H^2(H,\mu)$.
По упражнению~\ref{exer-cyclext}(iv) и условию~\eqref{eq:Shafarevich}
элемент
$\alpha|_H$ лежит в образе естественного отображения
$$H^2(H^{nr},\mu)\to H^2(H,\mu).$$
Теперь воспользуйтесь
упражнением~\ref{exercise:yadro-vycheta}.) Таким образом,
$\bar{\alpha}\in\Br^{nr}(K)$.
\item[(ii)]
Выведите из пункта (i), что многообразие $V/G$ не является стабильно
рациональным. (Указание: воспользуйтесь
упражнением~\ref{exercise:Br-nr-is-stable}.)
\end{itemize}
\end{prob}

В частности, если $V$ --- точное конечномерное представление группы
$G$, то фактор $V/G$ является унирациональным, но не стабильно
рациональным многообразием.

\subsection{Построение группы}

Отметим, что поскольку характеристика поля $k$ равна нулю,
то (неканонический) согласованный выбор
первообразных корней из единицы определяет
изоморфизм $G$-модулей $\Q/\Z\cong\mu$.
Построим группу $G$, для которой существует элемент
$$0\neq\alpha\in H^2(G,\Q/\Z)\cong H^2(G,\mu)\,,$$
удовлетворяющий введённому выше условию~\eqref{eq:Shafarevich}.
Пусть~$p\neq 2$~--- простое число, и $W\cong\F_p^4$ --- четырёхмерное
векторное пространство над конечным полем~$\F_p$. Положим
$$W'=\Lambda^2(W)\cong \F_p^6\,,$$
где $\Lambda^2$ обозначает внешний квадрат векторного пространства.
Построим группу $\widetilde{G}$ как
центральное расширение
$$0\to W'\to \widetilde{G}\to W\to 0\,,$$
взяв в качестве определяющего $2$-коцикла спаривание
$$\omega\colon (w_1,w_2)\mapsto w_1\wedge w_2\in W'\,,$$
где $W'$
рассматривается как тривиальный $W$-модуль (см. упражнение~\ref{exer-012}(ii)).
Группа~$\widetilde{G}$ имеет порядок $p^{10}$,
класс нильпотентности $2$ и
период $p$, и при этом выполнено
$$W'=[\widetilde{G},\widetilde{G}]=Z(\widetilde{G})\subset\widetilde{G}$$
(все эти утверждения
непосредственно следуют из определения группы $\widetilde{G}$).

\bigskip

\begin{prob}{\textbf{Свойства группы $\widetilde{G}$}}
\label{exercise:kommutatornoe-sparivanie}
\hspace{0cm}
\begin{itemize}
\item[(i)] Проверьте, что $\omega$ действительно является $2$-коциклом.
\item[(ii)] Покажите, что в группе $\widetilde{G}$ выполнено равенство
$$[g_1,g_2]=2w_1\wedge w_2\,,$$
где $w_i\in\nlb W$, и $g_i\in\widetilde{G}$ является (произвольным)
поднятием элемента $w_i$.
\item[(iii)] Докажите, что для любой нетривиальной подгруппы $M\subset W'$ расширение
$$0\to M\to\widetilde{G}\to\widetilde{G}/M\to 1$$
нетривиально. (Указание: выберите нетривиальный элемент $m\in M$ и такие
элементы $w_i, v_i\in W$, что
$$2\sum_i w_i\wedge v_i=m\,.$$ Тогда для любых поднятий
$g_i, f_i\in\nlb\widetilde{G}$ элементов $w_i$ и $v_i$ выполнено
равенство
$$\mbox{$\prod\limits_i [g_i,f_i]=m\in\widetilde{G}\,.$}$$
С другой стороны, если обсуждаемое расширение тривиально, то в $\widetilde{G}$ есть подгруппа $\widetilde{G}_M\cong\widetilde{G}/M$,
и $[\widetilde{G}_M, \widetilde{G}_M]\subset\widetilde{G}_M$.
Остаётся заметить, что можно выбрать $g_i, f_i\in\widetilde{G}_M\subset\widetilde{G}$, и получить противоречие.)
\end{itemize}
\end{prob}

\bigskip
\bigskip
\bigskip
\bigskip

\begin{prob}{\bf Нетривиальность расширений}\label{exer-nontrivext}

Пусть $\Gamma$ --- произвольная конечная группа.
\begin{itemize}
\item[(i)]
Рассмотрим элемент
$$\varphi\in\Hom(\Gamma,\Q/\Z)\cong H^1(\Gamma,\Q/\Z)$$
и кограничное отображение
$$\delta\colon H^1(\Gamma,\Q/\Z)\to H^2(\Gamma,\Z/p\Z)\,,$$
возникающее из точной последовательности тривиальных $\Gamma$-модулей
\begin{equation}\label{eq-exseqtriv}
0\to \Z/p\Z\to\Q/\Z\stackrel{p}\to\Q/\Z\to 0\,.
\end{equation}
Докажите, что элемент $\delta(\varphi)\in H^2(\Gamma,\Z/p\Z)$
соответствует расширению
$$
0\to \Z/p\Z\to \widetilde{\Gamma}\to \Gamma\to 1\,,
$$
возникающему как обратный образ расширения~\eqref{eq-exseqtriv}
относительно гомоморфизма $\varphi\colon\Gamma\to \Q/\Z$.
\item[(ii)]
Пусть ненулевой элемент $\beta\in H^2(\Gamma,\Z/p\Z)$ соответствует
расширению
$$
0\to \Z/p\Z\to \widetilde{\Gamma}\to \Gamma\to 1\,,
$$
для которого группа $\widetilde{\Gamma}$ имеет период $p$.
Покажите, что образ элемента $\beta$
относительно естественного отображения
$$H^2(\Gamma,\Z/p\Z)\to H^2(\Gamma,\Q/\Z)$$
нетривиален.
(Указание: предположите обратное
и воспользуйтесь длинной точной последовательностью когомологий.
Далее воспользуйтесь пунктом (i), чтобы построить в $\widetilde{\Gamma}$ элемент порядка строго больше,
чем $p$.)
\end{itemize}
\end{prob}

\bigskip
Выберем нетривиальный элемент $z\in W'\subset\widetilde{G}$ (до
последнего момента нам не потребуется от него никаких специальных
свойств) и породим им подгруппу
$$\Z/p\Z\cong\nlb\langle z\rangle\subset\widetilde{G}\,.$$
Определим группу $G$ как
$\widetilde{G}/\langle z\rangle$, и рассмотрим элемент $\beta\in
H^2(G, \Z/p\Z)$, соответствующий расширению
$$0\to \langle z\rangle \to \widetilde{G}\to G\to 1\,.$$
Пусть $H\subset G$ --- произвольная подгруппа. Рассмотрим элемент $u\in Z(H)$
и порождённую им (центральную циклическую) подгруппу $N\subset H$.
Пусть
$$I\subset H^2(H,\Z/p\Z)$$ является образом естественного отображения
$$H^2(H/N, \Z/p\Z)\to H^2(H,\Z/p\Z)\,,$$
а $\gamma\in H^2(H,\Z/p\Z)$
является ограничением $2$-коцикла $\beta$ с $G$ на $H$.

\begin{prob}{\textbf{Свойство~\eqref{eq:Shafarevich} для группы $G$}}

\begin{itemize}
\item[(o)]
Рассмотрим центральное расширение
\begin{equation}\label{eq:rasshirenie-H}
0\to \langle z\rangle \to \widetilde{H}\to H\to 1\,,
\end{equation}
соответствующее $2$-коциклу $\gamma$.
Покажите, что $\gamma\in I$ тогда и только тогда, когда у расширения~\eqref{eq:rasshirenie-H}
существует (теоретико-групповое) сечение над $N$, образ которого является нормальной подгруппой в $\widetilde{H}$.
\item[(i)] Докажите, что
$$Z(G)=Z(\widetilde{G})/\langle z\rangle=W'/\langle z\rangle\,,$$
и имеется точная последовательность
$$0\to Z(G)\to G\to W\to 0\,.$$
\item[(ii)] Предположим, что $u\in Z(G)$. Докажите, что
$\gamma\in I$. (Указание: используйте пункты~(o) и~(i), а также то, что $W'$ является абелевой группой периода $p$.)
\item[(iii)] Предположим,
что $u\not\in Z(G)$, и образ $v\in W$ элемента $u$ порождает
образ группы $H$ при гомоморфизме $G\to W$. Докажите, что
$\gamma=0$. (Указание: в этом случае группа $H$ и её прообраз
$\widetilde{H}$ в $\widetilde{G}$ --- абелевы группы периода $p$, так что
$\langle z\rangle\subset\widetilde{H}$ выделяется прямым слагаемым.)
\item[(iv)] Предположим,
что $u\not\in Z(G)$, и $u'\in H$ --- такой элемент,
что его образ~$v'$ в~$W$ не лежит в
подгруппе, порождённой элементом $v$. Докажите, что в этом случае
$$v\wedge v'=lz$$
при некотором $l\not\equiv 0 \pmod{p}$. (Указание:
используйте упражнение~\ref{exercise:kommutatornoe-sparivanie}(ii).)
\item[(v)] Покажите, что при подходящем выборе $z$ описанная
в пункте~(iv) ситуация не реализуется ни для каких
$H$ и $u\in Z(H)$. (Указание: выберите в качестве
$$z\in W'=\Lambda^2 W$$
неразложимый бивектор.)
\item[(vi)] Докажите,
что элемент $\beta\in H^2(G, \Z/p\Z)$ ненулевой. (Указание: воспользуйтесь
упражнением~\ref{exercise:kommutatornoe-sparivanie}(iii).)
\item[(vii)]
Покажите, что образ $\alpha$ элемента $\beta$
относительно естественного отображения
$$H^2(G,\Z/p\Z)\to H^2(G,\Q/\Z)$$
нетривиален. (Указание:
воспользуйтесь пунктом~(vi) и упражнением~\ref{exer-nontrivext}(ii)).
Используя пункты~(ii), (iii) и (v), докажите, что
при подходящем выборе $z$ элемент
$\alpha\in H^2(G,\Q/\Z)$ удовлетворяет условию~\eqref{eq:Shafarevich}.
Таким образом,
для любого точного представления~$V$
группы~$G$ над алгебраически замкнутым полем~$k$ характеристики нуль
многообразие~\mbox{$V/G$} не является стабильно рациональным согласно
упражнению~\ref{exercise:factor-po-hitrovyebnutoj-gruppe}.
\end{itemize}
\end{prob}

\newpage
\section{Арифметика двумерных квадрик}
\label{section:2-dim-quadrics}

\subsection{Инварианты квадрик}
Пусть $K$ --- произвольное поле характеристики, не равной $2$.
Выберем четырёхмерное векторное пространство $V$ над $K$ и
невырожденную квадратичную форму $q$ на пространстве~$V$.
Пусть $Q\subset\Pb(V)\cong\Pb^3$ --- двумерная квадрика, заданная формой $q$.

Как и раньше, для расширения полей $K\subset L$ через $X_L$ будет
обозначаться расширение скаляров объекта $X$ с $K$ на $L$. В
частности, определены квадратичная форма~$q_L$ на векторном
пространстве $V_L=L\otimes_K V$ над $L$, а также квадрика $Q_L$ над
полем~$L$.

\begin{defin}\label{defin:discrquadr}
Выберем базис в $V$ над $K$ и рассмотрим определитель матрицы,
соответствующей форме $q$. Заметим, что его класс в факторгруппе
$K^*/(K^*)^2$ не зависит от выбора базиса, а также не зависит от
выбора $q$ при фиксированной квадрике $Q$ из-за чётности
размерности. Этот класс в $K^*/(K^*)^2$ будем называть {\it
дискриминантом} квадрики $Q$ и будем обозначать через $\disc(Q)$.
\end{defin}

Из определения~\ref{defin:discrquadr} видно, что дискриминант
корректно определён для любой гладкой
чётномерной квадрики.

\bigskip

\begin{prob}{\bf Квадрики с тривиальным дискриминантом}
\label{prob-section}

Предположим, что $\disc(Q)=1$.
\begin{itemize}
\item[(i)]
Покажите, что
квадрика $Q$ {\it пфистерова}, то есть
задаётся уравнением
$$x^2-ay^2-bz^2+abt^2=0\,$$
для некоторых элементов~$a,b\in\nlb K^*$.
\item[(ii)]
Пусть $C_0\subset \Pb^2$ --- коника над полем $K$, заданная уравнением
$$u^2-av^2-bw^2=0\,,$$
где $a,b\in K^*$ определяются по квадрике $Q$ так, как описано
в пункте~(i). Другими словами, $C_0$ является пересечением
квадрики $Q$ с плоскостью~\mbox{$t=0$}.
Докажите, что
существует определённая над полем $K$
проекция~\mbox{$\pi_0\colon Q\to C_0$}, слои которой являются прямыми
относительно вложения~\mbox{$Q\subset\Pb^3$}
(ср. с упражнением~\ref{exer-quaternion-alg}(ii)),
и при этом ограничение отображения $\pi_0$ на
плоское сечение~\mbox{$C_0\subset Q$} является тождественным отображением.
(Указание: пусть $\pi_0$ переводит точку~\mbox{$(x:y:z:t)$}
в~\mbox{$(u:v:w)$}, где $w=z^2-at^2$ и
$$u+\sqrt{a}\,v=(x+\sqrt{a}\,y)(z+\sqrt{a}\,t)\,.$$
Для проверки корректности этого определения можно
воспользоваться
тем, что отображение нормы $x+\sqrt{a}\,y\mapsto x^2-ay^2$ является
гомоморфизмом мультипликативных групп.
Слой над точкой $(u:v:w)$ на конике $C_0$ является прямой на $Q$,
проходящей через точки $(u:v:w:0)$ и $(av:u:0:-w)$,
ср. с упражнением~\ref{prob-quatconic}(i).)
\item[(iii)]
Рассмотрим инволюцию $\iota_0\colon Q\to Q$,
заданную формулой
$$(x:y:z:t)\mapsto (x:y:z:-t)\,.$$
Проверьте, что композиция $\pi_0'=\pi_0\circ\iota_0$
обладает всеми свойствами проекции~$\pi_0$, перечисленными
в пункте~(ii). Таким образом, на квадрике $Q$ имеется два семейства
прямых, определённых над $K$. Пользуясь тем, что на
$Q$ нет других семейств прямых (даже и над алгебраическим замыканием поля $K$),
покажите, что для данной коники $C_0\subset Q$
проекции $\pi_0$ и $\pi_0'$ --- единственные морфизмы,
обладающие этими свойствами.
В частности,
выбор одного из двух семейств прямых на $Q$ (или, что то же самое,
на $Q_{\bar K}$)
равносилен выбору одного из двух морфизмов
$\pi_0$ и $\pi_0'$.
\item[(iv)]
Зафиксируем одно из двух семейств прямых на $Q$.
Докажите, что для любого гладкого плоского сечения $C$ квадрики $Q$
это определяет изоморфизм~\mbox{$C\cong C_0$}. (Указание:
по пункту (iii) мы можем считать, что выбрана одна из двух
проекций $\pi_0$ и $\pi_0'$. Её ограничение на $C$
даёт требуемый изоморфизм.)
\item[(v)]
Пусть $C$ --- гладкое плоское сечение квадрики $Q$.
Докажите, что
существуют ровно две
тождественные на $C$ проекции $\pi, \pi'\colon Q\to C$,
слои которых являются прямыми (ср. с пунктом~(iii)).
\item[(vi)]
Зафиксируем одно из двух семейств прямых на $Q$.
Докажите, что это определяет изоморфизм $Q\cong C\times C$ многообразий
над $K$, где $C$ --- произвольное гладкое плоское сечение
квадрики $Q$.
(Указание: соответствующие проекции~\mbox{$Q\to C$}
задаются как в пункте~(v).)
\item[(vii)] Докажите непосредственно, что
слои проекции $C\times C\to C$ изоморфны $\P^1$.
(Указание: диагональ $C\to C\times C$ определяет сечение
семейства многообразий Севери--Брауэра $C\times C\to C$.
Теперь можно воспользоваться упражнением~\ref{exer-SBforms}(i).)
На самом деле, можно показать, что $\P^1$-расслоение $C\times C\to C$ является проективизацией $\P(V)\to C$ расслоения $V$ ранга $2$ на $C$, определенного следующим образом:
$$
V={\rm Ker}\big(H^0(C,T_C)\otimes \OO_C\to T_C\big)\,,
$$
где $T_C$ является касательным расслоением на $C$.
\item[(viii)]
Докажите, что для произвольной гладкой квадрики
$Q\subset \P^3$ над $K$ следующие условия равносильны:
дискриминант квадрики $Q$ тривиален; квадрика $Q$ пфистерова;
имеется изоморфизм $Q\cong C\times C$ многообразий над $K$,
где $C$~--- некоторая гладкая коника над $K$.
(Указание: равносильность первого и второго условий
следует из пункта~(i).
Равносильность второго и третьего условий следует из пункта~(vi).)
\end{itemize}
\end{prob}

\bigskip
Основные утверждения из упражнения~\ref{prob-section}
можно также интерпретировать в терминах кватернионных алгебр.
А именно, пусть дана кватернионная алгебра~$A(a,b)$,
см. упражнение~\ref{exer-quaternion-alg}, где $a,b\in K^*$.
Тогда гиперповерхность в $\P(A)$, заданная уравнением
$$
\Nrd(\alpha)=0,\quad \alpha\in A\,,
$$
см. упражнение~\ref{exer-rednorm},
является пфистеровой квадрикой $Q$.
При этом коника $C_0$ канонически соответствует алгебре $A$,
см. упражнение~\ref{prob-quatconic}(i), а два семейства прямых на $Q$
с базой $C_0$ являются семействами правых и левых идеалов в $A$
(морфизм~$\pi_0$ соответствует правым идеалам, а $\pi_0'$ --- левым).

\bigskip
Из упражнения~\ref{prob-section}(iv) следует, что если
дискриминант $\disc(Q)$ тривиален, то корректно определена
коника $C$ над полем $K$, являющаяся гладким плоским сечением квадрики $Q$.

\begin{defin}\label{defin:Clifford}
В этом случае класс $b(C)$ коники $C$ в группе Брауэра $\Br(K)$
(см. определение~\ref{defin-b}) будем называть
{\it инвариантом Клиффорда} квадрики $Q$ и будем обозначать через $\cl(Q)$.
\end{defin}

На самом деле, инвариант Клиффорда корректно определён для
любой гладкой чётномерной квадрики с тривиальным дискриминантом.
Таким образом, возникает последовательность инвариантов
квадрик (чётность размерности, дискриминант,
инвариант Клиффорда),
в которой каждый следующий инвариант определён
при условии тривиальности предыдущих.
Эти инварианты
принимают значения в группах~\mbox{$\Z/2\Z$},~\mbox{$K^*/(K^*)^2$} и
$\Br(K)_2$, соответственно.
По упражнениям~\ref{exer-Kummer}(ii)
и~\ref{exer-brtor}(iii)
данные группы изоморфны группам
когомологий Галуа $H^n(G_K,\Z/2\Z)$ при~\mbox{$n=0,1,2$}.

Существует обобщение этих инвариантов на случай произвольных
значений $n$. При $n=3$ такой инвариант был построен
Дж.\,Арасоном в~\cite{Arason}.
Для некоторых конкретных квадрик инварианты строятся явным образом,
см.~\cite[\S18]{Karpenko}.
Построение инвариантов в общем случае основано
на двух доказанных теперь гипотезах Милнора. Одна из них
связывает квадратичные формы над полем $K$ с~\mbox{$\KM$-группами} Милнора
$\KM^M_n(K)$ (см. определение~\ref{defin:K-group} для случая $n=2$)
и была доказана А.\,Вишиком, В.\,Воеводским и Д.\,Орловым
в~\cite{OVV}.
Другая гипотеза утверждает,
что естественный гомоморфизм
$$\KM^M_n(K)/2\to H^n(G_K,\Z/2\Z)$$
является изоморфизмом для любого натурального числа $n$
и была доказана В.\,Воеводским в~\cite{Voevodsky}.
Отметим, что и данное утверждение,
и теорема Меркурьева--Суслина (см. обсуждение после
упражнения~\ref{exer-norm-symbol})
являются частными случаями более общей гипотезы
Блоха--Като, доказанной недавно В.\,Воеводским и М.\,Ростом.

\bigskip
\begin{prob}{\bf Нетривиальность инварианта Клиффорда}
\label{prob-nonsection}
\hspace{0cm}
\begin{itemize}
\item[(i)]
Пусть $F=K\big(\sqrt{a}\big)$ --- произвольное квадратичное расширение
поля~$K$. Предположим, что
квадрика $Q_F$ содержит прямую, и при этом $Q(K)=\varnothing$.
Покажите, что тогда
форма $q$ является тензорным произведением
двух квадратичных форм на двумерных
векторных пространствах. В частности, дискриминант~$\disc(Q)$
тривиален. (Указание: рассмотрим
в $V_F$ двумерное подпространство~$U$ над $F$,
соответствующее прямой на $Q_F$.
Тогда $U$ изотропно относительно формы $q_F$, и
$U\cap\sigma(U)=0$, где~$\sigma$ обозначает нетривиальный
автоморфизм поля $F$ над $K$.
Покажите, что в $U$ существует такой
базис $\{e,f\}$ над $F$, что
$$b\big(e,\sigma(f)\big)=0\,,$$
где $b$ ---
симметрическая билинейная форма на $V_F$,
ассоциированная с $q_F$. Затем рассмотрите базис
$$
\left\{e+\sigma(e),\,(e-\sigma(e))/\sqrt{a},\,f+\sigma(f),\,(f-\sigma(f))/\sqrt{a}\,\right\}
$$
в $V$ над $K$ и соответствующую матрицу Грама.)
Отметим, что точно так же доказывается обобщение этого факта
для квадратичных форм от большего числа переменных.
Пусть $q$ --- невырожденная квадратичная форма
от $2n+2$ переменных, и $Q$ --- соответствующая квадрика
в $\P^{2n+1}$.
Предположим, что
квадрика $Q_F$ содержит подпространство $\P^n$,
и при этом $Q(K)=\varnothing$.
Тогда форма $q$ (тензорно) делится на квадратичную
форму на двумерном векторном пространстве,
заданную диагональной матрицей
$\mathrm{diag}(1,-a)$.
\item[(ii)]
Рассмотрим поле
$$L=K\big(\sqrt{\disc(Q)}\big)\,,$$
совпадающее с полем $K$ в случае тривиального дискриминанта
$\disc(Q)$ и являющееся
его квадратичным расширением в случае нетривиального дискриминанта.
Пусть коника $C$ над полем $L$ --- гладкое плоское сечение
квадрики $Q_L$.
Докажите, что
$C(L)\neq\varnothing$ тогда и только тогда,
когда~\mbox{$Q(K)\neq\varnothing$}.
(Указание:
в случае тривиального дискриминанта воспользуйтесь
упражнением~\ref{prob-section}(ii),(iv).
Если дискриминант $\disc(Q)$ нетривиален и $C(L)\neq\varnothing$,
воспользуйтесь
упражнением~\ref{prob-section}(v) и докажите, что квадрика $Q_L$
содержит прямую. После этого примените пункт~(i).
Для доказательства обратной импликации заметьте, что
если $Q(K)\neq\varnothing$, то и $Q_L(L)\neq\varnothing$,
так что всё сводится к уже разобранному случаю тривиального
дискриминанта.)
\item[(iii)]
В обозначениях пункта~(ii) докажите, что $\cl(Q_L)=0$
в группе~$\Br(L)$ тогда и только тогда, когда~\mbox{$Q(K)\ne\varnothing$}.
(Указание:
воспользуйтесь пунктом~(ii) и упражнением~\ref{prob:SB-s-tochkoi}.)
\end{itemize}
\end{prob}

\subsection{Геометрический смысл инвариантов квадрик}

\begin{prob}{\bf Геометрический смысл дискриминанта}
\label{prob-discr-geom-interpretation}
\hspace{0cm}
\begin{itemize}
\item[(i)] Рассмотрим произвольную прямую $l$ на $Q_{K^{sep}}$
и элемент $g$ группы Галуа $G_K$.
Докажите, что пересечение $g(l)\cap l$ является точкой,
то есть $l$ и $g(l)$ принадлежат разным семействам
прямых на $Q_{K^{sep}}$,
тогда и только тогда, когда
$$g\big(\sqrt{\disc(Q)}\big)=-\sqrt{\disc(Q)}\,.$$
(Указание: зафиксируем однородные координаты $(x_0:x_1:x_2:x_3)$ на $\P^3$,
в которых $Q$ задаётся уравнением
$$\sum\limits_{i=0}^3 a_ix_i^2=0\,.$$
Поскольку группа (проективных) автоморфизмов квадрики $Q_{K^{sep}}$
действует транзитивно на множестве прямых на $Q_{K^{sep}}$, мы можем
считать, что прямая $l$ проходит через точки
\mbox{$(\sqrt{a_1}:\sqrt{-a_0}:0:0)$} и~\mbox{$(0:0:\sqrt{a_3}:\sqrt{-a_2})$}.
Заметим, что~\mbox{$\disc(Q)=(-a_0)a_1(-a_2)a_3$}.
Таким образом,
равенство $g\big(\sqrt{\disc(Q)}\big)=-\sqrt{\disc(Q)}$ выполняется
тогда и только тогда, когда
$g$ меняет знаки у нечётного количества элементов из
$\sqrt{-a_0}$, $\sqrt{a_1}$, $\sqrt{-a_2}$ и $\sqrt{a_3}$.
Простой перебор случаев показывает, что это условие равносильно тому,
что пересечение $g(l)\cap l$ является точкой.)
\item[(ii)] Воспользуйтесь пунктом (i) для того, чтобы
решить упражнение~\ref{prob-nonsection}(ii), не используя
упражнение~\ref{prob-nonsection}(i).
(Указание: предположим, что дискриминант~$\disc(Q)$ нетривиален
и $C(L)\neq\varnothing$. Тогда по
упражнению~\ref{prob-section}(v) квадрика~$Q_L$
содержит некоторую прямую $l$. Пусть $\sigma$ --- нетривиальный
элемент группы Галуа $\Gal(L/K)$.
Тогда
$$\sigma\big(\sqrt{\disc(Q)}\big)=-\sqrt{\disc(Q)}\,.$$
По пункту~(i) это даёт $K$-точку $l\cap\sigma(l)$ на $Q$.)
\end{itemize}
\end{prob}

\bigskip
Теория Куммера (см. упражнение~\ref{exer-Kummer}(ii))
отождествляет группы~$K^*/(K^*)^2$ и~$\Hom(G_K,\Z/2\Z)$.
По упражнению~\ref{prob-discr-geom-interpretation}(i)
гомоморфизм из $G_K$ в $\Z/2\Z$,
соответствующий дискриминанту~$\disc(Q)$,
определяется действием группы Галуа $G_K$ на двух
эффективных образующих группы
Пикара
$$\Pic(Q_{K^{sep}})\cong\Z\oplus\Z\,.$$
Аналогичный факт имеет
место для квадрик произвольной чётной размерности после замены
группы Пикара на группу Чжоу алгебраических циклов средней
размерности на квадрике.

\bigskip
Рассмотрим ортогональный грассманиан формы $q$.
Другими словами, рассмотрим подмногообразие ${\rm G}(q)$ в
грассманиане~${\rm G}(V,2)$,
параметризующее двумерные изотропные подпространства в $V$
относительно формы $q$, то есть прямым на квадрике~$Q$.
Над
сепарабельным замыканием $K^{sep}$ многообразие
${\rm G}(q)_{K^{sep}}$ изоморфно объединению двух непересекающихся
копий $\Pb^1$.
Напомним, что задать морфизм многообразий $B\to \G(V,2)$
над~$K$~--- это то же самое,
что задать подмноогообразие
$F\subset\P(V)\times B$ над $K$,
для которого соответствующая проекция в $B$ плоская
и её слои являются прямыми в~$\P(V)\cong\P^3$.
В частности, тождественный морфизм грассманиана
соответствует многообразию инцидентности $I\subset \P(V)\times\G(V,2)$.

\begin{prob}{\bf Геометрический смысл инварианта Клиффорда}
\label{prob-Clifford-geom-interpretation}
\hspace{0cm}
\begin{itemize}
\item[(i)]
Покажите, что дискриминант $\disc(Q)$ тривиален
тогда и только тогда, когда грассманиан $\G(q)$
приводим над $K$. В этом случае имеется
изоморфизм
$$\G(q)\cong C\sqcup C$$
многообразий над $K$,
где коника $C$ является произвольным гладким плоским
сечением квадрики $Q$ над $K$.
(Указание: при условии, что дискриминант тривиален,
воспользуйтесь упражнением~\ref{prob-section}(v).
Рассмотрите графики морфизмов $\pi$ и $\pi'$,
вложенные в $Q\times C\subset \P^3\times C$, чтобы отождествить
$C\sqcup C$ с грассманианом.
Для обратной импликации, используя изоморфизм
$$\G(q)_{K^{sep}}\cong \P^1\sqcup\P^1\,,$$
покажите, что
$$\G(q)\cong C'\sqcup C''\,,$$
где $C'$ и $C''$ --- гладкие коники над $K$.
Далее рассмотрите ограничение $I_{C'}$ многообразия
инцидентности $I$ на $C'\subset \G(q)$.
Поскольку $I_{C'}$ изоморфно проецируется на $Q\subset \P(V)$,
это задаёт морфизм $Q\to C'$.
Аналогично, ограничение~$I$ на $C''$ задаёт проекцию из $Q$ на $C''$.
Поскольку слои обеих проекций являются прямыми, мы видим,
что коники $C'$ и $C''$ изоморфны друг другу,
а также изоморфны произвольному гладкому плоскому сечению $C$ квадрики $Q$,
ср. с упражнением~\ref{prob-section}.
Таким образом, возникает изоморфизм $Q\cong C\times C$
(то, что это изоморфизм, достаточно проверить над алгебраическим замыканием поля~$K$).
Теперь остаётся воспользоваться упражнением~\ref{prob-section}(viii).)
\item[(ii)] Выведите из пункта~(i)
упражнение~\ref{prob-discr-geom-interpretation}(i).
(Указание: чтобы доказать, что два гомоморфизма из какой-либо
группы в $\Z/2\Z$ совпадают, достаточно проверить, что у них одинаковые ядра.
Для этого
можно воспользоваться пунктом~(i)
для соответствующих квадратичных расширений поля~$K$.)
\item[(iii)]
Предположим, что дискриминант $\disc(Q)$ нетривиален.
Пусть $L=K\big(\sqrt{\disc(Q)}\big)$, и~$C$~--- коника над
$L$, представляющая инвариант Клиффорда $\cl(Q_L)\in\Br(L)$.
Любое многообразие $Y$ над $L$
имеет структуру многообразия над $K$, возникающую
из композиции морфизмов
\begin{equation}
\label{eq:tupoe-ogranichenie}
Y\to\Spec(L)\to\Spec(K)\,.
\end{equation}
Рассмотрим на $C$ такую структуру многообразия над $K$.
Докажите, что имеется изоморфизм
$C\cong \mathrm{G}(q)$ многообразий над $K$.
(Указание:
по упражнению~\ref{prob-section}(v) проекция $\pi\colon Q_L\to C$
задаёт морфизм $j\colon C\to {\mathrm G}(q)_L$
многообразий над $L$.
Поскольку ограничение скаляров, происходящее из
морфизма~\eqref{eq:tupoe-ogranichenie}, является левым сопряжённым функтором
к расширению скаляров, это задаёт морфизм $f\colon C\to G(q)$
многообразий над $K$. Для того, чтобы доказать, что~$f$ является
изоморфизмом, достаточно доказать,
что изоморфизмом является его расширение скаляров
$$f_L\colon C_L\to G(q)_L\,.$$
Пусть $\sigma$ обозначает
нетривиальный элемент группы Галуа $\Gal(L/K)$.
Тогда
$$C_L\cong C\sqcup \sigma_*C\,,$$
где $\sigma_*C$ определяется так же, как в примере~\ref{examp-spusk}(v).
При этом отображение~$f_L$ соответствует паре
морфизмов~$j$ и~$\sigma_*j$, а морфизм
$$
\sigma_*j\colon\sigma_*C\to \sigma_*\big(G(q)_L\big)\cong G(q)_L
$$
задаётся проекцией
$$
\sigma_*\pi\colon \sigma_*(Q_L)\cong Q_L\to \sigma_*C\,.
$$
Так как слои проекции $\sigma_*\pi$ получаются применением $\sigma$
к слоям исходной проекции $\pi$,
то из упражнения~\ref{prob-discr-geom-interpretation}(i) следует,
что слои проекций $\pi$ и $\sigma_*\pi$ принадлежат разным семействам.
Таким образом, $f_L$ является изоморфизмом.)
\end{itemize}
\end{prob}

\subsection{Вырождения квадрик}

Пусть $K$ является полным полем дискретного нормирования $v$
с кольцом дискретного нормирования $\OO_K$ и
совершенным полем вычетов $\kappa$ характеристики, не равной~$2$.
Пусть задана гладкая двумерная квадрика $Q\subset\P^3$ над полем~$K$.

\bigskip

\begin{prob}{\bf Неразветвлённость инварианта Клиффорда}
\label{exer-Clifunr}

Предположим, что у квадрики $Q$ существует проективная модель
$\mathcal Q\subset\P^3$ над~$\OO_K$,
согласованная с исходным проективным вложением квадрики $Q$
(см. раздел~\ref{subsection:models}),
замкнутый слой которой над полем $\kappa$
является либо гладкой квадрикой, либо конусом
над гладкой коникой (ср. с понятием хорошей модели коники,
введенным перед упражнением~\ref{prob-conic-models}).
Положим
$$L=K(\sqrt{\disc(Q)})\,.$$
Докажите, что
$$\res\big(\cl(Q_{L})\big)=0\,.$$
(Указание: рассмотрите класс~\mbox{$\delta\in\Br(K)$}
произвольного гладкого плоского сечения квадрики~$Q$. Покажите,
что~\mbox{$\cl(Q_{L})$} является
образом элемента $\delta$ при естественном
отображении
$$\Br(K)\to\Br(L)\,.$$
Далее, основная идея рассуждения заключается в том,
что при данных условиях на вырождение модели $\QQ$ у нее
существует гладкое плоское сечение с гладким замкнутым слоем.
А именно, так как $2$ обратимо в $\OO_K$,
форму $q$ можно записать в диагональном виде:
$$q(x)=\sum_{i=0}^3a_i x_i^2, \quad a_i\in \OO_K\,.$$
При этом по условию на вырождение можно считать, что $v(a_i)=0$
при $1\leqslant i\leqslant 3$.
Теперь в качестве $\delta$ рассмотрите класс коники,
заданной уравнениями $q(x)=0$, $x_0=0$.
Наконец, воспользуйтесь упражнениями~\ref{prob-quatconic}(i)
и~\ref{exer-vychet-dlya-cycl}(iv)
для доказательства того, что $\res(\delta)=0$,
и примените коммутативность диаграммы
из упражнения~\ref{exer-Brcompldv}(iv).)
\end{prob}

\newpage

\section{Нерациональные двойные накрытия $\Pb^3$}
\label{section:Artin-Mumford}

В этой главе мы применяем
неразветвлённую группу Брауэра для доказательства нерациональности
некоторых унирациональных многообразий размерности~$3$.
Эти многообразия получаются как
двойные накрытия проективного
пространства~$\P^3$ с ветвлением в (особых) квартиках. Первоисточником для
соответствующего
примера является статья М.\,Артина и Д.\,Мамфорда~\cite{ArtinMumford},
однако наше изложение основывается на более простом подходе
М.\,Гросса (см.~\mbox{\cite[Appendix]{Aspinwall-i-td}}).

Как и раньше, препятствием к стабильной рациональности многообразия
$X$ будет нетривиальность неразветвлённой группы Брауэра
$\Br^{nr}\big(k(X)\big)$ поля $k(X)$. При этом в нашем случае
для упрощения вычислений мы воспользуемся
(см. упражнение~\ref{prob:netrivialnost-AM}) тем, что группа
$\Br^{nr}\big(k(X)\big)$ совпадает
с неразветвлённой группой Брауэра $\Br^{nr}(X)$ многообразия $X$.
Группа $\Br^{nr}(X)$ удобнее группы $\Br^{nr}\big(k(X)\big)$
тем, что позволяет работать только с нормированиями поля $k(X)$,
соответствующими дивизорам $D\subset X$.

\subsection{Ещё раз о неразветвлённой группе Брауэра}

Пусть $X$ --- многообразие размерности не меньше $2$ над полем характеристики нуль.
Напомним, что изолированная особая точка $x$ многообразия $X$ называется
{\it обыкновенной двойной}, если некоторая её открытая окрестность $U$ в топологии Зарисского
изоморфна гиперповерхностной особенности с невырожденным гессианом.
Другими словами, имеется изоморфизм
$$(U, x)\cong (\{f(t)=0\},0)\subset \A^{n+1}\,,$$
где $f(0)=0$, линейная часть многочлена $f$ нулевая,
а квадратичная часть $f$ является невырожденной квадратичной формой.
Tакую особенность можно разрешить одним раздутием $\widetilde{X}\to X$
с центром
в $x$, при котором вместо $x$ в~$X$ ``вклеивается'' гладкая
квадрика~$D$. По критерию Серра
(см.~\cite[Предложение~II.8.23]{Hartshorne}) многообразие, имеющее
лишь обыкновенные двойные особенности, является нормальным.

В дальнейшем нам понадобится следующий общий факт о неразветвлённой
группе Брауэра многообразия (см. определение~\ref{defin:unramified-Br-X}).
Пусть $D$
является гладким неприводимым приведённым дивизором в гладком неприводимом
многообразии~$\widetilde{X}$ над полем $k$ характеристики нуль. Тогда для
любого элемента
$$\alpha\in\Br^{nr}\big(\widetilde{X}\setminus D\big)\subset
\Br\big(k(\widetilde{X})\big)$$
его вычет
$$\res_D(\alpha)\in\nlb\Hom(G_{k(D)},\Q/\Z)$$
соответствует неразветвлённому циклическому накрытию дивизора $D$.
Это можно доказать, используя интерпретацию неразветвлённой группы Брауэра
в терминах этальных когомологий пучка корней из единицы
(см. разделы~\ref{subsection:etale-Kummer} и~\ref{subsection:etale-Brauer}),
а также последовательность локализации для этальных когомологий с конечными
коэффициентами (см.~\cite[VI.5.4(b)]{Mil}).

\bigskip
\begin{prob}{\bf Неразветвлённая группа Брауэра
и обыкновенные двойные особенности}\label{exer-unrBrsing}

Пусть $X$ --- неприводимое проективное многообразие над
полем~$k$ характеристики нуль, причём все особые точки многообразия~$X$
обыкновенные двойные. Докажите, что естественное вложение
$$\Br^{nr}\big(k(X)\big)\hookrightarrow \Br^{nr}(X)$$
является
изоморфизмом. (Указание: рассмотрите разрешение особенностей
$\widetilde{X}$
многообразия $X$ и воспользуйтесь аналогичным утверждением для
гладких многообразий, а также тем, что у гладких квадрик нет
нетривиальных неразветвлённых накрытий.)
\end{prob}

\subsection{Семейства двумерных квадрик}\label{subsec:famquadr}

Пусть теперь $B$ --- неприводимое нормальное многообразие над полем
$k$ характеристики нуль, и пусть $K$ обозначает поле $k(B)$
рациональных функций на $B$. Пусть~$\QQ$~--- семейство квадрик в
$\Pb^3$ с базой $B$, слой которого над схемной общей
точкой в~$B$ является гладкой квадрикой $Q$ над $K$.
Таким образом,
$$\mathcal{Q}\subset\P^3_B=\P^3\times B$$
является проективной моделью квадрики $Q$,
согласованной с её проективным вложением над~$K$
(см. раздел~\ref{subsection:models}).
Можно показать, что $\QQ$ ---
подмногообразие в~$\Pb^3\times B$, являющееся нулями сечения
линейного расслоения вида
$$\OO_{\Pb^3}(2)\boxtimes\LL\,,$$
где $\LL\in\Pic(B)$, причём морфизм $\QQ\to B$
плоский, а его общий слой гладкий. Здесь ${}\boxtimes{}$ обозначает
тензорное произведение на $\Pb^3\times B$ обратных образов
линейных расслоений относительно проекций $\Pb^3\times B\to\Pb^3$
и $\Pb^3\times B\to B$.

Пусть $\Delta^i\subset B$ --- подмногообразие, состоящее из точек,
над которыми слой $\QQ$ является квадрикой коранга не меньше $i$. В
частности, $\Delta^1$ является дивизором на~$B$,
поскольку $\Delta^1$ локально задаётся одним уравнением (что это за
уравнение?). Предположим, что дивизор $\Delta^1$ непуст (для
проективного $B$ это равносильно тому, что семейство квадрик $\QQ$
нетривиально), и что ${\rm codim}_B(\Delta^2)\geqslant 2$.
Предположим также, что дискриминант $\disc(Q)$ (см. определение~\ref{defin:discrquadr}) нетривиален.

Положим
$$L=K\big(\sqrt{\disc(Q)}\big)\,,$$
и пусть $X$
является нормализацией $B$ в $L$. Таким образом, $X$ ---
неприводимое нормальное многообразие над полем $k$ с полем
рациональных функций $k(X)\cong L$,
и имеется сюръективный конечный морфизм $X\to B$
степени $2$, соответствующий расширению полей $K\subset L$ над общей точкой.

\bigskip
\bigskip
\bigskip
\bigskip

\begin{prob}{\bf Особенности двойного накрытия}\label{exer-sing}
\hspace{0cm}
\begin{itemize}
\item[(i)]
Докажите, что $\Delta^1$ является дивизором ветвления сюръективного
конечного морфизма $X\to B$.
\item[(ii)]
Предположим, что многообразие $B$ гладкое, а все особенности
дивизора~$\Delta^1$ обыкновенные двойные. Покажите, что все
особенности многообразия~$X$ также обыкновенные двойные. (Указание:
рассмотрите уравнение, локально задающее $X$ как дивизор в
$B\times{\mathbb A}^1$.)
\end{itemize}
\end{prob}

\bigskip

\begin{prob}{\bf Нетривиальность неразветвлённой группы Брауэра}
\label{prob:netrivialnost-AM}

Предположим, что $B$ --- гладкое проективное многообразие,
все особенности дивизора $\Delta^1$ --- обыкновенные двойные,
и $Q(K)=\varnothing$. Положим (см. определение~\ref{defin:Clifford})
$$\alpha=\cl(Q_L)\in \Br(L)\,.$$
\begin{itemize}
\item[(i)]
Покажите, что $\alpha\ne 0$.
(Указание: воспользуйтесь упражнением~\ref{prob-nonsection}(iii).)
\item[(ii)]
Покажите, что для любого неприводимого приведённого дивизора
$D$ на $X$ имеется равенство $\res_D(\alpha)=0$,
то есть
$$\alpha\in\Br^{nr}(X)\,.$$
(Указание: рассмотрите нормирование поля~$K$,
задаваемое образом дивизора~$D$ в~$B$,
и воспользуйтесь тем, что имеется оценка
$$\mathrm{codim}_B(\Delta^2)\geqslant 2\,$$
на коразмерность $\Delta^2$ в $B$,
а также упражнением~\ref{exer-Clifunr}.)
\item[(iii)]
Покажите, что
$$\alpha\in\Br^{nr}(L)\,.$$
(Указание: воспользуйтесь
упражениями~\ref{exer-sing}(ii),~\ref{exer-unrBrsing} и пунктом (ii).)
\item[(iv)]
Покажите, что многообразие $X$ не является стабильно рациональным.
\end{itemize}
\end{prob}

\subsection{Построение геометрического примера}

Пусть $V$ --- векторное пространство размерности $n\geqslant 3$ над
полем $k$ характеристики, не равной $2$. Для каждого
$i\geqslant 0$ определим
$$\Sigma^i\subset\P\big(\mathrm{Sym}^2(V^{\vee})\big)$$
как подмножество,
состоящее из точек,
соответствующих квадратичным формам на $V$
коранга не меньше $i$.
Здесь $V^{\vee}$ обозначает пространство, двойственное к~$V$,
а~\mbox{$\mathrm{Sym}^2$}
обозначает симметрический квадрат векторного пространства;
таким образом,~\mbox{$\mathrm{Sym}^2(V^{\vee})$} является пространством
всех симметрических билинейных форм на~$V$.
В частности, $\Sigma^1$ является дивизором,
заданным уравнением $\det(q)=0$, где~\mbox{$q\in\mathrm{Sym}^2(V^{\vee})$}.

\bigskip
\begin{prob}{\textbf{Особенности многообразия вырожденных
квадрик}}
\label{prob:Sigma-1-2}
\hspace{0cm}
\begin{itemize}
\item[(o)]
Рассмотрим естественное действие группы $\GL(V)$ на
$\P\big(\mathrm{Sym}^2(V^{\vee})\big)$. Подмногообразия $\Sigma^i$ инвариантны
относительно этого действия. Проверьте, что
при любом $i\geqslant 0$ множество
$\Sigma^{i}\setminus\Sigma^{i+1}$ состоит из одной
$\GL(V)$-орбиты, плотной в~$\Sigma^i$.
(Указание: приведите квадратичную форму к диагональному виду.)
\item[(i)]
Проверьте, что $\Sigma^1$ является (приведённым)
неприводимым дивизором степени~$n$.
(Указание: подмногообразие $\Sigma^1$
задаётся одним уравнением~\mbox{$\det(-)=0$}.
Неприводимость~$\Sigma^1$ следует из пункта~(o) и
неприводимости многообразия~\mbox{$\GL(V)$}.
Наконец, уравнение $\det(-)=0$ задаёт приведённый
дивизор, что можно проверить, рассмотрев его дифференциал
в удобной точке $q\in\Sigma^1$.)
\item[(ii)]
Докажите, что
$$\codim_{\Sigma^1}(\Sigma^2)=2\,,$$
и $\Sigma^2$ является множеством особенностей многообразия~$\Sigma^1$.
(Указание: пусть
$$I_2=\mathrm{diag}(0,0,\underbrace{1,\ldots,1}_{n-2})\,.$$
Согласно пункту~(o) достаточно найти размерность $\Sigma^2$
в точке~$I_2$ и проверить, что $\Sigma^1$ особо в этой точке.
Первое делается при помощи явного описания
стабилизатора матрицы~$I_2$. Второе следует из явного вычисления
частных производных определителя.)
\item[(iii)]
Пусть $\mathfrak{H}_q$ --- матрица вторых производных
функции $\det(-)$ в точке
$$q\in\P\big(\mathrm{Sym}^2(V^{\vee})\big)\,.$$
Докажите, что $\mathfrak{H}_q$ имеет ранг $3$, если
$q$ --- общая точка многообразия $\Sigma^2$.
(Указание: достаточно проверить это в точке на $\Sigma^2$,
соответствующей матрице~$I_2$.)
\end{itemize}
\end{prob}

\bigskip
Пусть $V$ --- четырёхмерное векторное пространство над
полем комплексных чисел~$\Cb$. Пусть
$B$~--- общее трёхмерное проективное подпространство в
проективном пространстве
$$\P\big(\mathrm{Sym}^2(V^{\vee})\big)\cong\P^9\,,$$
пятимерное многообразие $\QQ$ --- соответствующее расслоение на
квадрики в $\P(V)$ с базой $B$, а $Q$ --- слой $\QQ$ над схемной общей
точкой базы $B$. В частности, $Q$ является гладкой квадрикой над
полем~$K=\Cb(B)$.
Для алгебраического многообразия~$U$ над полем $\Cb$
через $H^i\big(U(\Cb),\Z\big)$
мы будем обозначать $i$-ые когомологии Бетти
соответствующего топологического пространства
$U(\Cb)$ с классической топологией. Напомним, что имеется произведение
$$
H^i\big(U(\Cb),\Z\big)\times H^j\big(U(\Cb),\Z\big)\to H^{i+j}\big(U(\Cb),\Z\big)\,.
$$
Если $U$ неприводимое и гладкое, то для каждого (неприводимого) подмногообразия~\mbox{$Z\subset U$} коразмерности $i$ определён его класс $[Z]\in H^{2i}\big(U(\Cb),\Z\big)$. Наконец, если~$U$ ещё и проективно, то имеется канонический изоморфизм $H^{2d}\big(U(\Cb),\Z\big)\cong \Z$, где~$d$~--- размерность многообразия $U$.

\bigskip

\begin{prob}{\textbf{Стабильная нерациональность специального
двойного накрытия $\P^3$ с ветвлением в квартике}}
\label{prob:Artin-Mumford}
\hspace{0cm}
\begin{itemize}
\item[(o)] Покажите, что
уравнение $\det(Q)=0$ задаёт в $B$ неприводимый приведённый дивизор
$\Delta^1$ степени~$4$.
Проверьте, что $\Delta^2$ ---
конечный набор точек, совпадающий с множеством
особых точек дивизора~$\Delta^1$,
и особенности дивизора~$\Delta^1$~--- обыкновенные двойные.
Докажите, что дискриминант $\disc(Q)$ нетривиален.
(Указание: все утверждения, кроме последнего,
следуют из упражнения~\ref{prob:Sigma-1-2}
и теоремы Бертини. Для доказательства нетривиальности дискриминанта
посчитайте значение на $\disc(Q)$ нормирования, соответствующего
дивизору~$\Delta^1$.)
\item[(i)] Докажите, что многообразие $\QQ$ является
гладким дивизором бистепени~\mbox{$(1,2)$}
в гладком проективном многообразии
$$W=B\times\P(V)\,.$$
(Указание: воспользуйтесь теоремой Бертини
в следующей форме (см., например,~\mbox{\cite[теорема~17.16]{Harris}}).
Если дан морфизм $\psi\colon M\to\Pb^n$ из гладкого многообразия
$M$, то прообраз общей гиперплоскости
относительно $\psi$ гладок. Эту теорему надо применить
к проекции многообразия инцидентности
$$I\subset\P\big(\mathrm{Sym}^2(V^{\vee})\big)\times\P(V)$$
на первый сомножитель.)
\item[(ii)]
Пусть $\QQ_b$ ---
слой расслоения $\QQ$ над произвольной замкнутой точкой $b\in\nlb B$.
Докажите, что для любого элемента $c\in H^4\big(\QQ(\Cb),\Z\big)$ произведение
$$c\cdot [\QQ_b]\in H^{10}\big(\QQ(\Cb),\Z\big)\cong\Z$$
чётно. (Указание: сначала
докажите аналогичное утверждение для произвольного элемента
$$c'\in H^4\big(W(\Cb),\Z\big)$$
и класса~$\QQ_b$ в $H^8\big(W(\Cb),\Z\big)$,
используя пункт (i).
Затем воспользуйтесь теоремой Лефшеца о гиперплоском сечении, чтобы
найти такой элемент~\mbox{$c'\in H^4\big(W(\Cb),\Z\big)$},
что его ограничение на
$\QQ$ равно $c$.)
\item[(iii)] Докажите, что $Q(K)=\varnothing$.
(Указание: для точки из $Q(K)$
рассмотрите замыкание по Зарисскому графика соответствующего
рационального отображения~\mbox{$B\dasharrow\QQ$} и
воспользуйтесь пунктом~(ii).)
\item[(iv)] Пусть двойное накрытие $X\to B$ построено так же, как выше в разделе~\ref{subsec:famquadr}.
Покажите, что многообразие $X$ не является стабильно рациональным.
(Указание: воспользуйтесь упражнением~\ref{prob:netrivialnost-AM}(iv).)
\item[(v)] Замените в предыдущих рассмотрениях $B\cong\P^3$ на $B\cong\P^2$.
Что изменится?
\end{itemize}
\end{prob}

\subsection{Некоторые конструкции унирациональности}

Данный раздел предполагает некоторые дополнительные знания из алгебраической
геометрии по сравнению с остальными частями книги: в частности,
предполагается, что читатель знаком с такими
понятиями, как
вложение, соответствующее линейной системе,
и индекс пересечения кривых на поверхности.
Сведения об этом можно найти, например,
в учебнике~\cite{Hartshorne}.
Отметим, что в дальнейшем
мы не будем использовать результаты этого раздела.

Пусть $F$ --- поле характеристики, не равной двум,
а $\bar{F}$~--- алгебраическое замыкание поля~$F$.
Пусть $S$ --- гладкая поверхность дель Пеццо степени $2$ над $F$,
то есть гладкая поверхность с обильным антиканоническим
дивизором $-K_S$, для которой~\mbox{$K_S\cdot K_S=2$}.
Тогда поверхность $S_{\bar{F}}$ изоморфна раздутию
семи $\bar{F}$-точек на~$\P^2$;
в частности, поверхность~$S_{\bar{F}}$ рациональна.
Антиканоническая линейная система задаёт морфизм
$$\varphi=\varphi_{|-K_S|}\colon S\to\P^2\,,$$
который является двулистным накрытием
$\P^2$ с ветвлением в гладкой квартике~\mbox{$T\subset\Pb^2$}.
Напомним, что \emph{прямыми} на
поверхности дель Пеццо степени~$2$ (или, в более общем случае,
на поверхности дель
Пеццо степени не меньше~$2$)
называются кривые антиканонической степени $1$.
На $S_{\bar{F}}$ имеется конечное число прямых,
причём они являются гладкими кривыми рода нуль
(чего всякий и ожидает от прямых). Отметим, что на поверхности дель Пеццо степени $1$ общая антиканоническая кривая является эллиптической кривой антиканонической степени $1$; с этим связано наше ограничение на степень поверхности~$S$.
Подробнее про основные свойства поверхностей дель Пеццо
можно прочитать, например, в главе~IV книги~\cite{Manin-KubFormy},
в главе~8 книги~\cite{Dolgachev} или в разделе III.3
книги~\cite{Kollar-RatCurves}.

\bigskip
Следующее упражнение в основном воспроизводит рассуждение из доказательства
теоремы~IV.7.7 в книге~\cite{Manin-KubFormy}.
Более общее и подробное изложение можно найти в~\cite{Testa}.
Пусть $R$ обозначает прообраз на~$S$ кривой
ветвления~\mbox{$T\subset\Pb^2$}.

\begin{prob}{\textbf{Унирациональность поверхности дель Пеццо степени~$2$
с точкой}}
\label{prob:DP2}

Пусть $U\subset S$~---
дополнение к объединению кривой $R$ и всех прямых на $S_{\bar{F}}$.
Согласно сказанному выше, $U$ является открытым по Зарисскому
подмножеством в~$S$.

\begin{itemize}
\item[(o)]
Пусть $l\subset\Pb^2$ --- прямая (определённая над $\bar{F}$).
Проверьте, что полный прообраз
$$\tilde{l}=\varphi^{-1}(l)\subset S_{\bar{F}}$$
приводим
тогда и только тогда, когда $\tilde{l}$ распадается
в объединение двух прямых на $S_{\bar{F}}$.
\item[(i)]
Докажите, что для каждой точки $x\in S(F)$, не лежащей ни на одной
прямой на $S_{\bar{F}}$,
существует кривая
$$C\in |-2K_S-3x|\,,$$
то есть кривая из линейной системы
$|-2K_S|$, имеющая в точке~$x$
особенность кратности не меньше~$3$.
(Указание: непустоту линейной системы достаточно доказать над замыканием
основного поля. Поэтому с самого начала предположим, что
поле $F$ алгебраически замкнуто.
Тогда имеется бирациональный морфизм $\pi\colon S\to\P^2$,
являющийся раздутием семи точек $P_1$, \ldots, $P_7$ на $\P^2$
с исключительными дивизорами $E_1$, \ldots, $E_7$.
При этом точка $x$ не лежит ни на одном из исключительных
дивизоров $E_i$.
Напомним, что имеется линейная эквивалентность
$$K_S\sim \pi^* K_{\P^2}+\sum\limits_{i=1}^7E_i\,.$$
Кроме того, дивизор $D$ на $\P^2$ имеет кратность не меньше $2$
в точке $P_i$ тогда и только тогда, когда
дивизор $\pi^*D$ имеет кратность
не меньше $2$ вдоль $E_i$.
Наконец, так как точка $x$ не содержится ни в одном из дивизоров $E_i$,
морфизм~$\pi$ является изоморфизмом в окрестности точки $x$.
Таким образом, имеется отображение
$$|-2K_{\P^2}-2\sum\limits_{i=1}^7P_i-3\pi(x)|\longrightarrow
|-2K_S-3x|\,,$$
заданное формулой
$$D\mapsto\pi^{-1}(D)+\sum\limits_{i=1}^7\big(\mult_{P_i}(D)-2\big)E_i\,.$$
Можно проверить, что это отображение является биекцией,
но мы не будем это использовать.
Класс $-2K_{\P^2}$ является классом линейной эквивалентности
кривых степени $6$ на $\P^2$.
Требование, что кривая имеет кратность
по крайней мере $2$ (соответственно, по крайней мере $3$)
в данной точке, накладывает на уравнение кривой
$3$ линейных условия (соответственно, $6$ линейных условий).
Таким образом, размерность линейной системы
$$|-2K_{\P^2}-2\sum_{i=1}^7 P_i-3\pi(x)|$$
больше или равна
$$27-7\cdot 3-6=0\,.$$
Следовательно, эта линейная система
непуста, а значит, непуста и линейная система
$|-2K_S-3x|$.) На самом деле, это утверждение можно доказать и в том
случае, если точка $x$ лежит на одной из прямых на $S$, однако
это требует чуть больше усилий.
\item[(ii)]
Покажите, что если $x\in S(F)$ не лежит ни на одной прямой
на $S_{\bar{F}}$, то
никакая кривая, проходящая через $x$, не может содержать неприводимую
компоненту (определённую над $\bar{F}$), которая при антиканоническом
морфизме $\varphi\colon S\to\P^2$ изоморфно отображается на
прямую, проходящую через точку $\varphi(x)$.
(Указание: воспользуйтесь пунктом~(o).)
\item[(iii)]
Покажите, что если $x\in U(F)$, то
любая кривая $C\in |-2K_S-3x|$ неприводима (причём даже
над $\bar{F}$)
и $F$-рациональна, и при этом кратность $\mult_x(C)$ кривой~$C$ в точке~$x$
равна~$3$.
(Указание: антиканонический
образ $\varphi(C)\subset\Pb^2$ кривой~$C$ является кривой степени
не больше~$4$. Поскольку $x\not\in R$, в окрестности~$x$ отображение $\varphi$
неразветвлено и, следовательно,
$$\mult_{\varphi(x)}\big(\varphi(C)\big)\geqslant\mult_x(C)\geqslant 3\,.$$
Если $\mult_{\varphi(x)}\big(\varphi(C)\big)=4$ или
$\deg\big(\varphi(C)\big)=3$, то $\varphi(C)$ является
объединением четырёх или трёх прямых над $\bar{F}$, проходящих
через $\varphi(x)$, соответственно.
По пункту~(ii) это невозможно. Итак,
$\mult_{\varphi(x)}\big(\varphi(C)\big)=3$ и
$\deg\big(\varphi(C)\big)=4$. Следовательно, $\mult_x(C)=3$,
и $C$ бирационально отображается на $\varphi(C)$.
Предположим, что $C$ приводима. Тогда $\varphi(C)$ также
приводима, и из условий на кратность и степень следует,
что $\varphi(C)$ содержит прямую над $\bar{F}$,
проходящую через $\varphi(x)$. Это снова противоречит
пункту~(ii). Поэтому $C$ неприводима и бирациональна $\varphi(C)$.
Так как $\mult_{\varphi(x)}\big(\varphi(C)\big)=3$ и
$\deg\big(\varphi(C)\big)=4$, то
кривая~$\varphi(C)$, а значит, и кривая $C$, рациональна.)
\item[(iv)]
Покажите, что если $x\in U(F)$, то кривая $C\in |-2K_S-3x|$ единственна.
(Указание: предположим, что найдутся две различные кривые
$$C_1, C_2\in |-2K_S-3x|\,.$$
Тогда $C_1$ и $C_2$ неприводимы по пункту~(iii), так что
$$8=C_1\cdot C_2\geqslant\mult_x(C_1)\cdot\mult_x(C_2)=9\,,$$
противоречие.)
В дальнейшем мы будем обозначать эту кривую через $C_x$.
Если точка $x$ определена не над полем $F$, а над
полем $F'\supset F$, то и кривая $C_x$ также определена над полем
$F'$.
\item[(v)]
Пусть на поверхности $S$ есть $F$-точка
$$x\in U(F)\subset S(F)\,.$$
Построим по ней кривую $C_x$ (см. пункт~(iv)). Далее, для общей точки
$y\in C_x$ построим кривую~$C_y$ тем же способом. Рассмотрим многообразие
инцидентности
$$\widetilde{S}=\{(y, z) \mid y\in C_x\cap U, z\in C_y\}.$$
Докажите, что $\widetilde{S}$~--- рациональная поверхность над полем
$F$. (Указание: проекция поверхности $\widetilde{S}$ на $C_x\cap U$
задаёт на $\widetilde{S}$ структуру расслоения.
Общий слой этого расслоения изоморфен кривой $C_{\eta}$
над полем $F(C_x)$,
построенной по общей точке $\eta$ кривой $C_x$ как в пункте~(iv).
Из пункта~(iii) следует, что кривая~$C_{\eta}$
рациональна над полем $F(C_x)$, то есть над
полем рациональных функций от одной переменной над $F$.)
\item[(vi)] Предположим, что поверхность $S$ имеет $F$-точку,
не лежащую ни на одной прямой на $S_{\bar{F}}$, а также не лежащую
на кривой ветвления
$R\subset S$.
Докажите, что $S$ унирациональна над полем $F$.
(Указание: покажите, что морфизм
$$\widetilde{S}\to S,\quad (y,z)\mapsto z\,,$$
доминантен. Это следует из того, что поверхность $\widetilde{S}$ неприводима
и её образ содержит как минимум две неприводимые кривые,
$C_x$ и $C_y$ для общей точки~\mbox{$y\in C_x$}.)
\item[(vii)]
Докажите $F$-унирациональность поверхности дель Пеццо степени $d\geqslant 2$,
определённой над $F$ и имеющей достаточно много
(например, плотное по Зарисскому множество) $F$-точек. (Указание:
раздуйте точки.)
\end{itemize}
\end{prob}

Читатель, имеющий более основательное
знакомство с геометрией поверхностей дель Пеццо, может провести
конструкцию кривой $C$ из упражнения~\ref{prob:DP2}(i),
исходя из геометрических соображений. Пусть
$\pi\colon S'\to S$ --- раздутие точки $x$
с исключительным дивизором $E$.
Так как точка $x$ не лежит на
прямых на $S_{\bar{F}}$, можно проверить, что
поверхность $S'$ является поверхностью дель Пеццо степени~$1$.
Линейная система $|-2K_{S'}|$ задаёт морфизм
$$\psi\colon S'\to Q\,,$$
который является двулистным накрытием квадратичного
конуса $Q$ (см., например,~\cite[\S8.8.2]{Dolgachev}). Рассмотрим инволюцию
Галуа
$$\iota\colon S'\to S'$$
двойного накрытия~$\psi$.
Искомая кривая $C$ получается как $\pi\big(\iota(E)\big)$.
При этом неприводимость и
рациональность~$C$ очевидны из конструкции.
Отметим также, что условия общности, которые накладываются
на точку~$x\in S$, можно значительно ослабить (см.~\cite{Testa}).

\bigskip

Теперь мы воспользуемся упражнением~\ref{prob:DP2} для того, чтобы доказать
унирациональность многообразия из упражнения~\ref{prob:Artin-Mumford}.

\begin{prob}{\textbf{Унирациональность двойного накрытия $\P^3$
с ветвлением в квартике}}

Пусть $k$ --- алгебраически замкнутое поле характеристики нуль, и
$V\subset\P^3$ --- определённая над $k$ поверхность
степени $4$ с изолированными особенностями (в частности, эта поверхность неприводима).
Пусть~$\pi\colon X\to\P^3$ --- двойное накрытие
проективного пространства $\P^3$, разветвлённое в поверхности~$V$; многообразие $V$ неприводимо и нормально.

\begin{itemize}
\item[(o)] Докажите, что на $X$ имеется двухпараметрическое семейство прямых, то есть таких кривых,
которые изоморфно отображаются на прямые в $\P^3$ при отображении $\pi$.
(Указание: образ прямой на $X$
при отображении $\pi$ является бикасательной к поверхности $V\subset\P^3$.
Размерность пространства бикасательных к $V$ можно посчитать, например,
рассмотрев
его как подмногообразие в грассманиане~\mbox{$\mathrm{G}(4,2)$}
двумерных подпространств в четырёхмерном векторном пространстве.)
\item[(i)] Рассмотрим проекцию $\phi\colon\P^3\dasharrow\P^1$
из общей прямой $L$ в $\P^3$. Она задаёт рациональное отображение
$\theta\colon X\dasharrow\P^1$.
Пусть $\psi\colon \widetilde{X}\to X$ является раздутием
прообраза прямой $L$ на $X$. Тогда $\psi$ даёт разрешение неопределённостей
отображения $\theta$. Композиция $\theta\circ\psi$ является
расслоением $\widetilde{X}\to\Pb^1$.
Пусть $S$ --- поверхность над полем
$$F=k(\Pb^1)\cong k(t)\,,$$
являющаяся общим (схемным) слоем морфизма $\theta\circ\psi$.
Покажите, что $S$ является гладкой
поверхностью дель Пеццо степени $2$.
\item[(ii)]
Покажите, что множество $F$-точек плотно по Зарисскому
на поверхности~$S$. (Указание: воспользуйтесь
пунктом~(o).) Выведите отсюда и из
упражнения~\ref{prob:DP2}, что поверхность $S$ унирациональна.
\item[(iii)] Докажите, что многообразие $\widetilde{X}$
(а значит, и многообразие $X$) унирационально. (Указание:
воспользуйтесь пунктом~(ii).)
\end{itemize}
\end{prob}

\newpage
\section{Ограничение скаляров по Вейлю и алгебраические торы}
\label{section:Weil}

\subsection{Ограничение скаляров по Вейлю}
Пусть $k\subset K$ --- конечное сепарабельное расширение полей
степени $d$, а $X$ --- квазипроективное многообразие над полем $K$.
Вложим поле $K$ в $k^{sep}$ над $k$.
Пусть $L$ является нормальным замыканием поля $K$
в $k^{sep}$ над $k$. Тогда
существует ровно $d$ вложений
$$\sigma_1,\ldots,\sigma_d\colon K\hookrightarrow L$$
над полем $k$.
При помощи вложений
$\sigma_i$ из многообразия $X$ получается набор многообразий $X_{\sigma_i}$
над полем~$L$: если $X$ задано над $K$ набором уравнений, то,
применяя~$\sigma_i$ к их коэффициентам, получаем набор уравнений для
$X_{\sigma_i}$. Заметьте, что многообразия~$X_{\sigma_i}$ не обязательно
изоморфны друг другу, даже над сепарабельным замыканием~$L^{sep}=k^{sep}$.

Рассмотрим многообразие $X_{\sigma_1}\times\ldots\times
X_{\sigma_d}$ над полем $L$. Для любого
элемента~\mbox{$g\in\Gal(L/k)$} имеется изоморфизм
$$g_*\big(X_{\sigma_1}\times\ldots\times X_{\sigma_d}\big)\cong
X_{g\sigma_1}\times\ldots\times X_{g\sigma_d}$$
многообразий над $L$.
Действуя далее перестановкой сомножителей,
получаем изоморфизм
$$g_*\big(X_{\sigma_1}\times\ldots\times X_{\sigma_d}\big)\cong
X_{\sigma_1}\times\ldots\times X_{\sigma_d}\,.$$
Это задаёт канонические данные спуска
на многообразии $X_{\sigma_1}\times\ldots\times
X_{\sigma_d}$ относительно расширения полей
$k\subset L$.
Согласно упражнению~\ref{exer-generalspusk}(iv), это однозначно определяет
такое квазипроективное многообразие $Y$ над $k$, что расширение скаляров
$Y_L$ изоморфно~\mbox{$X_{\sigma_1}\times\ldots\times X_{\sigma_d}$}.

\begin{defin}
Построенное выше многообразие $Y$
называется \emph{ограничением скаляров по Вейлю} многообразия $X$ относительно расширения полей $k\subset K$ и
обозначается $R_{K/k}(X)$.
\end{defin}

Отметим, что ограничение скаляров по Вейлю
является алгебраическим аналогом аналитической
процедуры, позволяющей рассматривать
$n$-мерное комплексное многообразие как $2n$-мерное
вещественное многообразие, ср. с
упражнениями~\ref{prob:gruppa-ogranichivaetsya-v-gruppu}(iv),~\ref{prob:geom-svojstva-Weil}(i), и~\ref{prob:ogranichenie-proektivnogo}(iii).

\bigskip
\begin{prob}{\bf Функториальные свойства ограничения скаляров по Вейлю}
\label{prob:gruppa-ogranichivaetsya-v-gruppu}
\hspace{0cm}
\begin{itemize}
\item[(o)] Проверьте, что ограничение скаляров по Вейлю
определяет функтор $R_{K/k}$ из категории квазипроективных
многообразий над полем~$K$ в категорию квазипроективных
многообразий над полем~$k$.
\item[(i)]
Пусть $X_1$ и $X_2$ --- многообразия, определённые над $K$.
Докажите, что
$$R_{K/k}(X_1\times_K X_2)\cong R_{K/k}(X_1)\times_k R_{K/k}(X_2)\,.$$
\item[(ii)]
Докажите, что ограничение скаляров по Вейлю алгебраической группы
над~$K$ является алгебраической группой над $k$. (Указание:
воспользуйтесь пунктом~(i).)
\item[(iii)]
Докажите, что функтор ограничения скаляров по Вейлю является правым
сопряжённым к функтору расширения скаляров, то есть для любого
\mbox{$K$-многообразия} $X$ и $k$-многообразия $Y$ есть каноническая
биекция
$$
\Hom_K\big(Y_K, X\big)\cong\Hom_k\big(Y, R_{K/k}(X)\big)\,,
$$
где $\Hom_F(U,V)$ обозначает множество морфизмов
между алгебраическими многообразиями $U$ и $V$ над полем $F$.
\item[(iv)]
Покажите, что имеется каноническая биекция $X(K)\cong
R_{K/k}(X)(k)$.
\item[(v)]
Пусть многообразие $Y$ определено над полем $k$.
Действие группы Галуа~\mbox{$\Gal(L/k)$} на множестве вложений
$K$ в $L$ над $k$ задаёт гомоморфизм
$$\varphi\colon \Gal(L/K)\to\mathrm{S}_d\subset\Aut\big(Y^{{}\times d}\big)\,.$$
Докажите, что многообразие
$R_{K/k}(Y_K)$ является скруткой многообразия $Y^{{}\times d}$
при помощи $\varphi$, см. упражнения~\ref{exer-spusk}(i),(iii).
\end{itemize}
\end{prob}

\bigskip
Напомним (см. раздел~\ref{subsection:Br-i-SB}), что многообразие $X$,
определённое над полем $k$,
называется геометрически неприводимым,
если неприводимо многообразие $X_{\bar{k}}$
(и геометрически приводимым в противном случае).

\bigskip
\begin{prob}{\bf Геометрические свойства ограничения скаляров по Вейлю}
\label{prob:geom-svojstva-Weil}
\hspace{0cm}
\begin{itemize}
\item[(i)]
Докажите, что
$$\dim_k\big(R_{K/k}(X)\big)=[K:k]\cdot\dim_K(X)\,.$$
\item[(ii)]
Пусть $Z\subset X$ --- замкнутое подмногообразие в $X$, определённое над
$K$, а~\mbox{$U=X\setminus Z$}~--- открытое дополнение.
Покажите, что $R_{K/k}(Z)$ является замкнутым подмногообразием
в $R_{K/k}(X)$, а $R_{K/k}(U)$ является открытым подмножеством в
$R_{K/k}(X)$. Верно ли, что
$$R_{K/k}(X)\cong R_{K/k}(Z)\cup R_{K/k}(U)\,?$$
(Указание: рассмотрите пример $X=\Ab^1_K$, $Z=\{0\}$, и
расширьте скаляры с $k$ на~$k^{sep}$.)
\item[(iii)] Докажите, что ограничение скаляров по Вейлю
геометрически неприводимого многообразия геометрически неприводимо.
\item[(iv)] Верно ли, что ограничение скаляров по Вейлю
неприводимого многообразия неприводимо?
(Указание: примените
упражнение~\ref{prob:gruppa-ogranichivaetsya-v-gruppu}(v)
к $Y=\Spec(E)$, где~$E$~---
такое конечное расширение Галуа поля $k$, что $K\otimes_k E$ не имеет
делителей нуля.)
\item[(v)]
Покажите, что имеется каноническая биекция
$$R_{K/k}(X)(k^{sep})\cong\Map_{G_K}\!\big(G_k, X(k^{sep})\big)\,,$$
коммутирующая с действием группы
$G_k$, где действие~$G_k$ на правой части задаётся по формуле
$$\big({}^g\varphi\big)(g')=\varphi(g'g),\quad \varphi\in
\Map_{G_K}\!\big(G_k, X(k^{sep})\big),
\quad g,g'\in G_k\,,$$
а действие группы $G_K$ на $G_k$ задаётся левыми сдвигами.
(Указание: по определению имеется каноническая $G_k$-эквивариантная
биекция
$$
R_{K/k}(X)(k^{sep})\cong X_{\sigma_1}(k^{sep})\times\ldots\times X_{\sigma_d}(k^{sep})\,.
$$
Пусть для определённости $\sigma_1$ соответствует
выбранному вложению поля $K$ в поле $k^{sep}$. Тогда элементу
$$
(x_{\sigma_1},\ldots, x_{\sigma_d})\in X_{\sigma_1}(k^{sep})\times\ldots\times X_{\sigma_d}(k^{sep})
$$
соответствует отображение
$$
g\mapsto g\big(x_{g^{-1}\sigma_1}\big),\quad
g\in G_k\,,
$$
а отображению
$$\varphi\in\Map_{G_K}\big(G_k,X(k^{sep})\big)$$
соответствует набор точек
$$\big(g_1\varphi(g_1^{-1}),\ldots,g_d\varphi(g_d^{-1})\big)\,,$$
где $g_i\in G_k$ --- такие элементы, что $g_i\sigma_1=\sigma_i$.)
\item[(vi)]
Предположим, что $X$ является коммутативной алгебраической
группой над полем~$K$. Покажите, что тогда имеется изоморфизм \mbox{$G_k$-модулей}
$$R_{K/k}(X)(k^{sep})\cong i_*X(k^{sep})\,,$$
где $i\colon G_K\hookrightarrow G_k$ обозначает естественное вложение,
а $i_*$~--- коиндуцирование, см. определение~\ref{defin-coind}.
(Указание: воспользуйтесь пунктом (v) и замечанием после
определения~\ref{defin-coind}.)
\end{itemize}
\end{prob}

\bigskip
\bigskip
\bigskip

\begin{prob}{\bf Ограничение по Вейлю аффинных многообразий}
\label{prob:ogranichenie-affinnogo}
\hspace{0cm}

\begin{itemize}
\item[(i)] Докажите, что имеется (неканонический) изоморфизм
$$R_{K/k}(\A^1)\cong\A^d\,.$$
(Указание: по упражнению~\ref{prob:gruppa-ogranichivaetsya-v-gruppu}(iii)
для любого $k$-многообразия
$Y$ имеется каноническая биекция между
множеством~\mbox{$\Hom_k\big(Y,R_{K/k}(\Ab^1)\big)$} и кольцом регулярных
функций
$$K[Y_K]\cong K\otimes_k k[Y]\,.$$
Далее, выбор базиса
$K$ над $k$ определяет биекцию $K\otimes_k k[Y]$
c $k[Y]^{\oplus d}$, то есть с множеством
$\Hom_k(Y,\Ab^d)$.)
\item[(ii)] Докажите, что $R_{K/k}(\A^n)\cong\A^{nd}$.
(Указание: воспользуйтесь пунктом (i) и
упражнением~\ref{prob:gruppa-ogranichivaetsya-v-gruppu}(i).)
\item[(iii)] Покажите, что если многообразие
$X\subset\A^n$ аффинно, то $R_{K/k}(X)$ тоже аффинно и
вложено в~$\A^{nd}$. (Указание: воспользуйтесь пунктом~(ii)
и упражнением~\ref{prob:geom-svojstva-Weil}(ii).)
\item[(iv)]
Пусть $X\subset \Ab^n$ --- аффинная гиперповерхность,
заданная уравнением
$$f(x_1,\ldots,x_n)=0\,.$$
Докажите, что аффинное
многообразие $R_{K/k}(X)\subset \Ab^{nd}$ задаётся следующими~$d$
уравнениями. Выберем базис $e_1,\ldots,e_d$ в $K$ над $k$.
Рассмотрим формальные переменные $x_{ij}$,
где $1\leqslant i \leqslant n$ и $1\leqslant j\leqslant d$.
Подставим в многочлен $f$ вместо $x_i$ выражения $\sum_{j=1}^d x_{ij}e_j$,
используя при этом соотношения
$$
e_{\alpha}\cdot e_{\beta}=\sum\limits_{\gamma=1}^{d} c_{\alpha\beta}^{\gamma} e_{\gamma}\,,
$$
выполненные в поле $K$, где $c_{\alpha\beta}^{\gamma}\in k$. Таким образом возникает равенство
$$
f\left(\sum\limits_{j=1}^d x_{1j}e_j,\ldots,
\sum\limits_{j=1}^d x_{nj}e_j\right)=\sum\limits_{l=1}^df_l(x_{ij})e_l\,,
$$
задающее $d$ уравнений $f_l(x_{ij})=0$, $1\leqslant l\leqslant d$.
(Указание: рассуждайте так же, как в пункте (i).)
\item[(v)] Пусть $X\subset \Ab^{n}$ --- аффинное многообразие. Проверьте, что уравнения, задающие $R_{K/k}(X)$, получаются из уравнений, задающих $X$, способом, описанным в пункте~(iv).
\end{itemize}
\end{prob}

\bigskip
\begin{prob}{\bf Ограничение по Вейлю рациональных многообразий}
\label{prob:ogranichenie-racionalnogo}

Пусть $X$ является $K$-рациональным многообразием.
Докажите, что многообразие $R_{K/k}(X)$ является $k$-рациональным.
(Указание: воспользуйтесь упражнениями~\ref{prob:ogranichenie-affinnogo}(ii)
и~\ref{prob:geom-svojstva-Weil}(ii),(iii).)
\end{prob}

\bigskip
\begin{prob}{\bf Ограничение по Вейлю проективных многообразий}
\label{prob:ogranichenie-proektivnogo}
\hspace{0cm}

\begin{itemize}
\item[(i)]
Пусть $X\subset\P^n=\P(W)$ --- проективное многообразие над $K$.
Постройте замкнутое вложение
$$R_{K/k}(X)\subset\P^{(n+1)^d-1}$$
над $k$.
(Указание:
каждое вложение $\sigma\colon K\hookrightarrow L$ над $k$
определяет $L$-векторное пространство
$W\otimes_{K,\sigma} L$ полученное из $W$ расширением скаляров
относительно $\sigma$ (то есть расширением
скаляров с $\sigma(K)$ на $L$).
Рассмотрите $L$-векторное пространство
$$U=(W\otimes_{K,\sigma_1} L)\otimes_L\ldots
\otimes_L (W\otimes_{K,\sigma_d} L)$$
и постройте замкнутое вложение
$$X_{\sigma_1}\times\ldots\times X_{\sigma_d}\hookrightarrow\P(U)\,.$$
Дальше действуйте аналогично тому, что объясняется в указании
к упражнению~\ref{prob:spusk-projective}(iii).)
\item[(ii)]
Предположим, что характеристика поля $k$ не равна $2$.
Покажите, что для квадратичного расширения
$$k\subset K=k(\sqrt{a}),\quad a\in k^*\,,$$
многообразие $R_{K/k}(\P^n)$ задаётся в проективизации пространства
матриц размера~\mbox{$(n+1)\times (n+1)$} условием
$$\mathrm{rk}\big((A+A^T)+\sqrt{a}(A-A^T)\big)\leqslant 1\,,$$
где $\mathrm{rk}$ обозначает ранг матрицы, а $A^T$ является транспонированной матрицей для $A$.
(Указание: пусть $g$ обозначает нетривиальный элемент
группы~\mbox{$G=\Gal(K/k)$}. Для произвольного $(n+1)$-мерного $k$-векторного
пространства~$V$ отождествите $K$-векторное пространство
$V_K\otimes_K g_*V_K$ вместе с естественными данными спуска
(см. упражнение~\ref{prob:spusk-projective}(iii))
с пространством матриц размера $(n+1)\times (n+1)$ над $K$ вместе с
$K$-полулинейной инволюцией, заданной по формуле
$$M\mapsto g(M)^T\,.$$
Таким образом, $k$-векторное пространство $(V_K\otimes_K g_*V_K)^G$ отождествляется с пространством
эрмитовых матриц над $K$. В свою очередь, пространство эрмитовых матриц отождествляется с пространством
всех матриц размера~\mbox{$(n+1)\times (n+1)$} над $k$:
матрице $A$ над $k$
сопоставляется эрмитова матрица $(A+A^T)+\sqrt{a}(A-A^T)$ над~$K$.
Далее воспользуйтесь упражнениями~\ref{prob:spusk-projective}(i),(ii).)
\item[(iii)]
В обозначениях пункта~(ii) покажите, что
поверхность $R_{K/k}(\P^1)$ изоморфна квадрике в $\P^3$,
заданной уравнением
$$x^2-y^2+z^2-au^2=0\,.$$
В частности, при $k=\mathbb{R}$ и $a=-1$ имеем
$K=\mathbb{C}$, и $\mathbb{R}$-точки вещественного
алгебраического многообразия $R_{\mathbb{C}/\mathbb{R}}(\mathbb{P}^1)$
образуют сферу Римана.
\item[(iv)]
Пусть $Q\subset\P^3$ --- гладкая двумерная квадрика над $k$,
для которой $\disc(Q)\neq 1$ (см. определение~\ref{defin:discrquadr}).
Пусть $K=k\big(\sqrt{\disc(Q)}\big)$, а $C$ --- гладкая коника над $K$,
соответствующая $\cl(Q_K)$ (см. определение~\ref{defin:Clifford}).
Докажите, что $R_{K/k}(C)\cong Q$.
(Указание:
по упражнению~\ref{prob-section}(vi)
имеется изоморфизм $Q_K\cong C\times C$ многообразий над $K$.
По упражнению~\ref{prob-discr-geom-interpretation}(i)
нетривиальная инволюция поля $K$ над $k$ переставляет сомножители в этом
разложении.)
\item[(v)] Выведите из пункта~(iv) упражнение~\ref{prob-nonsection}(ii).
\end{itemize}
\end{prob}

\subsection{Алгебраические торы}

Всюду в дальнейшем $\Gb_m$ обозначает определённую над $\Z$ аффинную
групповую схему $\Spec(\Z[t, t^{-1}])$.
Для любого поля $F$ мы также будем обозначать через $\Gb_m$
соответствующий одномерный тор над полем $F$, то есть
расширение скаляров схемы~$\Gb_m$ с~$\Z$ на $F$.

\begin{defin}
\emph{Алгебраическим тором} над полем $k$ называется алгебраическая
группа $T$ над $k$, для которой
существует расширение полей $k\subset F$ и изоморфизм
$F$-многообразий
$T_{F}\cong \Gb_{m}^{{}\times n}$
при некотором натуральном~$n$.
В этом случае будем говорить, что $T$ \emph{расщепим} над $F$.
\end{defin}

Отметим, что любой алгебраический тор расщепляется
над некоторым конечным сепарабельным расширением поля определения
(см.~\cite[Теорема~3.30]{Voskresenskii-rus}).
Таким образом, алгебраические торы --- то же самое, что формы
$\Gb_m^{{}\times n}$.

Как и выше, пусть $k\subset K$ является конечным сепарабельным расширением полей степени $d$. Ограничение скляров по Вейлю~$R_{K/k}(\Gb_{m})$ является алгебраическим тором размерности $d$ над~$k$ (это следует из
упражнения~\ref{prob:gruppa-ogranichivaetsya-v-gruppu}(ii),(v)).
По упражнению~\ref{prob:gruppa-ogranichivaetsya-v-gruppu}(iv)
имеется биекция
$$R_{K/k}(\Gb_{m})(k)\cong K^*\,.$$
Для краткости мы
будем обозначать этот алгебраический тор над $k$ символом $K^*$.
Поскольку поле $k$ в дальнейшем будет зафиксировано, такое
обозначение не должно привести к путанице. В частности, $k^*$ обозначает
тор $\Gb_{m}$.

Всюду далее $K'$ обозначает конечное сепарабельное расширение поля $K$.

\bigskip
\begin{prob}{\bf Общие свойства алгебраических торов}
\label{prob:restr-tori}
\hspace{0cm}

\begin{itemize}
\item[(i)] Докажите, что ограничение скаляров по Вейлю
алгебраического тора
является алгебраическим тором.
\item[(ii)]
Покажите, что тор $K^*$ является скруткой многообразия $\Gb_m^{{}\times d}$,
соответствующей гомоморфизму
$$\varphi\colon G_k\to\mathrm{S}_d\subset\Aut\big(\Gb_m^{{}\times d}\big)\,.$$
(Указание: воспользуйтесь упражнением~\ref{prob:gruppa-ogranichivaetsya-v-gruppu}(v).)
\item[(iii)]
Докажите, что вложение $k^*\hookrightarrow K^*$
является морфизмом алгебраических торов.
(Указание: воспользуйтесь
упражнением~\ref{prob:gruppa-ogranichivaetsya-v-gruppu}(iii).)
\item[(iv)]
Докажите, что норма
$$\Nm_{K/k}\colon K^*\to k^*$$
является морфизмом алгебраических торов.
(Указание:
покажите, что $\Nm_{K/k}$ является скруткой
морфизма
$$\Gb_m^{{}\times d}\to\Gb_m\,, \quad (z_1,\ldots,z_d)\mapsto
z_1\cdot\ldots\cdot z_d$$
при помощи гомоморфизма $\varphi$ из пункта~(ii).)
\item[(v)]
Найдите уравнения, задающие алгебраический тор
$R_{K/k}(\Gb_m)$.
(Указание: рассмотрите норму
$$\Nm_{K/k}\colon K^*\to k^*$$
как форму степени
$d$ над $k$. Из пунктов~(ii) и~(iv) следует, что
алгебраический тор $R_{K/k}(\Gb_m)$ является дополнением
к соответствующей гиперповерхности в аффинном пространстве~$\Ab^d$.)
\item[(vi)]
Докажите, что вложение $K^*\hookrightarrow K'^*$ и норма
$\Nm_{K'/K}\colon K'^*\to K^*$ являются морфизмами алгебраических торов.
(Указание: сначала воспользуйтесь пунктами~(iii) и~(iv), заменив
$k$ и $K$ на $K$ и $K'$, соответственно. Потом примените ограничение
скаляров по Вейлю с $K$ на $k$.)
\end{itemize}
\end{prob}

\bigskip

\begin{prob}{\bf Рациональность некоторых алгебраических торов}
\label{prob:tor-F/E-racionalen} \hspace{0cm}

\begin{itemize}
\item[(o)] Покажите, что алгебраический тор $K^*$
рационален над $k$. (Указание: воспользуйтесь
упражнением~\ref{prob:ogranichenie-racionalnogo}.)
\item[(i)] Докажите, что алгебраический тор $K^*/k^*$
рационален над $k$. (Указание: этот алгебраический тор является
открытым подмножеством в определённом над~$k$ проективном
пространстве~$\Pb(K)\cong\Pb^{d-1}$.)
\item[(ii)] Докажите, что алгебраический тор $K'^*/K^*$ рационален над $k$.
(Указание: сначала воспользуйтесь пунктом~(i), заменив
$k$ и $K$ на $K$ и $K'$, соответственно. Потом примените ограничение
скаляров по Вейлю с $K$ на $k$ и воспользуйтесь
упражнением~\ref{prob:ogranichenie-racionalnogo}.)
\end{itemize}
\end{prob}

\bigskip
\begin{prob}{\bf Ядро отображения нормы}
\label{prob:Ker-Nm} \hspace{0cm}

\begin{itemize}
\item[(i)]
Докажите, что ядро нормы $\Nm_{K/k}\colon K^*\to k^*$
является $(d-1)$-мерным алгебраическим тором над $k$.
(Указание: при помощи указания к упражнению~\ref{prob:restr-tori}(iv)
проверьте,
что над~$k^{sep}$ это ядро изоморфно тору $\Gb_m^{{}\times (d-1)}$.)
\item[(ii)]
Докажите, что ядро нормы $\Nm_{K'/K}\colon K'^*\to K^*$
является алгебраическим тором над $k$.
(Указание: сначала воспользуйтесь пунктом~(i), заменив
$k$ и $K$ на~$K$ и $K'$, соответственно. Потом примените ограничение
скаляров по Вейлю с $K$ на $k$ и воспользуйтесь
упражнением~\ref{prob:restr-tori}(i).)
\end{itemize}
\end{prob}

\bigskip
\begin{prob}{\bf Одномерные торы}
\label{prob:1-dim-tori} \hspace{0cm}

Пусть $k\subset K$ --- сепарабельное квадратичное расширение полей.
Согласно упражнению~\ref{prob:Ker-Nm}(i),
ядро отображения нормы $\Nm_{K/k}\colon K^*\to k^*$
является одномерным алгебраическим тором над $k$.
Докажите, что это определяет биекцию между множеством
сепарабельных квадратичных расширений поля $k$
и множеством одномерных алгебраических торов над $k$.
(Указание: оба этих множества биективны множеству квадратичных
характеров группы Галуа $G_k$. Для множества одномерных алгебраических
торов это следует из того, что $\Aut(\Gb_m)\cong\Z/2\Z$.)
\end{prob}

\subsection{Алгебраические торы и модули Галуа}

Как и раньше, $G_k$ обозначает группу Галуа $\Gal\big(k^{sep}/k\big)$.

\begin{defin}\label{defin-Pic-lattice}
Для алгебраического тора $T$, определённого над полем $k$,
\emph{двойственный модуль Галуа} над группой $G_k$ определяется формулой
$$T^{\vee}=\Hom(T_{k^{sep}}, \Gb_m)\,.$$
Здесь действие группы Галуа $G_k$ на $T^{\vee}$ задаётся по формуле
$$
\big({}^g\chi\big)(t)=g\big(\chi(g^{-1}t)\big),\quad
g\in G_k,\quad \chi\in T^{\vee},\quad t\in T(k^{sep})\,.
$$
\end{defin}

\bigskip
\begin{prob}{\bf Данные спуска для расщепимых торов}
\label{prob:split-tori-spusk} \hspace{0cm}

Пусть $M$ --- свободная конечно порожденная абелева группа. Рассмотрим расщепимый алгебраический
тор $T=\Spec\big(k[M]\big)$ над полем $k$,
где $k[M]$ является групповой алгеброй группы $M$,
а групповая структура на $T$ соответствует гомоморфизму алгебр
$$
k[M]\to k[M]\otimes_k k[M],\quad m\mapsto m\otimes m, \quad m\in M\,.
$$
\begin{itemize}
\item[(i)]
Покажите, что имеется канонический изоморфизм абелевых групп
$M\cong T^{\vee}$.
\item[(ii)]
Покажите, что имеется естественный изоморфизм колец
$$\End(T)\cong \End(T^{\vee})\,,$$
где $\End$ в левой части
обозначает кольцо
эндоморфизмов алгебраического тора~$T$,
а $\End$ в правой части обозначает кольцо
эндоморфизмов абелевой
группы~$T^{\vee}$. (Указание: пользуясь пунктом (i),
постройте взаимно обратные гомоморфизмы этих колец.)
\item[(iii)]
Пусть $k\subset K$ --- конечное расширение Галуа.
Постройте каноническую биекцию между множеством данных спуска
(см.~пример~\ref{examp-spusk}(iv))
на алгебраическом торе~\mbox{$\Spec\big(K[M]\big)$} и множеством действий группы Галуа $\Gal(K/k)$ на $M$.
(Указание: воспользуйтесь пунктом~(ii)
и упражнением~\ref{exer-spusk}(i).)
\end{itemize}
\end{prob}

\bigskip
\begin{defin}[{ср. с определением~\ref{defin-Pic-lattice}}]
\label{defin:dual-Galois-module}
\emph{Двойственный тор} $M^{\vee}$ для \mbox{$G_k$-модуля}~$M$,
конечно порождённого и свободного как $\Z$-модуль,
определяется как спуск с $K$ на $k$ алгебраического тора
$\Spec\big(K[M]\big)$ над $K$, где $k\subset K$ является
таким конечным расширением Галуа, что действие $G_k$ на $M$
пропускается через факторгруппу $\Gal(K/k)$
(см. упражнения~\ref{prob:split-tori-spusk}(iii)
и~\ref{exer-generalspusk}(vi)).
\end{defin}

Легко проверить, что определение~\ref{defin:dual-Galois-module}
корректно, то есть не зависит от выбора поля~$K$.

\bigskip
\begin{prob}{\bf Двойственность между алгебраическими торами и модулями Галуа}
\label{prob:Phi-Psi-quasiobratnye}

Проверьте, что функторы $\Phi\colon T\mapsto T^{\vee}$
и $\Psi\colon M\mapsto M^{\vee}$
являются квазиобратными контравариантными
функторами между категорией алгебраических торов
над $k$ и категорией \mbox{$G_k$-модулей}, конечно порождённых и свободных
как $\Z$-модули
(то есть функторы~\mbox{$\Phi\circ\Psi$} и~$\Psi\circ\Phi$ изоморфны
тождественным функторам).
\end{prob}

\bigskip
Естественно ожидать, что алгебраические торы, двойственные
пермутационным или стабильно пермутационным $G_k$-модулям,
обладают интересными свойствами.

\begin{defin}\label{defin:perm-torus}
Пусть $T$ --- алгебраический тор, определённый над полем $k$.
Говорят, что $T$ является \emph{пермутационным}
(соответственно, \emph{стабильно пермутационным}) тором,
если $T^{\vee}$ является пермутационным
(соответственно, стабильно пермутационным)
$G_k$-модулем (см. определения~\ref{defin:perm-module}
и~\ref{defin:stab-perm-module}).
\end{defin}

\bigskip
\begin{prob}{\bf Строение пермутационных торов}
\label{prob:stab-perm-torus-easy} \hspace{0cm}
\begin{itemize}
\item[(i)]
Постройте изоморфизм $G_k$-модулей
$$(K^*)^{\vee}\cong\Z[G_k/G_K]$$
для конечного сепарабельного расширения $k\subset K$, где $G_k/G_K$
обозначает множество классов
смежности.
\item[(ii)] Докажите, что всякий пермутационный тор $T$ над $k$ имеет вид
$$\mbox{$T\cong\prod\limits_i K_i^*\,,$}$$
где $K_i$ --- некоторые конечные
сепарабельные расширения поля $k$.
(Указание: воспользуйтесь пунктом~(i),
упражнением~\ref{exercise:permutation-module}(i)
и упражнением~\ref{prob:Phi-Psi-quasiobratnye}.)
\end{itemize}
\end{prob}

\bigskip
\begin{prob}{\bf Тривиальность торсоров над стабильно пермутационным тором}
\label{prob:stab-perm-torus-H1} \hspace{0cm}

\begin{itemize}
\item[(i)] Пусть $T$ --- стабильно пермутационный тор, определённый над
полем $k$.
Докажите, что
$$H^1\big(G_k,T(k^{sep})\big)=0\,.$$
(Указание: для пермутационного тора воспользуйтесь
упражнением~\ref{prob:stab-perm-torus-easy}(ii),
упражнением~\ref{prob:geom-svojstva-Weil}(vi),
леммой Шапиро из упражнения~\ref{exer-ind-inj}(iii)
и теоремой Гильберта~90. Случай стабильно пермутационного тора легко сводится
к случаю пермутационного тора.)
\item[(ii)]
Докажите, что любой торсор над стабильно пермутационным тором тривиален.
(Указание: воспользуйтесь пунктом~(i)
и упражнением~\ref{prob:torsor-vs-H1}(i).)
\end{itemize}
\end{prob}

\bigskip
\begin{prob}{\bf Стабильная рациональность стабильно пермутационных торов}
\label{prob:stab-perm-stab-rational} \hspace{0cm}

\begin{itemize}
\item[(i)] Докажите, что любой пермутационный тор рационален.
(Указание: воспользуйтесь упражнениями~\ref{prob:stab-perm-torus-easy}(ii)
и~\ref{prob:tor-F/E-racionalen}(o).)
\item[(ii)]
Воспользовавшись пунктом~(i), покажите,
что любой стабильно пермутационный тор стабильно рационален
(см. определение~\ref{defin:stably-rational}).
\end{itemize}
\end{prob}

\bigskip
\begin{prob}{\bf Торсоры над стабильно пермутационными торами}
\label{prob:torsor-over-stab-perm-torus}

Пусть $T$ --- стабильно пермутационный алгебраический тор, определённый над
полем $k$, и $X$ --- неприводимое
многообразие над~$k$. Пусть $\phi\colon V\to X$ --- торсор над~$T$
(см. определение~\ref{definition:torsorbase}).
Докажите, что многообразия $V$ и $X\times T$ стабильно бирационально
эквивалентны
в смысле определения~\ref{defin:stably-rational}.
(Указание: воспользуйтесь
упражнением~\ref{exercise:stab-permutation-module-restriction},
а затем упражнениями~\ref{prob:stab-perm-torus-H1}(ii),
\ref{prob:torsor-birational}(i) и~\ref{prob:stab-perm-stab-rational}(ii).)
\end{prob}

\subsection{Универсальный торсор}

В этом разделе предполагается, что $X$ --- гладкое многообразие, определённое
над полем $k$, и модуль Галуа
$\Pic(X_{k^{sep}})$ конечно порождён и свободен как $\Z$-модуль.

\begin{defin}
\emph{Тором Нерона--Севери} многообразия $X$ называется
алгебраический тор
$$T_{\NS}(X)=\Pic(X_{k^{sep}})^{\vee}$$
над полем~$k$, где $\Pic(X_{k^{sep}})^{\vee}$ обозначает
двойственный тор для $G_k$-модуля $\Pic(X_{k^{sep}})$ в смысле
определения~\ref{defin:dual-Galois-module}.
\end{defin}

\bigskip
\begin{prob}{\bf Конструкция универсального торсора над $X_{k^{sep}}$}
\label{prob:konstrukciya-universalnogo-torsora-nad-ksep} \hspace{0cm}

\begin{itemize}
\item[(i)] Пусть $L_1,\ldots, L_n$ --- линейные расслоения,
классы $l_1, \ldots, l_n$ которых образуют базис в~$\Pic(X_{k^{sep}})$.
Пусть $L_i^{\circ}$ обозначает дополнение к нулевому сечению в тотальном
пространстве расслоения $L_i$.
Рассмотрим
расслоенное произведение
$$W=L_1^{\circ}\times_X L_2^{\circ}\times_X\ldots\times_X L_n^{\circ}\,.$$
Выберем отождествление
$$T_{\NS}(X)_{k^{sep}}\cong\Gb_m^{{}\times n}\,,$$
которое соответствует
при эквивалентности категорий из упражнения~\ref{prob:Phi-Psi-quasiobratnye}
отождествлению
$$\mbox{$\Pic(X_{k^{sep}})\cong\bigoplus\limits_{i=1}^n\Z\cdot l_i$}\,.$$
Покажите, что это задаёт на $W$ структуру $T_{\NS}(X)_{k^{sep}}$-торсора
над $X_{k^{sep}}$.
\item[(ii)] Докажите, что класс изоморфизма
$T_{\NS}(X)_{k^{sep}}$-торсора $W$ не зависит от выбора расслоений $L_i$.
\end{itemize}
\end{prob}

\bigskip
\begin{defin}\label{defin:universal-torsor-over-ksep}
Торсор $W$, построенный
в упражнении~\ref{prob:konstrukciya-universalnogo-torsora-nad-ksep},
называется \emph{универсальным торсором} над многообразием $X_{k^{sep}}$.
\end{defin}

\bigskip
\begin{prob}{\bf Спуск для универсального торсора}
\label{prob:spusk-universalnogo-torsora} \hspace{0cm}

Предположим, что на многообразии $X$ имеется $k$-точка $p$, и
все обратимые функции на~$X_{k^{sep}}$ постоянны.
Как и в определении~\ref{defin:universal-torsor-over-ksep}, символ
$W$ обозначает универсальный торсор над многообразием~$X_{k^{sep}}$.

\begin{itemize}
\item[(i)] Покажите, что для каждого элемента $g\in G_k$
имеется канонический изоморфизм
$$g_*W\cong g_*L_1^{\circ}\times_X\ldots
\times_X g_*L_n^{\circ}\,,$$
где $g_*W$ и $g_*L_i^{\circ}$
определяются как в примере~\ref{examp-spusk}(v).
\item[(ii)]
Докажите, что имеется естественный изоморфизм между группой автоморфизмов $T_{\NS}(X)_{k^{sep}}$-торсора $W$
над многообразием~$X_{k^{sep}}$ и группой автоморфизмов слоя $W\vert_p$ торсора~$W$
над точкой $p$, который рассматривается
как~\mbox{$T_{\NS}(X)_{k^{sep}}$-торсор} над полем $k^{sep}$.
(Указание: воспользуйтесь тем, что все обратимые
функции на $X_{k^{sep}}$ постоянны.)
\item[(iii)]
Пользуясь пунктом~(ii), докажите, что классы изоморфизма
данных спуска на торсоре $W$ над~$X_{k^{sep}}$ находятся в естественном
взаимно-однозначном соответствии с классами изоморфизма
данных спуска на торсоре
$W\vert_p$, определённом над полем~$k^{sep}$.
(Указание: воспользуйтесь пунктом~(ii) и упражнением~\ref{exer-spusk}(i).)
В частности, на торсоре $W$ существуют данные спуска.
\item[(iv)] Предположим, что
$$H^1\big(G_k, T_{\NS}(X)_{k^{sep}}\big)=\{1\}\,.$$
Докажите, что в этом случае существует
единственный $T_{\NS}(X)$-торсор $V_X$ над~$X$, для которого
имеется изоморфизм $T_{\NS}(X)_{k^{sep}}$-торсоров
$$\big(V_X\big)_{k^{sep}}\cong W\,.$$
\end{itemize}
\end{prob}

\bigskip
\begin{defin}\label{defin:universal-torsor}
Торсор $V_X$, построенный
в упражнении~\ref{prob:spusk-universalnogo-torsora}(iv),
называется \emph{универсальным торсором} над многообразием $X$.
\end{defin}

\bigskip
\begin{prob}{\bf Стабильная бирациональная эквивалентность $X$ и~$V_X$}
\label{prob:X-vs-UX} \hspace{0cm}
Предположим, что на многообразии $X$ имеется $k$-точка,
и все обратимые функции на~$X_{k^{sep}}$ постоянны.
Предположим также, что $G_k$-модуль
$\Pic(X_{k^{sep}})$ является стабильно пермутационным.
Докажите, что в этом случае многообразия $X$ и~$V_X$ стабильно бирационально
эквивалентны. (Указание:
воспользуйтесь упражнением~\ref{prob:torsor-over-stab-perm-torus}.)
\end{prob}

\subsection{Поверхности Шатле и стабильно пермутационные модули}
\label{subsection:Chatelet}

Цель следующих упражнений ---
построить происходящий из геометрии
пример (дискретного)
стабильно пермутационного модуля Галуа,
не сводящийся
к упражнению~\ref{exercise:stab-permutation-module-from-Pic}.
А именно, мы построим поверхность~$S$,
группа Пикара которой является стабильно пермутационным модулем Галуа,
но при этом будет по крайней мере не очевидно,
что такая поверхность рациональна, а в общем случае это даже не верно.

Более того, уже сейчас мы постепенно начинаем приближаться к нашей
следующей ``большой'' цели ---
построению стабильно рационального, но
при этом нерационального многообразия (см.~главу~\ref{section:Polietilen}).
Поэтому нам интересно найти такую поверхность, которая
``имеет шансы'' оказаться стабильно рациональной.
В частности,~$S$ должна быть по крайней мере унирациональной.
При небольших дополнительных условиях из этого следует
рациональность~$S$ над сепарабельным
замыканием основного поля (а для поля нулевой характеристики
для этого даже не надо накладывать дополнительных условий).
Поэтому нас не интересуют тривиальные примеры
типа того, который получается, если в качестве~$S$ взять
поверхность общего типа с группой Пикара, изоморфной~$\Z$.

По тем же причинам мы будем работать именно с поверхностями,
и не удовлетворимся примерами других сортов~---
типа многомерных нерациональных многообразий Фано
с группой Пикара, изоморфной~$\Z$.
В качестве такого ``неинтересного'' примера можно было бы взять
гладкую трёхмерную квартику, то есть гладкую гиперповерхность степени~$4$
в~$\P^4$. Всякая такая гиперповерхности нерациональна
(см.~\cite{IskovskikhManin})
и имеет группу Пикара, изоморфную~$\Z$; более того,
для таких гиперповерхностей в некоторых случаях
известны конструкции унирациональности. С другой стороны,
по крайней мере очень общая трёхмерная квартика
не является стабильно рациональной (см.~\cite{CT-P}),
в то время как именно стабильной рациональностью
мы будем интересоваться
в главе~\ref{section:Polietilen}.

Мы будем следовать конструкции
А.\,Бовиля, Ж.-Л.\,Кольо-Телена,
Ж.-Ж.\,Сан\-сюка и П.\,Свиннертон-Дайера
из~\cite[\S2]{CT-i-drugie}.

\medskip
Пусть поле $k$ имеет характеристику, отличную от~$2$.
Как обычно, $G_k$ обозначает группу Галуа $\Gal(k^{sep}/k)$.
Выберем сепарабельный многочлен $P$ над $k$.
По определению сепарабельность многочлена $P$ означает, что он не
имеет кратных корней в алгебраическом замыкании поля $k$.
В частности, все корни многочлена $P$ лежат в $k^{sep}$.
Пусть $R\subset\Ab^1$ обозначает множество корней
многочлена $P$ в поле~$k^{sep}$. До конца этой главы мы будем предполагать,
что степень~\mbox{$d=2r-1\geqslant 1$} многочлена~$P$ нечётна.

Рассмотрим аффинное пространство $\Ab^3$ над $k$ с координатами
$v_1$, $v_2$ и $u$.
Зафиксируем элемент $a\in k^*$ и зададим поверхность~$S^0$
в~$\Ab^3$ уравнением
\begin{equation}\label{eq:Chatelet-affine}
v_1^2-av_2^2=P(u)\,.
\end{equation}
Аффинная поверхность $S^0$ является открытым подмножеством поверхности,
заданной в~\mbox{$\Pb^2\times\Ab^1$} уравнением
$$v_1^2-av_2^2-P(u)v_3^2=0\,,$$
где $v_i$ --- однородные координаты на~$\Pb^2$, а $u$ --- координата на
$\Ab^1$.

Пусть $\Pb$ --- проективизация расслоения $\OO(r)\oplus\OO(r)\oplus\OO$
на $\P^1$.
Рассмотрим вложение $\Pb^2\times\Ab^1\hookrightarrow\Pb$,
при котором координаты $v_1$ и $v_2$ на $\Pb^2$ соответствуют
расслоениям~$\OO(r)$, а координата $v_3$ соответствует расслоению~$\OO$.
Пусть~$S$ --- замыкание~$S^0$~в~$\Pb$. Поверхность~$S$
называют \emph{поверхностью Шатл\'e}.
Естественная проекция $\Pb\to\Pb^1$ определяет на $S$
структуру расслоения
на коники $f\colon S\to\Pb^1$.
На~$S^0$ морфизм $f$ задаётся формулой
$$(v_1,v_2,u)\mapsto u\,.$$
Пусть $F_{p}\subset S$
обозначает слой расслоения~$f$ над точкой $p\in\Pb^1$;
в частности, $F_{\infty}$ обозначает
слой~$f$ над точкой $\infty\in\Pb^1$.

\bigskip
\begin{prob}{\bf Дивизоры на поверхности Шатле}
\label{prob:Chatelet-easy} \hspace{0cm}

\begin{itemize}
\item[(o)]
Покажите, что поверхность $S\subset\Pb$
в подходящих координатах может быть задана ``на бесконечности''
уравнением
$$V_1^2-aV_2^2-Q(U)V_3^2=0\,,$$
где
$$Q(U)=U^{2r}P(U^{-1})\,.$$
(Указание: вспомните, как устроены функции перехода
для расслоения $\OO(r)$, а именно, $V_1=v_1u^{-r}$.
Имейте в виду, что в соответствии с нашими соглашениями
проективизация расслоения~--- это многообразие прямых в слоях расслоения.)
Убедитесь, что поверхность $S$ гладкая.
\item[(i)] Проверьте, что при всех $p\in R\cup\{\infty\}$
слои $F_{p}$ геометрически приводимы
и~\mbox{$F_{p}=D_{p}\cup D_{p}'$}
для некоторых (различных) неприводимых дивизоров~$D_{p}$ и~$D_{p}'$
над полем~$k^{sep}$. Более того,
все остальные слои расслоения $f$ гладкие.
\item[(ii)]
Покажите, что над полем $k^{sep}$ имеет место разложение
$$S\setminus S^0=F_{\infty}\cup E\cup E'\,,$$
причём дивизоры $E$ и $E'$ не пересекаются и являются
сечениями расслоения~$f$ (определёнными, вообще говоря,
над~$k^{sep}$).
В дальнейшем мы будем считать, что дивизоры $E$, $E'$,
а также $D_{p}$ и $D_{p}'$ при всех~\mbox{$p\in R\cup\{\infty\}$},
выбраны таким образом, что
кривая $E$ пересекает все кривые~$D_{p}$
и не пересекает ни одну кривую
$D_{p}'$, а кривая $E'$, наоборот,
пересекает все кривые~$D_{p}'$ и не пересекает ни одну кривую
$D_{p}$. Проверьте, что пересечения устроены так, что
мы всегда можем выбрать такие обозначения.
\item[(iii)]
Проверьте, что над полем $k^{sep}$
имеет место линейная эквивалентность
дивизоров
$$\sum\limits_{p\in R} D_p+D_{\infty}+E-E'-rF_{\infty}\sim 0\,.$$
(Указание: рассмотрите функцию $(v_1-\sqrt{a}v_2)/v_3$.)
\item[(iv)] Покажите, что группа Пикара
$\Pic\big(S_{k^{sep}}\big)$ свободно порождена
над $\Z$ дивизорами $E$, $E'$, $F_{\infty}$
и $D_{p}$, $p\in R$.
(Указание: стягивая над $k^{sep}$ компоненты приводимых слоёв расслоения $f$,
получите $\Pb^1$-расслоение над $\Pb^1$. Выведите отсюда,
что группа Пикара
$\Pic\big(S_{k^{sep}}\big)$
свободно порождена
над $\Z$ дивизорами $D_{\infty}$,~$E$,~$F_{\infty}$
и~$D_{p}$, $p\in R$. После этого воспользуйтесь
пунктом~(iii).)
\item[(v)]
Покажите, что имеется естественная точная последовательность
$G_k$-модулей
$$0\to\Z[R]\to
\Z\big[\{D_{p},D_{p}'\}_{p\in R}\big]\oplus
\Z\big[\{E,E'\}\big]\oplus
\Z\big[\{F_{\infty}\}\big]\to
\Pic\big(S_{k^{sep}}\big)\to 0\,.$$
(Указание: определите первое отображение так,
чтобы оно переводило $p\in R$ в дивизор
$D_{p}+D_{p}'-F_{\infty}$, и воспользуйтесь пунктом~(iv).)
\end{itemize}
\end{prob}

\bigskip
Напомним, что при $n>1$
\emph{группой диэдра}~$\mathrm{D}_{n}$ порядка~$2n$
называется группа симметрий плоскости, сохраняющих правильный $n$-угольник.
В частности, есть изоморфизмы $\mathrm{D}_2\cong\Z/2\Z\times\Z/2\Z$
и $\mathrm{D}_3\cong\mathrm{S}_3$.
Для любого $n$ имеется (неканонически расщепимая) точная последовательность
групп
$$
1\to \Z/n\Z\to\mathrm{D}_n\to\Z/2\Z\to 1\,.
$$

Будем говорить, что многочлен степени $n>1$ над полем~$k$
\emph{диэдрален}, если он сепарабелен, неприводим над $k$,
и группа Галуа его поля разложения изоморфна
группе диэдра~$\mathrm{D}_{n}$.
Будем говорить, что пара $(P, a)$ является \emph{диэдральной},
если~$P$~--- диэдральный многочлен,
и при этом квадратичное расширение поля $k$, соответствующее
каноническому сюръективному
гомоморфизму из группы диэдра на группу~$\Z/2\Z$,
совпадает с расширением~\mbox{$k\subset k(\sqrt{a})$}.
В частности, если пара~\mbox{$(P,a)$} диэдральна, то~$a$ не является
квадратом в поле~$k$.
Примером диэдральной пары является пара~\mbox{$(P, a)$},
где $P$ --- неприводимый многочлен степени $3$, и $a$ ---
его дискриминант, при условии, что $a\not\in (k^*)^2$.

\bigskip
\begin{prob}{\bf Диэдральные пары}
\label{prob:dihedral-pairs}

Допустим, что пара $(P,a)$ диэдральна. Напомним, что в соответствии с нашими предположениями степень $d$ многочлена~$P$ нечётна.

\begin{itemize}
\item[(i)]
Отождествим группу Галуа $G$ поля разложения многочлена $P$
с группой диэдра $\mathrm{D}_{d}$, действующей
симметриями правильного $d$-угольника.
Докажите, что существует
\mbox{$G$-эквивариантная}
биекция между множеством $R$ корней многочлена $P$
и множеством вершин $d$-угольника.
(Указание:
поскольку $G$ действует на $R$ транзитивно, стабилизатор
$\mathrm{Stab}_G(\theta)$
любого корня
$\theta\in R$ имеет порядок~$2$.
Так как $d$ нечётно, то
$\mathrm{Stab}_G(\theta)$ сохраняет единственную вершину $v$ в $d$-угольнике.
Теперь можно воспользоваться транзитивностью
действия $\mathrm{D}_d$ на вершинах и проверить, что отображение
$\theta\mapsto v$ даёт требуемую биекцию.)
\item[(ii)]
Пусть~\mbox{$\theta\in k^{sep}$}~--- корень многочлена $P$.
Покажите, что поле $k(\theta, \sqrt{a})$ является полем разложения
многочлена $P$ над $k$.
(Указание: пусть $K$ обозначает поле разложения
многочлена $P$ над $k$. Из определения диэдральной пары
следует, что группа $\Gal\big(K/k(\sqrt{a})\big)$
является циклической подгруппой
$\Z/d\Z$ в группе~\mbox{$G\cong \mathrm{D}_{d}$}.
По теории Галуа группа $\Gal\big(K/k(\theta, \sqrt{a})\big)$
является стабилизатором в группе~\mbox{$\Z/d\Z$} корня $\theta$.
Из пункта~(i) следует, что этот стабилизатор тривиален.)
\item[(iii)]
Пусть~\mbox{$\theta\in k^{sep}$}~--- корень многочлена $P$.
Покажите, что над полем $k(\theta)$ имеется разложение
$$P(u)=\alpha (u-\theta)\big(P_1(u)^2-aP_2(u)^2\big)\,,$$
где $P_1$ и $P_2$ --- многочлены степени $(d-1)/2$, а $\alpha$~---
элемент поля~$k$.
(Указание:
из пункта~(ii) мы знаем, что
$$G\cong\Gal\big(k(\theta,\sqrt{a})/k\big)\,.$$
По пункту~(i) подгруппа
$$\Gal\big(k(\theta, \sqrt{a})/k(\theta)\big)\cong\Z/2\Z$$
группы $G$
действует свободно на множестве $R\setminus\{\theta\}$.
Рассмотрим два различных корня
$$\eta=x+\sqrt{a}y,\quad \bar{\eta}=x-\sqrt{a}y, \quad x,y\in k(\theta)$$
многочлена $P$,
переставляемые этим действием. Тогда
$$(u-\eta)(u-\bar{\eta})=(u-x)^2-ay^2.$$
Теперь из мультипликативности нормы для расширения полей
$$k(u,\theta)\subset k(u,\theta, \sqrt{a})$$
следует, что
над полем $k(\theta)$ многочлен $P$ имеет требуемый вид.
Здесь в качестве $\alpha$ можно взять, например, старший коэффициент
многочлена~$P$.)
\end{itemize}
\end{prob}

\bigskip
В дальнейшем, как и в начале этого раздела, $P$ является произвольным
сепарабельным многочленом над~$k$ нечётной степени~$d$, и $a\in k^*$.
При этом~$S$
обозначает соответствующую поверхность Шатле.

\begin{prob}{\bf Рациональность поверхности Шатле над расширениями поля $k$}
\label{prob:Chatelet-rational} \hspace{0cm}

\begin{itemize}
\item[(i)]
Покажите, что над полем $k(\sqrt{a})$ поверхность
$S$ рациональна.
(Указание: проверьте, что сечения $E$ и $E'$ определены над полем
$k(\sqrt{a})$.)
\item[(ii)] Предположим, что пара $(P,a)$ диэдральна.
Пусть $\theta\in k^{sep}$ --- корень многочлена~$P$.
Докажите, что над полем $k(\theta)$ поверхность
$S$ рациональна.
(Указание:
по упражнению~\ref{prob:dihedral-pairs}(iii)
над полем $k(\theta)$ имеется разложение
$$P(u)=\alpha (u-\theta)\big(P_1(u)^2-aP_2(u)^2\big)\,.$$
Пользуясь мультипликативностью нормы,
после замены координат уравнение~\eqref{eq:Chatelet-affine}
аффинного открытого подмножества
$S^0\subset S$ можно переписать в виде
$$\hat{v}_1^2-a\hat{v}_2^2=\hat{u}\,.$$
Проекция поверхности $S_0$ на плоскость $\Ab^2$ с координатами $\hat{v}_1$ и
$\hat{v}_2$
является изоморфизмом.)
\end{itemize}
\end{prob}

\bigskip
Упражнение~\ref{prob:Chatelet-rational} показывает, что поверхность
Шатле рациональна над многими просто устроенными
конечными расширениями поля $k$.
Оказывается, что над самим полем $k$ такие поверхности часто бывают
нерациональны.
Чтобы это проверить, мы будем пользоваться теоремой
В.\,А.\,Исковских
о нерациональности расслоений на коники (см.~\cite[Теорема~1.6]{Is1},
\cite[Теорема~1]{Is2} и~\cite[Теорема~2]{Is3};
несколько более современное изложение можно найти
в~\cite[\S4]{Is-UMN}).

Предположим, что поле $k$ совершенно.
Пусть на гладкой проективной поверхности $\Sigma$
над полем $k$ задана структура расслоения на коники
$f\colon\Sigma\to\P^1$. Предположим, что данное расслоение на коники
\emph{относительно минимально} над $k$,
то есть на $\Sigma$ над полем
$k^{sep}$ не существует
$G_k$-инвариантного набора попарно непересекающихся
\mbox{$(-1)$-кривых}, стягиваемых отображением $f$.
\emph{Теорема Исковских} утверждает, что в этом случае
рациональность поверхности $\Sigma$ над $k$ равносильна тому,
что на $\Sigma$ есть $k$-точка и количество геометрически приводимых слоёв
расслоения $f$ над $k^{sep}$ не больше трёх.
(Обратите внимание, что из гладкости поверхности $\Sigma$
следует, что все слои морфизма $f$ приведены!)

В частности, теорема Исковских даёт другой (более геометрический)
способ решить упражнение~\ref{prob:Chatelet-rational}(ii).
Действительно, над полем $k(\theta)$ можно стянуть компоненты геометрически приводимых
слоёв расслоения $f$ так, что у получившегося расслоения
на коники останется не более двух геометрически приводимых слоёв,
каждый из которых определён над~$k(\theta)$.
На соответствующей поверхности есть точка
над полем~$k(\theta)$, откуда следует её рациональность.

\begin{prob}{\textbf{Нерациональность поверхности Шатле над полем $k$}}
\label{prob:X-nonrational} \hspace{0cm}

\begin{itemize}
\item[(i)]
Докажите, что расслоение на коники $f\colon S\to\P^1$ относительно минимально над~$k$,
если и только если
для любого корня $\theta$ многочлена $P$ поле $k(\theta)$ не содержит
$\sqrt{a}$.
(Указание:
предположим сначала, что для любого $\theta$ поле $k(\theta)$ не содержит
$\sqrt{a}$. Тогда
для каждого геометрически приводимого слоя $F_p$ расслоения $f$
в группе $G_k$ найдётся элемент, сохраняющий точку $p\in\Pb^1$
и при этом меняющий местами компоненты этого слоя.
Это означает, что
поверхность $S$ относительно минимальна над $k$.
Предположим теперь, что для некоторого корня $\theta$
поле $k(\theta)$ содержит $\sqrt{a}$.
Проверьте, что $G_k$-орбита любой из двух компонент слоя $F_{\theta}$
состоит из непересекающихся $(-1)$-кривых.)
\item[(ii)]
Пусть поле $k$ совершенно. Предположим, что
для любого корня $\theta$
многочлена~$P$ поле $k(\theta)$ не содержит
$\sqrt{a}$ и степень многочлена $P$ не меньше~$3$.
Покажите, что в этом случае
поверхность Шатле $S$ нерациональна над $k$.
(Указание: воспользуйтесь пунктом~(i) и теоремой Исковских.)
В частности, это выполнено, если многочлен~$P$ неприводим,
так как по предположению его степень нечётна.
\end{itemize}
\end{prob}

\bigskip
Цель следующего
упражнения --- показать, что для диэдральных пар
группа Пикара
$\Pic\big(S_{k^{sep}}\big)$
является стабильно пермутационным $G_k$-модулем.
Мы будем пользоваться результатом, доказанным Эндо и Миятой
(см.~\cite[Theorem~1.5]{EndoMiyata}
и~\mbox{\cite[Lemma~1.2]{EndoMiyata}}).
Пусть $G$ --- конечная группа, у которой все силовские подгруппы
циклические.
Пусть дана короткая точная последовательность $G$-модулей
$$0\to M'\to M\to M''\to 0\,,$$
где $M'$ и $M''$ удовлетворяют
условию~\eqref{eq:Endo-Miyata} из упражнения~\ref{exercise:stab-permutation-module}(ii) и являются свободными конечно порождеными абелевыми группами.
\emph{Теорема Эндо--Мияты} утверждает,
что такая точная последовательность расщепляется.

\begin{prob}{\textbf{Группа Пикара и универсальный торсор для поверхности Шатле}}
\label{prob:Chatelet-Pic} \hspace{0cm}

Предположим, что пара $(P,a)$ диэдральна.
Пусть~\mbox{$G\cong \mathrm{D}_{d}$} обозначает группу
Галуа поля разложения многочлена $P$.

\begin{itemize}
\item[(o)]
Проверьте, что действие группы $G_k$ на все модули, входящие в точную последовательность из упражнения~\ref{prob:Chatelet-easy}(v),
пропускается через факторгруппу~$G$.
\item[(i)]
Пусть $H\subset G$ --- циклическая подгруппа порядка $d$.
Покажите, что группа~$\Pic\big(S_{k^{sep}}\big)$
является стабильно пермутационным
$H$-модулем.
(Указание: поверхность $S$ рациональна
над полем $k(\sqrt{a})$
по упражнению~\ref{prob:Chatelet-rational}(i).
Теперь можно воспользоваться
упражнением~\ref{exercise:stab-permutation-module-from-Pic}(iv).
Обратите внимание, что в случае поверхностей использовавшееся
в упражнении~\ref{exercise:stab-permutation-module-from-Pic}
предположение о том, что характеристика поля $k$ равна нулю, необязательно,
так как в размерности~$2$ имеется разрешение особенностей
над произвольным полем.)
\item[(ii)]
Пусть $\theta\in k^{sep}$ --- корень многочлена $P$.
Рассмотрим его стабилизатор
$$\Z/2\Z\cong H_{\theta}\subset G\,.$$
Покажите, что $\Pic\big(S_{k^{sep}}\big)$
является стабильно пермутационным
$H_{\theta}$-модулем.
(Указание: поверхность $S$ рациональна
над полем $k(\theta)$
по упражнению~\ref{prob:Chatelet-rational}(ii), так что
можно воспользоваться
упражнением~\ref{exercise:stab-permutation-module-from-Pic}(iv).)
\item[(iii)]
Пусть $l$ --- простое число, и $G_l\subset G$ --- силовская
$l$-подгруппа.
Проверьте, что~$\Pic\big(S_{k^{sep}}\big)$ является стабильно пермутационным
$G_l$-модулем.
(Указание: так как степень многочлена
$P$ нечётна, каждая нетривиальная $2$-подгруппа в~$G$ является одной из групп~$H_{\theta}$.
Поэтому при $l=2$ утверждение следует из пункта~(ii).
Таким образом, можно считать, что $l$ нечётно. Тогда $G_l$
содержится в подгруппе~$H$.
В~этом случае утверждение следует из пункта~(i)
и упражнения~\ref{exercise:stab-permutation-module-restriction}.)
\item[(iv)]
Покажите, что $\Pic\big(S_{k^{sep}}\big)$
является стабильно пермутационным
$G$-модулем, а следовательно, и стабильно пермутационным
$G_k$-модулем.
(Указание: примените к точной последовательности
из упражнения~\ref{prob:Chatelet-easy}(v)
теорему Эндо--Мияты, воспользовавшись пунктом~(iii)
и упражнением~\ref{exercise:stab-permutation-module}(iii)
для проверки её условий.)
\item[(v)]
Докажите, что для поверхности $S$ существует
универсальный торсор $V_S$, и при этом $S$ и $V_S$ стабильно
бирационально эквивалентны.
(Указание: поверхность~$S$ проективна, так что
все обратимые функции на ней постоянны. Кроме того,
поверхность~$S$ содержит $k$-точку --- её нетрудно найти на геометрически приводимом
слое $D_{\infty}$ расслоения на коники $f$.
Таким образом, универсальный торсор $V_S$ существует
по упражнению~\ref{prob:spusk-universalnogo-torsora}(iv).
То, что $S$ и $V_S$
стабильно бирационально
эквивалентны, следует из упражнения~\ref{prob:X-vs-UX}.)
\end{itemize}
\end{prob}

\newpage
\section{Пример нерационального стабильно рационального многообразия}
\label{section:Polietilen}

\subsection{План построения примера}
В этой главе мы построим, следуя
статье А.\,Бовиля, Ж.-Л.\,Кольо-Телена,
\mbox{Ж.-Ж.\,Сансюка} и П.\,Свиннертон-Дайера~\cite{CT-i-drugie},
многообразие $X$ над совершенным полем~$k$ характеристики, не равной~$2$,
которое будет нерациональным,
но стабильно рациональным (в частности, унирациональным).
Если $k$ алгебраически замкнуто, то такое многообразие
$X$, очевидно, не может быть кривой.
На самом деле, то же верно и для произвольного поля $k$
(это следует из того, что коника без точек не может быть
унирациональной).
Из критерия рациональности Кастельнуово следует,
что над алгебраически замкнутым полем характеристики~$0$ всякая
унирациональная поверхность рациональна
(см.~\cite[Corollary~V.5]{Beauville}), так что
в случае алгебраически замкнутого поля~$k$ нулевой характеристики
многообразие $X$ не может быть поверхностью.
Таким образом, естественно искать пример среди нерациональных
поверхностей над алгебраически незамкнутым полем $k$, которые
рациональны над алгебраическим замыканием $\bar{k}$
поля $k$.
Один из простейших примеров поверхностей
такого сорта можно найти среди расслоений на коники над $\P^1$,
пользуясь теоремой Исковских, сформулированной
в разделе~\ref{subsection:Chatelet}.
Точнее говоря, мы будем
работать с квазипроективными поверхностями,
и $X$ будет открытым подмножеством в некотором нерациональном
(проективном) расслоении
на коники~$\widetilde{X}$.

Результаты, собранные в разделе~\ref{subsection:Chatelet},
подсказывают нам направление, в котором можно пытаться найти интересующий
нас пример. А именно, мы хотели бы построить
проективную поверхность $\widetilde{X}$,
которая была бы нерациональной, но стабильно рациональной.
Согласно упражнению~\ref{exercise:stab-permutation-module-from-Pic}(v),
группа Пикара $\Pic\big(\widetilde{X}_{k^{sep}}\big)$ в этом случае является
стабильно пермутационным $G_k$-модулем, что сильно сужает область
поиска. Более того, серия примеров, обладающих такими свойствами,
у нас уже заготовлена~--- это (многие) поверхности Шатле,
как мы знаем из упражнений~\ref{prob:X-nonrational}(ii)
и~\ref{prob:Chatelet-Pic}(iv).
Далее,
из стабильной пермутационности $G_k$-модуля
$\Pic\big(\widetilde{X}_{k^{sep}}\big)$
и упражнения~\ref{prob:X-vs-UX}
(см. также упражнение~\ref{prob:Chatelet-Pic}(v))
следует, что поверхность $\widetilde{X}$ стабильно бирационально
эквивалентна своему универсальному торсору $V_{\widetilde{X}}$.
Таким образом, для решения исходной задачи
достаточно установить стабильную рациональность
универсального торсора для подходящей поверхности Шатле
(на самом деле, мы будем проверять их рациональность, а не просто
стабильную рациональность).
Оказывается, что такие торсоры можно задать
явными уравнениями (см.~\cite{CT-S-SDa} и~\cite{CT-S-SDb}).
Вопрос об их рациональности
сводится к вопросу о рациональности
пересечений двух квадрик.

Нам будет удобнее работать даже не с универсальным
торсором~$V_{\widetilde{X}}$,
а с некоторым его фактормногообразием по подходящим
образом выбранному алгебраическому тору, которое также является
рациональным, и при этом тоже стабильно бирационально эквивалентно
поверхности $X$ (это фактормногообразие мы будем обозначать через $V$,
чтобы напомнить о его связи с универсальным торсором).
Смысл замены универсального торсора на это многообразие заключается
в понижении размерности. В конечном итоге мы опустим
промежуточные детали рассуждения с универсальным торсором,
и окажется, что
$X\times\Pb^3$ бирационально эквивалентно
пересечению двух квадрик в~$\Pb^7$,
а рациональность этого пересечения можно
установить при помощи проекции из прямой.

Отметим, что для алгебраически замкнутого поля $k$ можно построить пример
трёхмерного нерационального стабильно рационального многообразия.
Эта конструкция основана
на двумерном примере над незамкнутым полем; нерациональность
в этом случае устанавливается при помощи теории промежуточного якобиана
(подробности см. в~\cite[\S3]{CT-i-drugie}).

\bigskip
Выше мы объяснили, ``где искать'' нужный нам пример.
Теперь мы
объясним чуть более детально
план построения примера (как было уже упомянуто,
всё происходящее можно воспринимать,
вовсе не думая про универсальный торсор, но помнить о нём
всё же полезно для понимания общей картины).
Будет рассматриваться коммутативная
диаграмма многообразий

\begin{equation}\label{diag:plan-Polietilena}
\xymatrix{
X\ar@{^(->}[r]\ar@{->}[dd]&A\ar@{->}[d]^{f} &
A\times D\ar@{->}[l]_{\mathrm{pr}_A}\ar@{->}[d]^{{}\cong{}} \\
&B & C\ar@{->}[l]_{g}\\
\Gamma\ar@{^(->}[ru]&& V\ar@{->}[ll]\ar@{^(->}[u]
}
\end{equation}
В этой диаграмме $\Gamma$ будет рациональной кривой на
многообразии $B$. При этом~\mbox{$X=f^{-1}(\Gamma)$}
и~\mbox{$V=g^{-1}(\Gamma)$}.
Многообразие $X$ является неприводимой нерациональной поверхностью,
а многообразия~$V$ и $D$ рациональны.
Так как $C\cong A\times D$, то $V\cong X\times D$.
Таким образом,~$X$ будет давать искомый пример.
Более того,
в диаграмме~\eqref{diag:plan-Polietilena}
многообразия~\mbox{$A$,~$B$,~$C$} и~$D$
будут алгебраическими торами, а все морфизмы будут морфизмами алгебраических
торов. При этом после подходящих компактификаций морфизм~$f$
окажется ограничением расслоения
на коники,~$g$ окажется ограничением квадратичного рационального
отображения из $\P^7$ в $\P^3$, а $\Gamma$~--- открытым подмножеством прямой в $\P^3$. Это позволит рассматривать $V$ как
открытое подмножество в пересечении двух квадрик.

\subsection{Поля $K$, $k'$ и $K'$}

Везде в дальнейшем мы будем предполагать, что основное поле $k$ совершенно
и имеет характеристику, отличную от~$2$.
Отметим, что совершенность будет нужна для применения теоремы Исковских
в доказательстве нерациональности поверхности~$X$,
а в остальных местах без неё можно обойтись, предполагая только сепарабельность
рассматриваемых расширений поля~$k$. Предположение о характеристике
потребуется в доказательстве рациональности многообразия~$V$,
точнее говоря, для того, чтобы можно было эффективно работать
с квадриками, задающими компактификацию многообразия~$V$.

Пусть $P\in k[x]$ --- неприводимый многочлен степени
$3$ над полем $k$ с дискриминантом $a$.
Предположим, что $a\not\in (k^*)^2$, и рассмотрим
поля~\mbox{$k'=k(\sqrt{a})$} и~\mbox{$K=k(\theta)$},
где $\theta\in\bar{k}$ --- некоторый
корень многочлена $P$. Наконец, пусть~$K'$~--- композит расширений $K$ и $k'$
поля $k$, то есть $K'$ --- поле разложения многочлена~$P$.

\bigskip
\begin{prob}{\bf Группы Галуа}
\label{prob:Galois-3}

Докажите, что имеют место изоморфизмы:
$$\Gal(K'/k')\cong\Z/3\Z\,,\quad
\Gal(K'/K)\cong\Gal(k'/k)\cong
\Z/2\Z\,,\quad\Gal(K'/k)\cong\mathrm{S}_3\,.
$$
В частности, пара $(P,a)$ диэдральна
(ср. с примером перед упражнением~\ref{prob:dihedral-pairs}).
Проверьте, что расширение $k\subset K$ не является расширением
Галуа.
\end{prob}

Как и в главе~\ref{section:Weil},
в дальнейшем через $K^*$ мы будем обозначать как алгебраический тор
$R_{K/k}(\Gb_m)$, так и мультипликативную группу поля $K$
(канонически изоморфную группе $k$-точек этого тора).
Аналогичным образом, через $K$ мы будем обозначать как
алгебраическую группу $R_{K/k}(\Gb_a)$, так и аддитивную
группу поля $K$. Те же соглашения будут использоваться и для
других полей над $k$ (включая само поле~$k$). В частности, символ~$k$ обозначает
алгебраическую группу $\Gb_{a}=\Spec(k[t])$.

\subsection{Нерациональное расслоение на коники}

Рассмотрим морфизмы алгебраических торов $\Nm_{K/k}\colon K^*\to k^*$ и
$\Nm_{k'/k}\colon k'^*\to k^*$.
Пусть $A=K^*\times_{k^*}k'^*$ --- расслоенное произведение
$K^*$ и $k'^*$ над $k^*$, а $B=K^*$.
Пусть~\mbox{$f\colon A\to B$}~--- естественная проекция.

\begin{prob}{\bf Явный вид многообразия $A$}
\label{prob:equation-of-A} \hspace{0cm}

Покажите, что многообразие $A$ является открытым подмножеством
в аффинной гиперповерхности в $\Ab^5$, заданной уравнением
$$F_3(u_1,u_2,u_3)=v_1^2-av_2^2\,,$$
где $u_1$, $u_2$, $u_3$, $v_1$ и $v_2$ --- координаты
на $\Ab^5$, а
$F_3$ --- однородный кубический многочлен.
(Указание: $F_3$ задаёт норму $\Nm_{K/k}$.)
\end{prob}

\bigskip
Элемент $\theta\in K$ определяет $k$-точку алгебраической
группы $K$ над $k$.
Групповой сдвиг на $-\theta$, применённый к подгруппе $k\subset K$,
определяет подмногообразие~${k-\theta}$ в $K$ (это аффинная прямая в трёхмерном
аффинном пространстве над $k$). Кроме того, в~$K$ имеется открытое
по Зарисскому подмножество $K^*$.
Рассмотрим алгебраическую кривую
\begin{equation}\label{eq:Gamma}
\Gamma=(k-\theta)\cap K^*\subset K\,.
\end{equation}
Пусть $X=f^{-1}(\Gamma)\subset A$.

\bigskip
\bigskip
\bigskip
\bigskip
\bigskip
\begin{prob}{\bf Нерациональность поверхности $X$}
\label{prob:equation-of-X} \hspace{0cm}
\nopagebreak
\\
\begin{itemize}
\item[(i)] Покажите, что
$\Gamma$ изоморфно дополнению в $\Ab^1$ к замкнутому
подмножеству, заданному многочленом $P$.
(Указание:
согласно упражнению~\ref{prob:restr-tori}(v),
многообразие $K^*$ является дополнением в $K$ к гиперповерхности
$\{\Nm_{K/k}=0\}$.
Покажите, что
ограничение нормы $\Nm_{K/k}$ на $k-\theta\subset K$ равно $P(t)$,
где $t$~---
координата на $k\cong\Gb_a$. Это можно сделать,
например, рассмотрев расширение скаляров
с $k$ на $k^{sep}$ (или на $K'$)
и воспользовавшись указанием к упражнению~\ref{prob:restr-tori}(iv).)
\item[(ii)] Покажите,
что многообразие $X$ является открытым подмножеством
в аффинной поверхности в $\Ab^3$, заданной уравнением
$$P(u)=v_1^2-av_2^2\,,$$
где $u$, $v_1$ и $v_2$ --- координаты
на $\Ab^3$.
Таким образом, $X$ является открытым
подмножеством в поверхности Шатле (см. главу~\ref{section:Weil}).
\item[(iii)] Докажите, что поверхность $X$ нерациональна.
(Указание: воспользуйтесь упражнениями~\ref{prob:X-nonrational}(ii)
и~\ref{prob:Galois-3}.
Отметим, что в нашем случае степень многочлена $P$
равна трём, то есть соответствующая поверхность
Шатле имеет ровно четыре
геометрически приводимых слоя
над~$\P^1$. Поэтому
вместо теоремы Исковских в её полной общности, которую мы использовали
в упражнении~\ref{prob:X-nonrational}(ii),
достаточно результата
работы~\cite{Is3}.)
\end{itemize}
\end{prob}

\subsection{Рациональное пересечение двух квадрик}

Мы начнём этот раздел с одной геометрической конструкции,
устанавливающей рациональность (достаточно хорошего)
пересечения двух квадрик, содержащего прямую.
Напомним, что многообразие $C\subset\P^n$ называется
\emph{конусом с вершиной в точке~\mbox{$v\in C$}}, если для любой точки
$c\in C$, отличной от $v$, прямая $\langle v, c\rangle$
содержится в $C$. Нетрудно показать, что
множество всех вершин конуса образует проективное
подпространство в $\P^n$, которое тоже часто называют вершиной конуса
(к счастью, обычно из контекста ясно, о чём идёт речь).

\bigskip
Пусть многообразие $M\subset\P^n$ является полным пересечением двух
квадрик $Q_1$ и $Q_2$, определённых над полем $k$.
Иными словами, $\dim(M)=n-2$, и однородный идеал
многообразия $M$ порождается уравнениями квадрик $Q_1$ и $Q_2$.
В частности, для любой точки $p\in M$ имеется равенство
$$\mathbb{T}_p(M)=\mathbb{T}_p(Q_1)\cap \mathbb{T}_p(Q_2)\,,$$
где $\mathbb{T}_p(-)$ обозначает
вложенное проективное касательное пространство к соответствующему проективному
многообразию в точке $p$.
Предположим также,
что многообразие $M$ неприводимо и
содержит определённую над $k$ прямую $l$, для которой
выполнены два условия:
\begin{itemize}
\item[(A)]
прямая~$l$ не содержится в множестве особых точек
$\Sing(M)$ многообразия~$M$;
\item[(B)] многообразие
$M$ не является конусом с вершиной, лежащей на прямой~$l$.
\end{itemize}
Цель следующего упражнения~---
доказать рациональность многообразия~$M$.

\begin{prob}{\textbf{Рациональность пересечения двух квадрик, содержащего
прямую}}\label{prob:two-quadrics}

Пусть $\Pi\subset\Pb^n$ --- общая двумерная плоскость
над $\bar{k}$,
проходящая через прямую $l$. Обозначим через $q_i$ ограничения
квадрик $Q_i$ на плоскость $\Pi$.
Заметим, что из общности~$\Pi$ следует, что $q_i$ не совпадают с $\Pi$,
то есть являются (возможно, приводимыми или неприведенными)
кониками на~$\Pi$, содержащими~$l$.

\begin{itemize}
\item[(o)] Покажите, что многообразие $M$ геометрически неприводимо
и не содержится ни в какой гиперплоскости.
(Указание: так как прямая $l$ определена над полем~$k$,
и многообразие $M$ неприводимо, то $l$ содержится в каждой
неприводимой компоненте многообразия~$M_{\bar{k}}$.
Так как при этом $l\not\subset\Sing(M)$ по условию~(A),
то такая компонента только одна, то есть $M$ геометрически неприводимо.
Предположим теперь, что
$M$ содержится в какой-то гиперплоскости~\mbox{$H\subset\P^n$}.
Тогда $M$ содержится в каждом из пересечений
$$\hat{Q}_i=Q_i\cap H,\quad i=1,2\,.$$
Заметим, что ни одно
из этих пересечений не совпадает со всей гиперплоскостью $H$ --- иначе
$M$ либо приводимо, либо совпадает с $H$ (в последнем случае~$M$
не является полным пересечением).
Так как
$$\dim(M)=n-2=\dim(\hat{Q}_i)\,,$$
то $M$ является некоторой общей неприводимой компонентой $\hat{Q}_1$ и
$\hat{Q}_2$. В частности,
$$\deg(M)\leqslant \deg(\hat{Q}_1)\leqslant 2<4=\deg(Q_1)\cdot\deg(Q_2)\,,$$
что невозможно, так как по предположению многообразие $M$
является полным пересечением $Q_1$ и $Q_2$.)
Эти утверждения не будут использоваться в дальнейшем доказательстве,
однако, возможно, помогут более полно осознать
геометрическую картину.
\item[(i)]
Покажите, что $q_i\neq 2l$ при $i=1,2$.
(Указание: если $q_i=2l$ для некоторого~\mbox{$i=1,2$},
то квадрика $Q_i$
особа в каждой точке прямой~$l$.
Чтобы проверить это, рассмотрите касательное пространство
к $Q_i$ в соответствующей точке~--- по предположению оно содержит
общую плоскость, проходящую через прямую~$l$.
В этом случае также $l\subset\Sing(M)$, что противоречит условию~(A).)
Таким образом, $q_i=l\cup l_i$, где $l_i\neq l$.
\item[(ii)]
Покажите, что прямые $l_1$ и $l_2$ из пункта~(i) различны.
(Указание: в противном случае квадрики $Q_1$ и $Q_2$ совпадают.)
\item[(iii)]
Пусть $l'\subset l$ --- (непустое) открытое подмножество,
состоящее из таких точек~\mbox{$p\in l$},
что $p\not\in\Sing(M)$. Покажите, что прямые $l_1$ и $l_2$ не могут
пересекаться в точке $p\in l'$.
(Указание: предположите обратное. Выведите из этого,
что плоскость $\Pi$ содержится в касательном пространстве
$\mathbb{T}_p(M)\cong\P^{n-2}$.
Далее, покажите, что для общей $\Pi$ такое включение невозможно,
найдя оценку сверху на размерность многообразия всех плоскостей,
содержащихся в $\bigcup\limits_{x\in l'}\mathbb{T}_x(M)$.)
\item[(iv)]
Докажите, что прямые $l_1$ и $l_2$ пересекаются в (единственной) точке~$p$,
не лежащей на прямой~$l$.
(Указание:
предположим, что $p\in l$.
Из пункта~(iii) следует, что
$p$ принадлежит конечному множеству $l\setminus l'$. Из неприводимости
многообразия всех плоскостей, проходящих через $l$,
следует, что точка $p$ --- одна и та же для всех $\Pi$.
Рассматривая касательные пространства к квадрикам~$Q_i$ в точке~$p$,
мы видим, что обе квадрики $Q_i$ особы в~$p$, а значит, являются конусами
с вершиной в~$p$. Отсюда следует, что многообразие~$M$ также
является конусом с вершиной в точке $p\in l$, что противоречит
условию~(B).)
\item[(v)]
Докажите, что многообразие
$M$ рационально над $k$. (Указание: рассмотрите проекцию
$$M\dasharrow\P^{n-2}$$
многообразия $M$ с центром в прямой $l$.
Пользуясь пунктами~(i), (ii) и~(iv), покажите, что
это бирациональное отображение.)
\item[(vi)]
Где в предыдущих рассуждениях использовалась неприводимость
многообразия~$M$?
(Указание: из пунктов~(i), (ii) и~(iv) на самом деле следует только то,
что одна из неприводимых компонент $M$ при проекции из $l$
отображается на свой образ бирационально. Без предположения неприводимости
это не исключало бы того, что есть какая-то другая компонента,
размерность которой при проекции из $l$ падает. Например,
если квадрика
$$Q_1=H\cup H'$$
приводима, квадрика $Q_2$ выбрана общим образом,
а прямая $l$ выбрана общим образом в $Q_2\cap H$, то
при проекции из $l$ неприводимая компонента~\mbox{$Q_2\cap H'$}
многообразия $M$ отображается на свой образ бирационально, а
неприводимая компонента $Q_2\cap H$ с точностью до бирационального
преобразования является $\P^1$-расслоением над своим образом.)
\end{itemize}
\end{prob}

\bigskip
Вернёмся к построению нашего примера.
Рассмотрим алгебраический тор
$$\hat{C}=K'^*\times k'^*$$
и морфизм алгебраических торов
$$\hat{g}\colon\hat{C}\to B,\quad
(x,y)\mapsto \Nm_{K'/K}(x)\cdot\Nm_{k'/k}(y)^{-1}\,.$$
Пусть
$$C=\hat{C}/k^*\,,$$
где $k^*$ вложено в $K'^*\times k'^*$ диагонально.
Морфизм алгебраических торов $\hat{g}$ пропускается через морфизм алгебраических торов
$$g\colon C\to B\,.$$

\bigskip

\begin{prob}{\bf Неприводимость слоёв морфизма $g$}
\label{prob:Ker-g} \hspace{0cm}

\begin{itemize}
\item[(i)]
Докажите, что имеется точная последовательность алгебраических
групп
$$
1\to \Ker(\Nm_{K'/K})\times\Ker(\Nm_{k'/k})\to \Ker(\hat{g})\to k^*\to 1\,.
$$
\item[(ii)]
Выведите из пункта~(i), что многообразие $\Ker(\hat{g})$
неприводимо.
(Указание: воспользуйтесь
упражнением~\ref{prob:Ker-Nm} и тем, что
многообразие, расслоенное над неприводимой базой
с геометрически неприводимыми слоями, само неприводимо.)
\item[(iii)]
Докажите, что ядро морфизма $g$ неприводимо.
(Указание: имеется сюръективный морфизм
$\Ker(\hat{g})\to\Ker(g)$.)
\item[(iv)] Докажите, что все слои морфизма $g$ неприводимы.
\end{itemize}
\end{prob}

\bigskip

Положим $V=g^{-1}(\Gamma)\subset C$ (см.~\eqref{eq:Gamma}).
Напомним, что $\theta\in\bar{k}$ --- корень многочлена~$P$.

\begin{prob}{\bf Компактификация многообразия $V$}
\label{prob:compactification-of-C} \hspace{0cm}

\begin{itemize}
\item[(i)]
Рассмотрим $B$ как открытое подмножество в
проективном пространстве~\mbox{$\Pb(K\oplus k)\cong\Pb^3$},
а $C$ --- как открытое подмножество в
проективном пространстве~\mbox{$\Pb(K'\oplus k')\cong\Pb^7$}.
Покажите, что при
этом морфизм $g$ продолжается до квадратичного рационального отображения
$$\tilde{g}\colon\Pb(K'\oplus k')\dasharrow \Pb(K\oplus k)\,, \quad
(x:y)\mapsto \big(\Nm_{K'/K}(x):\Nm_{k'/k}(y)\big)\,.$$
\item[(ii)]
Покажите, что $V$ является плотным
открытым подмножеством в проективном многообразии
$$\widetilde{V}=\{(x:y)\in\Pb(K'\oplus k')\mid
\pi\big(\Nm_{K'/K}(x)+\theta\cdot\Nm_{k'/k}(y)\big)=0\}\,,$$
где $\pi\colon K\to K/k$ --- естественная проекция.
\end{itemize}
\end{prob}

\bigskip
Таким образом, многообразие $\widetilde{V}$ из
упражнения~\ref{prob:compactification-of-C}(ii) является
компактификацией многообразия~$V$.
Из упражнения~\ref{prob:Ker-g}(iv) следует, что многообразие $V$
(а~значит, и $\widetilde{V}$) неприводимо.
Цель следующего упражнения --- доказать, что у
многообразия~$\widetilde{V}$
нет особых точек, определённых над $k$.

\begin{prob}{\textbf{Особенности многообразия $\widetilde{V}$}}
\label{prob:singularities-of-V}

\begin{itemize}
\item[(i)]
Пусть $E\subset F$ --- квадратичное сепарабельное расширение поля $E$
с нетривиальным автоморфизмом $\sigma\in\Gal(F/E)$. Рассмотрим $E$ и $F$
как аффинные пространства над $E$.
Пусть морфизм $N\colon F\to E$ задан формулой
$$N(s)=\Nm_{F/E}(s)\,.$$
Отождествим касательное пространство к $F$ в точке $s\in F$
с самим пространством $F$, а касательное пространство
к $E$ в точке $t\in E$ с пространством $E$.
Тогда дифференциал
отображения $N$ в точке~$s$ равен
$$dN\colon\alpha\mapsto\Tr_{F/E}\big(\alpha\cdot\sigma(s)\big),
\quad \alpha\in F\,,$$
где $\Tr_{F/E}\colon F\to E$ обозначает след в сепарабельном
расширении $E\subset F$.
(Указание: воспользуйтесь тем, что
$$\Nm_{F/E}\big(s\big)=s\cdot\sigma(s)\,,$$
и продифференцируйте.)
\item[(ii)]
Для каждой точки $(x,y)\in K'\oplus k'$
рассмотрим отображение
$$T_{x,y}\colon (\alpha,\beta)\mapsto
\pi\Big(\Tr_{K'/K}\big(\alpha\cdot\sigma(x)\big)+
\theta\cdot\Tr_{k'/k}\big(\beta\cdot\sigma(y)\big)\Big)\,,\quad
(\alpha,\beta)\in K'\oplus k'\,.$$
Покажите, что точка $(x:y)\in \widetilde{V}(k)$ особа тогда и только тогда,
когда отображение
$$T_{x,y}\colon K'\oplus k'\to K/k$$
не сюръективно.
(Указание: воспользуйтесь уравнением для $\widetilde{V}$ из
упражнения~\ref{prob:compactification-of-C}(ii).)
\item[(iii)]
Покажите, что для $(x:y)\in \widetilde{V}(k)$ выполняется $x\neq 0$.
(Указание: если $x=0$, то из упражнения~\ref{prob:compactification-of-C}(ii)
следует, что
$$\theta\cdot\Nm_{k'/k}(y)\in k\,,$$
что невозможно, так как $\theta\not\in k$.
Заметьте, что это единственное место во всем упражнении,
где используется какая-либо информация про число $\theta$!)
\item[(iv)]
Докажите, что все $k$-точки многообразия $\widetilde{V}$ неособы.
(Указание: отображение из~$K'$ в~$K$, заданное формулой
\begin{equation}\label{eq:otobrazhenie-Tr}
z\mapsto \Tr_{K'/K}\big(z\cdot\sigma(x)\big)\,,
\end{equation}
сюръективно, так как расширение $K\subset K'$
сепарабельно, и $x\neq 0$
согласно пункту~(iii).)
\item[(v)] Можно ли аналогичным образом доказать,
что все $\bar{k}$-точки на $\widetilde{V}$ неособы?
(Указание: нельзя. Действительно,
мы можем доопределить отображение~\eqref{eq:otobrazhenie-Tr}
до отображения из~$\bar{k}^6$ в~$\bar{k}^3$,
используя ограничение скаляров по Вейлю, однако это отображение
уже ни в каком смысле не будет ``следом'' для расширения полей.
В частности, невозможно будет воспользоваться невырожденностью следа,
которая была существенна для пункта~(iv).
На самом деле, как показано в~\cite{CT-i-drugie},
многообразие~$\widetilde{V}$ имеет восемь изолированных особых точек
над полем~$\bar{k}$.)
\end{itemize}
\end{prob}

\bigskip
Из упражнения~\ref{prob:singularities-of-V}
можно получить информацию о структуре многообразия~$\widetilde{V}$,
и доказать его рациональность с помощью упражнения~\ref{prob:two-quadrics}.
На замыкании кривой~\mbox{$\Gamma\subset K^*$}
(см.~\eqref{eq:Gamma})
в проективном пространстве $\Pb(K\oplus k)$
лежит точка~\mbox{$p=\Pb(k\oplus 0)$}. Отсюда следует, что
многообразие~$\widetilde{V}$ содержит прямую
\begin{equation}\label{lucky-line}
l=\Pb(k'\oplus 0)\subset\Pb(K'\oplus k')\,.
\end{equation}

\bigskip
\begin{prob}{\textbf{Рациональность многообразия $V$}}
\label{prob:rationality-of-V}

\begin{itemize}
\item[(o)] Задайте многообразие $\widetilde{V}$ как пересечение двух
квадрик в $\Pb^7$.
(Указание: воспользуйтесь упражнением~\ref{prob:compactification-of-C}(ii).)
\item[(i)]
Покажите, что прямая $l$ не содержится в множестве особых точек
многообразия $\widetilde{V}$.
(Указание: согласно упражнению~\ref{prob:singularities-of-V}(iv),
у многообразия $\widetilde{V}$
нет особых $k$-точек, а прямая $l$ определена над $k$.)
\item[(ii)]
Покажите, что многообразие $\widetilde{V}$ не является конусом.
(Указание: вершина конуса --- проективное подпространство, определённое над $k$.
В частности, оно имеет $k$-точку, которая является особой точкой конуса.)
\item[(iii)]
Докажите, что многообразие $\widetilde{V}$ (а значит, и~$V$)
рационально.
(Указание: воспользуйтесь неприводимостью $\widetilde{V}$,
пунктами~(i), (ii) и
упражнением~\ref{prob:two-quadrics}.)
\end{itemize}
\end{prob}

\subsection{Стабильная бирациональная эквивалентность $X$ и $V$}

Наша конечная цель --- достроить морфизмы $f$ и $g$ до коммутативной
диаграммы~\eqref{diag:plan-Polietilena}, выбрав подходящим образом
рациональный тор $D$. Для начала покажем, что $g$ пропускается
через $f$. Для простоты
мы будем одинаково обозначать
элементы тора $\hat{C}$ и их образы в $C=\hat{C}/k^*$.

\bigskip
\bigskip
\bigskip
\bigskip
\bigskip
\bigskip
\begin{prob}{\textbf{Морфизм $h\colon C\to A$}}
\label{prob:h-from-C-to-A}

\begin{itemize}
\item[(i)] Постройте морфизм $h'\colon C\to k'^*$, для которого выполняется
равенство
$$\Nm_{K/k}\circ g=\Nm_{k'/k}\circ h'\,.$$
(Указание: положите
$$h'(x,y)=\Nm_{K'/k'}(x)\cdot\Nm_{k'/k}(y)^{-1}\cdot y^{-1}\,,$$
где $x\in K'$, $y\in k'^*$.)
\item[(ii)] Покажите, что морфизмы $g$ и $h'$ задают такой морфизм
$$h\colon C\to A=K^*\times_{k^*} k'^*\,,$$
что $g=f\circ h$.
\item[(iii)] Пусть $\tau$ ---
образующая группы $\Gal(K'/k')\cong\Z/3\Z$.
Определим морфизм
$$\lambda\colon A\to C$$
формулой
$$\lambda(r,s)=\big(\tau(r)^{-1},s^{-1}\big)\,.$$
Проверьте, что $\lambda$ является сечением морфизма $h$,
то есть
$$h\circ\lambda=\mathrm{id}_A\,.$$
(Указание:
пусть $\sigma$ обозначает нетривиальную инволюцию поля
$K'$ над $K$. Воспользуйтесь тем, что
$$\sigma\tau=\tau^2\sigma$$
в группе $\Gal(K'/k)$,
а также тем, что
$$\Nm_{K/k}(z)=\Nm_{K'/k'}(z)$$
для каждого $z\in K^*$.)
\end{itemize}
\end{prob}

\bigskip

Цель следующего упражнения --- проанализировать, насколько композиция
$\lambda\circ h$ отличается от тождественного морфизма $\mathrm{id}_C$.

\begin{prob}{\textbf{Тор $D$}}
\label{prob:tor-D}

\begin{itemize}
\item[(i)]
Покажите, что для $(x,y)\in C$ элемент
$$(x,y)\cdot \big(\lambda\circ h\big)(x,y)^{-1}\in C$$
зависит только от $x\in K'^*$. Это определяет морфизм торов
$K'^*\to C$. Докажите, что $\tau(K^*)\subset K'^*$ лежит в его ядре.
(Указание: воспользуйтесь явной формулой
$$(x,y)\cdot \big(\lambda\circ h\big)(x,y)^{-1}=
\Big(x\cdot\tau\big(\Nm_{K'/K}(x)\big), \Nm_{K'/k'}(x)\Big)\in C,$$
а также указанием к упражнению~\ref{prob:h-from-C-to-A}(iii).)
\item[(ii)]
Положим $D=K'^*/\tau(K^*)$. Покажите, что тор $D$ рационален.
(Указание: воспользуйтесь изоморфизмом
$$K'^*/\tau(K^*)\stackrel{\tau^{-1}}{\longrightarrow} K'^*/K^*$$
и упражнением~\ref{prob:tor-F/E-racionalen}(ii).)
\end{itemize}
\end{prob}

Из упражнения~\ref{prob:tor-D}(i) возникает морфизм
$\mu\colon D\to C$.

\bigskip
\begin{prob}{\textbf{Разложение $A\times D\cong C$}}
\label{prob:A-times-D-is-C}

\begin{itemize}
\item[(i)]
Используя морфизмы $\lambda$ и $\mu$ из упражнений~\ref{prob:h-from-C-to-A}
и~\ref{prob:tor-D}, постройте изоморфизм
$A\times D\stackrel{{}\sim{}}{\longrightarrow} C$.
(Указание: обратное отображение $C\to A\times D$ задаётся формулой
$$(x,y)\mapsto \big(h(x,y),x\big)\,,$$
где $x\in K'^*$, $y\in k'^*$.)
\item[(ii)]
Проверьте коммутативность диаграммы~\eqref{diag:plan-Polietilena}.
\item[(iii)] Проверьте, что
многообразие $V$ бирационально эквивалентно $X\times D$.
В частности, многообразия $X$ и $V$ стабильно бирационально
эквивалентны.
(Указание: воспользуйтесь пунктом~(ii) и
упражнением~\ref{prob:tor-D}(ii).)
\end{itemize}
\end{prob}

Таким образом, указанный в начале главы план полностью реализован.
А именно, в упражнении~\ref{prob:equation-of-X}(iii) мы доказали, что
поверхность $X$ нерациональна, а
в упражнении~\ref{prob:rationality-of-V}(iii)
мы проверили рациональность многообразия $V$. Наконец,
в упражнении~\ref{prob:A-times-D-is-C}(iii) мы установили
стабильную бирациональную эквивалентность~$X$ и~$V$.
Другими словами, поверхность $X$ является примером
нерационального стабильно рационального многообразия.

\subsection{Ещё одна конструкция стабильной рациональности}

Шепард-Баррон заметил, что полученный выше результат можно немного
усилить: оказывается, что рационально даже многообразие $X\times\Pb^2$.
Объясним это во введённых выше обозначениях. При этом будет
использоваться \emph{теорема Воскресенского}
(см.~\cite[Теорема~4.74]{Voskresenskii-rus}),
которая утверждает, что
всякий двумерный алгебраический тор рационален.

Как было показано в упражнении~\ref{prob:A-times-D-is-C}(i), имеется
изоморфизм алгебраических торов $C\cong A\times D$.
В то же время тор $D$
содержит одномерный подтор $T=k'^*/k^*$.
Из упражнения~\ref{prob:tor-D} следует, что
$T$-орбита точки $(x,y)\in C$ является образом отображения
$$
T\to C\,,\quad
t\mapsto
(x\cdot t\cdot\Nm_{k'/k}(t),y\cdot t^3)\,,\quad t\in k'^*\,.
$$
Поскольку это отображение имеет степень $3$ по $t$, замыкания орбит тора~$T$
являются (рациональными) кубическими кривыми в
проективном пространстве~\mbox{$\Pb^7=\Pb(K'\oplus k')$}.
Тор $T$ открыто вложен в проективную
прямую $\Pb^1=\Pb(k')$, причём дополнение является
подмногообразием~\mbox{$Z\subset\Pb(k')$}, заданным квадратичным уравнением
$\Nm_{k'/k}(t)=0$, где $t\in k'$ (таким образом, $Z$ является парой точек над
алгебраическим замыканием поля $k$). Отсюда следует, что для
любых~\mbox{$x\in K'$}
и~\mbox{$y\in k'^*$} граница $T$-орбиты точки $(x:y)\in\Pb^7$
равна
$$
Z\subset \Pb(k'\oplus 0)=l\subset\Pb^7\,,
$$
ср. с~\eqref{lucky-line}.
Для доказательства этого утверждения можно воспользоваться
тем, что если $\Nm_{k'/k}(t)=0$,
$t\neq 0$ и $y\neq 0$, то
$\Nm_{k'/k}(y\cdot t^3)=0$, и~\mbox{$y\cdot t^3\neq 0$}.
Поэтому любая гиперплоскость $H$, содержащая прямую $l$,
пересекается с такой \mbox{$T$-орбитой} ровно по одной точке.

Из коммутативности диаграммы~\eqref{diag:plan-Polietilena}
следует, что подмногообразие $V\subset C$ инвариантно относительно
действия тора $T$.
Из предыдущего обсуждения вытекает,
что сечение многообразия~$V$ общей гиперплоскостью $H$,
проходящей через прямую~$l$, бирационально отображается на фактор
$V/T$ (такое сечение называется ``slice'' в геометрической теории
инвариантов). Многообразие $V\cap H$ по-прежнему является
неприводимым пересечением двух квадрик, содержащим прямую,
которая удовлетворяет условиям упражнения~\ref{prob:two-quadrics}.
Поэтому многообразие $V\cap H$ рационально, а значит,
рационально и многообразие~$V/T$.

\bigskip

Так как $V\cong X\times D$, морфизм
$\mathrm{pr}_A\colon V\to X$ имеет рациональное сечение.
Поскольку $T\subset D$,
морфизм $V/T\to X$
также имеет рациональное сечение.
Более того, многообразие $X$ является фактором
многообразия~$V/T$ по действию двумерного тора \mbox{$S=D/T$}.
Значит, многообразие $V/T$ бирационально эквивалентно
произведению~\mbox{$S\times X$}.
Тор $S$ рационален над $k$ по теореме Воскресенского,
так как имеет размерность~$2$. Таким
образом, многообразие $X\times\Pb^2$ также рационально над $k$.

\newpage

\section{Теорема Минковского--Хассе}
\label{section:Min-Hasse}

В этой главе мы разберём доказательство теоремы Минковского--Хассе.
При этом мы будем пользоваться
некоторыми фундаментальными фактами из теории полей
классов (а именно, точностью последовательности~\eqref{eq:Br-sequence}).
Доказательства этих фактов слишком техничны для нашего изложения;
с ними можно ознакомиться, например,
в~\cite[VII]{CasselsFrolich}.

\subsection{Предварительные сведения}
\label{subsection:Min-Hasse-preliminary}

Напомним (см., например, \cite[II.1]{CasselsFrolich}),
что {\it (мультипликативным) нормированием}
поля $K$ называется гомоморфизм
$$
K^*\to\Rb_{>0},\quad x\mapsto |x|\,,
$$
для которого существует такая положительная константа $C\in\Rb_{>0}$,
что для всех элементов~\mbox{$x\in K^*$} из неравенства $|x|\leqslant 1$
следует неравенство $|1+x|\leqslant C$.
Кроме того, по определению~\mbox{$|0|=0$}.
Нормирование $|\cdot|_1$ {\it эквивалентно}
нормированию $|\cdot|_2$, если существует такая
положительная константа $\lambda\in\Rb_{>0}$,
что $|\cdot|_1=|\cdot|_2^\lambda$.
Нормирование {\it неархимедово}, если можно положить
$C=1$. Это равносильно тому, что
$$|x+y|\leqslant\max\{|x|,|y|\}$$
для всех $x, y\in K$.
В противном случае нормирование называется \emph{архимедовым}.
Легко видеть,
что отношение эквивалентности на нормированиях уважает свойство
(не)архимедовости нормирований.

Отметим, что произвольное нормирование $|\cdot |$ поля $K$
задаёт метрику на $K$ по формуле
$$\rho(x,y)=|x-y|,\quad x,y\in K\,.$$
В свою очередь, метрика $\rho$ определяет топологию на поле $K$,
относительно которой операции сложения, вычитания, умножения и
деления непрерывны.

\begin{prob}{\textbf{Примеры нормирований}}
\label{prob:normirovaniya-Q}

\begin{itemize}
\item[(i)] Докажите, что абсолютная величина
$x\mapsto |x|$ задаёт архимедово нормирование
поля $\Q$. При этом пополнение поля $\Q$ по соответствующей метрике
изоморфно полю вещественных чисел $\Rb$.
\item[(ii)] Пусть $p$ --- простое число. Для всякого ненулевого
целого числа $r$ положим~\mbox{$\ord_p(r)=d$},
где $r=p^d\cdot r'$, и $r'$ --- целое число, взаимно простое с $p$.
Пусть $x\in\Q^*$ представлено в виде $x=m/n$. Положим
$$\ord_p(x)=\ord_p(m)-\ord_p(n)\,.$$
Докажите, что функция
$x\mapsto p^{-\ord_p(x)}$
определяет неархимедово нормирование поля $\Q$.
При этом пополнение поля $\Q$ по соответствующей метрике
изоморфно полю $p$-адических чисел $\Q_p$.
\item[(iii)] Покажите, что для любого дискретного нормирования
$$v\colon K^*\to\Z$$
(см.~раздел~\ref{subsection:unramified-Brauer-generalities})
и вещественного числа $\alpha>1$
определено неархимедово нормирование
$$|x|_v=\alpha^{-v(x)},\quad x\in K^*\,.$$
(Отметим, что
само дискретное нормирование $v$ не является
мультипликативным нормированием!)
\end{itemize}
\end{prob}
\bigskip

Ввиду упражнения~\ref{prob:normirovaniya-Q}(iii) мы будем обозначать произвольное нормирование поля~$K$ символом $|\cdot|_v$, даже в случае, когда оно не возникает ни из какого дискретного нормирования (например, когда оно архимедово). Более того, обычно мы будем сокращать это обозначение просто до $v$, используя, тем не менее, обозначение~$|\cdot|_v$ для самой функции из $K^*$ в $\Rb_{>0}$. Соответственно, пополнение поля $K$ относительно~$|\cdot |_v$ (или, что то же самое, относительно $v$) будет обозначаться через~$K_v$.

\medskip

\emph{Теорема Островского}
(см., например,~\mbox{\cite[Теорема~II.3.2]{CasselsFrolich}})
утверждает, что все нормирования поля
$\Q$ с точностью до эквивалентности
исчерпываются построенными в упражнении~\ref{prob:normirovaniya-Q}.

{\it Глобальным полем} называется
конечное расширение поля $\Q$ или конечное расширение поля~$\F_p(T)$.
Существует
обобщение теоремы Островского для произвольного глобального поля (см., например,~\cite[II.11, II.12]{CasselsFrolich}).
В частности, оно утверждает, что у глобального поля есть лишь конечное число
архимедовых нормирований (при этом у глобального поля положительной
характеристики таких нормирований нет вообще).

{\it Локальным полем} называется пополнение $K_v$ глобального
поля $K$ по некоторому нормированию $v$
(ср. с определением~\ref{defin-local-field}).
Таким образом, если $v$ архимедово, то локальное поле $K_v$~--- это
либо поле вещественных чисел~$\Rb$, либо поле комплексных чисел~$\Cb$.
Если же $v$ неархимедово, то локальное поле $K_v$
либо является конечным расширением
поля $p$-адических чисел $\Q_p$, либо является полем $\Fb_q((T))$,
где $q$~--- степень простого числа.
Для произвольного конечного множества
нормирований $S$ глобального поля $K$ рассмотрим диагональное вложение
$$\mbox{$K\to \prod\limits_{v\in S}K_v\,.$}$$
\emph{Теорема о слабой аппроксимации}
(см.~\cite[Лемма~II.6.1]{CasselsFrolich})
утверждает, что образ поля~$K$
относительно этого вложения плотен (для случая $K=\Q$ это
сводится к китайской теореме об остатках).

\medskip

Пусть дана гладкая
проективная квадрика $Q$, определённая над глобальным полем~$K$, характеристика которого не равна $2$.
\emph{Теорема Минковского--Хассе} утверждает, что если для каждого
нормирования~$v$ поля~$K$ квадрика~$Q$ имеет
точку над полем~$K_v$, то~$Q$ имеет точку над~$K$.
Цель следующих упражнений --- вывести эту теорему
из некоторых утверждений теории полей классов.

\subsection{Квадрики над локальными полями}

\begin{prob}{\textbf{Вариация коэффициентов}}
\label{prob:small-variation}

Пусть $F$ --- локальное поле. На $F$
имеется естественная топология, заданная
нормированием поля $F$.

\begin{itemize}
\item[(o)] Докажите, что множество квадратов в $F$ открыто.
(Указание: для неархимедова случая
воспользуйтесь леммой Гензеля или разложением
функции~\mbox{$\sqrt{1+t}$} в ряд Тейлора).
\item[(i)] Пусть гладкая проективная квадрика $Q\subset\P^n$ над $F$ задана уравнением
$$\sum\limits_{i=0}^n a_ix_i^2=0\,,$$
где $a_i\in F^*$. Докажите,
что существует число $\varepsilon>0$, обладающее следующим свойством:
если для всех $0\leqslant i\leqslant n$ выполнены неравенства \mbox{$|b_i-a_i|<\varepsilon$}, то квадрика~$Q'$,
заданная уравнением
$$\sum\limits_{i=0}^n b_ix_i^2=0\,,$$
изоморфна~$Q$ над~$F$.
(Указание: воспользовавшись пунктом~(o),
выберите такое $\varepsilon>0$, что
$b_i/a_i$
является квадратом в $F$ для всех $0\leqslant i\leqslant n$
при $|b_i-a_i|<\varepsilon$.)
\item[(ii)] Пусть $Q$ --- гладкая проективная
квадрика, определённая над полем $F$.
Докажите, что при малом изменении коэффициентов
уравнения квадрики $Q$ она остаётся изоморфной самой себе над $F$.
(Указание: докажите, что матрицу, задающую невырожденную билинейную форму, можно привести
к диагональному виду так, что коэффициенты соответствующей
диагональной матрицы непрерывно зависят от коэффициентов исходной.
После этого воспользуйтесь пунктом~(i).)
\end{itemize}
\end{prob}

\bigskip
\begin{prob}{\textbf{Точки на квадриках с гладкой редукцией}}
\label{prob:quadric-over-local-field}

Пусть $Q\subset\P^n_F$ ---
квадрика положительной размерности над неархимедовым
локальным полем $F$, характеристика которого не равна $2$, заданная многочленом с
коэффициентами в кольце
нормирования~\mbox{$\OO_F\subset F$},
причём определитель соответствующей
матрицы является
обратимым элементом в~$\OO_F$.
Цель этого упражнения --- доказать, что $Q(F)\ne\varnothing$.

\begin{itemize}
\item[(o)] Рассмотрим редукцию уравнения квадрики $Q$ по модулю
максимального идеала $\mathfrak{m}_F\subset\OO_F$. Это задаёт квадрику
$\overline{Q}$ над конечным полем вычетов $k\cong\OO_F/\mathfrak{m}_F$.
Покажите, что квадрика $\overline{Q}\subset\P^n_{k}$
является гладкой.
(Указание: какой определитель у квадратичной формы над полем~$k$,
задающей $\overline{Q}$?)
\item[(i)] Докажите, что гладкая проективная квадрика положительной размерности,
определённая над конечным полем, имеет над ним точку.
(Указание: воспользуйтесь теоремой Шевалле--Варнинга,
доказанной в упражнении~\ref{exer-Chevalley},
или докажите это непосредственно.)
\item[(ii)] Выберем точку
$$\bar{p}\in\overline{Q}(k)\subset\P^n_{k}\,.$$
Проведём через $\bar{p}$ прямую $\bar{l}$,
не лежащую во вложенном проективном  касательном
пространстве
$$\mathbb{T}_{\bar{p}}\big(\overline{Q}\big)\subset\P^n_{k}\,.$$
Поднимем уравнение
прямой $\bar{l}$ произвольным образом до уравнения прямой~$l$ над $\OO_F$.
Докажите, что квадрика $Q$ имеет точку, лежащую на прямой $l$ и определённую
над $F$. (Указание: воспользуйтесь леммой Гензеля для
квадратичного уравнения, задающего пересечение $l$ и $Q$.)
\end{itemize}
\end{prob}

\bigskip
\begin{prob}{\textbf{Точки квадрики над пополнениями глобального поля}}
\label{prob:almostall}

Пусть $Q$ --- гладкая проективная квадрика положительной размерности
над глобальным полем~$K$, характеристика которого не равна $2$.
Докажите, что квадрика $Q$ имеет точку над локальным полем~$K_v$
для почти всех нормирований~$v$ поля $K$. (Указание:
выберите однородный многочлен, задающий квадрику $Q$ в соответствующем
проективном пространстве. Заметьте, что для почти всех неархимедовых
нормирований~$v$
это дает многочлен с коэффициентами в кольце нормирования
поля~$K_v$. После этого
воспользуйтесь упражнением~\ref{prob:quadric-over-local-field}.)
\end{prob}

\subsection{Редукция к случаю $\dim(Q)=1$}

Пусть $Q\subset\Pb^n$ --- гладкая проективная квадрика над глобальным полем $K$,
и~$\dim(Q)\geqslant 3$, то есть~$n\geqslant 4$. Рассмотрим общее
проективное подпространство $\Lambda\subset\Pb^n$ коразмерности~$2$ и
проекцию
$$\pi_{\Lambda}\colon\Pb^n\dasharrow \Pb^{1}$$
с центром в $\Lambda$.
Для произвольной точки $x\in\Pb^1$ через~$Q_x$ будем обозначать
замыкание прообраза~$\pi_{\Lambda}^{-1}(x)$ в $Q$.
Отметим, что квадрика $Q\cap\Lambda$
гладкая из-за общности выбора подпространства
$\Lambda$. Ясно, что $Q\cap\Lambda$ содержится в $Q_x$ для каждой точки $x\in\Pb^1$.

\begin{prob}{\textbf{Особые слои отображения $\pi_{\Lambda}$}}
\label{prob:malo-sloev}

Проверьте, что для
всех точек~$x\in\P^1$, кроме двух, многообразие $Q_x$ является
гладкой квадрикой в~$\Pb^{n-1}$. (Указание:
гиперплоскости в $\P^n$, проходящие через~$\Lambda$, параметризуются
прямой в двойственном пространстве~$(\P^n)^{\vee}$, а касательные
гиперплоскости к~$Q$ параметризуются квадрикой в
двойственном пространстве~$(\P^n)^{\vee}$.)
\end{prob}

\bigskip
По условию теоремы Минковского--Хассе
для каждого нормирования $v$ поля $K$
существует точка $p_v\in Q(K_v)$. Проекция из неё определяет
бирациональную эквивалентность $Q\dasharrow \Pb^{n-1}$
над полем~$K_v$, и поэтому множество $K_v$-точек плотно по Зарисскому в~$Q$.
Таким образом, по упражнению~\ref{prob:malo-sloev}
можно предполагать, что
для любого $v\in S$
точка $p_v$ не лежит на особых слоях отображения
$\pi_{\Lambda}\colon Q\dasharrow\P^1$.

Пусть множество $S$ состоит из всех нормирований $v$ поля $K$,
для которых квадрика $Q\cap\Lambda$ не имеет $K_v$-точек.
Так как размерность $Q\cap\Lambda$ положительна,
по упражнению~\ref{prob:almostall},
применённому к квадрике~\mbox{$Q\cap\Lambda$}, множество $S$ конечно.

\begin{prob}{\textbf{Хорошие слои отображения $\pi_{\Lambda}$}}
\label{prob:zamechatelnyj-sloj}

Для каждого нормирования $v$ поля $K$
рассмотрим точку \mbox{$x_v\in \Pb^1(K_v)$}, равную~$\pi_{\Lambda}(p_v)$.
Докажите, что
существует такая точка $x\in \Pb^1(K)$, что для любого нормирования $v\in S$
квадрика $Q_x$ изоморфна над полем $K_v$ квадрике~$Q_{x_v}$.
(Указание: воспользуйтесь
теоремой о слабой аппроксимации применительно к множеству нормирований~$S$
и упражнением~\ref{prob:small-variation}(ii).)
\end{prob}

\bigskip
\begin{prob}{\textbf{Редукция к случаю $\dim(Q)=2$}}

Покажите, что для доказательства теоремы Минковского--Хассе
в случае~\mbox{$\dim(Q)\geqslant 3$}
достаточно доказать её для случая $\dim(Q)=2$.
(Указание: воспользовавшись введёнными выше обозначениями, заметьте, что
по упражнению~\ref{prob:zamechatelnyj-sloj} слой $Q_x$
имеет точку над полем $K_v$ для любого $v\in S$.  По
построению для любого нормирования
$v\notin S$ квадрика $Q\cap \Lambda$ имеет точку
над полем $K_v$, а значит, квадрика
$$Q_x\supset Q\cap \Lambda$$
также имеет точку над полем $K_v$. Таким образом, $Q_x$ имеет точку над
$K_v$ для любого~$v$. Поскольку
$$\dim(Q_x)=\dim(Q)-1\geqslant 2\,,$$
то индукция
по $\dim(Q)$ показывает, что квадрика $Q_x$ имеет точку над полем
$K$. Cледовательно, квадрика $Q\supset Q_x$ тоже имеет точку над
полем $K$.)
\end{prob}

\bigskip
\begin{prob}{\textbf{Редукция к случаю $\dim(Q)=1$}}
\label{prob:reduction-to-dimQ-1}

Пусть $Q\subset\Pb^3$ --- гладкая двумерная квадрика над глобальным
полем $K$. Рассмотрим поле
$$L=K\big(\sqrt{\disc(Q)}\big)\,,$$
где $\disc(Q)$ --- дискриминант квадрики $Q$
(см. определение~\ref{defin:discrquadr}).
Пусть коника $C$ над полем $L$ --- гладкое плоское сечение
квадрики $Q_L$.
Докажите, что теорема Минковского--Хассе для квадрики $Q$
над полем $K$ сводится к теореме
Минковского--Хассе для коники $C$ над полем $L$.
(Указание:
для любого нормирования $w$ поля $L$
пополнение $L_w$ является расширением пополнения $K_v$, где~$v$
определяется как ограничение $w$ на $K$, причём
$$L_w\cong K_v\big(\sqrt{\disc(Q)}\big)\,.$$
Предположим, что для квадрики $Q$ выполнены условия теоремы
Минковского--Хассе, то есть $Q$ имеет точки над всеми
полями $K_v$.
По упражнению~\ref{prob-nonsection}(ii)
коника $C$ имеет точки над всеми полями $L_w$.
Если для коники $C$ выполнена теорема Минковского--Хассе,
то $C$ имеет точку над полем $L$. По упражнению~\ref{prob-nonsection}(ii)
отсюда следует, что квадрика $Q$ имеет точку над полем $K$.)
\end{prob}

\subsection{Случай $\dim(Q)\leqslant 1$}

Пусть $K$ --- произвольное глобальное поле, и дан элемент $\alpha\in\Br(K)$.
Можно показать, что для всех нормирований $v$, кроме конечного числа,
образ $\alpha$ относительно естественного отображения
$$\Br(K)\to\Br(K_v)$$
равен нулю (ср. с упражнением~\ref{prob:pochti-vse}).
Имеет место {\it точная последовательность теории полей классов}
(см.~\cite[VII.9.6, VII.11.2(bis)]{CasselsFrolich}):
\begin{equation}\label{eq:Br-sequence}
0\to \Br(K)\to\oplus_{v}\Br(K_v)\to \Q/\Z\to 0\,.
\end{equation}
Напомним, что для неархимедовых нормирований
имеется изоморфизм
$$\res\colon \Br(K_v)\stackrel{\sim}\longrightarrow\Q/\Z\,,$$
задаваемый отображением вычета
(см. упражнение~\ref{exer-Brlocfield}(i)),
в то время как для
архимедовых нормирований имеем
$$\Br(\Rb)\cong \frac{1}{2}\Z/\Z$$
(см. упражнение~\ref{prop-realBr}) и $\Br(\Cb)=0$.
При этом правое
отображение в последовательности~\eqref{eq:Br-sequence} соответствует
суммированию в группе $\Q/\Z$.

Предположим, что характеристика поля $K$ не равна $2$. Благодаря соответствию между гладкими кониками и
кватернионными алгебрами (см. упражнение~\ref{prob-quatconic}(i)),
точность последовательности~\eqref{eq:Br-sequence}
в первом члене влечёт теорему Минковского--Хассе для
одномерного случая.
Тот факт, что композиция двух отображений
нулевая для подгрупп \mbox{$2$-кручения}, в случае $K=\Q$
равносилен квадратичному закону
взаимности (см.~\cite[4.5.5]{PanchishkinManin}).
Точность в среднем члене является важным результатом
теории полей классов и даёт явное описание группы Брауэра глобальных
полей.

Случай $\dim(Q)=0$ легко разбирается явным образом для поля $\Q$,
а для произвольного глобального поля $K$ является нетривиальным
следствием из теории полей классов: элемент из $K^*$ является квадратом
тогда и только тогда, когда он является квадратом в $K_v$ для каждого
нормирования $v$ поля $K$.
Например, это утверждение можно вывести из теоремы Чеботарёва о плотности
(см.~\cite[VIII]{CasselsFrolich}).

\medskip

Для поля $K=\Q$ можно дать непосредственное доказательство теоремы
Мин\-ков\-ского--Хассе в одномерном случае. Сформулируем её в
терминах элементов группы Брауэра. Для краткости через $(a,b)$ будем
обозначать класс в $\Br(K)$ кватернионной алгебры $A(a,b)$, где
$a,b\in\Q^*$ (см. упражнение~\ref{exer-quaternion-alg}).

Предположим, что для
любого нормирования $v$ поля $\Q$ образ класса $(a,b)$ относительно
естественного отображения
$$\Br(\Q)\to \Br(\Q_v)$$
тривиален. Мы хотим
доказать, что $(a,b)=0$. Умножая на квадраты и используя
билинейность класса $(a,b)$ по~$a$ и $b$
(см.~упражнение~\ref{exer-quaternion-alg}(iii),(iv)),
можно считать, что $a$ и $b$ ---
ненулевые целые числа, свободные от квадратов, и $|a|\leqslant |b|$. Будем
вести индукцию по $M=|a|+|b|$.

\begin{prob}{\textbf{Индукция по $M$}}
\label{prob:ind-in-M}\hspace{0cm}

\begin{itemize}
\item[(o)] Разберите случай $M=2$ явным образом.
(Указание: если $a=1$ или $b=1$, то можно воспользоваться
упражнением~\ref{exer-quaternion-alg}(iii).
Чтобы показать, что случай~\mbox{$a=b=-1$} невозможен,
примените упражнение~\ref{exer-quaternion-alg}(ii)
над полем $\Rb$.)
\item[(i)]
Если $M>2$, то $|b|>1$. Рассмотрим произвольный простой
делитель $p$ числа $b$. Докажите, что $a$ является квадратом по
модулю $p$. (Указание: воспользуйтесь тем, что
по упражнению~\ref{exer-quaternion-alg}(ii)
коника, заданная уравнением
$$x^2-ay^2-bz^2=0\,,$$
имеет точку над полем~$\Q_p$.)
Следовательно, $a$ является квадратом по модулю~$b$, и имеется равенство
\begin{equation}\label{eq:acbb}
a=c^2-bb'\,,
\end{equation}
где $c$ и $b'$ --- целые числа.
Более того, изменив, если потребуется, число $b'$,
можно заменить в равенстве~\eqref{eq:acbb}
число $c$ на целое число из интервала от~\mbox{$-b/2$} до~\mbox{$b/2$}.
Таким образом, можно считать, что
$$|c|\leqslant \frac{|b|}{2}\,.$$
Выведите отсюда, что $|b'|<|b|$.
\item[(ii)] Докажите, что если число $b'$ выбрано так, как описано
в пункте~(i), то имеется равенство
$(a,b)=(a,b')$.
(Указание:
покажите, что
$$(a,b)+(a,b')=0\,,$$
пользуясь
билинейностью класса $(-,-)$, а также соотношением
Стейн\-берга~\mbox{$(d,1-d)=0$}
для любого $d\ne 0,1$,
см. упражнение~\ref{exer-quaternion-alg}(iii).
После этого воспользуйтесь равенством $2(a,b)=0$.)
Это завершает шаг индукции.
\end{itemize}
\end{prob}

\bigskip

Отметим, что приведённые выше рассуждения не дают полного
доказательства теоремы Минковского--Хассе даже для случая $K=\Q$,
несмотря на то, что в упражнении~\ref{prob:ind-in-M}
мы провели рассуждение для коник над $\Q$.
Это объясняется тем, что при редукции к коникам
приходится переходить к расширениям поля $\Q$
(см. упражнение~\ref{prob:reduction-to-dimQ-1}).
Существует ``самозамкнутое'' доказательство теоремы
Минковского--Хассе над $\Q$, не использующее другие числовые
поля (см.~\mbox{\cite[IV.3.2]{Serr-kurs-arifmetiki}}).
Тем не менее, оно использует теорему Дирихле о простых числах в
арифметической прогрессии, что является тонким результатом, также
имеющим непосредственное отношение к теории полей классов.

\subsection{Другие примеры выполнения принципа Хассе}

Пусть $\mathcal{M}$ --- некоторое множество
многообразий над фиксированным глобальным полем $K$
(например, $\mathcal{M}$ состоит из всех
проективных гиперповерхностей данной степени над $K$).
Говорят, что для $\mathcal{M}$ \emph{выполняется
принцип Хассе}, если для любого многообразия
$X$ из $\mathcal{M}$ верно следующее:
многообразие $X$ имеет точку над полем~$K$ тогда и только тогда,
когда оно имеет точку над полем $K_v$
для всех нормирований~$v$ поля $K$. Например, теорема Минковского--Хассе
утверждает, что принцип Хассе выполнен для
гладких квадрик (для особых квадрик он тоже выполнен~---
по очевидным причинам).

\begin{prob}{\bf Принцип Хассе для многообразий Севери--Брауэра}

Покажите, что для многообразий Севери--Брауэра
(см. определение~\ref{defin:Severi-Brauer})
выполняется принцип Хассе.
(Указание: воспользуйтесь
упражнением~\ref{prob:SB-s-tochkoi},
а также точной последовательностью теории полей
классов~\eqref{eq:Br-sequence}.)
\end{prob}

\bigskip
Из точной последовательности теории полей
классов~\eqref{eq:Br-sequence}
следует, что для конечного расширения Галуа $K\subset L$ имеет место точная последовательность
\begin{equation}\label{eq:Br-sequence-extension}
0\to \Br(L/K)\to \oplus_{v}\Br(L_w/K_v)\to \Q/\Z\,,
\end{equation}
где для каждого нормирования $v$ поля $K$ выбирается
произвольное нормирование~$w$ поля $L$,
продолжающее $v$ с $K$ на $L$ (на самом деле, эти два утверждения равносильны).

\bigskip

\begin{prob}{\bf Теорема Хассе о нормах}

Рассмотрим конечное циклическое расширение Галуа глобальных полей
$K\subset L$ (то есть расширение Галуа с конечной циклической группой Галуа).
\begin{itemize}
\item[(i)]
Докажите \emph{теорему Хассе о нормах}:
элемент $a\in K^*$ является нормой некоторого элемента $b\in L^*$
тогда и только тогда, когда для каждого нормирования~$v$ поля $K$
и каждого нормирования $w$ поля $L$, продолжающего~$v$,
существует такой элемент $b_w\in L_w^*$, что
$$\Nm_{L_w/K_v}(b_w)=a\,.$$
(Указание: воспользуйтесь точной
последовательностью~\eqref{eq:Br-sequence-extension}
и явным видом групп когомологий для циклической группы, полученным
в упражнении~\ref{exer-cyclgrp}(i).)
\item[(ii)]
Докажите, что принцип Хассе выполняется для аффинных многообразий в~$R_{L/K}(\Gb_m)$ над~$K$, заданных уравнением вида
$$
\Nm_{L/K}(x)=a,\quad x\in R_{L/K}(\Gb_m)\,,
$$
где $a$ --- некоторый элемент из $K^*$. (Указание: воспользуйтесь пунктом~(i).)
\end{itemize}
\end{prob}

\bigskip
Существует множество других примеров выполнения принципа
Хассе. Например, Бёрч в~\cite{Birch} доказал следующее утверждение:
принцип Хассе выполняется для гладких полных
пересечений $r$ гиперповерхностей в $\Pb^n$ над $\Q$,
каждая из которых имеет степень $d$, если выполняется неравенство
$$
n\geqslant r(r+1)(d-1)2^{d-1}\,.
$$
В частности, принцип Хассе выполняется для гладких пересечений
двух квадрик размерности не меньше $10$,
а также для гладких кубических гиперповерхностей размерности не меньше $15$.

С другой стороны, во многих интересных случаях принцип Хассе не выполняется.
Например, В.\,А.\,Исковских доказал в~\cite{Iskovskikh-Hasse},
что он не выполняется для пересечений двух квадрик в~$\Pb^4$.
Другие примеры, для которых можно построить препятствие к выполнению принципа Хассе,
мы обсудим в следующей главе.

\newpage
\section{Препятствие Брауэра--Манина}
\label{section:Br-Manin}

В главе~\ref{section:Min-Hasse} мы обсуждали принцип Хассе
для некоторых классов многообразий и отметили, что
в общем случае он не выполняется.
В этой главе мы определим препятствие
Брауэра--Манина к выполнению принципа Хассе,
и при помощи него построим контрпример к соответствующему
утверждению для кривых рода $1$.

\subsection{Определение препятствия Брауэра--Манина}
\label{subsection:Br-Manin-def}

В дальнейшем мы будем использовать
группы Брауэра для многообразий и схем
(см.~определение~\ref{defin:Br-scheme}).
Для простоты будем предполагать, что глобальное
поле~$K$ имеет нулевую характеристику, то есть
является конечным расширением поля~$\Q$.
Пусть~$X$~--- гладкое проективное
многообразие над глобальным полем $K$.
Тогда
группа Брауэра $\Br(X)$ многообразия $X$ совпадает с его неразветвлённой
группой Брауэра (см. определения~\ref{defin:Br-scheme}
и~\ref{defin:unramified-Br-X}, а также обсуждение
в разделе~\ref{subsection:Br-of-a-variety}). Как и раньше, для нормирования $v$ поля $K$
через $K_v$ будем обозначать соответствующее пополнение.
Для элемента $\alpha\in \Br(X)$, нормирования~$v$ поля $K$
и точки $p_v\in X(K_v)$ определим \emph{локальное спаривание}
$$
\langle\alpha,p_v\rangle_v=\res\big(p_v^*(\alpha)\big)\in\Q/\Z\,,
$$
где гомоморфизм групп
$$p_v^*\colon\Br(X)\to \Br(K_v)$$
строится по морфизму $p_v\colon\Spec(K_v)\to X$
так, как описано в разделе~\ref{subsection:Br-of-a-variety}.
Отображение $\res$ для неархимедова $v$ определяется как в
упражнении~\ref{exer-Brlocfield}(i).
Если $K_v\cong\mathbb{R}$, то отображение $\res$ является естественным
вложением
$$\Br(K_v)\cong\Z/2\Z\hookrightarrow\Q/\Z$$
(см. упражнение~\ref{prop-realBr}(i)).

Предположим, что для любого $v$ многообразие $X$ имеет \mbox{$K_v$-точку}.
Рассмотрим \emph{глобальное спаривание}
$$
\mbox{$\langle -, -\rangle\colon\Br(X)\times \prod\limits_v X(K_v)\to
\Q/\Z,
\quad
\langle \alpha,(p_v)\rangle=\sum\limits_v\langle\alpha,p_v\rangle_v\,.$}
$$
Чтобы убедиться в корректности этого определения,
надо объяснить, почему сумма в правой части конечна
(ср. с упражнением~\ref{prob:pochti-vse}).
Пусть $\OO_K$ обозначает кольцо целых
в поле~$K$.
Поле $K$ есть поле частных кольца $\OO_K$,
и $\OO_K$ является одномерным нётеровым целозамкнутым кольцом,
то есть $\Spec(\OO_K)$ --- регулярная одномерная схема.
Каждый ненулевой простой идеал~$\p$ в~$\OO_K$ задаёт
дискретное нормирование поля~$K$.
Легко показать, что это соответствие является биекцией между множеством
ненулевых простых идеалов $\p$ в $\OO_K$
и множеством классов эквивалентности
неархимедовых нормирований $v$ поля~$K$.
При этом кольцо нормирования~$\OO_v$
полного поля~$K_v$ является пополнением
по максимальному идеалу локального кольца~$(\OO_K)_{\p}$.

Далее, выберем проективную модель
$$\pi\colon\mathcal{X}\to\Spec(\OO_K)$$
над $\OO_K$ многообразия $X$ (см. раздел~\ref{subsection:models}).
Для любого элемента $\alpha\in\Br(X)$ существует такое непустое открытое
подмножество
$$U\subset \Spec(\OO_K)\,,$$
что $\alpha$ является ограничением на общий слой некоторого
элемента $\alpha_U\in\Br(\Xc_U)$,
где
$$\Xc_U=\pi^{-1}(U)\,.$$
Это непосредственно следует из определения группы Брауэра схемы
в терминах
алгебр Адзумаи, см. определение~\ref{defin:Br-scheme}.
Последнее утверждение можно также доказать,
используя описание группы Брауэра схемы
в терминах этальных когомологий, данное
в разделе~\ref{subsection:etale-Brauer}.
Поскольку морфизм~$\pi$ проективен,
каждая точка
$$p_v\colon\Spec(K_v)\to X$$
продолжается до
морфизма
$$\tilde{p}_v\colon\Spec(\OO_v)\to \Xc$$
по валюативному критерию собственности
(см., например,~\cite[Теорема~II.4.7]{Hartshorne}).
Поэтому для всех неархимедовых нормирований~$v$,
соответствующих точкам из~$U$,
элемент~\mbox{$p_v^*(\alpha)$} равен образу элемента $\tilde p_v^*(\alpha_U)$
относительно отображения
$$\Br(\OO_v)\to \Br(K_v)\,.$$
Значит, для таких~$v$ имеем $\langle\alpha,p_v\rangle_v=0$
(на самом деле, из упражнения~\ref{exer-Brlocfield}(i)
следует, что~$\Br(\OO_v)=0$).

\bigskip
\begin{prob}{\bf Препятствие Брауэра--Манина}\label{prob:BrManin}
\hspace{0cm}
\begin{itemize}
\item[(i)]
Покажите, что образ естественного отображения
$$\Br(K)\to\Br(X)$$
лежит в ядре глобального
спаривания $\langle -, -\rangle$.
(Указание: воспользуйтесь точной последовательностью
теории полей классов~\eqref{eq:Br-sequence}.)
Таким образом, глобальное спаривание пропускается через факторгруппу
$$
\widetilde{\Br}(X)={\rm Coker}\big(\Br(K)\to \Br(X)\big)\,.
$$
\item[(ii)]
Докажите, что образ диагонального вложения
$$\mbox{$X(K)\to\prod\limits_v X(K_v)$}$$
лежит в ядре спаривания $\langle -, -\rangle$. (Указание:
снова воспользуйтесь точной последовательностью теории полей
классов~\eqref{eq:Br-sequence}.)
\item[(iii)]
Пусть существует такой элемент $\alpha\in\widetilde{\Br}(X)$,
что для любого набора локальных точек
$$\mbox{$(p_v)\in \prod\limits_v X(K_v)$}$$
глобальное спаривание $\langle\alpha, (p_v)\rangle$ не равно нулю.
Покажите, что тогда $X$ не имеет $K$-точек.
\end{itemize}
\end{prob}

\bigskip
Иными словами, имеются вложения множеств
$$
\mbox{$X(K)\subset \Br(X)^{\bot}\subset \prod\limits_v X(K_v)\,,$}
$$
где $\Br(X)^{\bot}$ обозначает аннулятор
группы $\Br(X)$ в $\prod_v X(K_v)$
относительно глобального спаривания $\langle\cdot,\cdot\rangle$.
Разумеется,
если для какого-нибудь нормирования $v$ на $X$ нет $K_v$-точек или,
что равносильно, если $\prod_v X(K_v)=\varnothing$,
то на $X$ нет и $K$-точек.
Препятствие Брауэра--Манина даёт более сильное достаточное
условие пустоты множества~$X(K)$.
А именно, предположим, что
$$\mbox{$\prod\limits_v X(K_v)\ne\varnothing\,,$}$$
но множество
$\Br(X)^{\bot}$ пусто. Тогда на $X$ всё так же нет $K$-точек.
Пустота множества~\mbox{$\Br(X)^{\bot}$} и называется
\emph{препятствием Брауэра--Манина} к существованию
\mbox{$K$-точек} на~$X$.

\subsection{Вычисление препятствия Брауэра--Манина}

Наша следующая цель --- описание группы $\widetilde{\Br}(X)$,
определённой в упражнении~\ref{prob:BrManin}, и способа построения
её нетривиальных элементов в некоторых интересных случаях.
Вообще, во многих аналогичных случаях
мы будем использовать следующее обозначение:
для группы вида $\Br(-)$ с естественным отображением
$$b\colon\Br(K)\to\Br(-)$$
будем обозначать факторгруппу
$\Br(-)/b\big(\Br(K)\big)$ через
$\widetilde{\Br}(-)$.
В дальнейшем мы будем считать, что многообразие $X$ геометрически
неприводимо (напомним, что это означает неприводимость многообразия~$X_{\bar{K}}$).
Для расширения полей $K\subset F$ будем обозначать поле рациональных функций
на $X_F$ через $F(X)$. Нам понадобится нетривиальный факт из теории
полей классов: имеется равенство
$$H^3\big(G_K, \bar{K}^*\big)=0$$
для любого глобального поля~$K$
(см.~\cite[VII.11.4]{CasselsFrolich}).

\bigskip
\begin{prob}{\bf Группа $\widetilde{\Br}(X)$}\label{prob:Br0X}
\hspace{0cm}
\begin{itemize}
\item[(i)]
Постройте канонический изоморфизм
$$\widetilde{\Br}\big(\bar{K}(X)/K(X)\big)\cong
H^2\big(G_K, \bar{K}(X)^*/\bar{K}^*\big)\,.$$
(Указание: рассмотрите точную последовательность $G_K$-модулей
$$1\to\bar{K}^*\to\bar{K}(X)^*\to\bar{K}(X)^*/\bar{K}^*\to 1\,,$$
ср. с упражнением~\ref{propb:Picarinvar}(i).
Воспользуйтесь тем, что $H^3(G_K, \bar{K}^*)=0$.)
\item[(ii)]
Покажите, что имеется точная последовательность
\begin{equation}\label{H1-Br-Hom}
0\to H^1\big(G_K,\Pic(X_{\bar{K}})\big)\to
\widetilde{\Br}\big(\bar{K}(X)/K(X)\big)
\stackrel{\oplus\res_D}\longrightarrow\mbox{$\bigoplus\limits_{D\subset X}
\Hom\big(G_{K(D)},\Q/\Z\big)$}\,,
\end{equation}
где $D$ пробегает все неприводимые приведённые дивизоры на $X$.
(Указание: рассмотрите точную последовательность
$G_K$-модулей
$$1\to \bar{K}(X)^*/\bar{K}^*\to \Div(X_{\bar{K}})\to\Pic(X_{\bar{K}})
\to 0\,,$$
ср. с упражнением~\ref{propb:Picarinvar}(ii).
Затем воспользуйтесь тем,
что вложение полей
$$K_D\subset K(D)$$
индуцирует сюръективное
отображение групп Галуа
$$G_{K(D)}\to G_{K_D}$$
и вложение
$$\Hom\big(G_{K_D},\Q/\Z\big)\hookrightarrow
\Hom\big(G_{K(D)},\Q/\Z\big)\,.$$
После этого примените пункт~(i)
и упражнение~\ref{divisorscohom}(iii).)
\item[(iii)]
Покажите, что имеется точная последовательность
\begin{equation}\label{H1-Br-Brbar}
0\to H^1\big(G_K,\Pic(X_{\bar{K}})\big)\to\widetilde{\Br}(X)\to
\Br(X_{\bar{K}})\,.
\end{equation}
(Указание: из определения~\ref{defin:unramified-Br-X}
и обсуждения после него следует,
что имеется точная последовательность
\begin{equation}\label{byvshee-Brnr}
0\to \widetilde{\Br}(X)\to \widetilde{\Br}\big(K(X)\big)
\stackrel{\oplus\res_D}\longrightarrow\mbox{$\bigoplus\limits_{D\subset X}
\Hom\big(G_{K(D)},\Q/\Z\big)$}\,.
\end{equation}
Из пункта~(ii) следует, что имеется вложение
\begin{equation}\label{H1-Br-Br}
H^1\big(G_K,\Pic(X_{\bar{K}})\big)\hookrightarrow
\widetilde{\Br}\big(\bar{K}(X)/K(X)\big)\,.
\end{equation}
Из точных последовательностей~\eqref{H1-Br-Hom}
и~\eqref{byvshee-Brnr} следует, что образ группы~${H^1\big(G_K,\Pic(X_{\bar{K}})\big)}$ при композиции вложения~\eqref{H1-Br-Br} с естественным вложением $\widetilde{\Br}\big(\bar{K}(X)/K(X)\big)\hookrightarrow
\widetilde{\Br}\big(K(X)\big)$ лежит в подгруппе
$$\widetilde{\Br}(X)\subset\widetilde{\Br}\big(K(X)\big)\,.$$
Кроме того, существует коммутативная диаграмма
$$
\xymatrix{
\widetilde{\Br}\big(X\big)\ar@{^(->}[r]\ar@{->}[d]  &
\widetilde{\Br}\big(K(X)\big) \ar@{->}[d]\\
\Br\big(X_{\bar{K}}\big) \ar@{^(->}[r] &
\Br\big(\bar{K}(X)\big)
}
$$
Из неё следует точность последовательности~\eqref{H1-Br-Brbar}
во втором члене.)
\item[(iv)]
Предположим, что $\Br(X_{\bar{K}})=0$. (Например, $X_{\bar{K}}$
является гладкой проективной кривой
или гладким проективным рациональным многообразием над $\bar{K}$.
Другой пример такой ситуации возникает, если \mbox{$H^2(X,\OO_X)=0$}
и для расширения скаляров $X_{\Cb}$ многообразия $X$
относительно какого-либо (или, что равносильно, для любого) вложения
$K\hookrightarrow\Cb$ выполнено равенство
$$H^3\big(X(\Cb),\Z\big)_{tors}=0\,,$$
см. раздел~\ref{subsection:etale-complex}.)
Покажите, что в этом случае имеется изоморфизм
$$H^1\big(G_K,\Pic(X_{\bar{K}})\big)\cong\widetilde{\Br}(X)\,.$$
(Указание: это моментально следует из пункта~(iii).)
\item[(v)]
Предположим, что $\Br(X_{\bar{K}})=0$, и
группа $\Pic(X_{\bar{K}})$ не имеет кручения (в частности,
группа алгебраически
тривиальных дивизоров $\Pic^0(X_{\bar{K}})$ равна нулю).
Покажите, что тогда группа~\mbox{$\widetilde{\Br}(X)$} конечна.
(Указание: воспользуйтесь
тем, что в этом случае группа Пикара $\Pic(X_{\bar{K}})$
изоморфна группе Нерона--Севери многообразия $X_{\bar K}$
и, следовательно, конечно порождена.
Таким образом,~$\Pic(X_{\bar{K}})$
является свободной конечно порождённой абелевой группой.
Теперь воспользуйтесь
упражнением~\ref{prob:H1-finite}(iii).)
\item[(vi)]
Предположим, что $\Br(X_{\bar{K}})=0$ и $\Pic(X_{\bar{K}})\cong \Z$.
(Например, $X$ является гладкой квадрикой размерности не меньше трёх или многообразием Севери--Брауэра.)
Покажите, что в этом случае
$$\widetilde{\Br}(X)=0\,.$$
(Указание: воспользуйтесь
пунктом~(iv) и упражнением~\ref{exer-profcohom}(iii).)
\end{itemize}
\end{prob}

Отметим, что точную последовательность из упражнения~\ref{prob:Br0X}(iii)
можно также вывести из
описания группы Брауэра схемы в терминах этальных когомологий
(см.~раздел~\ref{subsection:etale-Brauer}) и
спектральной последовательности Хохшильда--Серра
для этальных когомологий

$$H^i\big(G_K,H^j_{\acute{e}t}(X_{\bar{K}},\mathbb{G}_m)\big)\implies
H_{\acute{e}t}^{i+j}(X, \mathbb{G}_m)\,.$$

\bigskip
\begin{prob}{\bf Элемент в группе Брауэра, построенный по дивизору}
\label{prob-Brdiv}

Пусть $K\subset L$ является конечным расширением Галуа с
циклической группой Галуа $G$ порядка $n$. Зафиксируем образующую $s$
группы $G$.

\begin{itemize}
\item[(o)]
Пусть даны дивизор $D\in\Div(X_L)$ и рациональная функция
$f\in K(X)^*$, для которых
$$\sum\limits_{g\in G} g(D)=(f)\in\Div(X)\,,$$
где $(f)$ обозначает дивизор функции $f$.
Докажите, что $D$ определяет класс
$$[D]\in H^1\big(G, \Pic(X_L)\big)\,.$$
(Указание: воспользуйтесь упражнением~\ref{exer-cyclgrp}(i).)
\item[(i)]
Проверьте, что $K(X)\subset L(X)$ является расширением Галуа,
группа Галуа которого канонически изоморфна группе $G$.
Пусть $A(f)$ обозначает циклическую алгебру над полем $K(X)$,
соответствующую элементу $f\in K(X)^*$,
расширению~\mbox{$K(X)\subset L(X)$}
и образующей $s$ группы Галуа $G$ (см. упражнение~\ref{exer-cyclic-alg}(i)).
Докажите, что образ класса $[D]$ в группе $\widetilde{\Br}\big(K(X)\big)$
относительно композиции
$$H^1\big(G, \Pic(X_L)\big)\to H^1\big(G_K, \Pic(X_{\bar{K}})\big)\to
\widetilde{\Br}(X)\to\widetilde{\Br}\big(K(X)\big)$$
равен классу $[A(f)]$.
(Указание: воспользуйтесь конструкцией изоморфизма из
упражнения~\ref{prob:Br0X}(ii), а также
упражнением~\ref{exer-cobound-map}(iv).)
В частности, класс~$[A(f)]$ принадлежит подгруппе
$$\Br(X)=\Br^{nr}\big(K(X)\big)\subset\Br\big(K(X)\big)\,.$$
\item[(ii)] Покажите непосредственно, что класс
$[A(f)]$ принадлежит подгруппе
$$\Br^{nr}\big(K(X)\big)\subset\Br\big(K(X)\big)\,,$$
то есть проверьте явно, что вычеты класса $[A(f)]\in\Br\big(K(X)\big)$
на все дивизоры на многообразии~$X$ равны нулю.
(Указание: по упражнению~\ref{exer-vychet-dlya-cycl}(ii)
имеется равенство
$$\res_E\big([A(f)]\big)=0$$
для любого неприводимого приведённого
дивизора $E$ на $X$, носитель которого не содержится в носителе
$\Supp(f)$ главного дивизора $(f)$. Для вычисления вычета $\res_D\big([A(f)]\big)$ рассмотрите такую рациональную функцию~\mbox{$h\in L(X)^*$}, что
$$(h)=D-D'\,,$$
где (не обязательно эффективный)
дивизор $D'$ не имеет общих компонент с дивизором
$$(f)=\sum_{g\in G}g(D)\,.$$
Далее воспользуйтесь тем, что класс циклической алгебры
$A\big(\Nm(h)\big)$ равен нулю в группе~\mbox{$\Br\big(K(X)\big)$},
см. упражнение~\ref{exer-cyclic-alg}(o). Здесь~$\Nm$
обозначает норму Галуа для расширения полей
$K(X)\subset L(X)$. Также используйте то, что имеется равенство
$$[A(f)]=\big[A\big(\Nm(h)\big)\big]+
\big[A\big(f\cdot\Nm(h)^{-1}\big)\big]$$
в группе $\Br\big(K(X)\big)$.)
\end{itemize}
\end{prob}

\bigskip
Чтобы вычислять препятствие Брауэра--Манина, удобно для заданного
элемента $\alpha\in\Br(X)$ найти условия на нормирование $v$ поля $K$,
при которых локальное спаривание $\langle \alpha, p_v\rangle_v$ равно
нулю для любой точки $p_v\in X(K_v)$.
Мы решим эту задачу для случая $\alpha=[A(f)]$, где
$f$ и $A(f)$ такие, как в упражнении~\ref{prob-Brdiv}.
Отметим, что спаривание~\mbox{$\langle \alpha, p_v\rangle_v$}
определено именно для элемента
$\alpha\in\Br(X)$, а не его образа в факторгруппе $\widetilde{\Br}(X)$.
При этом свойство \emph{лежать в ядре глобального спаривания}
зависит только от образа элемента группы $\Br(X)$ в
факторгруппе, см. упражнение~\ref{prob:BrManin}(i).

Введём ряд обозначений.
Пусть схема $\Xc$ является проективной моделью
многообразия $X$ над $\OO_K$ (см. разделы~\ref{subsection:models}
и~\ref{subsection:Br-Manin-def}).
Для неархимедова нормирования $v$ определим конечное поле
$\kappa_v=\OO_K/\mathfrak{p}$, где $\mathfrak{p}$~---
простой идеал в $\OO_K$, соответствующий нормированию~$v$.
Через $\Xc_v$ обозначим
произведение
$$\Xc_v=\Spec(\OO_v)\times_{\Spec(\OO_K)}\Xc\,,$$
а через $\Xc(v)$ обозначим слой схемы $\Xc_v$ над замкнутой точкой
$$\Spec(\kappa_v)\hookrightarrow\Spec(\OO_v)\,.$$
Отметим, что $\Xc(v)$ также является слоем схемы $\Xc$ над замкнутой
точкой
$$\Spec(\kappa_v)\hookrightarrow\Spec(\OO_K)\,.$$
Для рациональной функции $f$ из $K(X)^*$ (соответственно, из $K_v(X)^*$)
через $\SSupp(f)$
в дальнейшем мы будем обозначать носитель дивизора функции $f$ на всей
схеме~$\Xc$ (соответственно, на $\Xc_v$). При этом мы продолжаем обозначать через $\Supp(f)$ носитель дивизора функции $f$ на многообразии $X$ над $K$, являющемся общим слоем схемы $\Xc$.

Пусть $K\subset L$ --- расширение глобальных полей. Для любого нормирования
$v$ поля $K$ имеется разложение
\begin{equation}\label{eq:proizv-polej}
\mbox{$L\otimes_K K_v\cong\prod\limits_{i=1}^r L_{w_i}$}
\end{equation}
в произведение $r$ локальных полей
(отметим, что все утверждения
из упражнения~\ref{prob-val-extension} остаются
верными для произвольных нормирований, не только дискретных).
Будем говорить, что нормирование~$v$ \emph{распадается в $L$},
если все пополненные поля~$L_{w_i}$ изоморфны~$K_v$.
Неархимедово нормирование~$v$ \emph{неразветвлено в $L$}, если
для всех~$i$ неразветвлено расширение локальных полей $K_v\subset L_{w_i}$
(см. раздел~\ref{subsection:unramified-Brauer-generalities}).
Отметим, что все нормирования поля $K$, кроме конечного числа,
неразветвлены в данном расширении $L$
(см.~\mbox{\cite[Предложение~III.2.8]{Lang-chisla}}).

\bigskip
\begin{prob}{\bf Тривиальность локального спаривания}
\label{prob:uslovie-na-normirovanie}

Будем использовать обозначения и предположения из
упражнения~\ref{prob-Brdiv}. Как и выше, $v$ является нормированием поля $K$.

\begin{itemize}
\item[(i)]
Покажите, что образ класса $[A(f)]$ относительно естественного
отображения
$$\Br\big(K(X)\big)\to\Br\big(K_v(X)\big)$$
является классом
циклической алгебры над $K_v(X)$, построенной по
расширению $K_v(X)\subset L_w(X)$, элементу
$$f\in K(X)^*\subset K_v(X)^*$$
и образующей $s^r$ подгруппы
$$\Gal(L_w/K_v)\subset\Gal(L/K)$$
(см. упражнение~\ref{prob-val-extension}(v)),
где $w$~--- любое из возникающих в~\eqref{eq:proizv-polej}
нормирований~\mbox{$w_i$, $1\leqslant i\leqslant r$}.
(Указание: воспользуйтесь упражнениями~\ref{exer-invimage}(iv)
и~\ref{exer-cyclic-alg}(o).)
\item[(ii)]
Предположим, что нормирование $v$ распадается в $L$.
Покажите, что в этом случае
класс алгебры $A(f)$ в $\Br\big(K_v(X)\big)$ равен нулю.
(Указание: воспользуйтесь пунктом~(i).)
\item[(iii)]
Предположим, что функция $f$ регулярна в точке $p_v\in X(K_v)$,
и $f(p_v)\neq 0$. Покажите, что
$$p_v^*\big([A(f)]\big)\in\Br(K_v)$$
является циклической алгеброй, построенной по расширению $K_v\subset L_w$,
элементу $f(p_v)\in K_v^*$ и образующей $s^r$ группы $\Gal(L_w/K_v)$,
где $r$ и $w$~--- такие как в пункте~(i).
(Указание: воспользуйтесь пунктом~(i),
а также явным описанием циклической алгебры из
упражнения~\ref{exer-cyclic-alg}(i) и определением обратного
образа~$p_v^*$ из конца главы~\ref{section:unramified-Brauer}.)
\item[(iv)]
Предположим, что $v$ является неархимедовым нормированием, неразветвлённым в $L$. Как уже обсуждалось
в разделе~\ref{subsection:Br-Manin-def}, любая точка
$p_v\in X(K_v)$ определяет морфизм
$$p_v\colon\Spec(K_v)\to \Xc_v\,,$$
продолжающийся до морфизма
$$\tilde{p}_v\colon\Spec(\OO_v)\to\Xc_v\,.$$
Предположим, что для точки~$p_v$ носитель $\SSupp(f)$ не пересекается
с образом~\mbox{$\mathrm{Im}(\tilde{p}_v)$} на $\Xc_v$. Докажите,
что
$$\langle [A(f)], p_v\rangle_v=0\,.$$
(Указание: сначала покажите, что из условия на $p_v$ следует, что
функция $f$ регулярна в $p_v$, и выполнено
$$f(p_v)\in \OO_v^*\subset K_v\,.$$
Далее воспользуйтесь пунктом~(iii) и
упражнением~\ref{exer-vychet-dlya-cycl}(ii).)
\item[(v)]
Предположим, что неархимедово нормирование
$v$ неразветвлено в $L$, и носитель $\SSupp(f)$ не содержит
слой $\Xc(v)\subset\Xc_v$. Докажите,
что для любой точки~\mbox{$p_v\in X(K_v)$} существует рациональная функция
$h\in L_w(X)^*$, для которой носитель
$\SSupp\big(f\cdot\Nm(h)^{-1}\big)$
не пересекается с образом $\mathrm{Im}(\tilde{p}_v)$ на $\Xc_v$,
где $\Nm$ обозначает норму Галуа в расширении полей
$K_v(X)\subset L_w(X)$.
(Указание: положим
$$\Xc_w=\Spec(\OO_w)\times_{\Spec(\OO_v)}\Xc\,.$$
Рассмотрите дивизор $\mathcal{D}$ на схеме $\Xc_w$,
являющийся замыканием дивизора~\mbox{$D\subset X_L$} из
упражнения~\ref{prob-Brdiv}(o). Примените метод из решения
упражнения~\ref{prob-Brdiv}(ii) к дивизору $\mathcal{D}$ на $\Xc_w$,
используя относительное проективное вложение схемы
$\Xc$ над $\Spec(\OO_K)$.)
\item[(vi)]
Пусть множество $S$ состоит из всех таких
неархимедовых нормирований~$v$ поля~$K$,
что либо~$v$ разветвлено в~$L$, либо носитель $\SSupp(f)$
содержит слой~$\Xc(v)$, а также из всех вещественных нормирований
поля $K$, не распадающихся в $L$.
Покажите, что множество $S$ конечно.
Докажите, что если $v\not\in S$, то для любой точки
$p_v\in X(K_v)$ имеется равенство
$$\langle [A(f)], p_v\rangle_v=0\,.$$
(Указание: для неархимедовых нормирований воспользуйтесь
пунктами~(iv) и~(v), а также тем, что
$$[A(f\cdot\Nm(h)^{-1})]=[A(f)]$$
в группе $\Br\big(K_v(X)\big)$,
ср. с решением упражнения~\ref{prob-Brdiv}(ii).
Для архимедовых нормирований воспользуйтесь пунктом~(ii).)
\end{itemize}
\end{prob}

Заметим, что если кольцо $\OO_K$ факториально
(например, если $K=\Q$), то всегда существует такой элемент
$c\in K^*$, что носитель $\SSupp(c\cdot f)$
не содержит слой~$\Xc(v)$ ни для какого неархимедова
нормирования $v$.

\subsection{Препятствие Брауэра--Манина на кривой рода один}

В этом разделе мы применим препятствие Брауэра--Манина,
чтобы построить контрпример к принципу Хассе
для кривых рода один.
Напомним, что \emph{эллиптическая кривая}~$E$ над полем~$F$~---
это кривая рода~один над~$F$ с отмеченной $F$-точкой.

\begin{prob}{\bf Кривые рода один}\label{prob-genusone}
\hspace{0cm}
\begin{itemize}
\item[(i)]
Покажите, что любая гладкая проективная кривая $E$ рода один является
торсором
(см. определение~\ref{defin-torsorpoint})
над некоторой эллиптической кривой.
(Указание: рассмотрите якобиан~$J(E)$ кривой $E$
и воспользуйтесь теоремой Римана--Роха,
чтобы определить действие $J(E)$ на $E$.)
\item[(ii)]
Покажите, что для любой эллиптической кривой $E$
над конечным полем~$\Fb_q$ имеется равенство
$$H^1\big(G\,_{\Fb_q}, E(\bar\Fb_q)\big)=0\,.$$
(Указание: сначала докажите сюръективность морфизма
$$\Phi_q-\mathrm{Id}\colon E\,_{\bar\Fb_q}\to E\,_{\bar\Fb_q}\,,$$
рассмотрев его дифференциал,
где~$\Phi_q$ обозначает $\bar{\mathbb{F}}_q$-линейный морфизм Фробениуса,
заданный возведением координат точек в степень~$q$.
Затем воспользуйтесь упражнением~\ref{prob:H1-ot-Z-s-kryshkoj}.)
\item[(iii)]
Докажите, что любая гладкая проективная кривая рода
один над конечным полем $\Fb_q$ имеет $\Fb_q$-точку.
(Указание: воспользуйтесь пунктами~(i) и~(ii), а также
упражнением~\ref{prob:torsor-vs-H1}(i).)
\item[(iv)]
Докажите, что любая гладкая плоская
проективная кривая рода один с гладкой редукцией
(по модулю максимального идеала в кольце нормирования)
над неархимедовым локальным полем имеет точку над этим полем.
(Указание: воспользуйтесь пунктом~(iii), а также леммой Гензеля, как
в упражнении~\ref{prob:quadric-over-local-field}(ii).)
На самом деле, аналогичное утверждение верно для произвольных
(не обязательно плоских) гладких проективных кривых рода один с гладкой редукцией.
\end{itemize}
\end{prob}

Отметим, что Ленг доказал сюръективность отображения
$$g\mapsto\Phi_q(g)\cdot g^{-1}$$
для любой связной алгебраической группы
над полем $\bar{\mathbb{F}}_q$ (см.~\cite{Lang},
а также более короткое
доказательство Стейнберга в~\cite{Steinberg}).

\bigskip

Пусть $X$ --- гладкая проективная кривая рода один над глобальным полем $K$
характеристики нуль (мы не предполагаем наличие $K$-точки на $X$).
Предположим, что~$X$ имеет эффективный приведённый дивизор
$F\in\Div(X)$ степени два,
то есть что существует морфизм степени
два $\psi\colon X\to\Pb^1$ над $K$ (морфизм $\psi$ задаётся
линейной системой~$|F|$).
Обозначим через $Z\in\Div(X)$ дивизор ветвления морфизма~$\psi$.
Предположим, что $Z$ неприводим над $K$.
Пусть существует такое квадратичное расширение
$K\subset L$, что
$$Z=E+\sigma(E)\,,$$
где $E\in \Div(X_L)$
является эффективным приведённым дивизором степени два,
а $\sigma$ обозначает нетривиальную инволюцию поля $L$ над $K$.

\bigskip

\begin{prob}{\bf Построение элемента в группе Брауэра}
\label{prob:element-dlya-BrM}
\hspace{0cm}
\begin{itemize}
\item[(i)]
Покажите, что существует рациональная функция $f\in K(X)^*$,
для которой
$$Z-2F=(f)\,,$$
где, как и раньше, $(f)$ обозначает дивизор функции $f$.
(Отметим, что здесь не используется предположение о существовании поля~$L$.)
\item[(ii)]
Пусть $a\in K$ --- такой элемент, что $L=K(\sqrt{a})$.
Докажите, что класс
$$\alpha\in\Br\big(K(X)\big)$$
кватернионной алгебры $A_{-1}(a,f)$ принадлежит
подгруппе
$$\Br(X)\subset\Br\big(K(X)\big)\,.$$
(Указание: воспользуйтесь
упражнением~\ref{prob-Brdiv} для дивизора $D=E-F$.)
\end{itemize}
\end{prob}

\bigskip
Теперь
мы в явном виде построим кривую рода один,
доставляющую контрпример к принципу Хассе.
(Впервые этот пример был найден
К.-Э.\,Линдом в~\cite{Lind} и Г.\,Рейхардом в~\cite{Reichardt}.)
Рассмотрим гладкую проективную кривую $X$ над $\Q$,
открытое аффинное подмножество $U$ которой задаётся в $\Ab^2$
с координатами $x$ и $y$ уравнением
$$2y^2=x^4-17\,.$$
Заметим, что $X$ является объединением аффинной карты $U$ и другой
аффинной карты $V$, которая естественно отождествляется с кривой,
заданной в $\Ab^2$ с координатами $z$ и $w$ уравнением
$$2z^2=1-17w^4\,.$$ На $X$ имеется эффективный
приведённый дивизор степени два $F=\{w=0\}$.
При этом во введённых перед упражнением~\ref{prob:element-dlya-BrM}
обозначениях имеем:
$$Z=\{y=0\}=\{z=0\}\subset U\cap V,\quad L=\Q(\sqrt{17}),\quad
f=y=\frac{z}{w^2}.$$

\bigskip

Пусть $v_{\infty}$ обозначает (единственное) архимедово
нормирование поля $\Q$.

\begin{prob}{\bf Препятствие Брауэра--Манина для кривой $X$}
\hspace{0cm}
\begin{itemize}
\item[(o)] Покажите, что единственное простое число $p$,
для которого нормирование $v_p$ поля $\Q$ разветвлено
в расширении $\Q\subset\Q(\sqrt{17})$,
есть $p=17$. (Указание: нас интересует расширение локальных полей
$$\Q_p\subset\Q_p(\sqrt{17})\,.$$
Сначала проверьте, что при $p\not\in\{2, 17\}$ это расширение
неразветвлено (см.~раздел~\ref{subsection:unramified-Brauer-generalities}).
Затем проверьте, что нормирование $v_2$
распадается в расширении~\mbox{$\Q\subset\Q(\sqrt{17})$},
воспользовавшись тем, что $17$ является
квадратом в поле~$\Q_2$.)
\item[(i)]
Покажите, что для любого нормирования
$$v\not\in\{v_{\infty}, v_2, v_{17}\}$$
поля $\Q$ существует точка $p_v\in X(\Q_v)$.
(Указание: воспользуйтесь упражнением~\ref{prob-genusone}(iv).)
\item[(ii)]
Покажите, что для любого нормирования
$$v\in\{v_{\infty}, v_2, v_{17}\}$$
существует точка $p_v\in X(\Q_v)$.
(Указание: для $v=v_{\infty}$ утверждение очевидно.
Для $v=v_2$ имеем $Z(\Q_2)\neq\varnothing$, так как $17$ является четвёртой
степенью в $\Q_2$. Для $v=v_{17}$
имеем $F(\Q_{17})\neq\varnothing$, так как $2$ является квадратом
в $\Q_{17}$.)
\item[(iii)]
Покажите, что для любого нормирования $v\neq v_{17}$
и любой точки~\mbox{$p_v\in X(\Q_v)$} имеется равенство
$$\langle \alpha, p_v\rangle_v=0\,,$$
где $\alpha$ определяется как
в упражнении~\ref{prob:element-dlya-BrM}(ii).
(Указание:
воспользуйтесь пунктом~(o) и
упражнением~\ref{prob:uslovie-na-normirovanie}(vi).)
\item[(iv)]
Покажите, что для любой точки $p_{17}\in X(\Q_{17})$
выполнено
$$\langle \alpha, p_{17}\rangle_{v_{17}}\neq 0\,.$$
(Указание:
любая точка $p_{17}\in X(\Q_{17})$ имеет координаты из $\Z_{17}$
на $U$ или на $V$. Для $p_{17}\in U(\Z_{17})$
локальное спаривание
$$\langle \alpha, p_{17}\rangle_{v_{17}}\in\Q/\Z$$
соответствует символу Гильберта
$\big(17, y(p_{17})\big)_{17}$
как описано в упражнении~\ref{exer-vychet-dlya-cycl}(iv).
Для вычисления символа Гильберта покажите, что
$$v_{17}\big(y(p_{17})\big)=0\,,$$
в то время как~$y(p_{17})$ не является квадратом по модулю $17$,
поскольку $2$ не является четвёртой
степенью в $\mathbb{F}_{17}$. Аналогично рассматриваются
точки~\mbox{$p_{17}\in V(\Z_{17})$}.)
\item[(v)]
Докажите, что $X(\Q)=\varnothing$. Таким образом,
$X$ является контрпримером к принципу Хассе.
(Указание: воспользуйтесь пунктами~(iii) и~(iv),
а также упражнением~\ref{prob:BrManin}(iii).)
\end{itemize}
\end{prob}

\newpage
\appendix

\section{Этальные когомологии}
\label{section:etale}

В данном приложении перечисляются некоторые общие понятия и
утверждения из теории этальных когомологий, касающиеся группы
Брауэра многообразий. Прекрасный обзор
теории этальных когомологий содержится в~\cite{Dan}, а
подробное изложение всех фактов и доказательств содержится в~\cite{Mil}.
Ниже мы приводим более детальные ссылки на оба источника.
Также мы рекомендуем читателю книгу~\cite{FrK}.

\subsection{Этальные покрытия}

{\it Этальный морфизм} между схемами  --- это гладкий морфизм
относительной размерности нуль. Также
есть несколько равносильных определений этального морфизма, \cite[4.2.1]{Dan},
\cite[I.3]{Mil}. Расширение полей соответствует этальному
морфизму между их спектрами тогда и только тогда, когда оно
является конечным сепарабельным расширением.
{\it Этальное покрытие} схемы~$X$~---
это набор этальных морфизмов~${\{U_i\to X\}}$,
объединение образов которых равно~$X$, \mbox{\cite[II.1]{Mil}}.
Если~$X$ квазикомпактна, то обычно достаточно рассматривать только
конечные этальные покрытия, то есть конечные наборы~${\{U_i\to X\}}$.
В этом случае удобнее перейти к дизъюнктному
объединению~${U=\coprod_i U_i}$ и (сюръективному)
этальному морфизму~${U\to X}$. Класс этальных покрытий замкнут относительно взятия
композиции и расслоенных произведений, то есть
этальные покрытия образуют топологию
Гротендика,~\cite[Замечание~II.1.1]{Mil}.

Одной из основных причин
успеха использования этальной топологии является следующая лемма, принадлежащая М.\,Артину: для любой точки $x$ на гладком комплексном алгебраическом многообразии существует окрестность $U$ в топологии Зарисского, у которой тривиальны все
гомотопические группы, начиная со второй,~\cite[3.4.3]{Dan}. Таким образом, единственным препятствием к тому, чтобы окрестность $U$ была стягиваемой в комплексной топологии, является её фундаментальная группа. Устранение этого препятствия с ``проконечной'' точки зрения достигается за счёт перехода к конечным неразветвлённым накрытиям окрестности $U$, то есть рассмотрения этальных окрестностей точки $x$.

\subsection{Пучки в этальной топологии}

{\it Этальный предпучок множеств $\Fc$} на схеме $X$--- это контравариантый
функтор~\mbox{$U\mapsto\Fc(U)$}
из категории этальных схем над $X$ в категорию множеств,
\mbox{\cite[4.3.1]{Dan}},
\cite[II.1]{Mil}. Этальный предпучок множеств является {\it этальным пучком},
если он переводит прямые пределы в обратные, \cite[4.3.2]{Dan}.
Равносильное определение получится, если потребовать, что
для любого этального покрытия~\mbox{$\{V_i\to U\}$}
этальной схемы $U$ над $X$ отображение
$$\Fc(V)\to\coprod_i\Fc(V_i)$$ является
уравнителем двух отображений из $\coprod_i\Fc(V_i)$
в $\coprod_{i,j}\Fc(V_i\times_U V_j)$, соответствующих проекциям из
$V_i\times_U V_j$ на $V_i$  и $V_j$, \cite[II.1]{Mil}.
Например, категория этальных пучков на спектре
поля $K$ эквивалентна категории дискретных множеств
с непрерывным действием группы Галуа $G_K$
поля~$K$, \cite[4.3.2]{Dan}.

Произвольная схема $S$ представляет этальный пучок
$$U\mapsto h_S(U)=\Hom(U,S)\,,$$
\cite[Пример~4.3.1.1]{Dan},
\cite[Следствие~II.1.7]{Mil}. В частности, с каждой групповой схемой
связан этальный пучок групп, который обычно обозначается так же,
как и исходная групповая схема.
Например, имеются абелевы
этальные пучки групп~$\Gb_m$ и $\mu_l$, представленные коммутативными
групповыми схемами~\mbox{$\Spec(\Z[T,T^{-1}])$}
и~\mbox{$\Spec\big(\Z[T]/(T^l-1)\big)$}, соответственно,
\cite[Пример~4.4.1.2]{Dan},
\cite[Пример~II.2.18(b)]{Mil}. Для любого гладкого морфизма
схем $X\to Y$ морфизм этальных пучков \mbox{$h_X\to h_Y$} сюръективен,
\cite[4.3.3]{Dan}.

Произвольный квазикогерентный
пучок $\Fc$ на $X$ определяет этальный пучок $\Fc_{\acute e t}$ на~$X$ по формуле
$\Fc_{\acute e t}(U)=\Gamma(U,\pi^*\Fc)$
для любого этального морфизма $\pi\colon U\to X$,
\cite[Пример~4.3.1.2]{Dan},
\cite[Следствие~II.1.6]{Mil}.

\subsection{Когомологии абелевых этальных пучков}

{\it Когомологии абелевых пучков} на схеме $X$ ---
это правые производные функторы точного
слева функтора взятия глобальных сечений $\Fc\mapsto\Fc(X)$,
\cite[4.4.2]{Dan},
\mbox{\cite[Определение~III.1.5(a)]{Mil}}.
Обычно когомологии абелевых этальных пучков обозначаются
$H^i_{\acute e t}(X,\Fc)$. Для этальных пучков можно также рассматривать когомологии Чеха,
построенные при помощи перехода к прямому пределу по всем этальным покрытиям,
\cite[III.2]{Mil}.
Для квазикомпактной схемы, каждое конечное подмножество которой содержится в некотором открытом аффинном подмножестве (например, для квазипроективного многообразия), когомологии Чеха совпадают с определёнными выше
когомологиями для любого абелевого этального пучка,
\mbox{\cite[Теорема~III.2.17]{Mil}}.
Такие схемы мы будем называть {\it хорошими}.

Когомологии абелевого этального пучка на спектре поля $K$ канонически изоморфны
когомологиям Галуа соответствующего абелевого $G_K$-модуля,
\cite[4.4.3]{Dan}.
При этом прямой предел комплексов
Чеха по всем этальным покрытиям соответствует стандартному
комплексу, связанному с когомологиями групп (см. упражнение~\ref{exer-stcoml}).

Из леммы Артина следует совпадение
этальных когомологий и когомологий Бетти с постоянными конечными
коэффициентами,
\cite[4.6.5]{Dan}, \mbox{\cite[Теорема~III.3.12]{Mil}}.

Для любого квазикогерентного пучка $\Fc$ на $X$ имеется
канонический изоморфизм
$$H^i_{Zar}(X,\Fc)\cong H^i_{\acute e t}(X,\Fc_{\acute e t})\,,$$
где индекс $Zar$ соответствует когомологиям в топологии Зарисского,
а индекс~\mbox{$\acute e t$}~--- когомологиям в этальной топологии.
Подробности можно найти в~\cite[4.4.4]{Dan},
\cite[Предложение~III.3.7]{Mil},
\mbox{\cite[Замечание~III.3.8]{Mil}}.

Этальные когомологии гладкого
квазипроективного многообразия с коэффициентами в постоянном пучке $\Q$ тривиальны в положительных степенях: это следует из
аналогичных утверждений о когомологиях Зарисского
и о когомологиях проконечных групп, из возможности
вычислять этальные когомологии при помощи
когомологий Чеха и из спектральной последовательности
Хохшильда--Серра, \mbox{\cite[Теорема~III.2.20]{Mil}}.
Поэтому часто рассматриваются этальные когомологии
с коэффициентами в пучках кручения.
Пусть простое число $l$ не делит характеристику схемы~$X$.
Для натурального $n$ положим
$$\mu_l^{\otimes -n}=\Hom(\mu_l^{\otimes n},\Z/l\Z)\,,$$
и определим
{\it $l$-адические когомологии} по формуле
$$
H^i_{\acute e t}\big(X,\Z_l(n)\big)=
\varprojlim_r H^i_{\acute e t}\big(X,\mu_{l^r}^{\otimes n}\big),\quad
H^i_{\acute e t}\big(X,\Q_l(n)\big)=
\Q_l\otimes_{\Z_l}H^i_{\acute e t}\big(X,\Z_l(n)\big)\,,
$$
см.~\cite[4.7.1]{Dan}.
Для $l$-адических когомологий определено отображение обратного образа.
Для гладкого проективного многообразия $X$ размерности $d$
над сепарабельно замкнутым полем есть двойственность Пуанкаре,
то есть определено невырожденное спаривание
$$
H^i_{\acute e t}\big(X,\Q_l(n)\big)\otimes_{\Q_l}
H^{2d-i}_{\acute e t}\big(X,\Q_l(d-n)\big)\to \Q_l\,,
$$
см.~\cite[4.7.5]{Dan}, \cite[Теорема~VI.11.1]{Mil}.
Имеется отображение
прямого образа (гомоморфизм Гизина) на $l$-адических когомологиях,
удовлетворяющее формуле проекции,~\cite[4.7.6]{Dan}.
Из этого следует аналог формулы следа Лефшеца
из топологии,~\cite[4.7.8]{Dan}.
Также, можно рассматривать когомологии с компактным носителем
и с коэффициентами
в произвольном конструктивном $l$-адическом пучке, для которых выполнено
много утверждений, аналогичных теоремам для пучков на многообразиях
в классической топологии, \cite[4.7]{Dan},
\cite[Глава~VI]{Mil}.
Всё это позволило П.\,Делиню доказать при помощи $l$-адических когомологий
гипотезы Вейля о дзета-функциях многообразий
над конечными полями.
Подробности см. в~\cite{Del74}, \cite{Del80},
\cite[4.7.9, 4.7.10, 4.8]{Dan},
 \cite[VI.12, VI.13]{Mil}.

\subsection{Первые этальные когомологии с неабелевыми коэффициентами}
\label{subsection:etale-nonab}

Рассмотрим следующий вопрос: пусть дано (конечное) этальное покрытие \mbox{$U\to X$} и пусть дан некоторый алгебро-геометрический объект $B$ над $U$; какие условия гарантируют то, что $B$ является обратным образом некоторого объекта~$A$ над~$X$? Как описать все такие $A$?
Оказывается, что все основные определения и конструкции из
главы~\ref{section:Galois-cohomology}
обобщаются со случая спектра поля на случай произвольной схемы~$X$.
Более подробно об этом можно прочитать в~\cite{Grothendieck1959},
а также в обзоре~\cite{Vistoli}.

Итак, пусть даны схема $X$ и категория $\Mcal$,
расслоенная над этальными схемами над $X$,
то есть для каждого этального морфизма $U\to X$ определена
категория~$\Mcal(U)$, а для каждого морфизма $q\colon V\to U$
между этальными схемами над~$X$ определён функтор обратного
образа
$$q^*\colon\Mcal(U)\to\Mcal(V)\,.$$
При этом требуется, чтобы функторы
обратного образа были согласованы с композициями морфизмов между
этальными схемами над $X$. Поскольку для поставленных выше вопросов играют роль лишь изоморфизмы между объектами, будем также предполагать, что для любого~$U$ категория $\Mcal(U)$ является {\it группоидом}, то есть все морфизмы в $\Mcal(U)$ являются изоморфизмами. Например,
$\Mcal(U)$ может быть категорией векторных расслоений над~$U$,
или категорией схем конечного типа над $U$ (в качестве морфизмов
в обоих случаях берутся изоморфизмы).

Для простоты будем рассматривать только конечные этальные покрытия
(это достаточно, если $X$ квазикомпактна).
Пусть $f\colon U\to X$ --- этальное покрытие, и
$$p_i\colon U\times_X U\to U,\quad
p_{ij}\colon U\times_X U\times_X U\to U\times_X U$$
обозначают естественные проекции. Будем говорить, что на объекте $B$
из $\Mcal(U)$ определены {\it данные спуска}, если задан такой изоморфизм $\rho\colon p_1^*B\stackrel{\sim}\longrightarrow p_2^*B$, что
$$
p_{13}^*(\rho)=p_{23}^*(\rho)\circ p_{12}^*(\rho)\colon p_1^*B\to p_3^*B\,.
$$
Например, для любого объекта $A$ из $\Mcal(X)$
имеются канонические данные спуска на~$f^*A$.
Будем говорить, что $\Mcal$ {\it удовлетворяет условию спуска},
если функтор обратного образа $f^*$ осуществляет эквивалентность между
$\Mcal(X)$ и категорией объектов
из~$\Mcal(U)$ с данными спуска. Если такое условие выполнено не только для $X$, но и для любой этальной схемы $V$ над $X$,
то $\Mcal$ называется {\it стеком}. В этом случае
можно сказать, что $\Mcal$ удовлетворяет условию пучка в этальной топологии.

Для объекта $A$ из $\Mcal(X)$ определён этальный
предпучок групп $\underline{\rm Aut}(A)$, заданный
по формуле
$V\mapsto{\rm Aut}(q^*A)$ для этального морфизма \mbox{$q\colon V\to X$}.
Как и в упражнении~\ref{exer-spusk}(i),
имеется биекция между множеством классов изоморфизма данных спуска на $f^*A$ и
множеством $\check{H}^1_{\acute e t}\big(U/X,\underline{\rm Aut}(A)\big)$,
где первые этальные когомологии Чеха
с неабелевыми коэффициентами обобщают определение из
упражнения~\ref{exer-nonabcohom}(i), см.~\mbox{\cite[III.4]{Mil}}.
Далее, как в упражнении~\ref{exer-spusk}(iii),
если~$\Mcal$ удовлетворяет условию спуска,
то множество $\check{H}^1_{\acute e t}\big(U/X,\underline{\rm Aut}(A)\big)$
биективно множеству классов изоморфизма
$U$-форм объекта $A$, то есть таких объектов $A'$
из $\Mcal(X)$, что
\mbox{$f^*A\cong f^*A'$}, \mbox{\cite[III.4]{Mil}}, \mbox{\cite[4.4.5]{Dan}}.

Заметим также, что если $\Mcal$ является стеком,
то $\underline{\rm Aut}(A)$ является пучком в этальной топологии (последнее свойство выполнено для большинства разумных примеров).
Аналогично упражнениям~\ref{exer-nonabcohom}(ii) и~\ref{exer-nonabcohom}(iii),
имеется точная последовательность отмеченных множеств,
ассоциированная с короткой точной последовательностью пучков неабелевых
групп на хорошей схеме, см. шаг~3 в доказательстве
теоремы~IV.2.5 из~\cite{Mil}.

Категория квазикогерентных пучков удовлетворяет условию
спуска и является стеком, см., например, \mbox{\cite[Замечание~I.2.21]{Mil}},
\cite[Предложение~I.2.22]{Mil}. Доказательство этого факта в целом следует методу из упражнения~\ref{exer-spuskvp}. В частности, категория векторных расслоений удовлетворяет условию спуска, откуда следует, что каноническое отображение
$$\check{H}^1_{Zar}(X,\GL_n)\to\check{H}^1_{\acute e t}(X,\GL_n)$$
является биекцией,~\cite[Предложение~III.4.9]{Mil}.
Например, имеется изоморфизм
$$\Pic(X)\cong H^1_{\acute e t}(X,\Gb_m)\,.$$

Следующие категории также удовлетворяют условию спуска (и являются стеками):
категория аффинных схем,~\cite[Теорема~I.2.23]{Mil},
категория торсоров относительно гладкой
аффинной групповой схемы, \cite[Следствие~III.4.7]{Mil},
\mbox{\cite[Замечание~III.4.8]{Mil}},
категория схем с относительно обильным
пучком, \mbox{\cite[Proposition~VIII.7.8]{SGA1}}.
В частности, из последнего примера следует, что
категория схем Севери--Брауэра, то есть семейств многообразий Севери--Брауэра, удовлетворяет
условию спуска, и имеется биекция между множеством классов
изоморфизма схем Севери--Брауэра
и множеством~$\check{H}^1_{\acute{e}t}(X,\PGL_n)$,
см.~\mbox{\cite[III.4]{Mil}}.
Тем не менее условие спуска выполнено далеко не для
всех естественных категорий схем, например, категория
эллиптических схем не является стеком,~\cite[Example~4.39]{Vistoli}
(это происходит потому, что у неквазипроективного многообразия
может не существовать фактора по свободному
действию конечной группы в категории схем).
Это обстоятельство приводит к понятию алгебраического пространства.

\subsection{Последовательность Куммера}
\label{subsection:etale-Kummer}

Пусть простое число $l$ не делит характеристику схемы $X$,
и $r$ --- натуральное число.
Имеется точная последовательность этальных пучков, называемая {\it
последовательностью Куммера},~\cite[III.4]{Mil},~\cite[4.4.6]{Dan}:
\begin{equation}\label{eq:etal-Kummer}
1\to
\mu_{l^r}\to\Gb_m\stackrel{l^r}\to\Gb_m\to 1\,.
\end{equation}
Эта последовательность является аналогом экспоненциальной
последовательности пучков на комплексных многообразиях
и обобщает последовательность Куммера из
упражнения~\ref{exer-Kummer}(i).
Как было сказано в разделе~\ref{subsection:etale-nonab}, имеется изоморфизм
$$\Pic(X)\cong H^1_{\acute e t}(X,\Gb_m)\,.$$
Поэтому из длинной точной когомологической последовательности получаем короткие точные последовательности
\begin{align}
\label{seq-tors-etale-1}
&1\to \Gamma(X,\OO_X^*)/l^r\to
H^1_{\acute e t}(X,\mu_{l^r})\to \Pic(X)_{l^r}\to 0\,,\\
\label{seq-tors-etale-2}
&0\to \Pic(X)/l^r\to
H^2_{\acute e t}(X,\mu_{l^r})\to H^2_{\acute e t}(X,\Gb_m)_{l^r}\to 0\,.
\end{align}
В частности, если $X$ --- полная схема над сепарабельно замкнутым полем,
то из~\eqref{seq-tors-etale-1} возникает изоморфизм
$$H^1_{\acute e t}(X,\mu_{l^r})\cong \Pic(X)_{l^r}\,.$$
Отметим также, что точная
последовательность~\eqref{seq-tors-etale-2}
является обобщением изоморфизма из упражнения~\ref{exer-brtor}(iii).

Для произвольной схемы~$X$,
переходя в~\eqref{seq-tors-etale-2}
к прямому пределу по~$r$, получаем точную последовательность
\begin{equation}\label{eq:posle-predela}
0\to \Pic(X)\otimes_{\Z}\Q_l/\Z_l\to H^2_{\acute e
t}(X,\mu_{l^{\infty}})\to H^2_{\acute e t}(X,\Gb_m)_{l^{\infty}}\to
0\,,
\end{equation}
где $A_{l^{\infty}}=\bigcup_{r\geqslant 1}A_{l^r}$ для произвольной абелевой группы $A$; в частности, имеется равенство пучков абелевых групп~\mbox{$\mu_{l^{\infty}}=\bigcup_{r\geqslant 1}\mu_{l^r}$}. Если характеристика схемы $X$ равна нулю, то
имеется изоморфизм
\begin{equation}
\label{eq:coker-vs-H2}
H^2_{\acute e t}(X,\Gb_m)_{tors}\cong
\Coker\Big(\Pic(X)\otimes_{\Z}\Q\to
H^2_{\acute e t}(X,\mu_{\infty})\Big)\,,
\end{equation}
где $\mu_{\infty}=\prod_{l}\mu_{l^{\infty}}$. Таким образом, можно сказать, что группа
$H^2_{\acute e t}(X,\Gb_m)_{tors}$
является группой трансцендентных классов во вторых этальных когомологиях.

Важным свойством группы $H^2_{\acute e t}(X,\Gb_m)_{tors}$ является
её стабильная бирациональная инвариантность для гладких
проективных многообразий $X$ над полем $k$ характеристики нуль.
Действительно, из~\eqref{eq:coker-vs-H2} легко видеть, что
группа $H^2_{\acute e t}(X,\Gb_m)_{tors}$
не меняется при раздутиях с гладкими центрами и при умножении на
$\P^1$. Пусть теперь дано бирациональное отображение
$f\colon X\dasharrow Y$ между гладкими проективными многообразиями
над $k$.
По теореме Хиронаки существует гладкое проективное многообразие
$\widetilde{X}$ и бирациональные морфизмы $\pi$ и $\tilde{f}$,
образующие коммутативную диаграмму
$$
\xymatrix{
&\widetilde{X}\ar@{->}[ld]_{\pi}\ar@{->}[rd]^{\tilde{f}}&\\
X\ar@{-->}[rr]^{f}&&Y
}
$$
При этом можно считать, что морфизм
$\pi$ является последовательностью раздутий
с гладкими центрами. В частности, имеется
изоморфизм
$$\pi^*\colon
H^2_{\acute e t}(X,\Gb_m)_{tors}\stackrel{\sim}\longrightarrow
H^2_{\acute e t}(\widetilde{X},\Gb_m)_{tors}\,,$$
а также отображение
$$
(\pi^*)^{-1}\circ \tilde{f}^*\colon
H^2_{\acute e t}(Y,\Gb_m)_{tors}\to
H^2_{\acute e t}(X,\Gb_m)_{tors}\,.$$
Аналогичным образом, регуляризуя
бирациональное отображение $f^{-1}$, мы получаем отображение
$$
H^2_{\acute e t}(X,\Gb_m)_{tors}\to
H^2_{\acute e t}(Y,\Gb_m)_{tors}\,.$$
Можно проверить, что
эти отображения взаимно обратны.

\subsection{Группа Брауэра}
\label{subsection:etale-Brauer}

Основные определения и конструкции из главы~\ref{section:Brauer}
обобщаются с полей на схемы.
Используя построенную выше теорию спуска и действуя аналогично
пунктам~(i),~(ii) и~(iii) упражнения~\ref{exer-cohomBr},
можно показать, что имеется каноническое вложение
групп (см.~\mbox{\cite[Теорема~IV.2.5]{Mil}})
$$\lambda\colon\Br(X)\hookrightarrow H^2_{\acute e t}(X,\Gb_m)\,,$$
где $\Br(X)$ определяется в терминах алгебр Адзумаи,
см. определение~\ref{defin:Br-scheme}
(для такого подхода надо ещё предполагать, что схема $X$ хорошая).

Подобно упражнению~\ref{exer-brtor}(ii)
можно показать, что группа $\Br(X)$ состоит из элементов
кручения,~\cite[Предложение~IV.2.7]{Mil}
(для случая алгебр, ранг которых не взаимно прост с характеристикой,
надо рассматривать плоскую топологию).
Важно отметить, что упражнение~\ref{exer-cohomBr}(iv)
не обобщается на случай произвольных схем,
то есть вложение $\lambda$ не обязательно является изоморфизмом.
В частности, в~\mbox{\cite[Remarque~1.11]{Grothendieck66b}}
дан пример особой поверхности $X$,
для которой группа~\mbox{$H^2_{\acute e t}(X,\Gb_m)$}
имеет элемент бесконечного порядка. Тем не менее, относительно недавно
де Йонгом~\cite{deJong} было доказано, что для любой схемы $X$ с обильным
обратимым пучком
упоминавшееся выше вложение $\lambda$ даёт изоморфизм
$$\lambda\colon
\Br(X)\stackrel{\sim}\longrightarrow H^2_{\acute e t}(X,\Gb_m)_{tors}\,.$$

Для регулярного неприводимого многообразия $X$ над полем $k$
несложно показать
инъективность естественного отображения
$$\nu\colon H^2_{\acute e t}(X,\Gb_m)\to
\Br\big(k(X)\big)\,,$$
см.~\cite[Пример~III.2.22]{Mil}.
В~частности, поскольку группа $\Br\big(k(X)\big)$ состоит из элементов
кручения (см.~упражнение~\ref{exer-brtor}(o)),
для таких $X$ возникает равенство
\begin{equation}
\label{eq:tolko-tors}
H^2_{\acute e t}(X,\Gb_m)_{tors}= H^2_{\acute e t}(X,\Gb_m)\,.
\end{equation}
Таким образом, из приведённых в разделе~\ref{subsection:etale-Kummer}
рассуждений следует стабильная бирациональная инвариантность группы
$H^2_{\acute e t}(X,\Gb_m)$ для гладких проективных многообразий
над полем нулевой характеристики.

Для гладкого многообразия $X$ над полем характеристики нуль
отображение $\nu$ осуществляет изоморфизм
$$
\nu\colon H^2_{\acute e t}(X,\Gb_m)\stackrel{\sim}\longrightarrow
\Br^{nr}(X)\,.$$
Это следует из чистоты для этальных когомологий с конечными коэффициентами
и последовательности Куммера из раздела~\ref{subsection:etale-Kummer}.
Подробности можно найти в~\mbox{\cite[\S3.4]{CT}}.
Отсюда следует стабильная бирациональная инвариантность
группы~$\Br^{nr}(X)$ для гладких проективных многообразий
над полем нулевой характеристики.

Знаменитая гипотеза М.\,Артина предполагает конечность группы
Брауэра $\Br(X)$
для любой собственной схемы $X$ над $\Z$, \cite[Вопрос~IV.2.19]{Mil}.

\subsection{Случай комплексного алгебраического многообразия}
\label{subsection:etale-complex}

Пусть $X$ --- гладкое комплексное проективное многообразие.
Мы будем рассматривать топологию на множестве его $\Cb$-точек
$X(\Cb)$, заданную классической топологией на $\Cb$.
Через $H^i\big(X(\Cb),-\big)$
будем обозначать группы когомологий
Бетти топологического пространства~$X(\Cb)$,
или группы когомологий с коэффициентами в каком-либо
пучке на этом топологическом пространстве.

Выбор корней из единицы определяет изоморфизм
\begin{equation}\label{eq:H2et-H2an}
H^2_{\acute e t}(X,\mu_{\infty})\cong H^2\big(X(\Cb),\Q/\Z\big)\,.
\end{equation}
Из~\eqref{eq:coker-vs-H2}, \eqref{eq:tolko-tors} и~\eqref{eq:H2et-H2an},
а также упомянутого в разделе~\ref{subsection:etale-Brauer}
результата де~Йонга,
получаем изоморфизм
$$
\Br(X)\cong\Coker\Big(\Pic(X)\otimes_{\Z}\Q\stackrel{r}\longrightarrow
H^2\big(X(\Cb),\Q/\Z\big)\Big)\,.
$$
С другой стороны, из точной последовательности постоянных пучков в классической топологии
$$
0\to\Z\to\Q\to\Q/\Z\to 0
$$
возникает точная последовательность
$$
H^2\big(X(\Cb), \Z\big)\to H^2\big(X(\Cb), \Q\big)\to
H^2\big(X(\Cb),\Q/\Z\big)\to
H^3\big(X(\Cb), \Z\big)_{tors}\to 0\,.
$$
Поскольку отображение $r$ пропускается через $H^2\big(X(\Cb), \Q\big)$,
имеется короткая точная последовательность
\begin{multline*}
0\to\Coker\Big(\big(\Pic(X)\otimes_{\Z}\Q\big)\oplus H^2\big(X(\Cb), \Z\big)\to
H^2\big(X(\Cb), \Q\big)\Big)\to{}\\
{}\to \Br(X)\to
H^3\big(X(\Cb),\Z\big)_{tors}\to 0\,.
\end{multline*}
Наконец, (доказанная) гипотеза Ходжа для дивизоров утверждает, что
$$
\Big(H^0(X, \Omega^2_X)\oplus
H^1(X, \Omega^1_X)\Big)\cap H^2\big(X(\Cb), \Q\big)=
\mathrm{Im}\Big(\Pic(X)\otimes_{\Z}\Q\to H^2\big(X(\Cb), \Cb\big)\Big)\,,
$$
где $\Omega^i_X$ обозначает пучок алгебраических $i$-форм на многообразии $X$.
Поэтому
естественная проекция $H^2\big(X(\Cb), \Cb\big)\to H^2(X, \OO_X)$
индуцирует изоморфизм
\begin{multline*}
\Coker\Big(\big(\Pic(X)\otimes_{\Z}\Q\big)\oplus H^2\big(X(\Cb), \Z\big)\to
H^2\big(X(\Cb), \Q\big)\Big)\cong{}\\
{}\cong
\Coker\Big(H^2\big(X(\Cb), \Z\big)\to H^2(X, \OO_X)\Big)_{tors}
\end{multline*}
Таким образом, имеется точная последовательность
\begin{equation}\label{eq:Br-vs-H3}
0\to
\Coker\Big(H^2\big(X(\Cb), \Z\big)\to H^2(X, \OO_X)\Big)_{tors}\to
\Br(X)\to
H^3\big(X(\Cb),\Z\big)_{tors}\to 0\,.
\end{equation}
В частности, если нетривиальна группа $H^3\big(X(\Cb),\Z\big)_{tors}$, то
нетривиальна и группа~$\Br(X)$.
Используя~\eqref{eq:Br-vs-H3} и экспоненциальную последовательность,
можно увидеть, что естественное отображение
$$\Br(X)\to H^2\big(X(\Cb), \mathcal{A}_X^*\big)_{tors}$$
является изоморфизмом,
где $\mathcal{A}_X^*$~--- пучок обратимых голоморфных
функций на комплексном аналитическом многообразии~\mbox{$X(\Cb)$}.

\end{document}